%% file: thesis.tex
\documentstyle[12pt,epsf,epsfig,amsfonts,twoside]{article}

\def\part#1{\frac{\partial\phantom{#1}}{\partial#1}}
\newtheorem{thm}{Theorem}
\newtheorem{prp}[thm]{Proposition}
\newtheorem{lem}[thm]{Lemma}
\newtheorem{cor}[thm]{Corollary}

\newenvironment{prf}{\begin{trivlist}\item[]{\bf Proof} }%
{\hfill $\Box$ \end{trivlist}}
\newenvironment{dfn}{\begin{trivlist}\item[]{\bf Definition}\em }%
{\end{trivlist}}
\newenvironment{rmk}{\begin{trivlist}\item[]{\bf Remarks} }%
{\end{trivlist}}

\newcommand{\twoVgraph}{\raisebox{0pt}{
                 \begin{picture}(18,18)(-9,-5)
                 \put(0,0){\circle{16}} \put(-8,0){\line(1,0){16}}
                 \end{picture}}}

\newcommand{\fourVgraph}{\raisebox{0pt}{
                 \begin{picture}(18,26)(-9,-9)
                 \put(0,0){\oval(16,24)} \put(-8,4){\line(1,0){16}}
                 \put(-8,-4){\line(1,0){16}}
                 \end{picture}}}

\newcommand{\sixVgraph}{\raisebox{0pt}{
                 \begin{picture}(26,24)(-13,-8)
                 \put(-9,8){\circle{6}} \put(9,8){\circle{6}}
                 \put(-6,8){\line(1,0){12}} \put(0,-8){\circle{6}} 
                 \put(-9,5){\line(2,-3){7}} \put(9,5){\line(-2,-3){7}}
                 \end{picture}}}

\newcommand{\eightVgraphI}{\raisebox{0pt}{
                 \begin{picture}(26,26)(-13,-9)
                 \put(-9,9){\circle{6}} \put(9,9){\circle{6}}
                 \put(-9,-9){\circle{6}} \put(9,-9){\circle{6}} 
	      	 \put(-6,9){\line(1,0){12}}
                 \put(-9,6){\line(0,-1){12}}
                 \put(-6,-9){\line(1,0){12}}
                 \put(9,6){\line(0,-1){12}}
                 \end{picture}}}

\newcommand{\eightVgraphII}{\raisebox{0pt}{
                 \begin{picture}(28,28)(-14,-10)
                 \put(-13,13){\line(1,0){26}}
                 \put(-13,-13){\line(1,0){26}}
                 \put(-13,-13){\line(0,1){26}}
                 \put(13,-13){\line(0,1){26}}
                 \put(-3,3){\line(1,0){6}}
                 \put(-3,-3){\line(1,0){6}}
                 \put(-3,-3){\line(0,1){6}}
                 \put(3,-3){\line(0,1){6}}
                 \put(-13,13){\line(1,-1){10}}
                 \put(-13,-13){\line(1,1){10}}
                 \put(13,13){\line(-1,-1){10}}
                 \put(13,-13){\line(-1,1){10}}
                 \end{picture}}}

\newcommand{\tenVgraphI}{\raisebox{0pt}{
                 \begin{picture}(26,29)(-13,-9)
                 \put(-9,9){\circle{6}} \put(9,9){\circle{6}}
                 \put(-9,-9){\circle{6}} \put(9,-9){\circle{6}}
                 \put(0,9){\circle{6}}
                 \put(-6,9){\line(1,0){3}}
                 \put(6,9){\line(-1,0){3}}
                 \put(-9,6){\line(0,-1){12}}
                 \put(-6,-9){\line(1,0){12}}
                 \put(9,6){\line(0,-1){12}}
                 \end{picture}}}

\newcommand{\tenVgraphII}{\raisebox{0pt}{
                 \begin{picture}(28,31)(-14,-10)
                 \put(-13,13){\line(1,0){10}}
                 \put(0,13){\circle{6}}
                 \put(13,13){\line(-1,0){10}}
                 \put(-13,-13){\line(1,0){26}}
                 \put(-13,-13){\line(0,1){26}}
                 \put(13,-13){\line(0,1){26}}
                 \put(-3,3){\line(1,0){6}}
                 \put(-3,-3){\line(1,0){6}}
                 \put(-3,-3){\line(0,1){6}}
                 \put(3,-3){\line(0,1){6}}
                 \put(-13,13){\line(1,-1){10}}
                 \put(-13,-13){\line(1,1){10}}
                 \put(13,13){\line(-1,-1){10}}
                 \put(13,-13){\line(-1,1){10}}
                 \end{picture}}}

\newsavebox{\DISK}
\savebox{\DISK}[8pt]{\begin{picture}(8,8)(0,0)
                     \put(-2.5,-3){$\bullet$}
                     \end{picture}}
\newcommand{\disk}{\usebox{\DISK}}

\newcommand{\knotcrossing}{\raisebox{0pt}{
                 \begin{picture}(26,20)(-13,-5)
                 \put(-12,-12){\vector(1,1){24}}
		 \put(12,-12){\vector(-1,1){24}}
                 \put(0,0){\disk}
                 \end{picture}}}
\newcommand{\overcrossing}{\raisebox{0pt}{
                 \begin{picture}(26,20)(-13,-5)
                 \put(-12,-12){\vector(1,1){24}}
		 \put(12,-12){\line(-1,1){10}}
		 \put(-2,2){\vector(-1,1){10}}
                 \end{picture}}}
\newcommand{\undercrossing}{\raisebox{0pt}{
                 \begin{picture}(26,20)(-13,-5)
                 \put(-12,-12){\line(1,1){10}}
		 \put(2,2){\vector(1,1){10}}
		 \put(12,-12){\vector(-1,1){24}}
                 \end{picture}}}

\newcommand{\wtwowfourI}{\raisebox{0pt}{
  		 \begin{picture}(60,30)(-30,-10)
		 \put(-20,0){\circle{20}}
                 \put(-10,10){\oval(20,20)[t]}
                 \put(-10,-10){\oval(20,20)[b]}
		 \put(0,10){\line(1,0){20}}
		 \put(0,-10){\line(1,0){20}}
		 \put(0,-10){\line(0,1){20}}
		 \put(20,-10){\line(0,1){20}}
		 \put(20,-10){\line(1,-1){10}}
                 \put(20,10){\line(1,1){10}}
                 \end{picture}}}
\newcommand{\wtwowfourII}{\raisebox{0pt}{
  		 \begin{picture}(40,34)(-20,-10)
		 \put(-13,16){\line(1,0){26}}
		 \put(-13,-10){\line(1,0){26}}
		 \put(-13,-10){\line(0,1){26}}
		 \put(13,-10){\line(0,1){26}}
		 \put(-13,-10){\line(-1,-1){10}}
                 \put(13,16){\line(1,1){10}}
		 \put(13,-10){\line(-1,1){10}}
		 \put(-13,16){\line(1,-1){10}}
		 \put(0,3){\circle{8}}
                 \end{picture}}}

\newcommand{\twoVchord}{\raisebox{0pt}{
  		 \begin{picture}(30,30)(-15,-10)
                 \put(-15,-15){\vector(1,0){30}}
		 \put(0,-15){\oval(16,16)[t]}
                 \end{picture}}}
\newcommand{\thetachord}{\raisebox{0pt}{
  		 \begin{picture}(30,30)(-15,-10)
                 \put(-15,-15){\vector(1,0){30}}
		 \put(0,0){\circle{16}} \put(-8,0){\line(1,0){16}}
                 \end{picture}}}
\newcommand{\fourVchord}{\raisebox{0pt}{
  		 \begin{picture}(60,30)(-30,-10)
                 \put(-30,-15){\vector(1,0){60}}
		 \put(-15,-15){\oval(16,16)[t]}
		 \put(15,-15){\oval(16,16)[t]}
                 \end{picture}}}
\newcommand{\thetatwoVchord}{\raisebox{0pt}{
  		 \begin{picture}(30,30)(-15,-10)
                 \put(-15,-15){\vector(1,0){30}}
		 \put(0,-15){\oval(16,16)[t]}
		 \put(0,5){\circle{16}} \put(-8,5){\line(1,0){16}}
                 \end{picture}}}
\newcommand{\thetasqchord}{\raisebox{0pt}{
  		 \begin{picture}(30,30)(-15,-10)
                 \put(-15,-15){\vector(1,0){30}}
		 \put(0,0){\circle{16}} \put(-8,0){\line(1,0){16}}
		 \put(7,7){2}
                 \end{picture}}}

\begin{document}

%\title{\bf Rozansky-Witten invariants of hyperk{\"a}hler manifolds}
%\author{{\bf Justin Sawon} \\ 
%	{\bf Trinity College} \\
%	\vspace*{100mm} \\
%	{\em Thesis submitted for the degree of} \\
%	{\em Doctor of Philosophy} \\
%        \vspace*{15mm} \\
%}
%\date{October, 1999}
%\maketitle
%\vspace*{-40mm}
%\begin{figure}[htpb]
%\epsfxsize=60mm
%\centerline{\epsfbox{Crest.eps}}
%\end{figure}
%\thispagestyle{empty}
%
%  Note that original Crest.eps file included the lines:
%
%%%BoundingBox: 150 384 358 428
%%%HiResBoundingBox: 150.0693 384.6326 357.9306 427.3674

\title{\bf Rozansky-Witten invariants of hyperk{\"a}hler manifolds}
\author{{\bf Justin Sawon} \\ 
	{\bf Trinity College} \\
	\vspace*{30mm} \\
	{\bf University of Cambridge} \\
	\vspace*{80mm} \\
	{\em Thesis submitted for the degree of} \\
	{\em Doctor of Philosophy}
}
\date{October, 1999}
\maketitle
\vspace*{-160mm}
\begin{center}
\includegraphics*[scale=1.6]{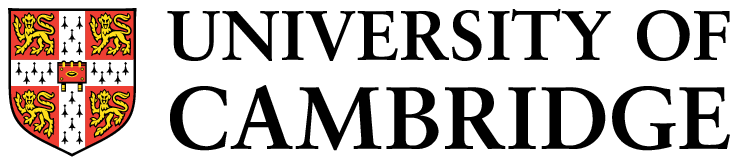}
\end{center}
\thispagestyle{empty}

\newpage
\pagenumbering{roman}
\thispagestyle{myheadings}
\markboth{}{}
\tableofcontents

\newpage
\thispagestyle{myheadings}
\markboth{}{}
\begin{abstract}
We investigate invariants of compact hyperk{\"a}hler manifolds
introduced by Rozansky and Witten: they associate an invariant to each
graph homology class. It is obtained by using the graph to perform
contractions on a power of the curvature tensor and then integrating
the resulting scalar-valued function over the manifold, arriving at a
number. For certain graph homology classes, the invariants we get are
Chern numbers, and in fact all characteristic numbers arise in this
way. 

We use relations in graph homology to study and compare these
hyperk{\"a}hler manifold invariants. For example, we show that the
norm of the Riemann curvature can be expressed in terms of the volume and
characteristic numbers of the hyperk{\"a}hler manifold. We also
investigate the question of whether the Rozansky-Witten invariants
give us something more general than characteristic numbers. Finally,
we introduce a generalization of these invariants which incorporates
holomorphic vector bundles into the construction.

\end{abstract}

\newpage
\addtocounter{page}{1}
\thispagestyle{myheadings}
\markboth{}{}
\vspace*{10mm}
\noindent
{\em The work contained in this dissertation is my own original
research carried out with the assistance of my PhD supervisor, Prof
N.\ Hitchin, and others mentioned in the
acknowledgements. Propositions 2, 5, 6, and Theorem 15 have
appeared in a joint preprint {\bf math.DG/9908114\/} with Prof Hitchin
entitled {\em ``Curvature and characteristic numbers of
hyperk{\"a}hler manifolds''\/}. Otherwise, to the best of my knowledge
and belief, this thesis contains no material previously published or
written by another person, except where due reference has been made in
the text.\/}

%\noindent
%{\em This work contains no material which has been accepted for the
%award of any other degree or diploma in any university or other
%tertiary institution and, to the best of my knowledge and belief,
%contains no material previously published or written by another
%person, except where due reference has been made in the text.\/} 

%\vspace*{5mm}
%\noindent
%{\em I give consent to this copy of my thesis, when deposited in the
%University Library, being available for loan and photocopying.\/} 

\vspace*{40mm}
\hfill{\em Justin Sawon\/} 

\hfill{\em October, 1999\/}

\newpage
\addtocounter{page}{1}
\thispagestyle{myheadings}
\markboth{}{}
\noindent
\begin{center}
{\bf Acknowledgements}
\end{center}

I would like to express my gratitude to my supervisor Professor Nigel
Hitchin for all the support and encouragement he has given me, and
especially for giving me an interesting and challenging topic to work
with, and then guiding me in all the right directions. His
mathematical insight and understanding of geometry has been nothing
short of inspirational.

I owe a great debt to Simon Willerton whose suggestion that the
Wheeling Theorem may play a part in this theory proved to be extremely
profitable, and whose helpful corrections of an early draft of my
thesis enlightened me on many of the subtler points of knot theory. I
am grateful to Stavros Garoufalidis who taught me a great deal about
graph homology in the early stages of this work, and to J{\o}rgen
Ellegaard Andersen for encouraging me to learn more about invariants
of knots and links, and for his and {\AA}rhus University's hospitality
during my visit to Denmark. My understanding of knots and
three-manifold invariants has also profited hugely from discussions
with Thang Le, Gregor Masbaum, Tomotada Ohtsuki, Dylan Thurston and
Pierre Vogel, all of whom I am indebted to.

I have benefited greatly from discussions with Manfred Lehn and Lothar
G{\"o}ttsche regarding the Chern numbers of Hilbert schemes. The data
they were able to provide me with was most helpful during my early
attempts to understand this theory. I am also grateful to Mikhail
Kapranov for discussing his work with me, and for sending me an old
version of his article on this topic which contained many useful
ideas. In the later stages of this work I received many interesting
correspondences from George Thompson, and I'm only sorry that I
haven't yet had the opportunity to explore all the suggestions he made
to me. Robbert Dijkgraaf, Nathan Habegger, and Lev Rozansky also
raised many interesting points in our conversations and I eagerly look
forward to exploring these speculations in the future.

Thanks are also due to Arnaud Beauville, Phil Boalch, Dorin Cheptea,
Gueo Grantcharov, Dominic Joyce, Maxim Kontsevich, Hiraku Nakajima,
Simon Salamon, Ivan Smith, Richard Thomas, Misha Verbitsky, and Pelham
Wilson, among others, for useful conversations, sharing their ideas,
and just generally taking an interest in my work. I would also like to
thank those people, especially Mike Eastwood, who taught me during my
years at the University of Adelaide, and without whom I would not have
had the education which led to the opportunity to study at the
University of Cambridge.

This research could not have been carried out without the considerable
financial generosity of Trinity College, who have amply supported me
throughout my studies in the UK. Support from the British government
in the form of an Overseas Research Student award is also gratefully
acknowledged, as is the hospitality of New College and the
Mathematical Institute during my lengthy visit to Oxford. Thanks also
to the Cambridge Philosophical Society and London Mathematical Society
for travel grants.

Lastly, but by no means least, thanks to my family and friends,
especially Anna, Christy, and my parents for encouragement and
unlimited support.

\newpage
\pagenumbering{arabic}
\addtocounter{section}{-1}
\include{intro}

\newpage
\include{ch1}

\addtocounter{page}{1}
\newpage
\include{ch2}

\addtocounter{page}{1}
\newpage
\include{ch3}

\newpage
\include{ch4}

\addtocounter{page}{1}
\newpage
\include{ch5}

\newpage
\include{ch6}

\addtocounter{page}{1}
\appendix
\include{app1}

\include{app2}

\include{app3}

\include{app4}

\include{app5}

\newpage
\thispagestyle{myheadings}
\markboth{}{}
%\nocite{*}
\bibliographystyle{amsplain}
\bibliography{refs}

\end{document}

%% file: intro.tex
\section{Introduction}

In this thesis we study new invariants of hyperk{\"a}hler
manifolds introduced by Rozansky and Witten in~\cite{rw97}. They
described a three-dimensional sigma model with target space a
hyperk{\"a}hler manifold, and showed that the partition function of
the theory is a three-manifold invariant of finite type. In a
``perturbative'' expansion the weights depend on the hyperk{\"a}hler
manifold and are indexed by trivalent graphs, and Rozansky and Witten give
a purely differential geometric construction of these weights. Roughly
speaking, we take a tensor power of the Riemann curvature tensor
of the hyperk{\"a}hler manifold; the trivalent graph tells us how to
contract indices (using the holomorphic symplectic form and its dual)
to get something we can integrate, resulting in a number.
$$
\begin{array}{ccccc}
\mbox{hyperk{\"a}hler manifold }X\hspace*{-10mm} & & & & \hspace*{-10mm}\mbox{trivalent graph
}\Gamma \\
 & \searrow & & \swarrow & \\
 & & \mbox{number }b_{\Gamma}(X) & &
\end{array}
$$

The weight $b_{\Gamma}(X)$ will be our main object of study. It
satisfies the following properties (Propositions $1$ to $4$
respectively in Chapter $1$):
\begin{itemize}
\item it is independent of the choice of complex structure on $X$,
\item it is invariant under deformations of the hyperk{\"a}hler metric
on $X$,
\item the dependence on $\Gamma$ is only through its graph homology
class,
\item taking products of hyperk{\"a}hler manifolds corresponds to
taking coproducts in graph homology.
\end{itemize}
Proposition $1$ follows from Rozansky and Witten's original
definition. While we use a complex geometric description of $X$, their
approach uses a description which is independent of choosing a complex
structure (which is compatible with the hyperk{\"a}hler metric) from
the start. Proposition $2$ was also known to Rozansky and
Witten, though we shall follow the approach of
Kapranov~\cite{kapranov99}, which was inspired by the ideas of 
Kontsevich~\cite{kontsevich99}. In this approach we pass from forms to
cohomology classes: the Riemann curvature tensor is replaced by the
Atiyah class~\cite{atiyah57}, which is represented by the curvature in
Dolbeault cohomology. We then work within the framework of complex
algebraic geometry. As a cohomology class the Atiyah class can be
described in many ways, most of which do not require us to choose a
hyperk{\"a}hler metric. This is useful for calculations as usually the
existence of a hyperk{\"a}hler metric is deduced from Yau's
theorem~\cite{yau78} and there is no explicit knowledge of what it
looks like. 
%In fact, Kapranov's approach allow us to define the
%weights $b_{\Gamma}(X)$ for any compact holomorphic symplectic
%manifold $X$.

Since the weights satisfy Proposition $2$ it makes sense to call them
{\em Rozansky-Witten invariants of hyperk{\"a}hler manifolds\/}, and
we adopt this terminology, though strictly speaking the
Rozansky-Witten invariants are the three-manifold invariants
constructed from these weights.

Graph homology is the space of linear combinations of oriented
trivalent graphs modulo the IHX and AS (anti-symmetry) relations. The
AS relations say that reversing the orientation of a graph gives us
minus the original graph. The IHX relations say that three graphs
which are identical except for small regions which look like $I$, $H$,
and $X$ respectively, are related by
\begin{center}
\begin{picture}(140,60)(-60,-30)
\put(-80,20){\line(0,-1){40}} 
\put(-80,20){\line(-2,1){20}} \put(-80,20){\line(2,1){20}}
\put(-80,-20){\line(-2,-1){20}} \put(-80,-20){\line(2,-1){20}}
\put(-40,0){$=$}
\put(-10,0){\line(1,0){20}}
\put(-10,0){\line(-1,3){10}} \put(-10,0){\line(-1,-3){10}}
\put(10,0){\line(1,3){10}} \put(10,0){\line(1,-3){10}}
\put(35,0){$-$}
%\put(40,0){and}
\put(60,30){\line(2,-3){40}} 
\put(60,-30){\line(2,3){19}} \put(100,30){\line(-2,-3){19}}
\put(70,-15){\line(1,0){20}}
\end{picture}
%Figure $1$: The IHX relation
\end{center}
Graph homology is graded by degree, which is given by half the number
of vertices of a graph.  Proposition $3$ says that the weights
$b_{\Gamma}(X)$ descend to graph homology.  In the ``perturbative''
expansion of the partition function it is the fact that the weights
satisfy the IHX relations that ensures we get a topological invariant
of the three-manifold. This does not rely at all on Proposition $2$;
indeed if the weights were not metric independent we would simply need
to choose a hyperk{\"a}hler manifold with a specified metric on it and
this would still lead to a three-manifold invariant.

On the other hand, the fact that the weights are hyperk{\"a}hler
manifold invariants does not depend at all on Proposition $3$. From this
point of view the IHX relations are really just a fancy way of writing
integration by parts. However, the power of this formalism should not
be under-estimated; whenever two graphs are homologous their
corresponding Rozansky-Witten invariants are identical.

It also follows from Proposition $3$ that hyperk{\"a}hler manifolds give
rise to elements of the dual of graph homology, ie.\ they can be
thought of as elements of graph cohomology. Conversely, we can regard
elements of graph cohomology as being {\em virtual\/} hyperk{\"a}hler
manifolds. We check that this idea makes sense in Subsection $2.3$,
before making use of it in Chapter $3$ by associating virtual
manifolds to the ${\frak{su}}(2)$ weight system arising in perturbative
${\mathrm SU}(2)$ Chern-Simons theory.

Many of the techniques we will use to evaluate the Rozansky-Witten
invariants will rely on the hyperk{\"a}hler manifold being
irreducible. Proposition $4$ allows us to extend these calculations to
reducible hyperk{\"a}hler manifolds.

Fundamental to most of our calculations is the fact that for certain
choices of trivalent graphs $\Gamma$ we get characteristic numbers of
$X$. More precisely, we take $\Gamma$ to be a {\em polywheel\/}, which
is defined to be the disjoint union of a collection of wheels
\begin{center}
\begin{picture}(150,70)(-140,15)
\put(-70,50){\circle{40}}
\put(-84,64){\line(-1,1){20}}
\put(-84,36){\line(-1,-1){20}}
\put(-56,36){\line(1,-1){20}}
\put(-56,64){\line(1,1){20}}
\put(-10,50){\circle{20}}
\put(-10,60){\line(0,1){25}}
\put(-10,40){\line(0,-1){25}}
\end{picture}
\end{center}
closed by summing over all ways of joining the spokes pairwise. Then
the corresponding Rozansky-Witten invariant $b_{\Gamma}(X)$ is a Chern
number (Proposition $5$ in Subsection $2.2$), and all characteristic
numbers arise in this way. Hence the Rozansky-Witten invariants can be
thought of as generalized characteristic numbers. The first question
which springs to mind is: how much more general are they, if at all?  The
corresponding question in graph homology would be: how many graphs are
there, if any, in the complement of the subspace spanned by
polywheels? In degree four and higher this complementary subspace is
always non-empty, and one might expect that the Rozansky-Witten
invariants corresponding to these graphs would be more general than
characteristic numbers, but this is by no means automatic. 

We will return to this question shortly. First let us mention that
it also pays to consider graphs which lie in the polywheel
subspace. The Rozansky-Witten invariants corresponding to these graphs
must necessarily be characteristic numbers, though this may not be at
all obvious from their structure. Let $\Theta$ be the unique trivalent
graph with two vertices (which looks just like theta); then the most
important example for us is the graph $\Theta^k$ consisting of $k$
disjoint copies of $\Theta$. Using the Wheeling Theorem~\cite{blt} we
can show that $\Theta^k$ lies in the polywheel subspace for all $k$
(Proposition $13$ in Subsection $4.2$). On the other hand, the
Rozansky-Witten invariant corresponding to $\Theta^k$ can be expressed
in terms of the ${\cal L}^2$-norm of the curvature and the volume of
our hyperk{\"a}hler manifold (Proposition $6$ in Subsection $2.4$). So
our theory says that this expression must be a characteristic number,
which is a surprising result. We can be more precise than this: it is
given by the multiplicative polynomial $\hat{A}^{1/2}$, where
$\hat{A}$ is the multiplicative polynomial corresponding to the
$\hat{A}$-genus (Theorem $15$ in Subsection $5.1$). Using the
positivity of the norm, we can also deduce some universal inequalities
for the characteristic numbers of hyperk{\"a}hler manifolds (Corollary
$16$ in Subsection $5.1$).

We can evaluate the Rozansky-Witten invariants for those graphs lying
in the polywheel subspace by writing then in terms of Chern numbers.
There are two main families of compact hyperk{\"a}hler manifolds, due
to Beauville~\cite{beauville83}: the Hilbert schemes of points on a
K$3$ surface and the generalized Kummer varieties. For these examples,
we begin with some limited information about the Chern numbers which
we obtain using the Riemann-Roch formula and the known values of the
Hirzebruch $\chi_y$-genus. Combining this with some direct reasoning
allows us to calculate the Rozansky-Witten invariants corresponding
to $\Theta^k$ for all $k$ (Propositions $19$ and $21$ in Chapter
$5$). This in turn helps us to evaluate more Chern numbers, and hence
more Rozansky-Witten invariants, leading to complete tables of
these values up to degree four, ie.\ real-dimension sixteen (see
Appendices D and E$.1$). We can also compute the values of the
Rozansky-Witten invariants in low degrees for products of these
manifolds, by virtue of Proposition $4$.

Armed with these calculations, we can now return to the fundamental
question posed above. To show that a given Rozansky-Witten invariant
is not a characteristic number all we need to do is exhibit two
hyperk{\"a}hler manifolds with the same Chern numbers which are
distinguished by the invariant. We can allow the more general
situation of disconnected manifolds since the invariants are additive
on disjoint sums. In degree four (real-dimension sixteen) we construct
two such manifolds by taking disjoint sums of the examples mentioned
above. This shows that the Rozansky-Witten invariant corresponding to
the disconnected graph
$$\fourVgraph\fourVgraph$$ 
is not a characteristic number (Theorem $22$ in Subsection $5.5$).
The interesting thing about this result is that all our calculations
were based upon knowledge of the characteristic numbers. In
particular, if the hyperk{\"a}hler manifold $X$ is reducible then the
invariant can be written in terms of the characteristic numbers of the
irreducible factors of $X$, and if $X$ is irreducible then it can be
expressed as a {\em rational\/} function of the characteristic numbers
of $X$. There exists a large family of graphs whose Rozansky-Witten
invariants have this property, ie.\ are given by a rational function
of characteristic numbers on irreducible manifolds (Theorem $17$ in
Subsection $5.1$). Indeed, together with the polywheels, these graphs
span graph homology up to (and including) degree five. On the other
hand, there are certainly higher degree graphs which are neither in
the polywheel subspace nor in this other family of graphs. The precise
nature of their Rozansky-Witten invariants, and how they are related
to the characteristic numbers (if at all) is still a mystery.

The Rozansky-Witten theory shares many of the properties of
Chern-Simons theory. In the ``perturbative'' models, the weights
$b_{\Gamma}(X)$ constructed from a hyperk{\"a}hler manifold $X$ are
analogous to those constructed from a semi-simple Lie algebra in
perturbative Chern-Simons theory. We can extend this analogy to
weights on chord diagrams $D$, which can roughly be described as
unitrivalent graphs whose univalent vertices lie on a collection of
oriented circles. The construction requires some additional
information, which in the Lie algebra case is given by attaching a
finite-dimensional representation to each oriented circle, and in the
Rozansky-Witten case by attaching to each circle a holomorphic vector
bundle $E_a$ over our hyperk{\"a}hler manifold. This gives us weights
$b_D(X;E_a)$ which in the simplest case (no vector bundles) reduce to
the Rozansky-Witten invariants $b_{\Gamma}(X)$.

The weights $b_D(X;E_a)$ satisfy many properties similar to those of
$b_{\Gamma}(X)$. For example, we can show (Propositions $23$, $24$,
and $25$ respectively in Chapter $6$):
\begin{itemize}
\item they do not depend on the choice of Hermitian structures (part
of the construction) on the vector bundles $E_a$,
\item the dependence on $D$ is only through its equivalence class in
the space of chord diagrams,
\item certain chord diagrams give rise to weights which are integrals
of the Chern classes of $E_a$.
\end{itemize}
Proposition $24$ implies that these weight systems lead to
finite-type invariants of knots and links, just as the weight system
$b_{\Gamma}(X)$ gives us a finite-type three-manifold
invariant. Proposition $25$ enables us to use the Wheeling Theorem to
derive a formula for the integral of the top component of the Mukai
vector of a vector bundle in terms of weights corresponding to simple
chord diagrams. This generalizes our earlier result which expressed
the Rozansky-Witten invariant corresponding to $\Theta^k$ as a
characteristic number given by the $\hat{A}^{1/2}$-polynomial.

We give here a rough guide to the contents of Chapters $1$ to $6$. In
the first chapter we give our definition of the Rozansky-Witten
invariant $b_{\Gamma}(X)$, as well as reviewing Rozansky and Witten's
original definition and Kapranov's cohomological interpretation (all
equivalent). These alternative descriptions are useful in proving
Propositions $1$ to $4$, which occupies us for the remainder of this
chapter.

In Chapter $2$ we prove the fundamental relation between Chern numbers
and Rozansky-Witten invariants corresponding to polywheels. We then
show that regarding elements of graph cohomology as virtual
hyperk{\"a}hler manifolds makes sense from the point of view of taking
products and its effect on characteristic numbers. We also establish
an expression for the invariant $b_{\Theta^k}(X)$ in terms of the norm
of the curvature and the volume of $X$.

Chapters $3$ and $4$ contain the results in graph homology which we
will later translate into results on hyperk{\"a}hler manifold by
taking the corresponding Rozansky-Witten invariants. The first of
these two chapters involves the use of the ${\frak{su}}(2)$ weight
system which arises in perturbative Chern-Simons theory with gauge
group ${\mathrm SU}(2)$; in the second we introduce some ideas and
results from the theory of finite-type invariants of knots. Together
with some direct fooling around with graph homology, these give us all
the relations we need.

Chapter $5$ begins by reinterpreting the graph homology
relations of the previous two chapters in terms of Rozansky-Witten
invariants. This gives us some general results which we combine with
some direct reasoning to calculate the Rozansky-Witten invariants for
some specific examples of compact hyperk{\"a}hler manifolds, namely
the Hilbert schemes of points on a K$3$ surface and the generalized
Kummer varieties.

In the final chapter we introduce the weight systems $b_D(X;E_a)$
constructed from collections of holomorphic vector bundles $E_a$ over
a hyperk{\"a}hler manifold, and prove that they satisfy the properties
outlined above.

Appendices A to E contain lists of graph homology relations and tables
of Chern numbers and Rozansky-Witten invariants, as well as tables of
other relations and data required in the calculations of Chapter $5$.

\begin{rmk}

We will make a few general comments here to save from repeating
them many times throughout the text. Firstly, in the topological
sigma-model introduced by Rozansky and Witten the hyperk{\"a}hler
manifold $X$ is taken to be either compact or asymptotically
flat. Certain non-compact examples are believed to give special
three-manifold invariants, in particular the generalized Casson-Walker
invariants. However, in the asymptotically flat case one still has to
be careful with the decay conditions on the curvature in order for the
integrals to converge. Instead we choose to avoid these difficulties
by working exclusively with compact hyperk{\"a}hler manifolds.

We wish to use the techniques of complex geometry but there is no
natural choice of complex structure on a hyperk{\"a}hler manifold,
which has an entire $S^2$ family of compatible complex structures.
So we pick one at random. We have already mentioned Proposition $1$,
which says that the Rozansky-Witten invariants are independent of this
choice. Characteristic classes also play an important role in this
theory. On a complex manifold we can define Chern classes, though in
the hyperk{\"a}hler case the odd Chern classes vanish and the even
ones can be related to the Pontryagin classes. Therefore the Chern
numbers are equivalent to the Pontryagin numbers. So although we will
tend to write expressions in terms of Chern numbers, one should
remember that these numbers are really topological invariants and
hence independent of the complex structure. For example, the Todd
genus is really the $\hat{A}$-genus, if we wished to write it in a
manifestly complex structure independent form.

Finally, after we introduce graph homology we abuse notation by simply
using the notation for a graph to denote the graph homology class it
represents. Indeed, we will even use the terminology ``graph'' when
referring to its class. Thus when we equate two different graphs we
mean that they are homologous, ie.\ they represent the same graph
homology class. We expect this should not cause any confusion.

\end{rmk}

%% file: ch1.tex
\section{Definitions and basic properties}

\subsection{Topological sigma-models}

In~\cite{rw97} Rozansky and Witten introduced a
$3$-dimensional topological sigma-model whose target space is a
(compact or asymptotically flat) hyperk{\"a}hler manifold $X$. The
partition function of this theory is a finite-type invariant of the
three-manifold $M$. A Feynman diagram calculation shows that it can be
written in the form 
$$\sum_{\Gamma}b_{\Gamma}(X)I_{\Gamma}^{\mathrm RW}(M),$$
where we sum over all trivalent graphs $\Gamma$, the {\em weights\/}
$b_{\Gamma}(X)$ depend on the hyperk{\"a}hler manifold, and the terms 
$I_{\Gamma}^{\mathrm RW}(M)$ depend on the three-manifold. This is not
really a perturbative expansion as for a given hyperk{\"a}hler
manifold $X$ of real-dimension $4k$ the weights $b_{\Gamma}(X)$ vanish
except when $\Gamma$ has $2k$ vertices, ie.\ we only get degree $k$
terms. There is evidence to suggest that $I_{\Gamma}^{\mathrm RW}(M)$
is the LMO invariant of Le, Murakami, and Ohtsuki~\cite{lmo98}. More
precisely, if we write the LMO invariant of our three-manifold $M$ as
$$Z^{\mathrm LMO}(M)=\sum_{\Gamma}I_{\Gamma}^{\mathrm LMO}(M)\Gamma$$
then it is believed that the $I_{\Gamma}^{\mathrm RW}(M)$ agree with the
coefficients $I_{\Gamma}^{\mathrm LMO}(M)$. This is shown by
Habegger and Thompson~\cite{ht99} at a ``physical level of rigour''
when the first Betti number of $M$ is greater than zero. The LMO
invariant takes values in the space of graph
homology, the space of linear combinations of oriented trivalent
graphs moduli the IHX relations and anti-symmetry. An element of graph
cohomology gives a weight system on graph homology, ie.\ a weight
system on trivalent graphs compatible with the IHX relations and
anti-symmetry. In general, composing the LMO invariant with a weight
system gives us a (scalar-valued) finite-type invariant. Although this
statement is true for all three-manifolds, only for rational homology
spheres can we say that all finite-type invariants arise in this way,
as the LMO invariant is only universal in this special case.

The most well-known weight systems are those arising from quadratic
Lie algebras and super-algebras, which includes all semi-simple Lie
algebras. Some of these weight systems arise naturally in perturbative
Chern-Simons theory, first studied by Axelrod and
Singer~\cite{as92,as94} (see also Bar-Natan~\cite{barnatan95I} and
Kontsevich~\cite{kontsevich94}). In fact, there is a close analogy
with the Rozansky-Witten theory. Expanding around the trivial
connection we get a finite-type invariant of three-manifolds which
looks like 
$$\sum_{\Gamma}c_{\Gamma}({\frak g})I_{\Gamma}^{\mathrm CS}(M),$$
where now the weight system $c_{\Gamma}({\frak g})$ depends on the
Lie algebra $\frak g$ of the gauge group, and the terms $I_{\Gamma}^{\mathrm
CS}(M)$ depend on the three-manifold. The latter have a configuration
space interpretation (see Bott and Taubes~\cite{bt94}), and when $M$
is a rational homology sphere they are believed to coincide with the
coefficients of the Aarhus integral of Bar-Natan, Garoufalidis,
Rozansky, and Thurston~\cite{bgrt97}. Thus we again recover the LMO
invariant, as the Aarhus integral is an alternative definition in this
case.

These perturbative expansions of the partition functions are perhaps
only true at a ``physical level of rigour''. We can also say that the
theory of the three-manifold terms $I_{\Gamma}(M)$ is not completely
understood. However, we can make precise mathematical sense of the
weight systems $b_{\Gamma}(X)$ and $c_{\Gamma}({\frak g})$. The main
object of our study will be the weights $b_{\Gamma}(X)$ arising from a
(compact) hyperk{\"a}hler manifold $X$, and we investigate their basic
properties in this chapter.

\subsection{Hyperk{\"a}hler geometry}

Let $X^{4k}$ be a compact hyperk{\"a}hler manifold of real-dimension
$4k$. This means that $X$ has holonomy contained in ${\mathrm
Sp}(k)$ (with holonomy equal to ${\mathrm Sp}(k)$ if and only if $X$
is irreducible). Then $X$ has the following structures:
\begin{itemize}
\item complex structures $I$, $J$, and $K$, acting like the
quaternions on the tangent bundle $T=TX$,
\item a metric $g$ which is K{\"a}hlerian with respect to each of $I$,
$J$, and $K$,
\item corresponding K{\"a}hler forms $\omega_1$, $\omega_2$, and
$\omega_3$, which are $d$-closed skew-symmetric two-forms.
\end{itemize}
More generally, if $a$, $b$, and $c$ are real numbers satisfying
$a^2+b^2+c^2=1$ then $aI+bJ+cK$ gives a complex structure on $X$ with
respect to which $g$ is K{\"a}hlerian, and hence we obtain a
two-sphere of compatible complex structures. 

Fix a specific complex structure $I$, and regard $X$ as a complex
manifold. Then $\omega=\omega_2+i\omega_3$ is a complex symplectic
form which is holomorphic with respect to $I$. Let $T$ and $T^*$ be
the complex tangent and cotangent bundles respectively. We can use
$\omega$ to identify them, and thus $T\cong T^*$. 

Take the Levi-Civita connection associated to the metric $g$. The
Riemann curvature tensor of this connection is a section
$$K\in\Omega^{1,1}(X,{\mathrm End}T)=\Omega^{0,1}(X,T^*\otimes{\mathrm
End}T)$$
with components $K^i_{\phantom{i}jk\bar{l}}$ relative to local complex
coordinates (with respect to the complex structure $I$), where the
$ij$ indices refer to ${\mathrm End}T$ and the $k\bar{l}$ 
indices refer to $\Omega^{1,1}$. Using $\omega$ to identify $T$ and
$T^*$, we get a section
$$\Phi\in\Omega^{0,1}(X,T^*\otimes T^*\otimes T^*).$$
defined by
$$\Phi_{ijk\bar{l}}=\sum_m\omega_{im}K^m_{\phantom{i}jk\bar{l}}.$$
Since $X$ is hyperk{\"a}hler, there is an ${\mathrm Sp}(2k,{\Bbb C})$
reduction of the frame bundle. Therefore the curvature takes values in
the Lie algebra of ${\mathrm Sp}(2k,{\Bbb C})$ consisting of matrices
of the form $A^k_j$ such that $S_{ij}=\sum_k\omega_{ik}A^k_j$ is
symmetric. In other words, $\Phi_{ijk\bar{l}}$ is symmetric in
$ij$. It is also symmetric in $jk$ as the Levi-Civita connection is
torsion-free and preserves the complex structure. Thus
$$\Phi\in\Omega^{0,1}(X,{\mathrm Sym}^3T^*).$$
This will be one of the main ingredients in the construction of the
weights $b_{\Gamma}(X)$. The others will be the holomorphic symplectic
form
$$\omega\in{\mathrm H}^0(X,{\Lambda}^2T^*)$$
and its dual
$$\tilde{\omega}\in{\mathrm H}^0(X,{\Lambda}^2T).$$
The latter is the skew form on $T^*$ dual to $\omega$. Note that if
$\omega$ is represented by the matrix $\omega_{ij}$ relative to local
complex coordinates, then in a dual basis the matrix $\omega^{ij}$ of
$\tilde{\omega}$ is minus the inverse of $\omega_{ij}$.

\subsection{The Rozansky-Witten invariants}

Let $\Gamma$ be a (possibly disconnected) oriented trivalent graph
with $2k$ vertices. The orientation is an equivalence class of
orientations on the edges and an ordering of the vertices; if the
orderings differ by a permutation $\pi$ and $n$ edges are oriented in
the reverse manner, then we regard these as equivalent if ${\mathrm
sgn}\pi=(-1)^n$. This definition differs from the usual notion of
orientation on trivalent graphs, which is an equivalence class of
cyclic orderings of the outgoing edges at each vertex, with two
orderings being equivalent if they differ on an even number of
vertices. However, we shall see later that these two definitions are
in fact equivalent.

Let the vertices of $\Gamma$ be $v_1,\ldots,v_{2k}$. Place a copy of
$\Phi$ at each vertex 
and label the outgoing edges with the holomorphic indices of
$\Phi$. For example, at $v_t$ we label the outgoing edges by $i_t$,
$j_t$, and $k_t$. We can do this in any order, as
$\Phi_{i_tj_tk_t\bar{l}_t}$ is symmetric in these indices. We then
contract along each edge using the dual holomorphic symplectic form,
$\tilde{\omega}$. Thus if an edge is labelled by $i_t$
at one end and $i_s$ at the other, then we contract with
$\omega^{i_ti_s}$ if the edge is oriented from $v_t$ to $v_s$, or with
$\omega^{i_si_t}$ if it has the opposite orientation. This gives us a
section of $(\bar{T}^*)^{\otimes 2k}$ with components
$$\Phi_{i_1j_1k_1\bar{l}_1}\cdots\Phi_{i_{2k}j_{2k}k_{2k}\bar{l}_{2k}}\omega^{i_1*}\omega^{j_1*}\omega^{k_1*}\cdots$$
where summation over repeated indices is assumed. Projecting this to
the exterior product gives us an element 
$$\Gamma(\Phi)\in\Omega^{0,2k}(X).$$
We multiply this form by the $k$th power of the holomorphic symplectic
form $\omega^k$ in $\Omega^{2k,0}(X)$. This gives us an element in
$\Omega^{2k,2k}(X)$ which we can integrate over the manifold.
\begin{dfn}
The Rozansky-Witten invariant of the hyperk{\"a}hler manifold $X$
corresponding to the trivalent graph $\Gamma$ is
$$b_{\Gamma}(X)=\frac{1}{(8\pi^2)^kk!}\int_X\Gamma(\Phi)\omega^k.$$
\end{dfn}
Strictly speaking this is the weight corresponding to $\Gamma$
occurring in the ``perturbative'' expansion of the Rozansky-Witten
invariant of three-manifolds, but we shall adopt the above
terminology. Roughly speaking, $b_{\Gamma}(X)$ is given by taking a
product of curvature tensors, contracting according to the particular
trivalent graph $\Gamma$, and then integrating over the manifold. We
have chosen the additional factor for several reasons which will
become evident in due course. For the time being, let us just state
that dividing by $\pi^{2k}$ makes the invariants integral for the
examples of hyperk{\"a}hler manifolds we will study, the $k!$ factor
ensures that the invariants behave nicely with respect to products of
manifolds, and the $8^k$ factor allows us to make simple comparisons
with characteristic numbers. The overall normalization seems to give
integral invariants on compact hyperk{\"a}hler manifolds; at least
this is true for the examples we will look at. In general, the
Rozansky-Witten invariants of three-manifolds satisfy some integrality
properties, but this still leaves some freedom in choosing the
normalization of the weights. Apart from being the `right one' for
the various reasons outlined above, our choice also agrees with the
one made by Rozansky and Witten~\cite{rw97}.

In the rest of this chapter we shall look at some of the basic
properties of $b_{\Gamma}(X)$; in particular, we shall prove the
following results (which were at least known to Rozansky and Witten,
even if they were not explicitly proved in~\cite{rw97}):
\begin{itemize}
\item it is independent of the choice of complex structure on $X$
(ie.\ there is nothing special about using $I$ in our definition),
\item it is constant on connected components of the moduli space of
hyperk{\"a}hler metrics on $X$ (we shall follow Kapranov's
approach~\cite{kapranov99} in order to prove this),
\item the dependence on the graph $\Gamma$ is only through its graph
homology class,
\item products of hyperk{\"a}hler manifolds correspond to coproducts
in graph homology (which has the structure of a commutative
cocommutative Hopf algebra).
\end{itemize}

\subsection{The physicist's definition}

The original `physics' definition of $b_{\Gamma}(X)$ as given
in~\cite{rw97} is slightly different to the one given above. Indeed,
our definition is the same as that given in the Appendix
of~\cite{rw97} for complex symplectic manifolds $X$. These generalize
hyperk{\"a}hler manifolds, and of course the definition reduces to
one equivalent to the original definition in the case that $X$ is
hyperk{\"a}hler. From Rozansky and Witten's original definition,
however, it will be clear that there is no dependence on the complex
structure, so we review it here.

The holonomy of the Levi-Civita connection of an arbitrary Riemannian
metric on $X$ lies in ${\mathrm SO}(4k)$, but for a hyperk{\"a}hler
metric it lies in an ${\mathrm Sp}(k)$ subgroup. Let $TX$ be the
(real) tangent bundle of $X$. Then it is known (see Salamon's
book~\cite{salamon89}) that the complexified tangent bundle
of $X$ decomposes into the tensor product of a rank $2k$ complex
vector bundle $V$ with structure group ${\mathrm Sp}(k)$ and a
trivializable rank $2$ complex vector bundle $S$ with structure group
${\mathrm Sp}(1)$,
$$TX\otimes_{\Bbb R}{\Bbb C}=V\otimes S$$
Using indices $I$, $J$, $\ldots\in\{1,\ldots,2k\}$ on $V$ and $A$,
$B$, $\ldots\in\{1,2\}$ on $S$, we note that both bundles possess
non-degenerate skew-symmetric two-forms
$\epsilon_{IJ}$ and $\epsilon_{AB}$ with inverses $\epsilon^{IJ}$ and
$\epsilon^{AB}$ respectively. If $i$, $j$, $\ldots$ denote indices on
$TX$ coming from local real coordinates on $X$, then we use the
covariantly constant tensors $\gamma_{IA}^i$ and $\gamma_i^{IA}$ to
describe the maps
$$\gamma_{IA}^i:TX\otimes_{\Bbb R}{\Bbb C}\rightarrow V\otimes S$$
and 
$$\gamma_i^{IA}:V\otimes S\rightarrow TX\otimes_{\Bbb R}{\Bbb C}.$$
The Levi-Civita connection on $TX\otimes_{\Bbb R}{\Bbb C}$ reduces to
an ${\mathrm Sp}(k)$ connection on $V$ tensored with the trivial
connection on $S$. The Riemann curvature tensor can be written
$$R_{ijkl}=-\gamma_i^{IA}\gamma_j^{JB}\gamma_k^{CK}\gamma_l^{DL}\epsilon_{AB}\epsilon_{CD}\Omega_{IJKL}$$
and after some manipulations involving the symmetries of $R_{ijkl}$ we
find that $\Omega_{IJKL}$ is completely symmetric (see Proposition
$9.3$ in Salamon~\cite{salamon89}). Rozansky and
Witten~\cite{rw97} essentially define $b_{\Gamma}(X)$ using the bundle
$V$ instead of $TX$. Specifically, they use $\Omega_{IJKL}$ and
$\epsilon^{IJ}$ instead of $\Phi_{ijk\bar{l}}$ and $\omega^{ij}$.

A choice of complex structure on $X$ compatible with the
hyperk{\"a}hler structure corresponds to a choice of trivialization of
the bundle $S$. This allows us to identify $TX$ (which is now a
complex vector bundle with the chosen complex structure) and $V$, and
shows that our approach in Subsection $1.3$ gives the same result as
in Rozansky and Witten's original definition. Furthermore, in the
latter we do not assume any particular trivialization of $S$, which
means we do not need to specify a complex structure, and hence
\begin{prp}
The invariant $b_{\Gamma}(X)$ is independent of the choice of
compatible complex structure on $X$.
\end{prp}
From Rozansky and Witten's original definition we can also show that
$b_{\Gamma}(X)$ is a real number. This requires us to make use of the
quaternionic structure on the bundle $V$.

\subsection{Kapranov's approach}

A third (equivalent) definition due to Kapranov~\cite{kapranov99}
enables us to prove metric independence. The idea behind his
construction is to work with cohomology classes instead of
differential forms, and leads to a description which does not rely on
specific knowledge of the hyperk{\"a}hler metric. In most examples,
the existence of a hyperk{\"a}hler metric follows from Yau's
theorem~\cite{yau78}, but no explicit description is known. Kapranov's
approach enables us to avoid this obstacle and leads to a method of
calculation on specific examples of hyperk{\"a}hler manifolds. 

The form 
$$\Phi\in{\Omega}^{0,1}(X,{\mathrm Sym}^3T^*)$$
is $\bar{\partial}$-closed by the Bianchi identity for the Riemannian
curvature. Therefore it represents a Dolbeault cohomology class
$$[\Phi]\in{\mathrm H}^{0,1}_{\bar{\partial}}(X,{\mathrm Sym}^3T^*).$$
This is actually the Atiyah class $\alpha_T$ of the tangent
bundle of $X$ (see~\cite{atiyah57}). We shall say more about
this shortly, but first we complete the construction. As before we
construct
$$[\Gamma(\Phi)]\in{\mathrm H}^{0,2k}_{\bar{\partial}}(X)$$
which only depends on the cohomology class of $\Phi$. Now
${\mathrm H}^{0,2k}_{\bar{\partial}}(X)$ is one-dimensional and
generated by $[{\bar{\omega}}^k]$. Multiplying by $[\omega^k]$ as
before gives an element of ${\mathrm H}^{2k,2k}_{\bar{\partial}}(X)$
which we can integrate to get
$$b_{\Gamma}(X)=\frac{1}{(8\pi^2)^kk!}\int_X[\Gamma(\Phi)][\omega^k].$$
Here, and throughout, the integral of a top degree Dolbeault
cohomology class over the manifold means that we integrate some form
which represents the class. This definition of $b_{\Gamma}(X)$ is
clearly no different to before, since we can use the same forms as
before to represent the Dolbeault cohomology classes. However, since
$\Phi$ represents the Atiyah class we can use other representations of
$\alpha_T$ in order to get a quite different looking description.

In general, the Atiyah class of a holomorphic vector bundle $E$ over
$X$ is the obstruction to the existence of a global
holomorphic connection on $E$ (see Atiyah~\cite{atiyah57}). Suppose
that we have an open set $U\in X$ over which $E$ has a holomorphic
connection $\nabla_U$. We can think of $\nabla_U$ as a map from $E$ to
the first jet bundle $J_1(E)$ which gives the identity when composed
with the projection $J_1(E)\rightarrow E$. In other words, it is a
splitting of the short exact sequence
$$0\rightarrow E\otimes T^*\rightarrow J_1(E)\rightarrow E\rightarrow 0.$$
Tensoring with $E^*$ we get
$$0\rightarrow {\mathrm Hom}(E,E\otimes T^*)\rightarrow {\mathrm
Hom}(E,J_1(E))\rightarrow {\mathrm Hom}(E,E)\rightarrow 0$$
and $\nabla_U$ is a section of ${\mathrm Hom}(E,J_1(E))$ over $U$ which
maps to ${\mathrm Id}$ in ${\mathrm Hom}(E,E)$. Consider the
corresponding long exact sequence
$$0\rightarrow {\mathrm H}^0(X,{\mathrm Hom}(E,E\otimes
T^*))\rightarrow {\mathrm H}^0(X,{\mathrm
Hom}(E,J_1(E)))\hspace*{20mm}$$
$$\hspace*{25mm}\rightarrow {\mathrm H}^0(X,{\mathrm
Hom}(E,E))\rightarrow {\mathrm H}^1(X,{\mathrm Hom}(E,E\otimes
T^*))\rightarrow\ldots$$
A global connection $\nabla$ on $E$ is a section in ${\mathrm
H}^0(X,{\mathrm Hom}(E,J_1(E)))$ which maps to ${\mathrm Id}$ in
${\mathrm H}^0(X,{\mathrm Hom}(E,E))$. The image of ${\mathrm Id}$
$$\alpha_E\in{\mathrm H}^1(X,{\mathrm Hom}(E,E\otimes T^*))={\mathrm
H}^1(X,T^*\otimes{\mathrm End}E)$$
is the obstruction to the existence of such a $\nabla$. Clearly if
$\alpha_E\neq 0$ then $\nabla$ cannot exist by exactness, whereas if
$\alpha_E$ vanishes, then ${\mathrm Id}$ is the image of some
$\nabla$. We call $\alpha_E$ the Atiyah class of $E$, and it can be
described in the following ways:
\begin{itemize}
\item Take an open cover $\{U_i\}$ of $X$ and choose local holomorphic
connections $\nabla_i$ on $E|_{U_i}$ over each open set $U_i$ which
look like $d+A_i$ in local holomorphic coordinates, where
$A_i\in\Omega^{1,0}(U_i,{\mathrm End}E)$. Then
$$\nabla_i-\nabla_j=A_i-A_j\in\Omega^{1,0}(U_i\cap U_j,{\mathrm
End}E)={\mathrm H}^0(U_i\cap U_j,T^*\otimes {\mathrm End}E)$$
gives a {\v C}ech representative of $\alpha_E$.
\item Take a smooth global connection $\nabla$ of type $(1,0)$ on
$E$. Then the $(1,1)$ part of the curvature of $\nabla$ gives a
Dolbeault representative of $\alpha_E$. If $\nabla$ looks like $d+A$
over the open set $U$, with $A\in\Omega^{1,0}(U,{\mathrm End}E)$,
then locally the $(1,1)$ part of the curvature of $\nabla$ will be
$\bar{\partial}A$. Note that if $E$ is equipped with an Hermitian
metric $h$ then the curvature of the compatible $h$-connection is
automatically of type $(1,1)$.

In order to construct a smooth global connection of type $(1,0)$, we
could patch together local holomorphic connections $\nabla_i$ (as
above) using a partition of unity $\{\psi_i\}$. The terms
$\psi_i\nabla_i$ are only differential operators, but when we add them
we get
$$\nabla=\sum\psi_i\nabla_i$$
which since $\sum\psi_i\equiv 1$ is a connection (as can be seen from
looking at its symbol). Relative to a local trivialization $\nabla$
looks like $d+A$ with 
$$A=\sum\psi_i A_i.$$
Since the $A_i$ are holomorphic, the $(1,1)$ part of the curvature is  
$$\sum(\bar{\partial}\psi_i)A_i.$$ 
\item Under the right conditions, we can also represent $\alpha_E$ as
the residue of a meromorphic connection $\nabla$ on $E$. We require
$\nabla$ to have a simple pole along a smooth (not necessarily
connected) divisor
$D$. Let $L={\cal O}(D)$ be the line bundle corresponding to $D$ and
$s$ the canonical section, which vanishes along $D$. A meromorphic
section of a bundle with a simple pole along $D$ is the same as a
holomorphic section of the same bundle twisted by $L$; the
corresponding map is given by multiplying by $s$. For example, the
smooth differential operator $s\nabla$ on $X$ is a (non-singular)
global holomorphic section of ${\mathrm Hom}(E,J_1(E))\otimes L$. 

In general, if we have a bundle $F$ on $X$, then we get the short
exact sequence
$$F\stackrel{s}{\rightarrow}F\otimes L\rightarrow F\otimes
L|_D\rightarrow 0$$
where the first map is given by multiplying by the section $s$ and
the second map by restricting to the divisor $D$ (the zero set of
$s$). The corresponding long exact sequence is
$${\mathrm H}^0(X,F)\stackrel{s}{\rightarrow}{\mathrm H}^0(X,F\otimes
L)\rightarrow {\mathrm H}^0(D,F\otimes
L)\stackrel{\delta}{\rightarrow}{\mathrm H}^1(X,F)\rightarrow\ldots$$
The second group in this sequence is the space of meromorphic sections
of $F$ with simple poles along $D$, and restricting to $D$ gives us
the residue of the section. Let $F_1$, $F_2$, and $F_3$ be the bundles
${\mathrm Hom}(E,E\otimes T^*)$, ${\mathrm Hom}(E,J_1(E))$, and
${\mathrm Hom}(E,E)$ respectively. Then we get the following
commutative double long exact sequence.
$$\begin{array}{cccccccc}
{\mathrm H}^0(X,F_1) & \hspace*{-2mm}\rightarrow\hspace*{-2mm} &
{\mathrm H}^0(X,F_1\otimes L) &
\hspace*{-2mm}\rightarrow\hspace*{-2mm} & {\mathrm H}^0(D,F_1\otimes
L) & \hspace*{-2mm}\rightarrow\hspace*{-2mm} & {\mathrm H}^1(X,F_1) &
\hspace*{-2mm}\rightarrow\ldots \\
\downarrow & & \downarrow & & \downarrow & & \downarrow & \\
{\mathrm H}^0(X,F_2) & \hspace*{-2mm}\rightarrow\hspace*{-2mm} &
{\mathrm H}^0(X,F_2\otimes L) &
\hspace*{-2mm}\rightarrow\hspace*{-2mm} & {\mathrm H}^0(D,F_2\otimes
L) & \hspace*{-2mm}\rightarrow\hspace*{-2mm} & {\mathrm H}^1(X,F_2) &
\hspace*{-2mm}\rightarrow\ldots \\
\downarrow & & \downarrow & & \downarrow & & \downarrow & \\
{\mathrm H}^0(X,F_3) & \hspace*{-2mm}\rightarrow\hspace*{-2mm} &
{\mathrm H}^0(X,F_3\otimes L) &
\hspace*{-2mm}\rightarrow\hspace*{-2mm} & {\mathrm H}^0(D,F_3\otimes
L) & \hspace*{-2mm}\rightarrow\hspace*{-2mm} & {\mathrm H}^1(X,F_3) &
\hspace*{-2mm}\rightarrow\ldots \\
\downarrow & & \downarrow & & \downarrow & & \downarrow & \\
{\mathrm H}^1(X,F_1) & \hspace*{-2mm}\rightarrow\hspace*{-2mm} &
{\mathrm H}^1(X,F_1\otimes L) &
\hspace*{-2mm}\rightarrow\hspace*{-2mm} & {\mathrm H}^1(D,F_1\otimes
L) & \hspace*{-2mm}\rightarrow\hspace*{-2mm} & {\mathrm H}^2(X,F_1) &
\hspace*{-2mm}\rightarrow\ldots \\
\downarrow & & \downarrow & & \downarrow & & \downarrow & \\
\vdots & & \vdots & & \vdots & & \vdots & \\
\end{array}$$
Now $s\nabla$ in ${\mathrm H}^0(X,F_2\otimes L)$ maps down to
$s{\mathrm Id}$ in ${\mathrm H}^0(X,F_3\otimes L)$, which in turn is
the horizontal image of ${\mathrm Id}$ in ${\mathrm H}^0(X,F_3)$,
which maps down to the Atiyah class $\alpha_E$ in ${\mathrm
H}^1(X,F_1)$. In the other direction, since $s{\mathrm Id}$ must map
to zero horizontally, it follows that the horizontal image of
$s\nabla$, namely $s{\nabla}|_D$ in ${\mathrm H}^0(D,F_2\otimes L)$,
must map down to zero. Then by exactness, $s{\nabla}|_D$ must be the
vertical image of some element $\beta_E$ in ${\mathrm
H}^0(D,F_1\otimes L)$. Finally, the horizontal image of $\beta_E$,
which lies in ${\mathrm H}^1(X,F_1)$, must again be the Atiyah class
$\alpha_E$. In summary, we have the following pattern of maps and
elements.
$$\begin{array}{ccccccccc}
 & & & & \beta_E & \stackrel{\delta}{\rightarrow} & {\alpha}_E &
\rightarrow & 0 \\
 & & & & \downarrow & & \downarrow & & \\
 & & s\nabla & \rightarrow & s{\nabla}|_D & \rightarrow & 0 & & \\
 & & \downarrow & & \downarrow & & & & \\
{\mathrm Id} & \stackrel{s}{\rightarrow} & s{\mathrm Id} & \rightarrow
 & 0 &&&& \\
\downarrow & & \downarrow &&&&&& \\
{\alpha}_E & \rightarrow & 0 &&&&&& \\
\downarrow &&&&&&&& \\
0 &&&&&&&& \\
\end{array}$$
The point here is that $\beta_E$ in ${\mathrm
H}^0(D,T^*\otimes{\mathrm End}E\otimes L)$ is the residue of the
meromorphic connection $\nabla$, and $\alpha_E$ is completely
determined by this section over the divisor $D$, as
$\delta\beta_E=\alpha_E$.
\end{itemize}

Using any of these approaches on the tangent bundle $T$ gives us  
$$\alpha_T\in{\mathrm H}^1(X,T^*\otimes{\mathrm End}T).$$
We can use the holomorphic symplectic form to identify $T$ and $T^*$,
so that
$$\alpha_T\in{\mathrm H}^1(X,T^*\otimes T^*\otimes T^*).$$
In fact we get
$$\alpha_T\in{\mathrm H}^1(X,{\mathrm Sym}^3T^*)$$
because as before, we have an ${\mathrm Sp}(2k,{\Bbb C})$ reduction of
the frame bundle and the (local holomorphic or smooth global)
connections can be chosen to be torsion-free and complex structure
preserving. Taking the product of $2k$ copies of $\alpha_T$, one for
each vertex of the trivalent graph $\Gamma$, and using the graph to
contract indices with $\tilde{\omega}$, as before, we get  
$$\Gamma(\alpha_T)\in{\mathrm H}^{2k}(X,{\cal O}_X).$$
This cohomology space is one-dimensional, but we need to be careful as
there is no natural basis. Instead, to produce a number in a canonical
way we can pair $\Gamma(\alpha_T)$ with 
$$\omega^k\in{\mathrm H}^0(X,\Lambda^{2k}T^*)$$
using Serre duality. Up to the factor $\frac{1}{(8\pi^2)^kk!}$ this is
precisely $b_{\Gamma}(X)$.

Using the meromorphic connection approach, we can also write down a
residue formula for $b_{\Gamma}(X)$. Firstly let us assume that $D$ is
connected. Since $\beta_T$ maps to $\alpha_T$, it can be taken to have
the same symmetries, ie.\ 
$$\beta_T\in{\mathrm H}^0(D,{\mathrm Sym}^3T^*\otimes L).$$
Now instead of constructing $\Gamma(\alpha_T)$ as before, we can
replace one of the copies of $\alpha_T$ with $\beta_T$ and restrict
the entire construction to $D$. This gives us
$$\Gamma(\alpha_T,\beta_T)\in{\mathrm H}^{2k-1}(D,L)$$
and the map
$$\delta :{\mathrm H}^{2k-1}(D,L)\rightarrow{\mathrm H}^{2k}(X,{\cal
O}_X)$$
will take $\Gamma(\alpha_T,\beta_T)$ to $\Gamma(\alpha_T)$. By the
adjunction formula the canonical line bundle of $D$ is
\begin{eqnarray*}
{\cal K}_D & = & {\cal K}_X\otimes{\cal O}_X(D)|_D \\
   & = & {\cal K}_X\otimes L|_D
\end{eqnarray*}
where ${\cal K}_X$ is the canonical line bundle of $X$, and therefore
$$\Gamma(\alpha_T,\beta_T)\in{\mathrm H}^{2k-1}(D,{\cal K}_D\otimes
{\cal K}_X^*|_D).$$
Now as $X$ is a Calabi-Yau manifold, the bundle ${\cal
K}_X=\Lambda^{2k}T^*$ must be trivial, indeed $\omega^k$ is a
non-vanishing section. Therefore $\Gamma(\alpha_T,\beta_T)$ lies in
a one-dimensional cohomology group
$${\mathrm H}^{2k-1}(D,{\cal K}_D\otimes {\cal K}_X^*|_D)\cong{\mathrm
H}^{2k-1}(D,{\cal K}_D)$$
but again we must be careful as there is no natural isomorphism
between these two spaces. To obtain a number in a canonical way we
must pair $\Gamma(\alpha_T,\beta_T)$ with
$$\omega^k|_D\in{\mathrm H}^0(D,{\cal K}_X|_D)$$
using Serre duality. We claim that this gives the same result as the
Serre duality pairing of $\Gamma(\alpha_T)$ with $\omega^k$, up to a
factor of $2\pi i$. This is essentially just Cauchy's residue formula,
but we shall be a little more precise. 

Consider the Poincar{\'e} residue map, which takes sections of ${\cal
K}_X\otimes L$ to sections of ${\cal K}_D$. Suppose that $D$ is given in
local coordinates $z_1,\ldots,z_{2k}$ by $f(z)=0$. Write a section of
${\cal K}_X\otimes L$ (ie.\ a meromorphic section of ${\cal K}_X$) as
$$\nu=\frac{g(z)dz_1\wedge\ldots\wedge
dz_{2k}}{f(z)}=\frac{df}{f}\wedge\nu^{\prime}$$
then the map to sections of ${\cal K}_D$ takes $\nu$ to
$\nu^{\prime}$. By the adjunction formula, this map gives an
isomorphism for sections over $D$. Now consider
$$\Gamma(\alpha_T,\beta_T)\omega^k|_D\in{\mathrm H}^{2k-1}(D,{\cal
K}_X\otimes L|_D).$$
Using {\v C}ech cohomology, observe that this element can be described
as sections of ${\cal K}_X\otimes L$ on $2k$-fold intersections of
open sets in $D$. The Poincar{\'e} residue map takes these sections
isomorphically to sections of ${\cal K}_D$, and hence
$$\Gamma(\alpha_T,\beta_T)\omega^k|_D\in{\mathrm H}^{2k-1}(D,{\cal K}_D).$$
Under the map $\delta$, this maps to
$$\Gamma(\alpha_T)\omega^k\in{\mathrm H}^{2k}(X,{\cal K}_X)$$
and both these cohomology spaces are one-dimensional and can be
identified with ${\Bbb C}$ simply by integrating (either over contours
in the case of {\v C}ech cohomology or over the entire space if we use
Dolbeault cohomology). The difference between the two numbers we get
from the above elements is precisely a factor of $2\pi i$, as we have
`cancelled' off $df/f$ and the contour integral over this additional
variable would give
$$\int_{|f|=1}\frac{df}{f}=2\pi i$$
by Cauchy's residue formula. Lastly, in the case that $D$ is
disconnected we simply need to perform the above calculation on each of
the connected components and sum the results. Thus we have another
description of the Rozansky-Witten invariant 
$$b_{\Gamma}(X)=\frac{2\pi
i}{(8\pi^2)^kk!}\int_D\Gamma(\alpha_T,\beta_T)\omega^k|_D.$$
According to the description above, if locally we write
$$\frac{\omega^k}{f}|_D=\frac{df}{f}\wedge \nu^{\prime}$$
then the section $\omega^k|_D$ of ${\cal K}_X|_D$ should be replaced
by the section $f\nu^{\prime}$ of ${\cal K}_D\otimes L^*|_D$ in the
above integral.

The beauty of these approaches is that we no longer need to know the
hyperk{\"a}hler metric explicitly in order to perform the
construction. In fact, all we need is a compact complex manifold $X$
with a holomorphic symplectic form $\omega$, usually referred to as a
complex symplectic space. In the K{\"a}hler case $X$ will admit a
hyperk{\"a}hler metric, but there are non-K{\"a}hler examples as
well, such as Kodaira surfaces~\cite{kodaira64,kodaira66} (see also
Barth, Peters, and Van de Ven~\cite{bpv84}) and their Douady
spaces (as observed by Beauville~\cite{beauville83}). There are also
simply-connected examples due to Guan~\cite{guan95} (see
Bogomolov~\cite{bogomolov96i} for a clearer description of this
construction).

Returning to the hyperk{\"a}hler case, we can use Kapranov's approach
to prove the following result. 
\begin{prp}
The number $b_{\Gamma}(X)$ is invariant under a deformation of the
hyperk{\"a}hler metric. In particular, it is constant on connected
components of the moduli space of hyperk{\"a}hler metrics.
\end{prp}
\begin{prf}
To begin with, consider what happens when we rescale the metric $g$ by
a factor $\lambda$. The holomorphic symplectic form $\omega$ will also
rescale by $\lambda$, and its dual $\tilde{\omega}$ by
$\lambda^{-1}$. The Riemann curvature tensor $K$ does not change,
but since $\Phi$ is defined by contracting with $\omega$ it will
rescale by $\lambda$. Since $\Gamma(\Phi)$ is made from $2k$ copies of
$\Phi$ (one for each vertex) and $3k$ copies of $\tilde{\omega}$ (one for
each edge) it will rescale by $\lambda^{-k}$. Finally, we multiply
this by $\omega^k$, which rescales by $\lambda^k$ and hence the net
result is that $b_{\Gamma}(X)$ is invariant under such rescalings of
the metric. This corresponds to a volume-changing deformation of the
metric on $X$. Now consider deformations in the transverse direction.

The hyperk{\"a}hler metric is determined by the three closed
skew-forms $\omega_1$, $\omega_2$, and $\omega_3$
(see Hitchin, Karlhede, Lindstr{\"o}m, and Ro{\v
c}ek~\cite{hklr87}). Thus an arbitrary first order deformation of the  
metric will be given by a linear combination of the three variations
$\delta\omega_1$, $\delta\omega_2$, and $\delta\omega_3$, subject to
some additional algebraic conditions. Fixing the holomorphic
symplectic form $\omega=\omega_2+i\omega_3$ (with respect to $I$), we
see from Kapranov's approach that $b_{\Gamma}(X)$ doesn't depend on
the K{\"a}hler form $\omega_1$, and hence is invariant under the
deformation $\delta\omega_1$ (where it is implied that $\omega_2$ and
$\omega_3$ are kept fixed as we vary $\omega_1$). Of course, there is
nothing special about the complex structure $I$ (as we saw earlier),
and thus $b_{\Gamma}(X)$ must also be invariant under 
the deformations $\delta\omega_2$ and $\delta\omega_3$, and hence
under linear combinations of these deformations. 

Note that by fixing two of $\omega_1$, $\omega_2$, and $\omega_3$ and
varying the third we do not get the deformation described by rescaling
$g$, but we have already seen above that $b_{\Gamma}(X)$ is invariant
under this volume-changing deformation, and hence $b_{\Gamma}(X)$ is
invariant under arbitrary first-order deformations of the
hyperk{\"a}hler metric. Finally, it follows from the results of
Tian~\cite{tian87} and Todorov~\cite{todorov89} that the obstruction
to deforming a Calabi-Yau manifold vanishes, and hence the moduli
spaces of Calabi-Yau manifolds are smooth. In particular, the moduli
space of hyperk{\"a}hler metrics on $X$ is smooth. Indeed it was shown
earlier by Bogomolov~\cite{bogomolov78} that the Kuranishi family of
moduli of a compact hyperk{\"a}hler manifold is smooth at the base
point. If two hyperk{\"a}hler metrics can be joined by a path, then
integrating the variation of $b_{\Gamma}(X)$ (which we have shown to
be zero) along this path proves that $b_{\Gamma}(X)$ is constant on
connected components of the moduli space. 
\end{prf}

\subsection{Graph homology}

We defined an orientation of $\Gamma$ as being an equivalence class of
orientations on the edges and an ordering of the vertices; if the
orderings differ by a permutation $\pi$ and $n$ edges are oriented
in the reverse manner, then we regard these as equivalent if ${\mathrm
sgn}\pi=(-1)^n$. When associating $\tilde{\omega}$ to an edge, the
orientation of the edge tells us whether to use $\omega^{i_ti_s}$ or
$\omega^{i_si_t}=-\omega^{i_ti_s}$; thus changing the orientation on
one edge will reverse the sign of $\Gamma(\Phi)$. The ordering of the
vertices tells us in which order to multiply the copies of $\Phi$
associated to the vertices; since we later project to the exterior
product $\Omega^{0,2k}(X)$, changing the order of the vertices by an
odd permutation will also change the sign of $\Gamma(\Phi)$ (whereas
changing the order by an even permutation will have no effect). It
follows that reversing the orientation of $\Gamma$ will have the
effect of changing the sign of $\Gamma(\Phi)$, and hence
$$b_{\bar{\Gamma}}(X)=-b_{\Gamma}(X)$$
where $\bar{\Gamma}$ is $\Gamma$ with the opposite orientation. It is
clear that this is the right notion of orientation of trivalent graphs
to use when defining the Rozansky-Witten invariants (and hence we
shall call it the Rozansky-Witten orientation).

The standard notion of orientation of a trivalent graph is given by an
equivalence class of cyclic ordering of the outgoing edges at each
vertex, with two orderings being equivalent if they differ on an even
number of vertices. Thus a trivalent graph drawn in the plane inherits
a canonical orientation given by taking the anticlockwise cyclic
ordering at each vertex (whenever we draw a trivalent graph in the
plane we shall assume this orientation). We promised to show that
these two notions of orientation are equivalent; essentially we follow
Kapranov's proof~\cite{kapranov99}.

Given a trivalent graph $\Gamma$, let $V(\Gamma)$, $E(\Gamma)$, and
$F(\Gamma)$ be the sets of its vertices, edges, and flags respectively
(a {\em flag\/} being an edge together with a choice of vertex lying
on it). For a finite set $S$, let ${\mathrm det}S$ be the highest
exterior power of the vector space ${\Bbb R}^S$ (the real vector space
whose basis is given by the elements of $S$). The standard notion of
orientation on $\Gamma$ is the same as an orientation on the
one-dimensional vector space
$$\bigotimes_{v\in V(\Gamma)}{\mathrm det}F(v)$$
where $F(v)$ is the three-element set of flags whose chosen vertex is
$v$. Now ${\mathrm det}F(\Gamma)$ is the highest exterior power of the
vector space
$$\bigoplus_{f\in F(\Gamma)}{\Bbb R}^{\{f\}}=\bigoplus_{v\in
V(\Gamma)}{\Bbb R}^{F(v)}.$$
The right hand side is a direct sum of three-dimensional spaces
parametrized by $V(\Gamma)$. Since three-forms anticommute, this
implies
$${\mathrm det}F(\Gamma)\cong {\mathrm
det}V(\Gamma)\otimes(\bigotimes_{v\in V(\Gamma)}{\mathrm det}F(v)).$$
We can also write
$$\bigoplus_{f\in F(\Gamma)}{\Bbb R}^{\{f\}}=\bigoplus_{e\in
E(\Gamma)}{\Bbb R}^{F(e)}$$
where $F(e)$ is the two-element set of flags containing the edge
$e$. The right hand side is a direct sum of two-dimensional spaces
parametrized by $E(\Gamma)$. However, two-forms commute so
$${\mathrm det}F(\Gamma)\cong \bigotimes_{e\in E(\Gamma)}{\mathrm
det}F(e).$$ 
The Rozansky-Witten orientation on $\Gamma$ is given by an equivalence
class of orientations on the edges and ordering of the vertices, or in
other words an orientation on the one-dimensional vector space
$${\mathrm det} V(\Gamma)\otimes(\bigotimes_{e\in E(\Gamma)}{\mathrm
det}F(e))$$
since an orientation of the two-dimensional space ${\Bbb R}^{F(e)}$ is
clearly the same as an orientation of the edge $e$. This space is
isomorphic to
$${\mathrm det}V(\Gamma)\otimes{\mathrm det}F(\Gamma)\cong{\mathrm
det}V(\Gamma)^{\otimes 2}\otimes(\bigotimes_{v\in V(\Gamma)}{\mathrm
det}F(v)).$$
Being the square of a real line, the one-dimensional space ${\mathrm
det}V(\Gamma)^{\otimes 2}$ has a canonical orientation, and hence the
Rozansky-Witten orientation is equivalent to an orientation on
$$\bigotimes_{v\in V(\Gamma)}{\mathrm det}F(v)$$
ie.\ equivalent to the standard notion of orientation on $\Gamma$,
completing the proof.

\begin{figure}[htpb]
\epsfxsize=60mm
\centerline{\epsfbox{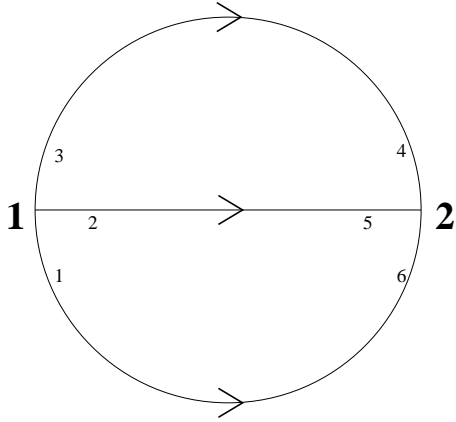}}
\caption{Compatibility of orientations}
\label{orientations}
\end{figure}

As an example, let $\Gamma$ be the two-vertex graph
$$\twoVgraph$$
with orientation induced from the planar embedding. Label the flags
$1,\ldots,6$ and the vertices $1,2$ as in
Figure~\ref{orientations}. 
\noindent
The corresponding orientation of ${\mathrm
det}F(\Gamma)$ is given by the element
$$f_1\wedge f_2\wedge f_3\wedge f_4\wedge f_5\wedge f_6$$
where the $f_i$'s are a basis for ${\Bbb R}^{F(\Gamma)}$. We can
rewrite this element as
$$f_1\wedge f_6\wedge f_2\wedge f_5\wedge f_3\wedge f_4$$
and this shows that the edges should be oriented from $1$ to $6$, from
$2$ to $5$, and from $3$ to $4$, as shown. More precisely, the edges
should be oriented in a manner {\em equivalent\/} to the one shown;
for example, the above element could also be rewritten
$$f_1\wedge f_6\wedge f_5\wedge f_2\wedge f_4\wedge f_3$$
which would correspond to orienting the edges from $1$ to $6$, from
$5$ to $2$, and from $4$ to $3$. This is equivalent as it differs by
reversing the orientation on an even number of edges.

We have shown that reversing the orientation of $\Gamma$ changes the
sign of $b_{\Gamma}(X)$. This is the first step in showing that the
Rozansky-Witten invariants descend to graph homology (which we shall
define shortly). The next step is to show compatibility with the IHX
relation.

Suppose that we are given three trivalent graphs $\Gamma_I$,
$\Gamma_H$ and $\Gamma_X$ which are identical except inside some ball
where they look like
\begin{center}
\begin{picture}(140,70)(-60,-35)
\put(-80,20){\line(0,-1){40}} 
\put(-80,20){\line(-2,1){20}} \put(-80,20){\line(2,1){20}}
\put(-80,-20){\line(-2,-1){20}} \put(-80,-20){\line(2,-1){20}}
%\put(-40,0){$=$}
\put(-10,0){\line(1,0){20}}
\put(-10,0){\line(-1,3){10}} \put(-10,0){\line(-1,-3){10}}
\put(10,0){\line(1,3){10}} \put(10,0){\line(1,-3){10}}
%\put(35,0){$-$}
\put(40,0){and}
\put(70,30){\line(2,-3){40}} 
\put(70,-30){\line(2,3){19}} \put(110,30){\line(-2,-3){19}}
\put(80,-15){\line(1,0){20}}
\end{picture}
%Figure $1$: The IHX relation
\end{center}
respectively. The orientations are induced from the planar embedding,
but let us just say that the corresponding Rozansky-Witten
orientations are given as follows: if the two vertices in each diagram
are labelled $t$ and $s$ then the connecting edge should be oriented
from $t$ to $s$ in $\Gamma_I$ and $\Gamma_X$ and from $s$ to $t$ in
$\Gamma_H$. We will show that the Rozansky-Witten invariants
corresponding to $\Gamma_I$, $\Gamma_H$, and $\Gamma_X$ are related by
$$b_{\Gamma_I}(X)=b_{\Gamma_H}(X)-b_{\Gamma_X}(X).$$
Reversing the orientation on $\Gamma_H$, this can be rewritten
$$b_{\Gamma_I}(X)+b_{\bar{\Gamma}_H}(X)+b_{\Gamma_X}(X)=0$$
where now the orientation of the edge joining the vertices $t$ and $s$
is the same in all three trivalent graphs.

The idea behind the proof is to show that 
$$\Gamma_I(\Phi)+\bar{\Gamma}_H(\Phi)+\Gamma_X(\Phi)\in\Omega^{0,2k}(X)$$
is a $\bar{\partial}$-exact form. Thus it is cohomologous to zero in
Dolbeault cohomology and hence using Kapranov's approach we obtain the
desired result. Consider $d_{LC}^2\Phi$, where $d_{LC}$ is the
exterior derivative corresponding to the Levi-Civita connection. Using
the facts that $d_{LC}^2$ gives the curvature and $\Phi$ is
$\bar{\partial}$-closed (by the Bianchi identity) we find
\begin{eqnarray*}
\bar{\partial}\partial_{LC}\Phi & = & d_{LC}^2\Phi \\
    & = & K\Phi\in\Omega^{1,2}(X,{\mathrm
    Sym}^3T^*)=\Omega^{0,2}(X,T^*\otimes{\mathrm Sym}^3T^*).
\end{eqnarray*}
Note that $\partial_{LC}$ also depends on the Levi-Civita connection,
whereas $\bar{\partial}$ does not. The product $K\Phi$ between
$K\in\Omega^{1,1}(X,{\mathrm End}T)$ and
$\Phi\in\Omega^{0,1}(X,{\mathrm Sym}^3T^*)$ is given by the induced
action of ${\mathrm End}T={\mathrm End}T^*$ on ${\mathrm Sym}^3T^*$
and wedging of forms. Note also that
$$\partial_{LC}\Phi\in\Omega^{1,1}(X,{\mathrm Sym}^3T^*)=\Omega^{0,1}(X,T^*\otimes{\mathrm Sym}^3T^*)$$
and so we can totally symmetrize $T^*\otimes{\mathrm Sym}^3T^*$ in
these terms to get forms taking values in ${\mathrm
Sym}^4T^*$. Denoting this symmetrizing operation by $S$, we obtain 
$$S(\partial_{LC}\Phi)\in\Omega^{0,1}(X,{\mathrm Sym}^4T^*)$$
and
$$S(K\Phi)\in\Omega^{0,2}(X,{\mathrm Sym}^4T^*)$$
such that $\bar{\partial}S(\partial_{LC}\Phi)=S(K\Phi)$. 

Consider the $I$ part of the graph $\Gamma_I$. Assume that the
vertices are labelled by $t$ and $s$ with $t<s$. We have seen that the
edge joining these two vertices is oriented from $t$ to $s$, and thus
when calculating $\Gamma_I(\Phi)$ this part contributes two copies of
$\Phi$ ``joined'' by a copy of $\tilde{\omega}$, namely the section
$$I(\Phi)\in C^{\infty}(X,(T^*)^{\otimes
4}\otimes(\bar{T}^*)^{\otimes 2})$$
with components
\begin{eqnarray*}
I(\Phi)_{j_tk_tj_sk_s\bar{l}_t\bar{l}_s} & = &
 \Phi_{i_tj_tk_t\bar{l}_t}\omega^{i_ti_s}\Phi_{i_sj_sk_s\bar{l}_s} \\
 & = &
 \omega_{i_tm}K^m_{\phantom{m}j_tk_t\bar{l}_t}\omega^{i_ti_s}\Phi_{i_sj_sk_s\bar{l}_s} \\
 & = &
 K^m_{\phantom{m}j_tk_t\bar{l}_t}\delta_m^{i_s}\Phi_{i_sj_sk_s\bar{l}_s} \\
 & = &
 K^m_{\phantom{m}j_tk_t\bar{l}_t}\Phi_{mj_sk_s\bar{l}_s}.
\end{eqnarray*}
The indices $j_tk_tj_sk_s$ refer to $(T^*)^{\otimes 4}$ and note that
this term is symmetric in $j_tk_t$ and $j_sk_s$. Similarly, the $\bar{H}$ part
of $\bar{\Gamma}_H$ contributes 
$$\bar{H}(\Phi)\in C^{\infty}(X,(T^*)^{\otimes 4}\otimes(\bar{T}^*)^{\otimes 2})$$
with components
$$\bar{H}(\Phi)_{j_tk_tj_sk_s\bar{l}_t\bar{l}_s}=K^m_{\phantom{m}j_tj_s\bar{l}_t}\Phi_{mk_tk_s\bar{l}_s}$$
to $\bar{\Gamma}_H(\Phi)$. This term is symmetric in $j_tj_s$ and
$k_tk_s$. Finally, the $X$ part of $\Gamma_X$ contributes the section 
$$X(\Phi)\in C^{\infty}(X,(T^*)^{\otimes
4}\otimes(\bar{T}^*)^{\otimes 2})$$
with components
$$X(\Phi)_{j_tk_tj_sk_s\bar{l}_t\bar{l}_s}=K^m_{\phantom{m}j_tk_s\bar{l}_t}\Phi_{mj_sk_t\bar{l}_s}$$
to $\Gamma_X(\Phi)$, and this term is symmetric in $j_tk_s$ and
$j_sk_t$. The sum of these three terms is totally symmetric in all
four holomorphic indices, ie.\
$$I(\Phi)+\bar{H}(\Phi)+X(\Phi)\in C^{\infty}(X,{\mathrm
Sym}^4T^*\otimes(\bar{T}^*)^{\otimes 2}).$$
In fact, projecting to the exterior product $\Omega^{0,2}(X,{\mathrm
Sym}^4T^*)$ the above sum gives $S(K\Phi)$, up to a factor. This is
because the symmetrization $S(K\Phi)$ of $K\Phi$ is given by the sum
of $24$ terms, but due to the symmetries already present in $K$ and
$\Phi$ there are actually just three distinct terms which correspond
to the $I(\Phi)$, $\bar{H}(\Phi)$, and $X(\Phi)$ terms above. So we
get eight copies of each of these three terms, and therefore
$$\Gamma_I(\Phi)+\bar{\Gamma}_H(\Phi)+\Gamma_X(\Phi)=\frac{1}{8}\Gamma_*(\Phi,S(K\Phi))\in\Omega^{0,2k}(X)$$
where $\Gamma_*$ is the graph which is identical to $\Gamma_I$ away
from the $I$ part (and hence identical to $\bar{\Gamma}_H$ away from
the $\bar{H}$ part, and $\Gamma_X$ away from the $X$ part), but
contains one tetravalent vertex instead of the $I$ part. The meaning
of the right hand side should be obvious; we use $\Phi$ at trivalent
vertices as before but now we substitute $S(K\Phi)$ for the
tetravalent vertex. Note that we can orient $\Gamma_*$ with the
Rozansky-Witten orientation, ie.\ an equivalence class of orientations
of the edges and ordering of the $2k-2$ trivalent vertices (the
tetravalent vertex is assumed to be labelled by $2k-1$). In this
situation the orientation is not so important because we only wish to
show that $\Gamma_*(\Phi,S(K\Phi))$ is exact, and as before, reversing
the orientation will merely change the sign of this term.

We have seen that $S(K\Phi)=\bar{\partial}S(\partial_{LC}\Phi)$, and so
we can substitute this into the above formula. Since $\Phi$ and
$\tilde{\omega}$ are $\bar{\partial}$-closed the right hand side
becomes 
$$\Gamma_*(\Phi,\bar{\partial}S(\partial_{LC}\Phi))=\bar{\partial}\Gamma_*(\Phi,S(\partial_{LC}\Phi))$$
where
$$\Gamma_*(\Phi,S(\partial_{LC}\Phi))\in\Omega^{0,2k-1}(X).$$
Therefore in Dolbeault cohomology
$$[\Gamma_I(\Phi)]+[\bar{\Gamma}_H(\Phi)]+[\Gamma_X(\Phi)]=0\in{\mathrm
H}_{\bar{\partial}}^{0,2k}(X),$$
and therefore
$$b_{\Gamma_I}(X)+b_{\bar{\Gamma}_H}(X)+b_{\Gamma_X}(X)=0$$
as claimed.

In the above argument we used a graph $\Gamma_*$ with one tetravalent
vertex and the other vertices trivalent. More generally, let
$\Gamma^{\prime}$ be a graph whose vertices are all of valency three
or greater. We can place copies of $\Phi$ at the trivalent vertices
and copies of $S(\partial_{LC}\Phi)$ at the tetravalent vertices. At
$b$-valent vertices we place copies of
$$\underbrace{S(\partial_{LC}\cdots
S(\partial_{LC}}_{b-3}\Phi)\cdots)\in\Omega^{0,1}(X,{\mathrm Sym}^bT^*).$$
Following the usual construction we arrive at a form
$$\Gamma^{\prime}(\Phi,S(\partial_{LC}\Phi),\ldots)\in\Omega^{0,V}(X)$$
where $V$ is the total number of vertices of $\Gamma^{\prime}$. Unlike
$\Gamma(\Phi)\in\Omega^{0,2k}(X)$, these more general forms are not
necessarily $\bar{\partial}$-closed, hence they do not represent
Dolbeault cohomology classes and we cannot obtain scalar-valued
invariants from them by integrating. However, they do enable us to
construct strong homotopy Lie algebra and operadic structures on the
spaces of forms and cohomology of hyperk{\"a}hler manifolds. For more
details see Kapranov~\cite{kapranov99}. Let us return to our
discussion of trivalent graphs.

Graph homology is the space of rational linear combinations of
oriented trivalent graphs modulo the AS and IHX relations. The AS
relations say that a graph with its orientation reversed is equivalent
to minus the graph, ie.\
$$\bar{\Gamma}\equiv -\Gamma.$$
The IHX relations say that 
$$\Gamma_I\equiv \Gamma_H-\Gamma_X.$$
We will continue to denote the equivalence class of a graph $\Gamma$
(its graph homology class) simply by $\Gamma$. The space of graph
homology is denoted ${\cal A}(\emptyset)$. It admits a grading given
by half the number of vertices; note that this grading is preserved by
the AS and  IHX relations. We denote the $k$th graded component (ie.\
graph homology classes with $2k$ vertices) by ${\cal
A}(\emptyset)^k$. The graphs in Figure $2$ can be taken as a basis for
graph homology in low degree. 
\begin{center}
$$
\begin{array}{ll}
k=1: & \twoVgraph \\
\phantom{0} & \\
k=2: & \twoVgraph^2 \qquad\fourVgraph \\
\phantom{0} & \\
k=3: & \twoVgraph^3 \qquad\twoVgraph\fourVgraph \qquad\sixVgraph \\
\phantom{0} & \\
k=4: &
\twoVgraph^4 \qquad\twoVgraph^2\fourVgraph \qquad\fourVgraph^2
\qquad\twoVgraph\sixVgraph \qquad\eightVgraphI \qquad\eightVgraphII
\\
\phantom{0} & \\
k=5: &
\twoVgraph^5 \qquad\twoVgraph^3\fourVgraph \qquad\twoVgraph\fourVgraph^2
\qquad\twoVgraph^2\sixVgraph \qquad\fourVgraph\sixVgraph
\\
 &
 \twoVgraph\eightVgraphI \qquad\twoVgraph\eightVgraphII
 \qquad\tenVgraphI \qquad\tenVgraphII \\
\phantom{k} & 
\end{array}
$$
Figure $2$: Basis for graph homology in low degree
\end{center}
\addtocounter{figure}{1}
\noindent
We will denote the two-vertex graph
$$\twoVgraph$$
by $\Theta$ and call it simply theta. Note that the graphs need not be
connected, and 
$$\twoVgraph^2$$
means to take the disjoint union of two copies of theta. We will
denote the graphs
$$\fourVgraph \qquad\sixVgraph \qquad\eightVgraphI \qquad\tenVgraphI \qquad\cdots$$
by $\Theta_2$, $\Theta_3$, $\Theta_4$, $\Theta_5$, etc.\ respectively,
and call them {\em necklace\/} graphs.

We can linearly extend the Rozansky-Witten invariants to arbitrary
rational linear combinations of trivalent graphs. Then from what we
have shown in this subsection the Rozansky-Witten invariants are
compatible with the AS and IHX relations (the graph homology
equivalence relations). Thus we have shown
\begin{prp}
The invariant $b_{\Gamma}(X)$ only depends on the trivalent graph
$\Gamma$ through its graph homology class. If $\Gamma_1$ and
$\Gamma_2$ are homologous then they define the same invariant
$b_{\Gamma_1}(X)=b_{\Gamma_2}(X)$.
\end{prp}
It follows that the Rozansky-Witten invariants allow us to associate
to a hyperk{\"a}hler manifold an element of graph cohomology ${\cal
A}(\emptyset)^*$, which is the space dual to the space of graph
homology. Conversely, we may at times choose to regard an element of
${\cal A}(\emptyset)^*$ as a ``virtual'' hyperk{\"a}hler manifold. In
the next chapter we shall relate certain Rozansky-Witten invariants to
characteristic numbers. For our virtual manifolds to make sense we
need to show that these characteristic numbers behave as they
should. For example, there will be certain integrality restraints,
though these can possibly be overcome by rescaling the graph
cohomology element. In the next chapter we will also check the
behaviour under taking products. The product on graph cohomology is
the dual of the coproduct on graph homology, which is given by the
formula 
$$\Delta(\Gamma)=\sum_{\gamma\sqcup{\gamma}^{\prime}=\Gamma}\gamma\otimes\gamma^{\prime}$$
where the sum is over all decompositions of $\Gamma$ into two disjoint
subgraphs. Graph homology also has a product given by disjoint union
of graphs. Together these structures make ${\cal A}(\emptyset)$ into a
commutative cocommutative Hopf algebra.  

Taking the coproduct of a graph is in some sense dual to taking the
product of hyperk{\"a}hler manifolds, as we shall now explain. Suppose
we have a reducible hyperk{\"a}hler manifold $X\times Y$, with
$X$ and $Y$ of real-dimensions $4k$ and $4l$ respectively. Let $p_1$
and $p_2$ be the projections from $X\times Y$ onto $X$ and $Y$
respectively. We wish to calculate the invariant $b_{\Gamma}(X\times
Y)$, where $\Gamma$ has $2k+2l$ vertices. Firstly, we know that the
holomorphic symplectic form, its dual, and the Riemann curvature
tensor of $X\times Y$ are all given by sums of the corresponding
objects pulled back from $X$ and $Y$ using $p_1^*$ and $p_2^*$
respectively. For example
$$\omega^{(X\times Y)}=p_1^*\omega^{(X)}+p_2^*\omega^{(Y)}.$$
If we take local coordinates on $X\times Y$ coming from local
coordinates on $X$ and $Y$, then this means that the matrix
$$\omega_{ij}^{(X\times Y)}=\left(\begin{array}{cc} \omega_{ij}^{(X)}
& 0 \\ 0 & \omega_{ij}^{(Y)} \end{array}\right)$$
of $\omega^{(X\times Y)}$ splits into block form, with the matrices of
$\omega^{(X)}$ and $\omega^{(Y)}$ along the diagonal. Similarly for
the dual symplectic form, the Riemann curvature tensor, and hence also
the section
$$\Phi^{(X\times Y)}=p_1^*\Phi^{(X)}+p_2^*\Phi^{(Y)}.$$
Let $\gamma$ be a connected component of $\Gamma$, with $2m$
vertices. Note that we can define
$$\gamma(\Phi^{(X\times Y)})\in\Omega^{0,2m}(X\times Y)$$ 
as before and it makes perfectly good sense. In fact, taking the wedge
product of all these forms for all connected components of $\Gamma$
gives us precisely 
$$\Gamma(\Phi^{(X\times Y)})\in\Omega^{0,2k+2l}(X\times Y).$$
Now when it comes to calculating $\gamma(\Phi^{(X\times Y)})$, the
indices must belong either all to $X$-coordinates or all to
$Y$-coordinates, for if they are mixed then at some stage we will need
to contract an $X$ index with a $Y$ index using
$\tilde{\omega}^{(X\times Y)}$, and this will give zero as
$\tilde{\omega}^{(X\times Y)}$ splits into block form in local
coordinates. Note that we have used the connectedness of $\gamma$
here. Therefore
$$\gamma(\Phi^{(X\times
Y)})=p_1^*\gamma(\Phi^{(X)})+p_2^*\gamma(\Phi^{(Y)}).$$
It follows that when calculating
$\Gamma(\Phi^{(X\times Y)})$, we must decompose $\Gamma$ into two
disjoint subgraphs $\gamma$ and $\gamma^{\prime}$, and then
$\Gamma(\Phi^{(X\times Y)})$ is given by the sum of
$$p_1^*\gamma(\Phi^{(X)})\wedge p_2^*\gamma^{\prime}(\Phi^{(Y)})$$
over all such decompositions. Now
\begin{eqnarray*}
(\omega^{(X\times Y)})^{k+l} & = &
(p_1^*\omega^{(X)}+p_2^*\omega^{(Y)})^{k+l} \\
 & = & {k+l \choose k}p_1^*(\omega^{(X)})^k\wedge
p_2^*(\omega^{(Y)})^l
\end{eqnarray*}
and therefore
\begin{eqnarray*}
b_{\Gamma}(X\times Y)\hspace*{-1mm} & = & \frac{1}{(8\pi^2)^{k+l}(k+l)!}\int_{X\times
Y}\Gamma(\Phi^{(X\times Y)})(\omega^{(X\times Y)})^{k+l} \\
 & = & \frac{1}{(8\pi^2)^{k+l}k!l!}\int_{X\times
Y}\sum_{\gamma\sqcup{\gamma}^{\prime}=\Gamma}p_1^*\gamma(\Phi^{(X)})\wedge p_2^*\gamma^{\prime}(\Phi^{(Y)})p_1^*(\omega^{(X)})^k\wedge
p_2^*(\omega^{(Y)})^l \\
 & = & \sum_{\gamma\sqcup{\gamma}^{\prime}=\Gamma}\left(\frac{1}{(8\pi^2)^kk!}\int_X\gamma(\Phi^{(X)})(\omega^{(X)})^k\right)\left(\frac{1}{(8\pi^2)^ll!}\int_Y\gamma^{\prime}(\Phi^{(Y)})(\omega^{(Y)})^l\right)
\\
 & = & \sum_{\gamma\sqcup{\gamma}^{\prime}=\Gamma}b_{\gamma}(X)b_{\gamma^{\prime}}(Y)
\end{eqnarray*}
If $\gamma$ has $2m$ vertices where $m>k$ then $\gamma(\Phi^{(X)})$ must
vanish as it would be given by projecting an element of $C^{\infty}(X,(\bar{T}^*)^{\otimes 2m})$ to the exterior product
$\Omega^{0,2m}(X)$, and the latter vanishes for ${\mathrm dim}_{\Bbb
R}X=4k$. Likewise if $\gamma^{\prime}$ has more than $2l$ vertices
then $\gamma^{\prime}(\Phi^{(Y)})$ must vanish. Hence on the
right hand side of the above formula, the terms in the sum vanish
unless $\gamma$ and $\gamma^{\prime}$ have precisely $2k$ and $2l$
vertices respectively. This formula can be rewritten in the following
way.
\begin{prp}
There is a correspondence between products of hyperk{\"a}hler
manifolds and coproducts of graphs. More precisely, if $X$ and $Y$ are
compact hyperk{\"a}hler manifolds then 
$$b_{\Gamma}(X\times Y)=b_{\Delta(\Gamma)}(X,Y).$$ 
\end{prp}
One important corollary of this proposition is that if $\Gamma$ is a
connected graph then its Rozansky-Witten invariant will vanish on
reducible hyperk{\"a}hler manifolds.

As another example, take $\Gamma$ to be the graph $\Theta^{k+l}$ given
by $k+l$ disjoint copies of theta. Since $\gamma$ must have $2k$
vertices, the only possibility is $\Theta^k$, but there are ${k+l
\choose k}$ different ways to choose which $k$ copies of theta make up
$\gamma$. Therefore
$$\frac{1}{(k+l)!}b_{\Theta^{k+l}}(X\times Y)=\frac{1}{k!}b_{\Theta^k}(X)\frac{1}{l!}b_{\Theta^l}(Y)$$
or in other words $\frac{1}{k!}b_{\Theta^k}$ is a multiplicative
invariant on hyperk{\"a}hler manifolds.

\subsection{Example calculation}

Let $S$ be a K$3$ surface, which is the unique compact irreducible
hyperk{\"a}hler manifold of real-dimension four. In this subsection we
shall explicitly
compute the Rozansky-Witten invariant $b_{\Theta}(S)$ using the
different approaches described in Subsection $1.5$. We shall assume
$S$ is the Kummer surface constructed as follows. Take a complex torus
$T^2={\Bbb C}^2/{\Bbb Z}^4$, and act on it with the involution
$$(z_1,z_2)\mapsto(-z_1,-z_2).$$
There are sixteen fixed-points: $(0,0)$, $(1/2,0)$, $(i/2,0)$, etc. We
blow-up the torus at these sixteen points to get $\widehat{T^2}$ and
then quotient by the involution. The resulting (smooth) surface $S$ is
a K$3$ surface. Denote the sixteen exceptional curves of the blow-up
in $\widehat{T^2}$ by $D_1,\ldots,D_{16}$; these are exceptional
curves of the first kind. They are fixed by the involution, and hence
we shall use the same notation $D_1,\ldots,D_{16}$ to denote the
images of these curves in $S$, which are really exceptional curves of
the second kind.

Near these exceptional curves $S$ looks like the cotangent bundle of
the projective line, $T^*{\Bbb P}^1$, with $D_i$ the zero section. The
space $T^*{\Bbb P}^1$ is a smooth resolution of
the double cone ${\Bbb C}^2/{\pm1}$. Away from the origin, the latter
can be described by coordinates $\pm(z_1,z_2)$. Suppose ${\Bbb P}^1$
has coordinates $\zeta$ and $\tilde{\zeta}$ on the complements of the
north and south poles respectively, with
$\tilde{\zeta}=-{\zeta}^{-1}$. A differential form $\eta d\zeta$
becomes $\tilde{\eta}d\tilde{\zeta}$ in the other coordinate patch,
where $\tilde{\eta}=\eta{\zeta}^2$. So $T^*{\Bbb P}^1$ can be covered
by two coordinate patches $U_1$ and $\tilde{U_1}$ whose coordinates
$(\zeta,\eta)$ and $(\tilde{\zeta},\tilde{\eta})$ are related in the
above way. The map onto ${\Bbb C}^2/{\pm1}$ is given by
$$(\zeta,\eta),(\tilde{\zeta},\tilde{\eta})\mapsto \pm(z_1,z_2)$$
where
$$(\eta,\eta\zeta,\eta{\zeta}^2)=(\tilde{\eta}\tilde{\zeta}^2,-\tilde{\eta}\tilde{\zeta},\tilde{\eta})=(z_1^2,z_1z_2,z_2^2).$$
An open cover for $S$ is given by
$$\{U_0,U_1,\tilde{U}_1,\ldots,U_{16},\tilde{U}_{16}\}$$
where each pair $\{U_i,\tilde{U}_i\}$ cover a small neighbourhood of
$D_i$ (so that $U_i\cup\tilde{U}_i$ does not intersect
$U_j\cup\tilde{U}_j$ for $i\neq j$), and the open set $U_0$ is the
complement of the exceptional curves.

To get a representation of the Atiyah class $\alpha_T$ of $S$ in {\v
C}ech cohomology we choose a holomorphic connection over each open set
in our cover. In fact, we may as well choose flat connections. For
example, on $U_0$ we choose the connection ${\nabla}_0$ characterized
by
$${\nabla}_0(\part{z_1})={\nabla}_0(\part{z_2})=0.$$
In $U_1$ coordinates this becomes
\begin{eqnarray*}
{\nabla}_0(\part{\zeta}) & = & \frac{1}{2\eta}(d\eta\otimes\part{\zeta}) \\
{\nabla}_0(\part{\eta}) & = &
\frac{1}{2\eta}(d\zeta\otimes\part{\zeta}-d\eta\otimes\part{\eta})
\end{eqnarray*}
So on $U_0\cap U_1$ we find
$$(\alpha_T)_{01}={\nabla}_0-{\nabla}_1=\frac{1}{2\eta}\left(\begin{array}{cc}
  d\eta & d\zeta \\
  0 & -d\eta
  \end{array}\right)\in{\mathrm H}^0(U_0\cap
U_1,T^*\otimes{\mathrm End}T)$$
in $U_1$ coordinates. Similarly for $(\alpha_T)_{0\tilde{1}}$,
$(\alpha_T)_{02}$, etc. The flat connection $\tilde{\nabla}_1$ on
$\tilde{U}_1$ is characterized by
$$\tilde{\nabla}_1(\part{\tilde{\zeta}})=\tilde{\nabla}_1(\part{\tilde{\eta}})=0.$$
In $U_1$ coordinates this becomes
\begin{eqnarray*}
\tilde{\nabla}_1(\part{\zeta}) & = &
\frac{2}{\zeta}(-d\zeta\otimes\part{\zeta}+d\eta\otimes\part{\eta})+\frac{6\eta}{\zeta^2}(d\zeta\otimes\part{\eta})
\\
\tilde{\nabla}_1(\part{\eta}) & = &
\frac{2}{\zeta}(d\zeta\otimes\part{\eta})
\end{eqnarray*}
So on $U_1\cap\tilde{U_1}$ we find
$$(\alpha_T)_{\tilde{1}1}=\tilde{\nabla}_1-{\nabla}_1=\frac{2}{\zeta^2}\left(\begin{array}{cc}
  -\zeta d\zeta & 0 \\
  \zeta d\eta+3\eta d\zeta & \zeta d\zeta
  \end{array}\right)\in{\mathrm
H}^0(U_1\cap\tilde{U_1},T^*\otimes{\mathrm End}T)$$
in $U_1$ coordinates. Similarly for $(\alpha_T)_{2\tilde{2}}$,
etc. Thus we have an element $\alpha_T$ in the {\v C}ech cohomology
group ${\mathrm H}^1(S,T^*\otimes{\mathrm End}T)$ which represents the
Atiyah class.

Next we rewrite $T^*\otimes{\mathrm End}T$ as $(T^*)^{\otimes 3}$ by
using the holomorphic symplectic form $\omega$ on $S$ to identify $T$
and $T^*$. Since Rozansky-Witten invariants do not change under
rescalings of the symplectic form, we can choose $\omega$ to look like
$2d\zeta\wedge d\eta$ on $U_1$. Then the matrices $\omega_{ij}$ and
$\omega^{ij}$ both equal
$$\left(\begin{array}{cc} 0 & 1 \\ -1 & 0
\end{array}\right)$$
and hence the identification of $T$ and $T^*$ is given by 
$$\part{\zeta}\leftrightarrow -d\eta$$
$$\part{\eta}\leftrightarrow d\zeta.$$ 
Thus using $\omega$ to rewrite $\alpha_T$, we find
\begin{eqnarray}
(\alpha_T)_{01} & = & \frac{-1}{2\eta}(d\eta\otimes d\zeta\otimes
d\eta+d\zeta\otimes d\eta\otimes d\eta+d\eta\otimes d\eta\otimes
d\zeta) \\
(\alpha_T)_{\tilde{1}1} & = & \frac{2}{\zeta}(d\zeta\otimes d\zeta\otimes
d\eta+d\eta\otimes d\zeta\otimes d\zeta+d\zeta\otimes d\eta\otimes
d\zeta)+\frac{6\eta}{\zeta^2}(d\zeta\otimes d\zeta\otimes
d\zeta).\nonumber
\end{eqnarray}
We get similar formula for $(\alpha_T)_{0\tilde{1}}$,
$(\alpha_T)_{02}$, etc.\ and in particular we see that $\alpha_T$ is in
${\mathrm H}^1(S,{\mathrm Sym}^3T^*)$.

To represent the Atiyah class $\alpha_T$ of $S$ in Dolbeault
cohomology we begin with a smooth global connection $\nabla$ of type
$(1,0)$ on $S$. To obtain such a connection we could take the flat
connections on each of the open sets
$$\{U_0,U_1,\tilde{U}_1,\ldots,U_{16},\tilde{U}_{16}\}$$
and then patch them together using a partition of unity
$$\{\psi_0,\psi_1,\tilde{\psi}_1,\ldots,\psi_{16},\tilde{\psi}_{16}\}.$$
For example, in $U_1$ coordinates this would look like
\begin{eqnarray*}
\nabla |_{U_1} & = &
 \psi_0{\nabla}_0+\psi_1{\nabla}_1+\tilde{\psi}_1\tilde{\nabla}_1 \\
 & = & \psi_0\left(d+\frac{1}{2\eta}\left(\begin{array}{cc} d\eta & d\zeta
 \\ 0 & -d\eta
 \end{array}\right)\right)+\psi_1d+\tilde{\psi}_1\left(d+\frac{2}{\zeta^2}\left(\begin{array}{cc}-\zeta
 d\zeta & 0 \\ \zeta d\zeta+3\eta d\zeta & \zeta d\zeta
 \end{array}\right)\right) \\
 & = & d+\psi_0\frac{1}{2\eta}\left(\begin{array}{cc} d\eta & d\zeta
 \\ 0 & -d\eta
 \end{array}\right)+\tilde{\psi}_1\frac{2}{\zeta^2}\left(\begin{array}{cc}-\zeta
 d\zeta & 0 \\ \zeta d\eta+3\eta d\zeta & \zeta d\zeta
 \end{array}\right) \\
\end{eqnarray*}
since $\psi_0+\psi_1+\tilde{\psi}_1$ is identically one on $U_1$. The
$(1,1)$ part of the curvature of this connection gives a Dolbeault
representative of the Atiyah class. If we write the above expression
$\nabla=d+A$, then this is precisely $\bar{\partial}A$, and
$$\alpha_T=\left[\frac{1}{2\eta}\left(\begin{array}{cc} d\eta & d\zeta \\ 0 & -d\eta \end{array}\right)\bar{\partial}\psi_0+\frac{2}{\zeta^2}\left(\begin{array}{cc}-\zeta d\zeta & 0 \\ \zeta d\eta+3\eta d\zeta & \zeta d\zeta \end{array}\right)\bar{\partial}\tilde{\psi}_1\right].$$
Using $\omega$ to rewrite our of representative of $\alpha_T$ as an
element of $\Omega^{0,1}(S,{\mathrm Sym}^3T^*)$, we get
\begin{eqnarray}
\alpha_T & = & \left[\frac{-1}{2\eta}(d\eta\otimes d\zeta\otimes
d\eta+d\zeta\otimes d\eta\otimes d\eta+d\eta\otimes d\eta\otimes
d\zeta)\bar{\partial}\psi_0\right. \\
 & & \hspace*{-10mm}+\left.\left(\frac{2}{\zeta}(d\zeta\otimes d\zeta\otimes d\eta+d\eta\otimes
d\zeta\otimes d\zeta+d\zeta\otimes d\eta\otimes
d\zeta)+\frac{6\eta}{\zeta^2}(d\zeta\otimes d\zeta\otimes
d\zeta)\right)\bar{\partial}\tilde{\psi}_1\right] \nonumber
\end{eqnarray}
in ${\mathrm H}^{0,1}_{\bar{\partial}}(S,{\mathrm Sym}^3T^*)$.

Our final description of the Atiyah class is as the residue of a
meromorphic connection. Let $D$ be the divisor on $S$ given by the sum
of the sixteen exceptional curves
$$D_1+\ldots +D_{16}$$
which are given locally by $\eta_i=\tilde{\eta}_i=0$. In $U_1$
coordinates, the flat connection $\nabla_0$ on $U_0$ contains a factor
of $\eta^{-1}$, and similarly for $\tilde{U}_1$, $U_2$, etc. Thus
$\nabla_0$ can be extended to a meromorphic connection on $S$ with a
simple pole along $D$, which we continue to denote by $\nabla_0$. Let
$L={\cal O}(D)$ be the line bundle associated to $D$, and $s$ the
canonical section, which is given locally by
$$s=\left\lbrace\begin{array}{ccc}
    1 & \mbox{in} & U_0 \\
    {\eta}_i & \mbox{in} & U_i \\
    \tilde{\eta}_i & \mbox{in} & \tilde{U}_i. 
    \end{array}\right.$$
Then $s{\nabla}_0$ is non-singular on $S$, and its restriction to $D$
is the image of the residue $\beta_T$ in ${\mathrm H}^0(D,{\mathrm
Sym}^3T^*\otimes L)$ as in Subsection $1.5$, where as usual we have
identified $T$ and $T^*$. A local calculation shows that
\begin{eqnarray}
(\beta_T)_1 & = & \frac{-1}{2}(d\eta\otimes d\zeta\otimes
d\eta+d\zeta\otimes d\eta\otimes d\eta+d\eta\otimes d\eta\otimes
d\zeta) \\
(\beta_T)_{\tilde{1}} & = & \frac{-1}{2}(d\tilde{\eta}\otimes
d\tilde{\zeta}\otimes d\tilde{\eta}+d\tilde{\zeta}\otimes
d\tilde{\eta}\otimes d\tilde{\eta}+d\tilde{\eta}\otimes
d\tilde{\eta}\otimes d\tilde{\zeta}) \nonumber
\end{eqnarray}
where we should remember that we have restricted to $D$ so that
$\eta=\tilde{\eta}=0$, and also that this is a section of a bundle
which has been twisted by $L$. We get similar formulae for
$(\beta_T)_2$, $(\beta_T)_{\tilde{2}}$, etc.\ and the Atiyah class is
given by $\alpha_T=\delta\beta_T$.

Now we calculate the Rozansky-Witten invariant $b_{\Theta}(S)$ using
the three different descriptions of the Atiyah class just given. We
begin with the {\v C}ech cohomology description.

We place a copy of the Atiyah class at each of the two vertices of
$\Theta$, and then contract along the edges using three copies of
$\tilde{\omega}$. The product of two {\v C}ech cohomology classes
$\{a_{ij}\}$ and $\{b_{ij}\}$ is given by
$$(ab)_{ijk}=\frac{1}{6}(a_{ij}b_{ik}-a_{ik}b_{ij}+a_{jk}b_{ji}-a_{ji}b_{jk}+a_{ki}b_{kj}-a_{kj}b_{ki}).$$
In the case of the Atiyah class~(1) this simplifies to
$$(\alpha_T\alpha_T)_{0\tilde{1}1}=\frac{1}{2}((\alpha_T)_{01}(\alpha_T)_{\tilde{1}1}-(\alpha_T)_{\tilde{1}1}(\alpha_T)_{01})$$
and similarly for $(\alpha_T\alpha_T)_{0\tilde{2}2}$, etc. In $U_1$
coordinates $\tilde{\omega}$ looks like
$$2\part{\zeta}\wedge\part{\eta}$$
and a calculation shows that using this to contract along the edges we
get
$$\Theta(\alpha_T)_{0\tilde{1}1}=\frac{-3}{\zeta\eta}.$$
Similar for $\Theta(\alpha_T)_{0\tilde{2}2}$, etc.\ and together these
give us $\Theta(\alpha_T)$ in ${\mathrm H}^2(S,{\cal O}_S)$. The Serre
duality pairing with $\omega$ in ${\mathrm H}^0(S,\Lambda^2T^*)$ can
be calculated by multiplying by this section and then integrating
along contours. For example, in $U_1$ coordinates we integrate along
the contour 
$$\{|\zeta |=|\eta |=1\}\subset U_0\cap\tilde{U}_1\cap U_1$$
to obtain
$$\int_{|\zeta |=1}\int_{|\eta |=1}\frac{-3}{\zeta\eta}2d\zeta\wedge
d\eta=-6(2\pi i)^2.$$
Summing over all sixteen exceptional curves and then dividing by the
factor $8\pi^2$ gives us the Rozansky-Witten invariant
$$b_{\Theta}(S)=48.$$

Next we do the same calculation with Dolbeault cohomology
classes. We work in $U_1$ coordinates. Taking two copies of our
representative of the Atiyah class from
Equation~(2), and contracting with $\tilde{\omega}$
gives us the section
$$\frac{-3}{\zeta\eta}(\bar{\partial}\psi_0\otimes\bar{\partial}\tilde{\psi}_1-\bar{\partial}\tilde{\psi}_1\otimes\bar{\partial}\psi_0)\in
C^{\infty}(S,(\bar{T}^*)^{\otimes 2})$$
and projecting to the exterior product gives us a Dolbeault
representative of
$$\Theta(\alpha_T)=\left[\frac{-6}{\zeta\eta}\bar{\partial}\psi_0\wedge\bar{\partial}\tilde{\psi}_1\right]\in{\mathrm
H}^{0,2}_{\bar{\partial}}(S,{\cal O}_S).$$
Multiplying by the symplectic form and integrating gives us
$$\int_{U_1}\frac{-12}{\zeta\eta}\bar{\partial}\psi_0\wedge\bar{\partial}\tilde{\psi}_1\wedge
d\zeta\wedge d\eta.$$
Since the integrand is supported inside $U_1$ anyway, we can think of
this as being an integral over all of $T^*{\Bbb P}^1$, and because
$d\zeta\wedge d\zeta$ and $d\eta\wedge d\eta$ vanish, we can write
\begin{eqnarray*}
\bar{\partial}\psi_0\wedge\bar{\partial}\tilde{\psi}_1\wedge
d\zeta\wedge d\eta & = &
\bar{\partial}\psi_0\wedge\bar{\partial}\tilde{\psi}_1\wedge
d\zeta\wedge d\eta
+\bar{\partial}\psi_0\wedge\partial\tilde{\psi}_1\wedge d\zeta\wedge
d\eta \\
 & = & \bar{\partial}\psi_0\wedge d\tilde{\psi}_1\wedge d\zeta\wedge
d\eta \\
 & = & \bar{\partial}\psi_0\wedge d\tilde{\psi}_1\wedge d\zeta\wedge
d\eta +\partial\psi_0\wedge d\tilde{\psi}_1\wedge d\zeta\wedge d\eta
\\
 & = & d\psi_0\wedge d\tilde{\psi}_1\wedge d\zeta\wedge d\eta.
\end{eqnarray*}
Substituting this back into the above integral we obtain 
$$\int_{T^*{\Bbb P}^1}\frac{-12}{\zeta\eta}d\psi_0\wedge d\tilde{\psi}_1\wedge
d\zeta\wedge d\eta.$$
Now $\psi_0$ vanishes on $U_1\backslash U_0$ and is identically one on
$U_0\backslash (U_1\cap\tilde{U}_1)$ (at least away from the 
other exceptional curves), and $\tilde{\psi}_1$ vanishes on
$U_1\backslash\tilde{U}_1$ and is identically one on
$\tilde{U}_1\backslash (U_0\cap U_1)$. In fact, we can choose these
functions to be invariant under the torus action on $T^*{\Bbb P}^1$.
In other words, converting to polar coordinates
$\zeta=re^{i\theta}$ and $\eta=se^{i\phi}$, we can choose
$\psi_0$ and $\tilde{\psi}_1$ to be functions of only $r$ and $s$,
with $\psi_0$ vanishing for $s=0$ and identically one for $s$ large,
and $\tilde{\psi}_1$ vanishing for $r=0$ or $s$ large, and identically
one for $r$ large and $s=0$. Hence the integral becomes
$$\int_{r=0}^{\infty}\int_{\theta=0}^{2\pi}\int_{s=0}^{\infty}\int_{\phi=0}^{2\pi}12d\psi_0\wedge 
d\tilde{\psi}_1\wedge d\theta\wedge d\phi=12(2\pi)^2\int_{r=0}^{\infty}\int_{s=0}^{\infty}d\psi_0\wedge
d\tilde{\psi}_1.$$
Writing the integrand $d(\psi_0d\tilde{\psi}_1)$ and applying Stokes'
theorem gives us an integral along the oriented boundary of the first
quadrant in ${\Bbb R}^2$, namely
$$48\pi^2\int_{\{r=0\}\cup\{s=\infty\}\cup\{\overline{r=\infty\phantom{0}\hspace*{-1.5mm}}\}\cup\{\overline{s=0}\}}\psi_0d\tilde{\psi}_1$$
where $\{\overline{r=\infty\phantom{0}\hspace*{-2mm}}\}$ and
$\{\overline{s=0}\}$ denote $\{r=\infty\}$ and $\{s=0\}$,
respectively, but with orientations opposite to the standard
ones. Observe that $\tilde{\psi}_1$ vanishes on $\{r=0\}$ and
$\{s=\infty\}$, and $\psi_0$ vanishes on $\{\overline{s=0}\}$. On
$\{\overline{r=\infty\phantom{0}\hspace*{-2mm}}\}$, $\psi_1$ vanishes
and therefore $\tilde{\psi}_1\equiv 1-\psi_0$ (as these functions form
a partition of unity, and our integral is supported in an open set
which is away from the other fifteen exceptional curves). Hence our
integral becomes
\begin{eqnarray*}
48\pi^2\int_{\{\overline{r=\infty\phantom{0}\hspace*{-1.5mm}}\}}\psi_0d(1-\psi_0) & = &
48\pi^2\int_{s=0}^{\infty}\psi_0d\psi_0|_{r=\infty} \\
 & = &
48\pi^2\left[\frac{1}{2}\psi_0^2\right]_{r=\infty,s=0}^{r=\infty,s=\infty} \\
 & = & 24\pi^2.
\end{eqnarray*}
Again we need to sum over all sixteen exceptional curves and divide by
the factor $8\pi^2$, obtaining
$$b_{\Theta}(S)=48.$$

Finally we calculate $b_{\Gamma}(S)$ using the residue description of
the Atiyah class~(3). First we compute
$\Theta(\alpha_T,\beta_T)$ in ${\mathrm H}^1(D,L)$ using a {\v C}ech
cohomology class for $\alpha_T$. In $U_1$ coordinates this is given by
multiplying $(\alpha_T)_{1\tilde{1}}$ by $(\beta_T)_1$, and then
contracting with three copies of $\tilde{\omega}$, and the result is
$$\Theta(\alpha_T,\beta_T)_{1\tilde{1}}=\frac{3}{\zeta}$$
where we have used the fact that $\eta$ vanishes on $D$. Next we need
to take the Serre duality pairing of this element with the section
$\omega|_D$ of $\Lambda^2T^*|_D$. Since $D$ is given locally by
$\eta=0$, we write
$$\omega|_D=d\eta\wedge (-2d\zeta)$$
and hence we should multiply $\Theta(\alpha_T,\beta_T)$ by the section
$-2\eta d\zeta$ of ${\cal K}_D\otimes L^*|_D$ and then perform an
integral along the contour 
$$\{|\zeta|=1\}\subset U_1\cap\tilde{U}_1.$$
Note that $\Theta(\alpha_T,\beta_T)_{1\tilde{1}}$ is a local section of
$L$, or in other words a local meromorphic function with a simple pole
along $D$. The latter description requires us to include a factor of
$\eta^{-1}$ which cancels with the $\eta$ in $-2\eta d\zeta$, and
hence we get 
$$\int_{|\zeta|=1}\frac{3}{\zeta}(-2d\zeta)=-6(2\pi i)$$
Alternatively, using a Dolbeault cohomology class for $\alpha_T$ we
obtain the expression
$$\Theta(\alpha_T,\beta_T)=\left[\frac{3}{\zeta}\bar{\partial}\tilde{\psi}_1|_D\right]\in{\mathrm 
H}^{0,1}_{\bar{\partial}}(D,L)$$
in $U_1$ coordinates. Here we regard $\tilde{\psi}_1|_D$ as simply a
function of $\zeta$ and $\bar{\zeta}$. In fact, converting to polar
coordinates $\zeta=re^{i\theta}$ as before, we may choose
$\tilde{\psi}_1|_D$ to be a function of $r$ only, vanishing for $r=0$
and identically one for $r$ large. Multiplying by $-2d\zeta$ as above
and integrating over $D_1$ gives us
\begin{eqnarray*}
\int_{D_1}\frac{3}{\zeta}\bar{\partial}\tilde{\psi}_1|_D\wedge
(-2d\zeta) & = & -6i\int_{r=0}^{\infty}\int_{\theta=0}^{2\pi}d\tilde{\psi}_1|_D\wedge d\theta \\
 & = & -6(2\pi i)\int_{r=0}^{\infty}d\tilde{\psi}_1|_D \\
 & = & -6(2\pi i)\left[\tilde{\psi}_1|_D\right]_{r=0}^{r=\infty} \\
 & = & -6(2\pi i)
\end{eqnarray*}
as before. Including the factor $2\pi i$ (from the Cauchy residue
formula), along with the usual factors $16$ and $(8\pi^2)^{-1}$, gives
us
$$b_{\Theta}(S)=48.$$

%% file: ch2.tex
\section {Characteristic numbers}

\subsection{Chern numbers and Chern-Weil theory}

By Chern-Weil theory we can write the Chern numbers of a manifold as
integrals of copies of the curvature, multiplied in some way. We will
make this precise below, but first let us note that Rozansky-Witten
invariants were defined in a similar way, and hence it is reasonable
to expect some relation between them. In fact all the Chern numbers
can be expressed in terms of Rozansky-Witten invariants, as we shall
show in this chapter. Whether or not the converse is true, ie.\ that
all Rozansky-Witten invariants can be expressed in terms of the Chern
numbers, is a fundamental question of this theory and one we shall
investigate further in Chapter $5$.

For a hyperk{\"a}hler manifold we shall see that all the odd Chern
classes vanish, and the remaining even Chern classes are really the
Pontryagin classes, which are topological invariants. So although we
find it easier to work with Chern classes, our results will be
completely independent of the choice of compatible complex structure
used to define them. So let $X$ be a compact hyperk{\"a}hler manifold
of real-dimension $4k$ and fix a complex structure $I$. By Chern-Weil
theory we can represent the Chern character of $X$ by traces of powers
of the Riemann curvature tensor
$$K^i_{\phantom{i}jk\bar{l}}\in\Omega^{1,1}(X,{\mathrm End}T).$$ 
More precisely, we have
\begin{eqnarray*}
ch(T) & = & 2k+ch_1(T)+ch_2(T)+\ldots +ch_{2k}(T) \\
      & = & \sum_{m=0}^{2k}\frac{(-1)^m}{m!(2\pi
      i)^m}[{\mathrm Tr}(K^m)]
\end{eqnarray*}
where we multiply $K$'s by composing ${\mathrm End}T$ and taking the
wedge product on forms, and then take the trace in ${\mathrm
End}T$. Note that since we can identify $T$ and $T^*$ using $\omega$,
we have
$$ch_m(T)=ch_m(T^*)=(-1)^mch_m(T)$$ 
and hence all odd components of the Chern character vanish for
hyperk{\"a}hler manifolds, as claimed above.

The Chern numbers of $X$ are obtained by taking components of the
Chern character and wedging them together to get something which can
be integrated over $X$. This means taking an even partition
$(\lambda_1,\ldots,\lambda_j)$ of $2k$ (ie.\ one for which every
$\lambda_i$ is even) and then calculating
$$ch_{\lambda_1}\cdots
ch_{\lambda_j}(X)=\int_Xch_{\lambda_1}(T)\wedge\cdots\wedge
ch_{\lambda_j}(T).$$
The integrand lies in $\Omega^{2k,2k}(X)$. Following Hirzebruch's
notation~\cite{hirzebruch78}, we shall work with a rescaling 
$$s_{\lambda}(T)={\lambda}!ch_{\lambda}(T)=\frac{1}{(2\pi
i)^{\lambda}}[{\mathrm Tr}(K^{\lambda})]\in{\mathrm
H}^{2\lambda}(X,{\Bbb Z})$$
of the components of the Chern character (note that we have dropped
the sign $(-1)^{\lambda}$ since $\lambda$ is even). Then we
get 
\begin{eqnarray*}
s_{\lambda_1}\cdots s_{\lambda_j}(X) & = &
\int_Xs_{\lambda_1}(T)\wedge\cdots\wedge s_{\lambda_j}(T) \\
   & = & \frac{1}{(2\pi i)^{2k}}\int_X{\mathrm
   Tr}(K^{\lambda_1})\wedge\cdots\wedge {\mathrm Tr}(K^{\lambda_j})
\end{eqnarray*}
which we shall still call Chern numbers.

\subsection{Relation to Rozansky-Witten invariants}

\begin{figure}[htpb]
\epsfxsize=60mm
\centerline{\epsfbox{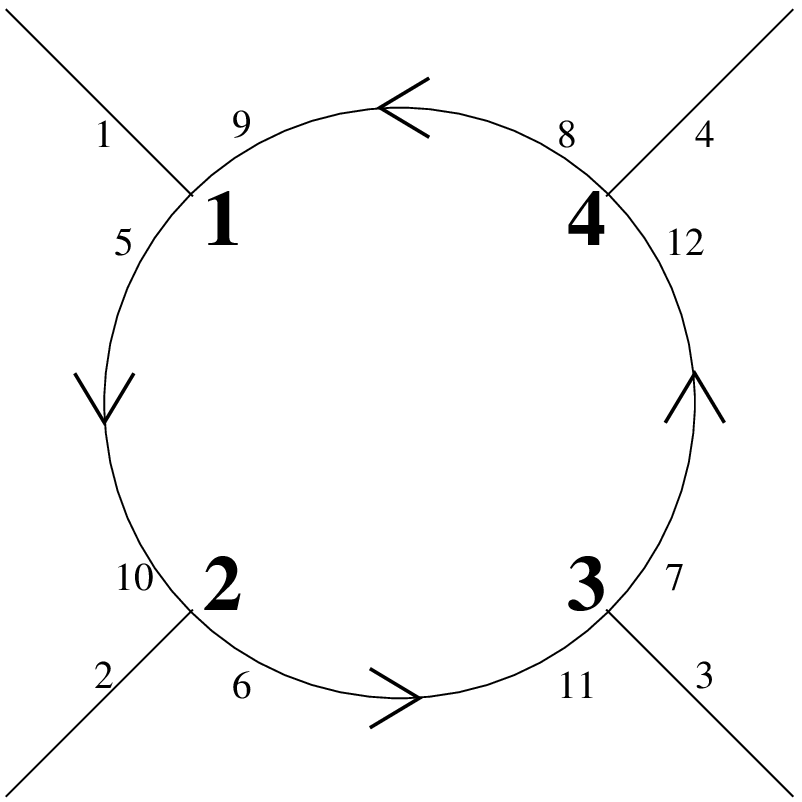}}
\caption{Orientation of a wheel}
\label{wheels_orient}
\end{figure}

We wish to express the Chern numbers of $X$ in terms of the
Rozansky-Witten invariants $b_{\Gamma}(X)$ for the appropriate choice
of graph homology class $\Gamma$. We will see that {\em wheels\/} play
an important role. Suppose that the graph $\Gamma$ contains a
$\lambda$-wheel $w_{\lambda}$, ie.\ a closed path with distinct
vertices and $\lambda$ edges. Each vertex has a third outgoing edge,
and we shall call these the spokes of the wheel. For example, an
$8$-wheel looks like:
\begin{center}
\begin{picture}(100,100) (50,50)
\put(100,100){\circle{40}}
\put(120,100){\line(1,0){30}}
\put(114,114){\line(1,1){20}}
\put(100,120){\line(0,1){25}}
\put(86,114){\line(-1,1){20}}
\put(80,100){\line(-1,0){30}}
\put(86,86){\line(-1,-1){20}}
\put(100,80){\line(0,-1) {25}}
\put(114,86){\line(1,-1) {20}}
\end{picture}
\end{center}
It is not too difficult to show using the AS and IHX relations that a
graph containing an odd wheel must be homologous to zero, so we shall
assume that $\lambda$ is even. We must be careful with our
orientation, which is induced from the planar embedding as
before. To see this as a Rozansky-Witten orientation let us assume,
without loss of generality, that the $\lambda$ vertices of the wheel
are ordered $1,\ldots,\lambda$ in an anticlockwise manner, and that
the flags at the $m$th vertex are labelled $m,\lambda+m,2\lambda+m$ in an
anticlockwise manner, with $m$ labelling the flag corresponding to the
spoke. The corresponding orientation of ${\mathrm
det}F(w_{\lambda})$ is given by the element
$$(f_1\wedge f_{\lambda+1}\wedge f_{2\lambda+1})\wedge (f_2\wedge
f_{\lambda+2}\wedge f_{2\lambda+2})\wedge\cdots\wedge (f_{\lambda}\wedge
f_{2\lambda}\wedge f_{3\lambda})$$
where the $f_i$'s are a basis for ${\Bbb R}^{F(w_{\lambda})}$ (recall
the example in the previous chapter). We can rewrite this element as
$$-(f_1\wedge f_2\wedge\cdots\wedge f_{\lambda})\wedge (f_{\lambda+1}\wedge
f_{2\lambda+2})\wedge (f_{\lambda+2}\wedge f_{2\lambda+3})\wedge\cdots\wedge (f_{2\lambda}\wedge f_{2\lambda+1})$$
and this shows that the edges of $\bar{w}_{\lambda}=-w_{\lambda}$
should be oriented in the anticlockwise direction, where we have
ordered the flags $1,2,\ldots,\lambda$ corresponding to the
spokes of the wheel in the same order as the vertices themselves. This
is best illustrated by an example, and Figure~\ref{wheels_orient}
shows the case $\lambda=4$.

When calculating $\Gamma(\Phi)$ this part of the graph will contribute
the section
$$w_{\lambda}(\Phi)\in C^{\infty}(X,(T^*\otimes\bar{T}^*)^{\otimes\lambda})$$
which with respect to local complex coordinates has components
\begin{eqnarray*}
w_{\lambda}(\Phi)_{k_1\bar{l}_1\cdots
k_{\lambda}\bar{l}_{\lambda}} & = &
-\Phi_{i_1j_1k_1\bar{l}_1}\omega^{j_1i_2}\Phi_{i_2j_2k_2\bar{l}_2}\omega^{j_2i_3}\cdots\Phi_{i_{\lambda}j_{\lambda}k_{\lambda}\bar{l}_{\lambda}}\omega^{j_{\lambda}i_1}
\\
 & = &
 -\omega_{i_1m_1}K^{m_1}_{\phantom{m_1}j_1k_1\bar{l}_1}\omega^{j_1i_2}\omega_{i_2m_2}K^{m_2}_{\phantom{m_2}j_2k_2\bar{l}_2}\omega^{j_2i_3}\cdots\omega_{i_{\lambda}m_{\lambda}}K^{m_{\lambda}}_{\phantom{m_{\lambda}}j_{\lambda}k_{\lambda}\bar{l}_{\lambda}}\omega^{j_{\lambda}i_1} \\
 & = &
 -K^{m_1}_{\phantom{m_1}j_1k_1\bar{l}_1}(-\delta^{j_1}_{m_2})K^{m_2}_{\phantom{m_2}j_2k_2\bar{l}_2}(-\delta^{j_2}_{m_3})\cdots(-\delta^{j_{\lambda-1}}_{m_{\lambda}})K^{m_{\lambda}}_{\phantom{m_{\lambda}}j_{\lambda}k_{\lambda}\bar{l}_{\lambda}}(-\delta^{j_{\lambda}}_{m_1}) \\
 & = &
 -(-1)^{\lambda}K^{m_1}_{\phantom{m_1}m_2k_1\bar{l}_1}K^{m_2}_{\phantom{m_2}m_3k_2\bar{l}_2}\cdots K^{m_{\lambda}}_{\phantom{m_{\lambda}}m_1k_{\lambda}\bar{l}_{\lambda}} \\
 & = & -{\mathrm
 Tr}(K^{\otimes\lambda})_{k_1\bar{l}_1\cdots
 k_{\lambda}\bar{l}_{\lambda}}
\end{eqnarray*}
where $k_1,\ldots,k_{\lambda}$ denote the indices attached to the
spokes, and we have used the fact that $\lambda$ is even. Note that in
$K^{\otimes\lambda}$ we take the tensor product of forms, whereas
$K^{\lambda}$ is obtained by taking the wedge product. The latter can
be recovered from the former by projecting to the exterior product.
Let us denote this projection by $S$ and $\bar{S}$, for projection to the
holomorphic and anti-holomorphic exterior products respectively; then
$$S\bar{S}(w_{\lambda}(\Phi))=-{\mathrm Tr}(K^{\lambda})\in\Omega^{0,\lambda}(X,{\Lambda}^{\lambda}T^*)=\Omega^{\lambda,\lambda}(X).$$
We will reorder the indices to separate the holomorphic and
anti-holomorphic forms; note that ${\mathrm
Tr}(K^{\lambda})\in\Omega^{\lambda,\lambda}(X)$ involves taking the
wedge product of $\lambda$ copies of $K$, so with respect to local
coordinates it will look like (some coefficient multiplied by)
$$dz_{k_1}\wedge d\bar{z}_{l_1}\wedge dz_{k_2}\wedge
d\bar{z}_{l_2}\wedge\cdots\wedge dz_{k_{\lambda}}\wedge
d\bar{z}_{l_{\lambda}}.$$
We can rewrite this as
$$(-1)^{\lambda/2}dz_{k_1}\wedge dz_{k_2}\wedge\cdots\wedge
dz_{k_{\lambda}}\wedge d\bar{z}_{l_1}\wedge
d\bar{z}_{l_2}\wedge\cdots\wedge d\bar{z}_{l_{\lambda}}$$
and hence 
$$S\bar{S}(w_{\lambda}(\Phi))_{k_1\cdots
k_{\lambda}\bar{l}_1\cdots\bar{l}_{\lambda}}=-(-1)^{\lambda/2}{\mathrm
Tr}(K^{\lambda})_{k_1\bar{l}_1\cdots
k_{\lambda}\bar{l}_{\lambda}}.$$

Suppose that $\Gamma$ is the necklace graph with $k$ beads
$\Theta_k$. For example, $\Theta_4$ looks like:
$$\eightVgraphI$$
This clearly contains an $8$-wheel. Indeed, the remainder of the graph
is simply four connecting edges:
\begin{center}
\begin{picture}(100,100) (50,50)
\put(100,100){\circle{40}}
\put(120,100){\line(1,0){10}}
\put(114,114){\line(1,1){16}}
\put(100,120){\line(0,1){10}}
\put(86,114){\line(-1,1){16}}
\put(80,100){\line(-1,0){10}}
\put(86,86){\line(-1,-1){16}}
\put(100,80){\line(0,-1) {10}}
\put(114,86){\line(1,-1) {16}}
\put(85,135){\oval(30,30)[t]}
\put(115,65){\oval(30,30)[b]}
\put(65,85){\oval(30,30)[l]}
\put(135,115){\oval(30,30)[r]}
\end{picture}
\end{center}
More generally, $\Theta_k$ is made up of a $2k$-wheel and $k$ edges
connecting adjacent spokes. To calculate the Rozansky-Witten invariant
$b_{\Theta_k}(X)$ we construct
$\Theta_k(\Phi)\in\Omega^{0,2k}(X)$. From our previous discussion, we
know that this is obtained from
$$w_{2k}(\Phi)=-{\mathrm Tr}(K^{\otimes 2k})\in C^{\infty}(X,(T^*\otimes\bar{T}^*)^{\otimes 2k})$$
by contracting with the $k$ copies of $\tilde{\omega}$ corresponding
to the $k$ connecting edges 
$$-{\mathrm Tr}(K^{\otimes 2k})_{k_1\bar{l}_1\cdots
k_{2k}\bar{l}_{2k}}\omega^{k_1k_2}\cdots\omega^{k_{2k-1}k_{2k}}d\bar{z}_{l_1}\otimes\cdots\otimes
d\bar{z}_{l_{2k}}$$
and then taking the projection $\bar{S}$ to the exterior product
$\Omega^{0,2k}(X)$. As always, we must be careful with the
orientation, but recall that the orientation was determined by the
element 
$$-(f_1\wedge f_2\wedge\cdots\wedge f_{2k})\wedge (f_{2k+1}\wedge
f_{4k+2})\wedge (f_{2k+2}\wedge f_{4k+3})\wedge\cdots\wedge (f_{4k}\wedge f_{4k+1})$$
which implies that the connecting edges should be oriented from $v_1$
to $v_2$, $\ldots,$ from $v_{2k-1}$ to $v_{2k}$ (the minus sign
already appears in our formula.) In other words, our use of
$\omega^{k_1k_2}$, $\ldots,$ $\omega^{k_{2k-1}k_{2k}}$ agrees with the
orientation.

Now suppose that instead of the connecting edges joining $v_1$ to
$v_2$, $\ldots,$ $v_{2k-1}$ to $v_{2k}$, they join $v_{\pi(1)}$ to
$v_{\pi(2)}$, $\ldots,$ $v_{\pi(2k-1)}$ to $v_{\pi(2k)}$ for some
permutation $\pi$ of $2k$ elements. Call the graph obtained in this
way $\Gamma_{\pi}$ (note that $\Gamma_{\mathrm Id}=\Theta_k$). Since
$$f_{\pi(1)}\wedge f_{\pi(2)}\wedge\cdots\wedge f_{\pi(2k)}=({\mathrm
sgn}\pi)f_1\wedge f_2\wedge\cdots\wedge f_{2k}$$
it follows that the orientation on $\Gamma_{\pi}$ induced from the
planar embedding is $-({\mathrm sgn}\pi)$ times the orientation given by
orienting the edges from $v_{\pi(1)}$ to $v_{\pi(2)}$, $\ldots,$ from
$v_{\pi(2k-1)}$ to $v_{\pi(2k)}$, where we assume that the vertices
are labelled $1$ to $2k$ in an anti-clockwise manner and the edges of
the wheel are oriented in an anti-clockwise manner, as before. It
follows that $\Gamma_{\pi}(\Phi)\in\Omega^{0,2k}(X)$ is given by
taking the projection $\bar{S}$ of
$$-({\mathrm sgn}\pi){\mathrm Tr}(K^{\otimes 2k})_{k_1\bar{l}_1\cdots
k_{2k}\bar{l}_{2k}}\omega^{k_{\pi(1)}k_{\pi(2)}}\cdots\omega^{k_{\pi(2k-1)}k_{\pi(2k)}}d\bar{z}_{l_1}\otimes\cdots\otimes
d\bar{z}_{l_{2k}}$$
to the exterior product $\Omega^{0,2k}(X)$.

Let 
$$\Gamma=\sum_{\pi\in{\cal S}_{2k}}\Gamma_{\pi}$$
where we have summed over all permutations of $2k$ objects. Each term
in this sum actually occurs $2^kk!$ times, because we can permute the
$k$ connecting edges in $k!$ ways and each one can be used to join two
spokes in two different ways (ie.\ the first spoke to the second, or
the second spoke to the first). We denote the graph homology class
given by summing over all graphs obtained by joining the spokes of
$w_{2k}$ in pairs by $\langle w_{2k}\rangle$, and call it the {\em
closure\/} of the wheel $w_{2k}$. Then $\Gamma=2^kk!\langle
w_{2k}\rangle$.

From what we have shown above, $\Gamma(\Phi)\in\Omega^{0,2k}(X)$ is
given by taking the projection $\bar{S}$ of
$$-\sum_{\pi\in{\cal S}_{2k}}({\mathrm sgn}\pi){\mathrm Tr}(K^{\otimes 2k})_{k_1\bar{l}_1\cdots
k_{2k}\bar{l}_{2k}}\omega^{k_{\pi(1)}k_{\pi(2)}}\cdots\omega^{k_{\pi(2k-1)}k_{\pi(2k)}}d\bar{z}_{l_1}\otimes\cdots\otimes
d\bar{z}_{l_{2k}}$$
to the exterior product. In constructing the Rozansky-Witten invariant
$b_{\Gamma}(X)$, we next multiply $\Gamma(\Phi)$ by $\omega^k$
to obtain an element of $\Omega^{2k,2k}(X)$. So we should consider
$$-\!\!\sum_{\pi\in{\cal S}_{2k}}\!({\mathrm sgn}\pi)\bar{S}({\mathrm Tr}(K^{\otimes 2k}))_{k_1\bar{l}_1\cdots
k_{2k}\bar{l}_{2k}}\omega^{k_{\pi(1)}k_{\pi(2)}}\cdots\omega^{k_{\pi(2k-1)}k_{\pi(2k)}}\omega^kd\bar{z}_{l_1}\wedge\cdots\wedge
d\bar{z}_{l_{2k}}.$$
Since this is a skewed-sum over all permutations of $2k$ objects, we
can insert the projection $S$ to the (holomorphic) exterior product in
front of $\bar{S}({\mathrm Tr}(K^{\otimes 2k}))$ and this will not
change the formula. It does mean, however, that $S\bar{S}({\mathrm
Tr}(K^{\otimes 2k}))$ can be replaced by ${\mathrm Tr}(K^{2k})$. We
also reorder the indices, introducing the appropriate sign adjustment,
to get
$$(\!-1\!)^{k+1}\!\!\!\sum_{\pi\in{\cal S}_{2k}}\!\!({\mathrm sgn}\pi){\mathrm
Tr}(K^{2k})_{k_1\cdots k_{2k}\bar{l}_1\cdots
\bar{l}_{2k}}\omega^{k_{\pi(1)}k_{\pi(2)}}\!\cdots\!\omega^{k_{\pi(2k-1)}k_{\pi(2k)}}\omega^kd\bar{z}_{l_1}\wedge\!\cdots\!\wedge
d\bar{z}_{l_{2k}}.$$

If $a_1,\ldots,a_{2k}$ are (local) holomorphic vector fields on $X$
then the algebra of exterior products tells us that
$$\omega^k(a_1,\ldots,a_{2k})=\frac{1}{(2k)!}\sum_{\pi\in{\cal
S}_{2k}}({\mathrm
sgn}\pi)\omega(a_{\pi(1)},a_{\pi(2)})\cdots\omega(a_{\pi(2k-1)},a_{\pi(2k)})$$
where the sum is over all permutations of $2k$ objects. Dually, if
$\alpha_1,\ldots,\alpha_{2k}$ are (local) holomorphic sections of the
cotangent bundle $T^*$ of $X$, then
$$\alpha_1\wedge\cdots\wedge\alpha_{2k}=\frac{1}{2^{2k}(k!)^2}\sum_{\pi\in{\cal
S}_{2k}}({\mathrm
sgn}\pi)\tilde{\omega}(\alpha_{\pi(1)},\alpha_{\pi(2)})\cdots\tilde{\omega}(\alpha_{\pi(2k-1)},\alpha_{\pi(2k)})\omega^k.$$
By writing ${\mathrm Tr}(K^{2k})$ in local coordinates, we can show that
\begin{eqnarray*}
2^{2k}(k!)^2{\mathrm Tr}(K^{2k})\hspace*{-22mm} & & \\
 & = & \sum_{\pi\in{\cal S}_{2k}}\!({\mathrm
sgn}\pi){\mathrm Tr}(K^{2k})_{k_1\cdots k_{2k}\bar{l}_1\cdots
\bar{l}_{2k}}\omega^{k_{\pi(1)}k_{\pi(2)}}\cdots\omega^{k_{\pi(2k-1)}k_{\pi(2k)}}\omega^kd\bar{z}_{l_1}\wedge\cdots\wedge
d\bar{z}_{l_{2k}}.
\end{eqnarray*}
Substituting this into the above formula we get
$$\Gamma(\Phi)\omega^k=(-1)^{k+1}2^{2k}(k!)^2{\mathrm
Tr}(K^{2k})\in\Omega^{2k,2k}(X)$$
and therefore
\begin{eqnarray*}
b_{\Gamma}(X) & = & \frac{1}{(8\pi^2)^kk!}\int_X\Gamma(\Phi)\omega^k \\
 & = & \frac{-(-1)^k2^{2k}(k!)^2}{(8\pi^2)^kk!}\int_X{\mathrm
 Tr}(K^{2k}) \\
 & = & -2^kk!\frac{1}{(2\pi i)^{2k}}\int_X {\mathrm Tr}(K^{2k})
 \\ 
 & = & -2^kk!s_{2k}(X).
\end{eqnarray*}
In particular
\begin{eqnarray*}
s_{2k}(X) & = & \frac{-1}{2^kk!}b_{\Gamma}(X) \\
 & = & -b_{\langle w_{2k}\rangle}(X)
\end{eqnarray*}
gives a way of expressing the characteristic number $s_{2k}(X)$ in
terms of the Rozansky-Witten invariants.

Suppose we wish to write the characteristic number
$$s_{\lambda_1}\cdots s_{\lambda_j}(X)=\frac{1}{(2\pi
i)^{2k}}\int_X{\mathrm Tr}(K^{\lambda_1})\wedge\cdots\wedge {\mathrm
Tr}(K^{\lambda_j})$$
in a similar way, where $(\lambda_1,\ldots,\lambda_j)$ is an even
partition of $2k$. Following the calculations presented above, we
arrive at the following formula.
\begin{prp}
The Chern numbers of a compact hyperk{\"a}hler manifold can
all be expressed in terms of the Rozansky-Witten invariants as
$$s_{\lambda_1}\cdots s_{\lambda_j}(X)=(-1)^jb_{\langle
w_{\lambda_1}\cdots w_{\lambda_j}\rangle}(X)$$
where $\langle w_{\lambda_1}\cdots w_{\lambda_j}\rangle$ denotes the
sum of all graphs obtained by taking the disjoint union of a
$\lambda_1$-wheel, a $\lambda_2$-wheel, $\ldots,$ and a
$\lambda_j$-wheel and joining their spokes in pairs. 
\end{prp}
Note that this sum of graphs includes examples where two spokes of the
same wheel may be joined. We call this sum the {\em closure\/} of the
disjoint union of the $j$ wheels $w_{\lambda_1},\ldots,w_{\lambda_j}$,
though more often we shall call it a {\em polywheel\/}. The only
possibly mysterious part in this formula is the sign $(-1)^j$. Recall
that our wheels have a canonical orientation given by their planar
embedding, but in constructing the Rozansky-Witten invariants it was
more convenient to reverse this, and hence the minus sign in our
formula for $s_{2k}(X)$. Similarly, with $\langle w_{\lambda_1}\cdots
w_{\lambda_j}\rangle$ we need to reverse the orientation of all $j$
wheels, and hence the sign $(-1)^j$.

Any characteristic number can be expressed as a linear combination of
the Chern numbers $s_{\lambda_1}\cdots s_{\lambda_j}(X)$, and
therefore as a linear combination of Rozansky-Witten invariants
corresponding to polywheels $\langle w_{\lambda_1}\cdots
w_{\lambda_j}\rangle$. In Appendix A$.1$ we expand out the polywheels
with up to ten vertices to see which graph homology classes they
give. For two, four, and six vertices it is possible to invert
these relations, so that it is possible to write all graph homology
classes as linear combinations of polywheels (see Appendix A$.2$). It
follows that for $k=1$, $2$, and $3$ all Rozansky-Witten invariants
can be expressed in terms of Chern numbers. 

For larger degree $k$ it is not possible to invert the polywheel
relations. To see why consider the dimension of ${\cal
A}(\emptyset)^k_{\mathrm conn}$ (the connected graph homology
classes with $2k$ vertices). For degree less than or equal to three
this dimension is one, and so the dimension of ${\cal A}(\emptyset)^k$
is precisely $p(k)$, the number of partitions of $k$. This is the same
as the number of polywheels. On the other hand, for degree greater
than or equal to four the dimension of ${\cal A}(\emptyset)^k_{\mathrm
conn}$ is greater than one, and so there are more than $p(k)$ distinct
graph homology classes, but still only $p(k)$ polywheels. In general,
the behaviour of ${\mathrm dim}{\cal A}(\emptyset)^k_{\mathrm conn}$
is not known, but for small $k$ we have the following sequences of
numbers.
$$
\begin{array}{|l|cccccccccc|}  
\hline
k & 1 & 2 & 3 & 4 & 5 & 6 & 7 & 8 & 9 & 10 \\
\hline
{\mathrm dim}{\cal A}(\emptyset)^k_{\mathrm conn} & 1 & 1 & 1 & 2 & 2
& 3 & 4 & 5 & 6 & 8 \\
{\mathrm dim}{\cal A}(\emptyset)^k & 1 & 2 & 3 & 6 & 9 & 16 & 25 & 42
& 65 & 105 \\
p(k) & 1 & 2 & 3 & 5 & 7 & 11 & 15 & 22 & 30 & 42 \\
\hline
\end{array}
$$
Thus the polywheels span a subspace of graph homology which is of
codimension   
$${\mathrm dim}{\cal A}(\emptyset)^k-p(k)$$
in degree $k$. We shall call this the {\em polywheel subspace\/}. In
degree four this codimension is one. It is somewhat surprising then
that $\Theta^4$ and $\Theta^2\Theta_2$ both lie in this subspace (see
Appendix A$.2$). In degree five the polywheel subspace has codimension
two, but again it contains both $\Theta^5$ and
$\Theta^3\Theta_2$. Later we will see that this is no accident, but is
a consequence of a deep result known as the Wheeling Theorem which
arose from the study of finite-type invariants of knots.

We have shown that the Rozansky-Witten invariants associated to
polywheels give the Chern numbers, and hence any graph homology class
in the polywheel subspace will give rise to a Rozansky-Witten invariant
which is a characteristic number. At this stage we are unable to
conclude anything about the graphs which lie outside the polywheel
subspace, though later we will prove that some of them give rise to
invariants which are not characteristic numbers. However, we will also
show that some of these Rozansky-Witten invariants, which appear to be
more general, are related in other ways to the characteristic numbers.

\subsection{Products of virtual hyperk{\"a}hler manifolds}

In the previous chapter we mentioned briefly that elements of graph
cohomology ${\cal A}(\emptyset)^*$ may be treated as ``virtual''
hyperk{\"a}hler manifolds, since actual hyperk{\"a}hler manifolds give
rise to such elements. We now wish to show that this association is
compatible with characteristic numbers under taking products. Firstly
let us see how the characteristic numbers behave under taking products
of actual hyperk{\"a}hler manifolds.

Let $X$ and $Y$ be compact hyperk{\"a}hler manifolds of
real-dimensions $4k$ and $4l$ respectively. The tangent bundle of
their product is given by
$$T(X\times Y)=p_1^*TX\oplus p_2^*TY$$
where $p_1$ and $p_2$ are the projections from $X\times Y$ onto $X$
and $Y$ respectively. Therefore
\begin{eqnarray*}
s(T(X\times Y)) & = & s(p_1^*TX\oplus p_2^*TY) \\
 & = & p_1^*s(TX)+p_2^*s(TY)
\end{eqnarray*}
where
$$s(T(X\times Y))=2(k+l)+s_2(T(X\times Y))+s_4(T(X\times Y))+\ldots$$
and similarly for $s(TX)$ and $s(TY)$. It follows that
\begin{eqnarray*}
s_{\lambda_1}\cdots s_{\lambda_j}(X\times Y)\hspace*{-10mm} & & \\
 & = & \int_{X\times
Y}s_{\lambda_1}(T(X\times Y))\wedge\cdots\wedge
s_{\lambda_j}(T(X\times Y)) \\
 & = & \int_{X\times
Y}(p_1^*s_{\lambda_1}(TX)+p_2^*s_{\lambda_1}(TY))\wedge\cdots\wedge
(p_1^*s_{\lambda_j}(TX)+p_2^*s_{\lambda_j}(TY)) \\
 & = & \int_{X\times Y}\sum_{S\subset
\{\lambda_1,\ldots,\lambda_j\}}\left(\Lambda_{\lambda_i\in
S}p_1^*s_{\lambda_i}(TX)\right)\wedge\left(\Lambda_{\lambda_i\not\in
S}p_2^*s_{\lambda_i}(TY)\right) \\
 & = & \sum_{S\subset\{\lambda_1,\ldots,\lambda_j\}}\int_{X\times
Y}p_1^*\left({\Lambda}_{\lambda_i\in 
S}s_{\lambda_i}(TX)\right)\wedge p_2^*\left(\Lambda_{\lambda_i\not\in
S}s_{\lambda_i}(TY)\right) \\
 & = & \sum_{S\subset
\{\lambda_1,\ldots,\lambda_j\}}\left(\int_X\Lambda_{\lambda_i\in
S}s_{\lambda_i}(TX)\right)\left(\int_Y\Lambda_{\lambda_i\not\in
S}s_{\lambda_i}(TY)\right) \\
 & = & \sum_{S\subset
\{\lambda_1,\ldots,\lambda_j\}}\left(\prod_{\lambda_i\in
S}s_{\lambda_i}(X)\right)\left(\prod_{\lambda_i\not\in
S}s_{\lambda_i}(Y)\right)
\end{eqnarray*}
where we regard $\prod_{\lambda_i\in S}s_{\lambda_i}(X)$ as being zero 
unless $\sum_{\lambda_i\in S}\lambda_i=2k$ (and similarly
$\sum_{\lambda_i\not\in S}\lambda_i=2l$). 

Let $B_k\in ({\cal A}(\emptyset)^k)^*$ and $B_l\in ({\cal
A}(\emptyset)^l)^*$ be (homogeneous) elements of graph cohomology. Recall
that the product in graph cohomology is dual to the coproduct in graph
homology. Thus if $\Gamma$ is a trivalent graph with $2(k+l)$ vertices
then
\begin{eqnarray*}
(B_kB_l)(\Gamma) & = & (B_k,B_l)(\Delta\Gamma) \\
 & = &
 \sum_{\gamma\sqcup\gamma^{\prime}=\Gamma}B_k(\gamma)B_l(\gamma^{\prime})
\end{eqnarray*}
where $B_k$ and $B_l$ give zero unless evaluated on graphs with $2k$
and $2l$ vertices, respectively. Thus $B_kB_l\in({\cal
A}(\emptyset)^{k+l})^*$. If these elements are to be regarded as
virtual hyperk{\"a}hler manifolds of real-dimensions $4k$ and $4l$
respectively, then their Rozansky-Witten invariants should be given by
$$b_{\Gamma}(B_k)=B_k(\Gamma)$$
for $\Gamma$ a graph with $2k$ vertices (and similarly for
$B_l$). Their product $B_k\times B_l$ as virtual hyperk{\"a}hler
manifolds should correspond to their product $B_kB_l$ as graph
cohomology elements. Therefore if $\Gamma$ is a trivalent graph with
$2(k+l)$ vertices then
\begin{eqnarray*}
b_{\Gamma}(B_k\times B_l) & = & (B_kB_l)(\Gamma) \\
 & = &
 \sum_{\gamma\sqcup\gamma^{\prime}=\Gamma}B_k(\gamma)B_l(\gamma^{\prime}) \\
 & = &
 \sum_{\gamma\sqcup\gamma^{\prime}=\Gamma}b_{\gamma}(B_k)b_{\gamma^{\prime}}(B_l) \\
 & = & b_{\Delta(\Gamma)}(B_k,B_l).
\end{eqnarray*}
Finally, we wish to show that this corresponds to the product formula
for characteristic numbers when we make the appropriate choice for
$\Gamma$, namely, a polywheel $\langle w_{\lambda_1}\cdots
w_{\lambda_j}\rangle$ where $(\lambda_1,\ldots,\lambda_j)$ is an even
partition of $2(k+l)$. Recall that for an actual hyperk{\"a}hler
manifold $X$ the Rozansky-Witten invariant corresponding to this
polywheel is the Chern number $(-1)^js_{\lambda_1}\cdots
s_{\lambda_j}(X)$. Suppose
$$\langle w_{\lambda_1}\cdots
w_{\lambda_j}\rangle=\gamma\sqcup\gamma^{\prime}$$
is a decomposition of $\Gamma$ into two disjoint graphs. If at least
one vertex of $w_{\lambda_i}$ lies in $\gamma$ then the whole wheel
must, and it would follow that no vertices of $w_{\lambda_i}$ could
lie in $\gamma^{\prime}$. Therefore we can find a subset $S$ of
$\{\lambda_1,\ldots,\lambda_j\}$ such that if $\lambda_i\in S$ then
the wheel $w_{\lambda_i}$ lies in $\gamma$, and otherwise it lies in
$\gamma^{\prime}$. It follows that $\gamma$ and $\gamma^{\prime}$ are
graphs occurring in 
$$\langle \cup_{\lambda_i\in S}w_{\lambda_i}\rangle\qquad\mbox{and}\qquad\langle \cup_{\lambda_i\not\in S}w_{\lambda_i}\rangle$$
respectively. Furthermore, every such pair of graphs occurring in the
above two polywheels give a decomposition of $\Gamma$, and this is
true for all possible subsets $S$ of
$\{\lambda_1,\ldots,\lambda_j\}$. Therefore
\begin{eqnarray*}
\Delta\langle w_{\lambda_1}\cdots w_{\lambda_j}\rangle & = &
\sum_{\gamma\sqcup\gamma^{\prime}=\Gamma} \gamma\otimes\gamma^{\prime}
\\
 & = & \sum_{S\subset\{\lambda_1,\ldots,\lambda_j\}}\langle \cup_{\lambda_i\in S}w_{\lambda_i}\rangle\otimes\langle \cup_{\lambda_i\not\in S}w_{\lambda_i}\rangle.
\end{eqnarray*}
Taking the corresponding Rozansky-Witten invariants and introducing
the additional signs we get
\begin{eqnarray*}
(-1)^jb_{\langle w_{\lambda_1}\cdots w_{\lambda_j}\rangle}(B_k\times
B_l)\hspace*{-2mm} & =\hspace*{-1mm} & (-1)^jb_{\Delta\langle
w_{\lambda_1}\cdots w_{\lambda_j}\rangle}(B_k,B_l) \\
 & =\hspace*{-1mm} &
\sum_{S\subset\{\lambda_1,\ldots,\lambda_j\}}(-1)^{|S|}b_{\langle
\cup_{\lambda_i\in S}w_{\lambda_i}\rangle}(B_k)(-1)^{j-|S|}b_{\langle
\cup_{\lambda_i\not\in S}w_{\lambda_i}\rangle}(B_l)
\end{eqnarray*}
where as usual we regard the terms on the right hand side as being
zero unless $\sum_{\lambda_i\in S}\lambda_i=2k$ and
$\sum_{\lambda_i\not\in S}\lambda_i=2l$. For actual hyperk{\"a}hler
manifolds this is precisely the product formula for characteristic
numbers, and therefore virtual hyperk{\"a}hler manifolds are well
behaved with respect to characteristic numbers under taking products.

\subsection{The disconnected graph $\Theta^k$}

In Chapter $4$ we will prove that $\Theta^k$ lies in the polywheel
subspace for all $k$. Thus the corresponding Rozansky-Witten
invariants $b_{\Theta^k}$ are characteristic numbers. This result
relies on the Wheeling Theorem, which arose in the study of knot
invariants. First though, let us look at some of the more basic
properties of this invariant.

We have seen that the invariant $\frac{1}{k!}b_{\Theta^k}$ is
multiplicative, so let us assume now that $X$ is irreducible. The
special property we shall need for irreducible $X$ is the following
(see Beauville~\cite{beauville83}, for example):
$${\mathrm H}^{0,2m}_{\bar{\partial}}(X)=\langle
[\bar{\omega}^m]\rangle$$
for $0\leq m\leq k$. In other words, this cohomology group is
one-dimensional and the Dolbeault class represented by the form
$\bar{\omega}^m$ gives a generator. 

If $\Gamma$ has $2m$ vertices where $m<k$ then we can still define 
$$[\Gamma(\Phi)]\in{\mathrm H}^{0,2m}_{\bar{\partial}}(X)$$
as before, and it will be some multiple $\beta_{\Gamma}$ of
$[\bar{\omega}^m]$. Note that the constant $\beta_{\Gamma}$ still
depends on $X$. Taking the disjoint union of several graphs will
simply correspond to multiplying the $\beta$'s, ie.\
$$[\gamma\gamma^{\prime}(\Phi)]=\beta_{\gamma}\beta_{\gamma^{\prime}}[\bar{\omega}^m]\in{\mathrm
H}^{0,2m}_{\bar{\partial}}(X)$$
where $2m$ is the total number of vertices in
$\gamma\gamma^{\prime}$. In particular, if $m=k$ we find
$$b_{\gamma\gamma^{\prime}}(X)=\frac{1}{(8\pi^2)^kk!}\beta_{\gamma}\beta_{\gamma^{\prime}}\int_X\bar{\omega}^k\omega^k.$$
For example
$$b_{\Theta^k}(X)=\frac{1}{(8\pi^2)^kk!}\beta_{\Theta}^k\int_X\bar{\omega}^k\omega^k.$$
Since $\omega_2$ is a K{\"a}hler form for $X$, a volume form is given
by 
\begin{eqnarray*}
\frac{1}{(2k)!}\omega_2^{2k} & = &
\frac{1}{(2k)!}(\frac{\omega+\bar{\omega}}{2})^{2k} \\
 & = & \frac{1}{2^{2k}(k!)^2}\omega^k\bar{\omega}^k.
\end{eqnarray*}
Therefore
$$\int_X\bar{\omega}^k\omega^k=2^{2k}(k!)^2{\mathrm vol}(X)$$
and hence
$$b_{\Theta^k}(X)=\frac{k!}{(2\pi^2)^k}\beta_{\Theta}^k{\mathrm vol}(X).$$
To determine $\beta_{\Theta}$ observe that
$$\beta_{\Theta}\int_X\bar{\omega}^k\omega^k=\int_X[\Theta(\Phi)][\bar{\omega}^{k-1}][\omega^k].$$
On the other hand, since $X$ is Ricci-flat we have (see Besse~\cite{besse87})
\begin{eqnarray*}
\|K\|^2 & = & \frac{8\pi^2}{(2k-2)!}\int_Xc_2\omega_2^{2k-2} \\
 & = &
 -\frac{4\pi^2}{(2k-2)!}\int_Xs_2(\frac{\omega+\bar{\omega}}{2})^{2k-2} \\
 & = &
 -\frac{4\pi^2}{2^{2k-2}((k-1)!)^2}\int_Xs_2\omega^{k-1}\bar{\omega}^{k-1}
\end{eqnarray*}
where $\|K\|^2$ is the ${\cal L}^2$-norm of the curvature, and $c_2$
is the second Chern class of the tangent bundle, which is equal to
$-s_2/2$. We wish to relate the two integrals
$$\int_Xs_2\omega^{k-1}\bar{\omega}^{k-1}\qquad\mbox{and}\qquad\int_X[\Theta(\Phi)][\bar{\omega}^{k-1}][\omega^k].$$
Observe that $\Theta$ can be obtained by joining the spokes of a
$2$-wheel. So consider the section
$$w_2(\Phi)\in C^{\infty}(X,(T^*\otimes\bar{T}^*)^{\otimes 2})$$
given by
$$w_2(\Phi)_{k_1\bar{l}_1k_2\bar{l}_2}=-{\mathrm Tr}(K^{\otimes 2})_{k_1\bar{l}_1k_2\bar{l}_2}.$$
Contracting with $\tilde{\omega}$ then taking the projection $\bar{S}$
to the exterior product we get (paying careful attention to the
orientation, as always)
$$\Theta(\Phi)=-\bar{S}{\mathrm Tr}(K^{\otimes
2})_{k_1\bar{l}_1k_2\bar{l}_2}\omega^{k_1k_2}d\bar{z}_{l_1}\wedge
d\bar{z}_{l_2}\in\Omega^{0,2}(X).$$
Since $\omega^{k_1k_2}$ is skew-symmetric, we can rewrite the right
hand side as
$$-S\bar{S}{\mathrm Tr}(K^{\otimes
2})_{k_1\bar{l}_1k_2\bar{l}_2}\omega^{k_1k_2}d\bar{z}_{l_1}\wedge
d\bar{z}_{l_2}=-{\mathrm
Tr}(K^2)_{k_1\bar{l}_1k_2\bar{l}_2}\omega^{k_1k_2}d\bar{z}_{l_1}\wedge
d\bar{z}_{l_2}.$$
Taking the corresponding cohomology classes, we find
$$[\Theta(\Phi)]=-[{\mathrm
Tr}(K^2)_{k_1\bar{l}_1k_2\bar{l}_2}\omega^{k_1k_2}d\bar{z}_{l_1}\wedge
d\bar{z}_{l_2}]\in{\mathrm H}^{0,2}_{\bar{\partial}}(X).$$
In general, for a $(2,0)$-form $\alpha=\sum\alpha_{ij}dz_i\wedge dz_j$
we have
$$\alpha\wedge\omega^{k-1}=\frac{1}{2k}(\sum\omega^{ij}\alpha_{ij})\omega^k.$$
Applying this to
\begin{eqnarray*}
{\mathrm Tr}(K^2) & = & {\mathrm
Tr}(K^2)_{k_1\bar{l}_1k_2\bar{l}_2}dz_{k_1}\wedge d\bar{z}_{l_1}\wedge
dz_{k_2}\wedge d\bar{z}_{l_2} \\
 & = & -{\mathrm Tr}(K^2)_{k_1\bar{l}_1k_2\bar{l}_2}dz_{k_1}\wedge
dz_{k_2}\wedge d\bar{z}_{l_1}\wedge d\bar{z}_{l_2}
\end{eqnarray*}
we find
$${\mathrm Tr}(K^2)\wedge\omega^{k-1}=\frac{-1}{2k}(\sum\omega^{k_1k_2}{\mathrm
Tr}(K^2)_{k_1\bar{l}_1k_2\bar{l}_2}d\bar{z}_{l_1}\wedge
d\bar{z}_{l_2})\omega^k$$
and taking cohomology classes
\begin{eqnarray*}
(2\pi i)^2s_2[\omega^{k-1}] & = & [{\mathrm Tr}(K^2)][\omega^{k-1}] \\
 & = & \frac{-1}{2k}[\sum\omega^{k_1k_2}{\mathrm
 Tr}(K^2)_{k_1\bar{l}_1k_2\bar{l}_2}d\bar{z}_{l_1}\wedge
 d\bar{z}_{l_2}][\omega^k] \\
 & = & \frac{1}{2k}[\Theta(\Phi)][\omega^k].
\end{eqnarray*}
Therefore
$$-8\pi^2ks_2[\omega^{k-1}][\bar{\omega}^{k-1}]=[\Theta(\Phi)][\omega^k][\bar{\omega}^{k-1}]$$
and hence
$$-8\pi^2k\int_Xs_2\omega^{k-1}\bar{\omega}^{k-1}=\int_X[\Theta(\Phi)][\bar{\omega}^{k-1}][\omega^k].$$
We can now determine $\beta_{\Theta}$, and we find
$$\beta_{\Theta}=\frac{1}{2k}\frac{\|K\|^2}{{\mathrm vol}(X)}$$
which gives the following proposition.
\begin{prp}
For an irreducible hyperk{\"a}hler manifold $X$ of real-dimension $4k$
$$b_{\Theta^k}(X)=\frac{k!}{(4\pi^2k)^k}\frac{\|K\|^{2k}}{{\mathrm
vol}(X)^{k-1}}.$$
\end{prp}
This formula is intricately related to the metric on the
hyperk{\"a}hler manifold; it involves both the curvature and the
volume of $X$. However, in Chapter $4$ we will see that $\Theta^k$ is
in the polywheel subspace. Thus $b_{\Theta^k}$ is a characteristic
number, which in the hyperk{\"a}hler case means that it is (somewhat
remarkably) a topological invariant of the manifold.

%% file: ch3.tex
\section{Relations in graph homology I}

\subsection{Perturbative Chern-Simons theory}

The Rozansky-Witten invariants allow us to associate an element of
graph cohomology ${\cal A}(\emptyset)^*$, the dual of graph homology,
to a hyperk{\"a}hler manifold. Other elements of ${\cal
A}(\emptyset)^*$ can tell us something about relations in graph
homology, ie.\ if trivalent graphs $\Gamma_1$ and $\Gamma_2$ are
homologous and $B\in{\cal A}(\emptyset)^*$, then  
$$B(\Gamma_1)=B(\Gamma_2).$$
This statement is completely tautologous, but we shall use it to
deduce some interesting relations. In particular, the elements of
${\cal A}(\emptyset)^*$ we shall use come from the weights
$c_{\Gamma}({\frak g})$ arising from a semisimple Lie algebra, as in
perturbative Chern-Simons theory.

The partition function of Chern-Simons theory gives us an invariant of
three-manifolds. A Feynman diagram calculation (as in Axelrod and
Singer~\cite{as92,as94}) shows that it is of
finite type and looks like
$$\sum_{\Gamma}c_{\Gamma}({\frak g})I_{\Gamma}^{\mathrm CS}(M)$$
where the weights $c_{\Gamma}({\frak g})$ depend on the Lie algebra
$\frak g$ of the gauge group and $I_{\Gamma}^{\mathrm CS}(M)$ depends
on the three-manifold $M$. More precisely, it is believed that this
sum occurs as the expansion of the contribution coming from the trivial
connection. In any case, we shall only be concerned here with the
weights $c_{\Gamma}({\frak g})$. For any quadratic Lie algebra $\frak
g$, ie.\ Lie algebra with an invariant inner product, these weights
can be rigorously defined as follows. 

Choose a basis $\{x_1,\ldots,x_n\}$ for $\frak g$ and let the
invariant inner product be $\sigma^{ij}$ with respect to this
basis. The structure constants $c_{ij}^{\phantom{ij}k}$ of $\frak g$
are defined by
$$[x_i,x_j]=c_{ij}^{\phantom{ij}k}x_k.$$
Using $\sigma^{-1}=\sigma_{ij}$ to lower indices we obtain
$$c_{ijk}=\sum_mc_{ij}^{\phantom{ij}m}\sigma_{mk}.$$
This is totally skew-symmetric since the Lie bracket is skew and the
inner product invariant.

Let $\Gamma$ be an oriented trivalent graph with $2k$ vertices. We use
the standard notion of orientation here, namely an equivalence class
of cyclic orderings of the outgoing edges at each vertex, with two
orderings being equivalent if they differ at even number of vertices.
Recall that this is equivalent to a Rozansky-Witten orientation.
Place a copy of $c_{ijk}$ at each vertex and attach the indices $ijk$
to the outgoing edges in accordance with the cyclic ordering given by
the orientation, ie.\ if we first attach $i$ to one of the edges, then
which way we attach $j$ and $k$ to the remaining two edges is
determined by the cyclic ordering. Note that this is well-defined as
$$c_{ijk}=c_{jki}=c_{kij}$$
so it does not depend on which edge we first attach $i$ to. Now use
$\sigma^{ij}$ to contract the indices along edges. Thus if two
vertices are connected by an edge, and the two ends of the edge are
labelled $i_t$ and $i_s$, then we contract with
$\sigma^{i_ti_s}$. Since $\sigma^{ij}$ is symmetric this does not
require an orientation of the edge. The resulting number is the weight
$c_{\Gamma}({\frak g})$, and it does not depend on the choice of basis
for $\frak g$.

For example, if $\Gamma=\Theta$ (with the canonical orientation coming
from drawing in the plane) then
$$c_{\Theta}({\frak g})=\sum_{i_1,j_1,k_1,i_2,j_2,k_2}c_{i_1j_1k_1}\sigma^{i_1i_2}\sigma^{j_1k_2}\sigma^{k_1j_2}c_{i_2j_2k_2}.$$
Note that we have used sub-indices $1$ and $2$ but this construction
does not depend on an ordering of the vertices.

Reversing the orientation of the graph by reversing the cyclic
ordering at one vertex changes the indices used to label the outgoing
edges at that vertex. This is equivalent to swapping two indices on a
$c_{ijk}$. Since
$$c_{ikj}=-c_{ijk}$$
this reverses the sign of $c_{\Gamma}({\frak g})$, ie.\
$$c_{\bar{\Gamma}}({\frak g})=-c_{\Gamma}({\frak g}).$$
The IHX relation also holds as a consequence of the Jacobi identity,
ie.\
$$c_{\Gamma_I}({\frak g})=c_{\Gamma_H}({\frak g})-c_{\Gamma_X}({\frak g}).$$
Thus if we linearly extend $c_{\Gamma}({\frak g})$ to rational
linear combinations of trivalent graphs, then the weights we get are
well defined on graph homology and so define an element of the dual space
${\cal A}(\emptyset)^*$. In the next subsection we shall calculate
some of these numbers explicitly when ${\frak g}={\frak{su}}(2)$.

\subsection{Gauge group ${\mathrm SU}(2)$}

In the case of Lie algebras coming from semisimple Lie groups an
invariant inner product is given by $-1/2$ times the Killing form,
which is given by taking the trace of the product in the adjoint
representation. Therefore 
$$\sigma_{ij}=-\frac{1}{2}{\mathrm Tr}({\mathrm ad}x_i{\mathrm ad}x_j).$$
The structure constants are given by the adjoint representation, ie.\
$$[x_i,x_j]=({\mathrm ad}x_i)x_j=c_{ij}^{\phantom{ij}k}x_k$$
so that
$$c_{ij}^{\phantom{ij}k}=({\mathrm ad}x_i)_{jk}$$
where the $jk$ coefficients indicate that we regard ${\mathrm ad}x_i$
as an $n\times n$ matrix. 

When ${\frak g}={\frak{su}}(2)$ we can take as a basis the $3\times 3$
matrices
$$x_1=\left(\begin{array}{ccc} 0 & 1 & 0 \\ -1 & 0 & 0 \\ 0 & 0 & 0
            \end{array}\right)\qquad 
x_2=\left(\begin{array}{ccc} 0 & 0 & -1 \\ 0 & 0 & 0 \\ 1 & 0 & 0
            \end{array}\right)\qquad 
x_3=\left(\begin{array}{ccc} 0 & 0 & 0 \\ 0 & 0 & 1 \\ 0 & -1 & 0
            \end{array}\right)$$
%\right)
written in the adjoint representation (we drop the `${\mathrm ad}$'
notation from here onwards). The Killing form is given by ${\mathrm 
Tr}(x_ix_j)$ and a calculation shows that this is $-2$ times the
identity matrix, so that the inner product is simply $\delta$. The
structure constants are given by the unique totally skew-symmetric
tensor in three-dimensions
$$c_{ijk}={\epsilon}_{ijk}$$
(ie.\ unique up to scale, but the scale here is one). More often we
shall use the $3\times 3$ matrix $x_i$ to represent
$c_{ij}^{\phantom{ij}k}$. 

Let us now calculate some of the weights arising from the
${\frak{su}}(2)$ Lie algebra, beginning with $\Gamma=\Theta$. Then
\begin{eqnarray*}
c_{\Theta}({\frak{su}}(2)) & = &
\sum_{i_1,j_1,k_1,i_2,j_2,k_2}{\epsilon}_{i_1j_1k_1}\delta^{i_1i_2}\delta^{j_1k_2}\delta^{k_1j_2}{\epsilon}_{i_2j_2k_2} \\
   & = & -\sum_{i,j,k}({\epsilon}_{ijk})^2 \nonumber \\
   & = & -6.
\end{eqnarray*}
For a disconnected graph in this theory, we simply multiply the
corresponding weights for the connected components. For example
$$c_{\Theta^k}({\frak{su}}(2))=c_{\Theta}({\frak{su}}(2))^k=(-6)^k.$$
Suppose that $\Gamma$ contains a $2$-wheel $w_2$. The contribution to
$c_{\Gamma}({\frak{su}}(2))$ from the $2$-wheel is
\begin{eqnarray*}
w_2({\frak{su}}(2)) & = & \sum_{j_1,k_1,j_2,k_2}{\epsilon}_{i_1j_1k_1}\delta^{k_1j_2}{\epsilon}_{i_2j_2k_2}\delta^{k_2j_1}
\\
 & = & -\sum_{j,k}{\epsilon}_{i_1jk}{\epsilon}_{i_2jk} \\
 & = & -2{\delta}_{i_1i_2}.
\end{eqnarray*}
Therefore
$$c_{\Gamma}({\frak{su}}(2))=-2c_{\Gamma^{\prime}}({\frak{su}}(2))$$
where $\Gamma^{\prime}$ is the graph constructed from
$\Gamma$ by removing the $2$-wheel and joining the two univalent
vertices (ie.\ the two loose ends created).

More generally, suppose that $\Gamma$ contains a $\lambda$-wheel
$w_{\lambda}$. The contribution to $c_{\Gamma}({\frak{su}}(2))$ from
this part of the graph would be 
\begin{eqnarray*}
w_{\lambda}({\frak{su}}(2)) & = &
\sum_{j_1,k_1,\ldots,j_{\lambda},k_{\lambda}}{\epsilon}_{i_1j_1k_1}\delta^{k_1j_2}{\epsilon}_{i_2j_2k_2}\delta^{k_2j_3}\cdots
{\epsilon}_{i_{\lambda}j_{\lambda}k_{\lambda}}\delta^{k_{\lambda}j_1} \\ 
 & = &
 \sum_{j_1,\ldots,j_{\lambda}}c_{i_1j_1}^{\phantom{i_1j_1}j_2}c_{i_2j_2}^{\phantom{i_2j_2}j_3}\cdots c_{i_{\lambda}j_{\lambda}}^{\phantom{i_{\lambda}j_{\lambda}}j_1} \\ 
 & = & \sum_{j_1,\ldots,j_{\lambda}}(x_{i_1})_{j_1j_2}(x_{i_2})_{j_2j_3}\cdots (x_{i_{\lambda}})_{j_{\lambda}j_1} \\
 & = & {\mathrm Tr}(x_{i_1}x_{i_2}\cdots x_{i_{\lambda}})
\end{eqnarray*}
where the indices $i_1,\ldots,i_{\lambda}$ label the spokes of the
wheel. Thus if $\Gamma$ is the polywheel $\langle w_{\lambda_1}\cdots
w_{\lambda_j}\rangle$ we find that
$$c_{\langle w_{\lambda_1}\cdots
w_{\lambda_j}\rangle}({\frak{su}}(2))=\sum {\mathrm
Tr}(x_{i_{11}}x_{i_{12}}\cdots
x_{i_{1{\lambda_1}}})\cdots{\mathrm Tr}(x_{i_{j1}}x_{i_{j2}}\cdots
x_{i_{j\lambda_j}})$$
where the sum is over all ways of contracting together the spokes
using the inner product, or since this is $\delta$, simply pairing the
indices and summing. We can assume all $\lambda$'s are even as
otherwise the polywheel vanishes. Let $\lambda_1+\ldots
+\lambda_j=2k$, so that the polywheel $\langle w_{\lambda_1}\cdots
w_{\lambda_j}\rangle$ has $2k$ vertices. By considering what happens
when we join two spokes, we will express the above weight in terms of
weights of polywheels with fewer vertices, thus arriving at a
recursive formula. In order to proceed we need the following two
identities.
\begin{lem}
Let $A$ be a $3\times 3$ matrix. Then
$$\sum_i x_iAx_i=A^t-({\mathrm Tr}A)I.$$
\end{lem}
\begin{lem}
Let $A$ and $B$ be $3\times 3$ matrices. Then
$$\sum_i{\mathrm Tr}(Ax_i){\mathrm Tr}(Bx_i)=-\frac{1}{2}{\mathrm
Tr}((A-A^t)(B-B^t)).$$
\end{lem}
\begin{prf}
Simply let 
$$A=\left(\begin{array}{ccc} a & b & c \\ d & e & f \\ g & h & i
      \end{array}\right)\qquad\mbox{and}\qquad 
B=\left(\begin{array}{ccc} m & n & o \\ p & q & r \\ s & t & u
      \end{array}\right)$$ 
then expand both sides of the relations.
\end{prf}
The first lemma will describe the effect of joining two spokes on the
same wheel, and the second will describe the effect of joining two
spokes on different wheels.

So first consider a factor in $c_{\langle w_{\lambda_1}\cdots
w_{\lambda_j}\rangle}({\frak{su}}(2))$ which looks like
$$\sum_k{\mathrm Tr}(x_{i_1}\cdots x_{i_{a-1}}x_kx_{i_{a+1}}\cdots
x_{i_{a+b-1}}x_kx_{i_{a+b+1}}\cdots x_{i_{\lambda}})$$
where $1\leq b\leq{\lambda}-1$ and $1\leq a\leq\lambda -b$. The first
lemma allows us to rewrite this as two terms
\begin{eqnarray*}
{\mathrm Tr}(x_{i_1}\cdots x_{i_{a-1}}(x_{i_{a+1}}\cdots
x_{i_{a+b-1}})^tx_{i_{a+b+1}}\cdots x_{i_{\lambda}})\hspace*{-40mm} & & \\
 & = & (-1)^{b-1}{\mathrm Tr}(x_{i_1}\cdots x_{i_{a-1}}x_{i_{a+b-1}}\cdots
x_{i_{a+1}}x_{i_{a+b+1}}\cdots x_{i_{\lambda}})
\end{eqnarray*}
and
$$-{\mathrm Tr}(x_{i_1}\cdots x_{i_{a-1}}x_{i_{a+b+1}}\cdots
x_{i_{\lambda}}){\mathrm Tr}(x_{i_{a+1}}\cdots x_{i_{a+b-1}}).$$
Substituting these two terms back into the expression for
$c_{\langle w_{\lambda_1}\cdots w_{\lambda_j}\rangle}({\frak{su}}(2))$
we find that the first gives us
$$(-1)^{b-1}c_{\langle w_{\lambda_1}\cdots w_{\lambda_s-2}\cdots
w_{\lambda_j}\rangle}({\frak{su}}(2))$$
and the second gives
$$-c_{\langle w_{\lambda_1}\cdots
w_{\lambda_{s-1}}w_{\lambda_s-b-1}w_{b-1}w_{\lambda_{s+1}}\cdots w_{\lambda_j}\rangle}({\frak{su}}(2))$$
where we originally joined two spokes of the $s$th wheel, and the
second expression vanishes unless $b$ is odd. Summing from $a=1$ to
$\lambda_s-b$ and $b=1$ to $\lambda_s-1$, the first expression gives
us
$$\frac{\lambda_s}{2}c_{\langle w_{\lambda_1}\cdots w_{\lambda_s-2}\cdots
w_{\lambda_j}\rangle}({\frak{su}}(2))$$
and summing from $a=1$ to $\lambda_s-b$ and rewriting $b$ as $l+1$,
the second expression gives us
$$-\sum_{l=0,\mbox{\footnotesize{even}}}^{\lambda_s-2}(l+1)c_{\langle w_{\lambda_1}\cdots
w_{\lambda_{s-1}}w_{\lambda_s-l-2}w_{l}w_{\lambda_{s+1}}\cdots w_{\lambda_j}\rangle}({\frak{su}}(2)).$$
We also need to sum from $s=1$ to $j$ but we shall do this later.

Now consider what happens when we join spokes from different
wheels. There are $\lambda_1$ choices for the spoke on the first wheel
and $\lambda_2$ choices for the spoke on the second wheel, and each
term looks like
$$\sum_k{\mathrm Tr}(x_{i_{11}}x_{i_{12}}\cdots
x_{i_{1(\lambda_1-1)}}x_k){\mathrm Tr}(x_{i_{21}}x_{i_{22}}\cdots
x_{i_{2(\lambda_2-1)}}x_k).$$
By the second lemma we can rewrite this as the sum of four terms, all
of which give the same expression (after some relabelling of indices)
$$-\frac{1}{2}{\mathrm Tr}(x_{i_{11}}x_{i_{12}}\cdots
x_{i_{1(\lambda_1-1)}}x_{i_{21}}x_{i_{22}}\cdots
x_{i_{2(\lambda_2-1)}}).$$
Substituting this back into the expression for $c_{\langle
w_{\lambda_1}\cdots w_{\lambda_j}\rangle}({\frak{su}}(2))$ we get
$$-2c_{\langle w_{\lambda_1}\cdots w_{\lambda_{s-1}}w_{\lambda_{s+1}}\cdots
w_{\lambda_{t-1}}w_{\lambda_{t+1}}\cdots
w_{\lambda_j}w_{\lambda_s+\lambda_t-2}\rangle}({\frak{su}}(2))$$
where we originally joined spokes coming from the $s$th and $t$th
wheels. Recall that there are $\lambda_s\lambda_t$ such terms, coming
from the different ways of choosing the spokes on each wheel.

Combining the above results we can obtain a recursion relation for the
weight $c_{\langle w_{\lambda_1}\cdots w_{\lambda_j}\rangle}({\frak{su}}(2))$. It remains to note that a given graph occurring in
$\langle w_{\lambda_1}\cdots w_{\lambda_j}\rangle$ has $k$ pairs of
spokes joined together; on the other hand, we've considered what
happens when we join just one pair of spokes, for which there are $k$
choices. Thus we've counted everything $k$ times, and hence the
relation is
\begin{eqnarray*}
kc_{\langle w_{\lambda_1}\cdots w_{\lambda_j}\rangle}({\frak{su}}(2))\hspace*{-2mm} &
=\hspace*{-2mm} & \sum_s\frac{\lambda_s}{2}c_{\langle w_{\lambda_1}\cdots w_{\lambda_s-2}\cdots
w_{\lambda_j}\rangle}({\frak{su}}(2)) \\
 & & -\sum_s\sum_{l=0,\mbox{\footnotesize{even}}}^{\lambda_s-2}(l+1)c_{\langle w_{\lambda_1}\cdots
w_{\lambda_{s-1}}w_{\lambda_s-l-2}w_{l}w_{\lambda_{s+1}}\cdots
w_{\lambda_j}\rangle}({\frak{su}}(2)) \\
 & & -\sum_{\mbox{\footnotesize{pairs}}\, s,t}2\lambda_s\lambda_tc_{\langle w_{\lambda_1}\cdots w_{\lambda_{s-1}}w_{\lambda_{s+1}}\cdots
w_{\lambda_{t-1}}w_{\lambda_{t+1}}\cdots
w_{\lambda_j}w_{\lambda_s+\lambda_t-2}\rangle}({\frak{su}}(2)).
\end{eqnarray*}
The initial conditions are
\begin{enumerate}
\item $c_{\emptyset}({\frak{su}}(2))=1$,
\item $c_{\langle w_{\lambda_1}\cdots
w_{\lambda_j}w_0\rangle}({\frak su}(2))=3c_{\langle
w_{\lambda_1}\cdots w_{\lambda_j}\rangle}({\frak{su}}(2))$.
\end{enumerate}
The second of these follows from the fact that we regard the wheel
with no spokes, $w_0$, as a closed loop; thus it contributes a factor
of
$${\mathrm Tr}({\mathrm Id}_3)=3.$$
Next we solve the recursion relation.
\begin{prp}
Let $\lambda_1,\ldots,\lambda_j$ be positive even integers whose sum
is $2k$. Then
\begin{enumerate}
\item $c_{\langle w_{\lambda_1}\cdots
w_{\lambda_j}\rangle}({\frak{su}}(2))=2^{j-1}c_{\langle
w_{2k}\rangle}({\frak{su}}(2))$,
\item $c_{\langle w_{2k}\rangle}({\frak{su}}(2))=\frac{(-1)^k(2k+1)!}{2^{k-1}k!}$.
\end{enumerate}
\end{prp}
\begin{prf}
We will use induction on $k$, the first case ($k=j=1$, $\lambda_1=2$)
being trivial. Suppose the results are true up to level $k$, and
consider the next level $k+1$. Suppose that all $\lambda$'s are greater
than $2$. Using the recursion relation and applying the inductive
hypothesis to the right hand side we find
\begin{eqnarray*}
(k+1)c_{\langle w_{\lambda_1}\cdots w_{\lambda_j}\rangle}({\frak{su}}(2))\hspace*{-10mm} & & \\
 & = & \sum_s\frac{\lambda_s}{2}2^{j-1}c_{\langle
w_{2k}\rangle}({\frak{su}}(2))-\sum_s\sum_{l=2,\mbox{\footnotesize{even}}}^{\lambda_s-4}(l+1)2^jc_{\langle
w_{2k}\rangle}({\frak{su}}(2)) \\  
 & & -\sum_sc_{\langle w_{\lambda_1}\cdots
w_{\lambda_{s-1}}w_{\lambda_s-2}w_0w_{\lambda_{s+1}}\cdots
w_{\lambda_j}\rangle}({\frak{su}}(2)) \\
 & & -\sum_s(\lambda_s-1)c_{\langle w_{\lambda_1}\cdots
w_{\lambda_{s-1}}w_0w_{\lambda_s-2}w_{\lambda_{s+1}}\cdots
w_{\lambda_j}\rangle}({\frak{su}}(2)) \\
 & & -\sum_{\mbox{\footnotesize{pairs}}\, s,t}2\lambda_s\lambda_t2^{j-2}c_{\langle
w_{2k}\rangle}({\frak{su}}(2)) \\
 & = & \sum_s\lambda_s2^{j-2}c_{\langle w_{2k}\rangle}({\frak{su}}(2))-\sum_s((\frac{\lambda_s-2}{2})^2-1)2^jc_{\langle w_{2k}\rangle}({\frak{su}}(2)) \\  
 & & -\sum_s3\lambda_s2^{j-1}c_{\langle w_{2k}\rangle}({\frak{su}}(2))-\sum_{s\neq t}\lambda_s\lambda_t2^{j-2}c_{\langle
w_{2k}\rangle}({\frak{su}}(2)) \\
 & = & (\sum_s(\lambda_s-\lambda_s^2+4\lambda_s-6\lambda_s)-\sum_{s\neq
t}\lambda_s\lambda_t)2^{j-2}c_{\langle w_{2k}\rangle}({\frak{su}}(2))
\\
 & = & -(\sum_s\lambda_s +(\sum_s\lambda_s)^2)2^{j-2}c_{\langle w_{2k}\rangle}({\frak{su}}(2))
\\
 & = & -(k+1)(2k+3)2^{j-1}c_{\langle w_{2k}\rangle}({\frak{su}}(2)).
\end{eqnarray*}
When $j=1$ this tells us
\begin{eqnarray*}
(k+1)c_{\langle w_{2(k+1)}\rangle}({\frak{su}}(2)) & = &
-(k+1)(2k+3)c_{\langle w_{2k}\rangle}({\frak{su}}(2)) \\
 & = & (k+1)\frac{(-1)^{k+1}(2k+3)!}{2^k(k+1)!}
\end{eqnarray*}
which proves the second part of the proposition. It also shows that
$$(k+1)c_{\langle w_{\lambda_1}\cdots w_{\lambda_j}\rangle}({\frak{su}}(2))=(k+1)2^{j-1}c_{\langle w_{2(k+1)}\rangle}({\frak{su}}(2))$$
which proves the first part, though we assumed that all $\lambda$'s
were greater than $2$. If this assumption is not satisfied there are
some additional $0$-wheels occurring in the formulae of the proof but
the proposition is still satisfied.
\end{prf}
These calculations of the ${\frak{su}}(2)$ weight system are of use in
deriving relations in graph homology, for example those occurring in
Appendix A. In the next subsection we shall describe how this process
works.

\subsection{Graph homology relations}

The ${\frak{su}}(2)$ weight system gives us an element of the dual
space ${\cal A}(\emptyset)^*$. Equivalently, by breaking this up into
separate degrees, we get homogeneous elements $C_k$ in each of the
dual spaces $({\cal A}(\emptyset)^k)^*$. Recall that since graph
homology has a Hopf algebra structure we can multiply these $C_k$'s,
the product in graph cohomology being the dual of the coproduct in
graph homology. If $\Gamma$ has $2(k+l)$ vertices, then $C_kC_l\in({\cal
A}(\emptyset)^{k+l})^*$ is given by
\begin{eqnarray*}
(C_kC_l)(\Gamma) & = & (C_k,C_l)(\Delta\Gamma) \\
 & = &
 \sum_{\gamma\sqcup\gamma^{\prime}=\Gamma}C_k(\gamma)C_l(\gamma^{\prime})
\end{eqnarray*}
where $C_k$ and $C_l$ give zero unless evaluated on graphs with $2k$
and $2l$ vertices, respectively. More generally we can form
$$C_{m_1}C_{m_2}\cdots C_{m_j}\in({\cal A}(\emptyset)^{m_1+\ldots +m_j})^*.$$

We know that the number of polywheels with $2k$ vertices is equal to
$p(k)$, the number of partitions of $k$. Suppose a graph $\Gamma$ with
$2k$ vertices belongs to the polywheel subspace. Then we can write
$$\Gamma=\sum_{(m_1,\ldots,m_j)=\,\mbox{\footnotesize{partition
of}}\,k}a_{(m_1,\ldots,m_2)}\langle w_{2m_1}\cdots w_{2m_j}\rangle.$$
To determine the coefficients $a_{(m_1,\ldots,m_j)}$ we simply need to
act on both sides with $C_{n_1}\cdots C_{n_j}$ for all partitions
$(n_1,\ldots,n_j)$ of $k$ then solve the set of linear equations. 

For example, take
$$\Gamma=\Theta^2=a_{(1,1)}\langle w_2w_2\rangle +a_{(2)}\langle
w_4\rangle.$$
Acting with $C_2$ and $C_1C_1$ gives
\begin{eqnarray*}
(-6)^2 & = & 60a_{(1,1)}+30a_{(2)} \\
2(-6)^2 & = & 2(-6)^2a_{(1,1)}
\end{eqnarray*}
respectively, and hence solving we find
$$\Theta^2=\langle w_2w_2\rangle -\frac{4}{5}\langle w_4\rangle.$$

An equivalent formulation of the above procedure is given by regarding
the elements $C_k$ as virtual hyperk{\"a}hler manifolds of
real-dimension $4k$, as described in the previous chapters. The
Rozansky-Witten invariants of $C_k$ are given by
$$b_{\Gamma}(C_k)=C_k(\Gamma)=c_{\Gamma}({\frak{su}}(2))$$
for trivalent graphs with $2k$ vertices. The polywheels correspond to
Chern numbers
\begin{eqnarray*}
s_{\lambda_1}\cdots s_{\lambda_j}(C_k) & = & (-1)^jc_{\langle
w_{\lambda_1}\cdots w_{\lambda_j}\rangle}({\frak{su}}(2)) \\
 & = & \frac{(-1)^{k+j}(2k+1)!}{2^{k-j}k!}
\end{eqnarray*}
and the product $C_kC_l\in({\cal A}(\emptyset)^{k+l})^*$ corresponds
to the product manifold $C_k\times C_l$ of dimension $4(k+l)$.

If $\Gamma$ lies in the polywheel subspace we know that the
corresponding Rozansky-Witten invariant will be a characteristic
number, ie.\ a linear combination of Chern numbers. Let us suppose
this is the case for $\Gamma=\Theta^k$, and hence
$\frac{1}{48^kk!}b_{\Theta^k}$ is a characteristic number. We have
included the extra factor of $k!$ because this makes the invariant
multiplicative (see Subsection $1.6$); the power of $48$ is simply to make
it more convenient to state our end result. Following Hirzebruch's
work on multiplicative sequences (see~\cite{hirzebruch78}), this means
that the generating sequence
$$g(x)=1+\frac{1}{48}b_{\Theta}x^2+\frac{1}{48^22!}b_{\Theta^2}x^4+\frac{1}{48^33!}b_{\Theta^3}x^6+\ldots$$
is determined by a single power series $f(x)$ which enables us to
write each $b_{\Theta^k}$ as a characteristic number. More
specifically, in $4k$ dimensions let $\gamma_1,\ldots,\gamma_{2k}$ be
Chern roots, so that $s_{\lambda}=\gamma_1^{\lambda}+\ldots
+\gamma_{2k}^{\lambda}$. Then
$$g(x)=f(x\gamma_1)\cdots f(x\gamma_{2k})$$
and by expanding this out and writing the coefficients in terms of
$s_{\lambda}$ we get precise expressions for $b_{\Theta^k}$ in terms
of Chern numbers. It is convenient to write
$${\mathrm ln}f(x)=a_0+a_1x+a_2x^2+a_3x^3+a_4x^4+\ldots$$
and hence
$${\mathrm
ln}g(x)=2ka_0+a_1s_1x+a_2s_2x^2+a_3s_3x^3+a_4s_4x^4+\ldots$$
Since we are in the hyperk{\"a}hler case, all odd coefficients vanish,
and $a_0$ must also vanish so that the series is independent of $k$. Then
\begin{eqnarray*}
g(x) & = & {\mathrm exp}(a_2s_2x^2+a_4s_4x^4+a_6s_6x^6+\ldots) \\
 & = & 1+a_2s_2x^2+(a_4s_4+\frac{1}{2}a_2^2s_2^2)x^4+(a_6s_6+a_2a_4s_2s_4+\frac{1}{6}a_2^3s_2^3)x^6+\ldots
\end{eqnarray*}
The coefficients $a_i$ can now be determined by evaluating on a {\em basic
sequence\/} of manifolds, which generate the cobordism ring (see
Hirzebruch~\cite{hirzebruch78}). 

In our case
$$s_{2k}(C_k)=\frac{(-1)^{k+1}(2k+1)!}{2^{k-1}k!}$$
is non-zero, so the virtual manifolds $C_k$ form a basic sequence. For
this sequence of manifolds
\begin{eqnarray*}
g(x) & = &
1+\frac{1}{48}b_{\Theta}(C_1)x^2+\frac{1}{48^22!}b_{\Theta^2}(C_2)x^4+\frac{1}{48^33!}b_{\Theta^3}(C_3)x^6+\ldots
\\
 & = &
 1-\frac{1}{48}6x^2+\frac{1}{48^22!}6^2x^4-\frac{1}{48^33!}6^3x^6+\ldots \\
 & = &
 1-\frac{1}{8}x^2+\frac{1}{8^22!}x^4-\frac{1}{8^33!}x^6+\ldots \\
 & = & {\mathrm exp}(-\frac{1}{8}x^2) 
\end{eqnarray*}
and we also have
\begin{eqnarray*}
g(x)\hspace*{-2mm} & =\hspace*{-2mm} & 1+a_2s_2(c_1)x^2+(a_4s_4\!
 +\!\frac{1}{2}a_2^2s_2^2)(c_2)x^4+(a_6s_6\! +\! a_2a_4s_2s_4\! +\! \frac{1}{6}a_2^3s_2^3)(c_3)x^6+\ldots \\
 & =\hspace*{-2mm} & 1+a_2s_2(c_1)x^2+(a_4-a_2^2)s_4(c_2)x^4+(a_6-2a_2a_4+\frac{4}{6}a_2^3)s_6(c_3)x^6+\ldots
\end{eqnarray*}
Comparing this to the sequence
\begin{eqnarray*}
f(x)^{-2} & = & {\mathrm exp}(-2{\mathrm ln}f(x)) \\
 & = & {\mathrm exp}(-2a_2x^2-2a_4x^4-2a_6x^6-\ldots ) \\
 & = & 1-2a_2x^2+(-2a_4+2a_2^2)x^4+(-2a_6+4a_2a_4-\frac{8}{6}a_2^3)x^6+\ldots
\end{eqnarray*}
we see that $f(x)^{-2}$ has coefficients equal to
$$-2\frac{(-1)^k}{8^kk!s_{2k}(C_k)}=\frac{1}{4^k(2k+1)!}$$
and therefore
\begin{eqnarray*}
f(x)^{-2} & = & 1+\sum_{k=1}^{\infty}\frac{1}{4^k(2k+1)!}x^{2k} \\
   & = & \frac{e^{x/2}-e^{-x/2}}{x} \\
   & = & \frac{{\mathrm sinh}(x/2)}{(x/2)}.
\end{eqnarray*}
The generating sequence for the Todd genus ${\mathrm Td}$ is determined by the
function $\frac{x}{1-e^x}$. In the hyperk{\"a}hler case the Chern
roots occur in plus/minus pairs, so we can take instead the generating
function
$$(\frac{x}{1-e^x}\times\frac{-x}{1-e^{-x}})^{1/2}=\frac{x}{e^{x/2}-e^{-x/2}}=\frac{(x/2)}{{\mathrm
sinh}(x/2)}.$$
Comparing this with the formula for $f(x)^{-2}$, we see that the
generating sequence for $\frac{1}{48^kk!}b_{\Theta^k}$ is
precisely ${\mathrm Td}^{1/2}$. Of course, as stated earlier all our
results are independent of the choice of compatible complex
structure. In particular, if we wish to write things in a form that is
manifestly independent of this choice then we can rephrase our results
in terms of the (topologically invariant) Pontryagin classes instead
of the Chern classes. This would also mean replacing the Todd genus by
the $\hat{A}$-genus, and so it is really $\hat{A}^{1/2}$ which
generates $\frac{1}{48^kk!}b_{\Theta^k}$.  

The above arguments were based upon the assumption that $b_{\Theta^k}$
should be a characteristic number, and so do not constitute a complete
proof. However, in the next chapter we shall prove 
\begin{thm}
A generating sequence for $\frac{1}{48^kk!}b_{\Theta^k}$
is
\begin{eqnarray*}
{\mathrm Td}^{1/2} & = &
1-\frac{1}{48}s_2+\frac{1}{48^22!}(s_2^2+\frac{4}{5}s_4)-\frac{1}{48^33!}(s_2^3+\frac{12}{5}s_2s_4+\frac{64}{35}s_6)+
\\
 & &
 +\frac{1}{48^44!}(s_2^4+\frac{24}{5}s_2^2s_4+\frac{48}{25}s_4^2+\frac{256}{35}s_2s_6+\frac{1152}{175}s_8)-\ldots
\end{eqnarray*}
\label{thm1}
\end{thm}
In fact, we shall prove the corresponding relation in graph homology
from which this result follows. All we need to show is that $\Theta^k$
lies in the polywheel subspace for all $k$, but actually our methods
will reproduce the above formula.

\subsection{The ${\mathrm SU}(2)$ virtual hyperk{\"a}hler manifolds}

Before ending this chapter let us say a few more words about the
virtual hyperk{\"a}hler manifolds $C_k$ arising from the perturbative
Chern-Simons theory with gauge group ${\mathrm SU}(2)$. One could ask
whether there is a family of hyperk{\"a}hler manifolds whose
Rozansky-Witten invariants realize the ${\frak{su}}(2)$ weight system.
Appendix E contains tables of values of Rozansky-Witten invariants 
in low degrees for the compact hyperk{\"a}hler manifolds $S^{[k]}$ and
$T^{[[k]]}$ and their products. The former is the Hilbert scheme of
$k$ points on a K$3$ surface $S$ and the latter is the generalized
Kummer variety. We will say much more about these manifolds in Chapter
$5$, where we will also perform the calculations which lead to the
values in Appendix E. The last part of that appendix contains values
for the ${\frak{su}}(2)$ weight system, or equivalently, values of the
Rozansky-Witten invariants for the virtual hyperk{\"a}hler manifolds
$C_k$. Let us compare this weight system to the hyperk{\"a}hler weight
system. Firstly, using Hilbert schemes of points on a K$3$ surface
$S$, we find
$$C_1\sim -\frac{1}{8}S$$
$$C_2\sim -\frac{1}{12}S^{[2]}+\frac{7}{96}S^2$$
$$C_3\sim -\frac{3}{64}S^{[3]}+\frac{5}{48}S\times
S^{[2]}-\frac{85}{1536}S^3$$
where $\sim$ means that both sides have the same Rozansky-Witten
invariants, ie.\ give the same weight system. Secondly, using
generalized Kummer varieties (note that $T^{[[1]]}$ is the Kummer
model of the K$3$ surface, so it is actually the same as $S$), we get
$$C_1\sim -\frac{1}{8}T^{[[1]]}$$
$$C_2\sim -\frac{1}{36}T^{[[2]]}+\frac{1}{32}(T^{[[1]]})^2$$
$$C_3\sim -\frac{3}{320}T^{[[3]]}+\frac{29}{1440}T^{[[1]]}\times
T^{[[2]]}-\frac{17}{1536}(T^{[[1]]})^3$$
Of course in degrees one, two, and three we only require one, two, and
three weights, respectively, to span graph cohomology, so it is no
surprise that we can express the ${\frak{su}}(2)$ weight system in
terms of these hyperk{\"a}hler weight systems. However, in degree four
we find that this is not possible using only the Hilbert schemes of
points on a K$3$ surface, or using only the generalized Kummer
varieties. Instead we must combine both these hyperk{\"a}hler weight
systems, and then we can write
$$C_4\sim
\frac{1}{32}S^{[4]}-\frac{7}{800}T^{[[4]]}+\frac{5}{256}S\times
S^{[3]}+\frac{1}{48}S^{[2]}\times S^{[2]}-\frac{73}{768}S^2\times
S^{[2]}+\frac{263}{6144}S^4.$$
There does not appear to be any nice pattern emerging in this
behaviour.

Instead let us suppose that the virtual manifolds $C_k$ are
represented by {\em actual\/} hyperk{\"a}hler manifolds, which are
connected and irreducible. We know that
\begin{eqnarray*}
s_{\lambda_1}\ldots s_{\lambda_j}(C_k) & = & (-1)^jc_{\langle
   w_{\lambda_1}\cdots w_{\lambda_j}\rangle}({\frak{su}}(2)) \\
   & = & (-1)^j2^{j-1}c_{\langle w_{2k}\rangle}({\frak{su}}(2)) \\
   & = & (-2)^{j-1}s_{2k}(C_k)
\end{eqnarray*}
where $\lambda_1+\ldots+\lambda_j=2k$. Only the relative values are
important here: an irreducible hyperk{\"a}hler manifold of
real-dimension $4k$ must have Todd genus $k+1$, so we need to rescale
$C_k$ to begin with in order for there to be any hope that it be
represented by an actual hyperk{\"a}hler manifold. Solving for the
Chern classes, we find that
$$c_{\lambda_1}c_{\lambda_2}\cdots c_{\lambda_j}(C_k)$$
is independent of the (even) partition $(\lambda_1,\ldots,\lambda_j)$
of $2k$. Furthermore, after rescaling so that ${\mathrm Td}(C_k)=k+1$
all of these Chern numbers remain integral. Hence there is no
immediate reason why there shouldn't exist actual hyperk{\"a}hler
manifolds to represent (the rescaled) virtual hyperk{\"a}hler
manifolds $C_k$.

%% file: ch4.tex
\section{Relations in graph homology II}

\subsection{Knot theory, wheels, and wheeling}

In this subsection we describe a result known as the Wheeling Theorem
which arose from the study of finite type invariants of knots, and
which has been recently proved by Bar-Natan, Le, and
Thurston~\cite{blt}. In the following subsection we shall use this
result to show that the graphs $\Theta^k$ lie in the polywheel
subspaces, and we shall obtain precise expressions for them in terms
of polywheels. Of course, these are precisely the expressions
described in the previous chapter.

The Wheeling Theorem is an isomorphism of algebras constructed from
certain diagrams. We begin by briefly describing the knot theory
motivation behind these algebras. Almost all of the results stated
here require us to be working over a field of characteristic zero, and
we shall assume this to be the case from the outset. Indeed, we are
really only concerned with the base field being the rational numbers
$\Bbb Q$. A good reference for this material is Bar-Natan's
article~\cite{barnatan95}. 

Given an invariant of framed oriented knots $V$, we can extend it to
an invariant of knots with singular points by the rule
$$V(\knotcrossing)=V(\overcrossing)-V(\undercrossing)$$
where the three knots are identical except in a small ball where they
are as indicated. We call $V$ a {\em Vassiliev invariant of type\/}
$m$ if the extended invariant vanishes on knots with more than $m$
singular points. Vassiliev invariants, also known as finite type
invariants, encompass the Alexander-Conway, Jones, HOMFLY, and
Kauffman polynomials, and the Reshetikhin-Turaev quantum
invariants. It is still an open problem whether they separate all
knots.

The set of all Vassiliev invariants of $\cal V$ is a linear space
filtered by type
$${\cal V}_0\subset {\cal V}_1\subset {\cal V}_2\subset\ldots$$
The fundamental theorem of Vassiliev invariants (due to
Kontsevich~\cite{kontsevich93}) says that the graded
space ${\mathrm gr}{\cal V}$ associated to $\cal V$ can be identified
with the dual ${\cal A}(S^1)^*$ of a certain space ${\cal
A}(S^1)$ of linear combinations of certain diagrams modulo certain
relations. This is not quite a precise statement, as there is an
additional subtlety concerning the parity of the framing number of the
knot. If we temporarily forget the framing, we can say that the graded
space associated to the set of all (unframed) Vassiliev invariants is
isomorphic to the dual of ${\cal A}(S^1)/1{\mathrm T}$, where we have
quotiented by an additional relation $1{\mathrm T}$ known as the {\em
framing independence relation\/}. However, let us stick to framed
knots and ignore the parity problem, as it won't affect what we have
to say here.

The diagrams of ${\cal A}(S^1)$ are unitrivalent graphs together with
an oriented circle $S^1$, upon which all of the univalent (or {\em
external\/}) vertices lie. We call this circle the {\em skeleton\/} of
the diagram, and shall usually draw these diagrams with it broken at
one point to produce a directed line (precisely where we break it is
not important). These diagrams are oriented by an equivalence class of
cyclic orderings of the trivalent (or {\em internal\/}) vertices, with
two orderings being equivalent if they differ on an even number of
vertices. We further require that every connected component of the
unitrivalent graph has at least one external vertex (so that the
overall diagram is connected). The relations are the IHX relations on
internal vertices, and the STU relations near the skeleton, and the AS
relation which simply says that reversing the orientation is the same
as changing the sign. In fact, it can be shown that the IHX and AS
relations follow from the STU relations. All of this is perhaps best
illustrated in the pictures below.
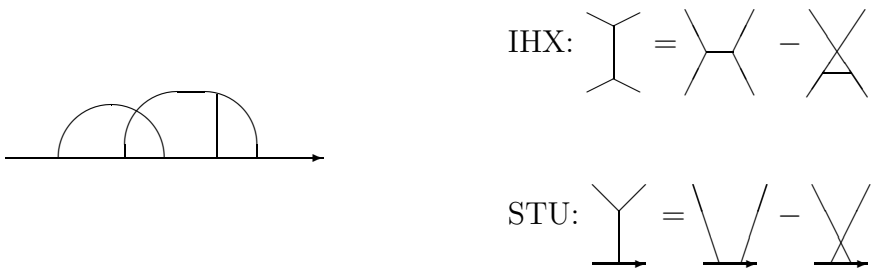
\begin{figure}[htpb]
\begin{picture}(150,130)(-140,0)
\put(-110,50){\vector(1,0){120}}
\put(-70,50){\oval(40,40)[t]}
\put(-40,50){\oval(50,50)[t]}
\put(-30,50){\line(0,1){24}}

\put(80,90){IHX:}
\put(120,100){\line(2,1){10}}
\put(120,100){\line(-2,1){10}}
\put(120,100){\line(0,-1){20}}
\put(120,80){\line(2,-1){10}}
\put(120,80){\line(-2,-1){10}}
\put(135,90){$=$}
\put(155,90){\line(1,0){10}}
\put(155,90){\line(-1,2){8}}
\put(155,90){\line(-1,-2){8}}
\put(165,90){\line(1,2){8}}
\put(165,90){\line(1,-2){8}}
\put(182,90){$-$}
\put(194,106){\line(2,-3){22}}
\put(215,106){\line(-2,-3){22}}
\put(199,82){\line(1,0){11}}

\put(80,25){STU:}
\put(112,10){\vector(1,0){20}}
\put(122,10){\line(0,1){20}}
\put(122,30){\line(1,1){10}}
\put(122,30){\line(-1,1){10}}
\put(138,25){$=$}
\put(154,10){\vector(1,0){20}}
\put(160,10){\line(-1,3){10}}
\put(168,10){\line(1,3){10}}
\put(182,25){$-$}
\put(196,10){\vector(1,0){20}}
\put(202,10){\line(1,2){15}}
\put(210,10){\line(-1,2){15}}
\end{picture}
\caption{A diagram in ${\cal A}(S^1)$, and the IHX and STU relations}
\end{figure}

We call the elements of ${\cal A}(S^1)$ {\em chord diagrams\/}, and
this space is graded by half the total number of vertices in the
diagram. We allow elements given by infinite series so long as there
are only a finite number of terms of any given degree.

The proof of ${\mathrm gr}{\cal V}\cong {\cal A}(S^1)^*$ relies on the
construction of a universal Vassiliev invariant of knots, taking
values in ${\cal A}(S^1)$. A {\em weight\/} of degree $m$ is an
element of the dual of ${\cal A}(S^1)^m$, and composing with the
universal invariant gives us a Vassiliev knot invariant of finite type
$m$. Conversely, all Vassiliev invariants of finite type $m$ can be
constructed from some system of weights of degrees $0$ to $m$. Taking
particular systems of weights (in all degrees) allows us to recover
the Jones polynomial, Alexander-Conway polynomial, and all of the
other invariants of knots mentioned above.

Up to normalization, all known examples of the universal Vassiliev
invariant are believed to give the same answer as the original
example, the Kontsevich integral $Z$~\cite{kontsevich93}. Actually,
this is really the framed Kontsevich integral, and composing with the
natural projection from ${\cal A}(S^1)$ to ${\cal A}(S^1)/1{\mathrm
T}$ will give us the original Kontsevich integral (of unframed
knots). Now $Z$ is extremely difficult to calculate, and it has only
recently been discovered how to calculate it for the most trivial
knot, the unknot $U$. In order to describe $Z(U)$ we first need to
introduce another description of the space ${\cal A}(S^1)$. 

Let $\cal B$ be the space of linear combinations of oriented
unitrivalent graphs, modulo the IHX and AS relations. In other words,
the diagrams in $\cal B$ are essentially the same as those in ${\cal
A}(S^1)$, except that we remove the skeletons. Thus the STU relations
are no longer relevant, but we retain the other relations. As before,
each connected component must have at least one univalent vertex,
though the overall diagram need not be connected. The space $\cal B$
is graded by half the total number of vertices. Some examples of
elements of $\cal B$ are wheels $w_{2m}$ and their disjoint unions. 
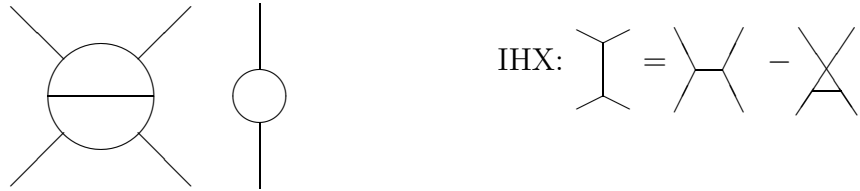
\begin{figure}[htpb]
\begin{picture}(150,110)(-140,0)

\put(-70,50){\circle{40}}
\put(-84,64){\line(-1,1){20}}
\put(-84,36){\line(-1,-1){20}}
\put(-56,36){\line(1,-1){20}}
\put(-56,64){\line(1,1){20}}
\put(-90,50){\line(1,0){40}}
\put(-10,50){\circle{20}}
\put(-10,60){\line(0,1){25}}
\put(-10,40){\line(0,-1){25}}

\put(80,60){IHX:}
\put(120,70){\line(2,1){10}}
\put(120,70){\line(-2,1){10}}
\put(120,70){\line(0,-1){20}}
\put(120,50){\line(2,-1){10}}
\put(120,50){\line(-2,-1){10}}
\put(135,60){$=$}
\put(155,60){\line(1,0){10}}
\put(155,60){\line(-1,2){8}}
\put(155,60){\line(-1,-2){8}}
\put(165,60){\line(1,2){8}}
\put(165,60){\line(1,-2){8}}
\put(182,60){$-$}
\put(194,76){\line(2,-3){22}}
\put(215,76){\line(-2,-3){22}}
\put(199,52){\line(1,0){11}}

\end{picture}
\caption{A diagram in ${\cal B}$, and the IHX relations}
\end{figure}

There is a natural map ${\cal B}\rightarrow {\cal A}(S^1)$
given by averaging over all ways of joining the univalent vertices to
a skeleton, and in fact this gives an isomorphism between $\cal B$ and
${\cal A}(S^1)$. The proof relies on using the equivalence relations
to rewrite an arbitrary element of ${\cal A}(S^1)$ like something
which is clearly in the image of the above map. We note here that
$\cal B$ has a product given by taking the disjoint union of
unitrivalent graphs. We denote this product by $\cup$. The space
${\cal A}(S^1)$ also has a product given by juxtaposition of
skeletons. This is why we choose to write skeletons as directed lines
rather than closed circles, but recall that it doesn't matter where we
break the circle in order to get a directed line (in other words, this
product is well-defined). Using the isomorphism ${\cal B}\cong{\cal
A}(S^1)$ we can transfer this product to $\cal B$, but note that it is
not the same as the disjoint union product. Instead we denote this
second product by $\times$.

The following theorem is known as the Wheels Theorem and was recently
proved by Bar-Natan, Le, and Thurston~\cite{blt}.
\begin{thm}
The framed Kontsevich integral $Z(U)\in{\cal A}(S^1)\cong{\cal B}$ of
the unknot is given by
$$\Omega ={\mathrm exp}_{\cup}\sum_{m=1}^{\infty}b_{2m}w_{2m}$$
where the `modified Bernoulli numbers' $b_{2m}$ are defined by
$$\sum_{m=0}^{\infty}b_{2m}x^{2m}=\frac{1}{2}{\mathrm
log}\frac{{\mathrm sinh}x/2}{x/2}$$
and ${\mathrm exp}_{\cup}$ means that we exponentiate using the
disjoint union product.
\end{thm}
The first few terms of $\Omega$ look like
$$1+\frac{1}{48}w_2+\frac{1}{48^22!}(w_2^2-\frac{4}{5}w_4)+\frac{1}{48^33!}(w_2^3-\frac{12}{5}w_2w_4+\frac{64}{35}w_6)+\ldots$$
From the definition of the coefficients $b_{2m}$, it is clear that
this formula can be obtained from ${\mathrm Td}^{-1/2}$ by replacing 
the characteristic classes $s_{2m}$ by wheels $w_{2m}$. 

Although we will not make use of this theorem, we do need to know
$\Omega$ which plays a part in the next theorem we will describe. We
wish to extend the space $\cal B$ by including diagrams which have
connected components with no univalent vertices, and we call this
larger space ${\cal B}^{\prime}$. Likewise, we extend ${\cal A}(S^1)$
by including diagrams with connected components consisting of
trivalent graphs with no univalent vertices, and we call this larger
space ${\cal A}(S^1)^{\prime}$. The isomorphism ${\cal B}\cong{\cal
A}(S^1)$ extends to ${\cal B}^{\prime}\cong{\cal A}(S^1)^{\prime}$.

We have seen that the space ${\cal B}$ admits two different products,
$\cup$ and $\times$. The disjoint union product $\cup$ clearly extends
to ${\cal B}^{\prime}$. Also, the juxtaposition product on ${\cal
A}(S^1)$ clearly extends to ${\cal A}(S^1)^{\prime}$, and then can be
transferred to ${\cal B}^{\prime}$ as these spaces are isomorphic. We
continue to denote these two products by $\cup$ and $\times$, and
denote the two algebras given by these products by ${\cal
B}^{\prime}_{\cup}$ and ${\cal B}^{\prime}_{\times}$ respectively.

Given a unitrivalent graph $C\in{\cal B}^{\prime}$, we get an operator
$$\hat{C}:{\cal B}^{\prime}\rightarrow{\cal B}^{\prime}$$ 
defined in the
following way. If $C$ has no more univalent vertices than
$C^{\prime}$, then $\hat{C}(C^{\prime})$ is given by summing over all
the ways of joining them to some (or all) of the univalent vertices of
$C^{\prime}$; otherwise $\hat{C}(C^{\prime})$ is zero. We can then
extend this definition linearly to any element $C\in{\cal
B}^{\prime}$, including infinite series provided that $C$ contains
only a finite number of terms of any given degree. For example, when
$C$ and $C^{\prime}$ are wheels we get 
$$\widehat{w_2}(w_4)=8\wtwowfourI +4\wtwowfourII,\qquad
\widehat{w_4}(w_2)=0.$$ 
\vspace*{1mm}

\noindent
Acting with the operator $\hat{\Omega}$ is known as {\em wheeling\/},
since $\Omega$ is made from wheels.

The following theorem is known as the Wheeling Theorem, and was
recently proved by Bar-Natan, Le, and Thurston~\cite{blt}. It is also
allegedly a corollary of Kontsevich's results~\cite{kontsevich97}.
\begin{thm}
The operator associated to $\Omega$ intertwines the two different
product structures on ${\cal B}^{\prime}$, ie.\ $\hat{\Omega}:{\cal
B}^{\prime}_{\cup}\rightarrow{\cal B}^{\prime}_{\times}$ is an
isomorphism of algebras.
\end{thm}
The Wheeling Theorem was conjectured by Bar-Natan, Garoufalidis,
Rozansky, and Thurston in~\cite{bgrt98}, and also independently by
Deligne~\cite{deligne96}. In the former article the Wheels Theorem was
also conjectured, and both were proved ``at the level of Lie
algebras'', though knowing this is not sufficient to deduce the
theorems in their full generality. It simply means verifying the
results when evaluated in a Lie algebra weight system, as in the
previous chapter where we evaluated trivalent graphs in the
${\frak{su}}(2)$ weight system. In a Lie algebra weight system, the
Wheeling Theorem reduces to the Duflo isomorphism~\cite{duflo70}. If
these weight systems generated graph cohomology (or more generally,
the dual space $({\cal B}^{\prime})^*$) then this would suffice, but
unfortunately it is believed that they do not. For example, it is
known that there is an element of ${\cal B}$ of degree $16$ which is
killed by all simple Lie algebra weight systems (see
Vogel~\cite{vogel96}), and it is suspected that this element will in
fact vanish in all Lie algebra weight systems. Of course, there should
also be corresponding elements of graph homology of degree $16$ which
satisfy the same vanishing properties in Lie algebra weight systems.
An interesting question, which we do not answer here, is whether these
graph homology classes give us Rozansky-Witten invariants which are
non-zero on some hyperk{\"a}hler manifold. Such a manifold would need
to have real-dimension $64$, which is considerably larger than in the
examples we will present in the next chapter. 

In the following subsection we shall use the Wheeling Theorem to prove
Theorem~\ref{thm1} from Chapter $3$.

\subsection{The proof of Theorem~\ref{thm1}}

The Wheeling Theorem says that wheeling with $\hat{\Omega}$
intertwines the two different product structures on ${\cal
B}^{\prime}$. Let $\ell$ be the unitrivalent graph given by a single
line, ie.\ with two univalent vertices connected by a single
edge. Then by the Wheeling Theorem
$$\hat{\Omega}({\ell}^{\cup k})=(\hat{\Omega}(\ell))^{\times k}$$
where the superscripts $\cup$ and $\times$ are to indicate that we are
using the two different products to calculate the $k$th powers. We
can consider the terms on each side of this equality which have no
univalent vertices, and equate these parts. This makes sense because
the number of univalent vertices is preserved under the IHX and AS
relations. 

The left hand side is given by 
$$\hat{\Omega}(\underbrace{\ell\cup\cdots\cup\ell}_k).$$
To get a term with no univalent vertices, we need to take the $k$th
term $\Omega_k$ in the series for $\Omega$ and join all of its
univalent vertices to those of $\ell^{\cup k}$. Up to a factor, this
is the same as taking the sum over all ways of joining the $2k$
univalent vertices of $\Omega_k$ to themselves, ie.\ the closure
$\langle\Omega_k\rangle$ of $\Omega_k$. The factor is precisely
$2^kk!$ since there are $k!$ ways to order the copies of $\ell$ and
two ways of joining each at a given location (ie.\ we can reverse the
ends). We already remarked that $\Omega$ looks like ${\mathrm
Td}^{-1/2}$ with the characteristic classes $s_{2m}$ replaced by
wheels $w_{2m}$. Therefore modulo terms with univalent vertices, the
left hand side equals $2^kk!$ times ${\mathrm Td}^{-1/2}_k$ with the
characteristic classes replaced by wheels, and summed over 
all the ways of joining the spokes of these wheels to produce a
polywheel. For example, the first few terms are
$$\begin{array}{lll}
  \hat{\Omega}(\ell)   & = & \frac{1}{24}\langle w_2\rangle \\
  \hat{\Omega}(\ell^{\cup 2}) & = &
  \frac{1}{24^2}(\langle w_2^2\rangle -\frac{4}{5}\langle w_4\rangle) \\
  \hat{\Omega}(\ell^{\cup 3}) & = &
  \frac{1}{24^3}(\langle w_2^3\rangle-\frac{12}{5}\langle
  w_2w_4\rangle +\frac{64}{35}\langle w_6\rangle) \\
  \hat{\Omega}(\ell^{\cup 4}) & = & \frac{1}{24^4}(\langle
  w_2^4\rangle-\frac{24}{5}\langle w_2^2w_4\rangle
  +\frac{48}{25}\langle w_4^2\rangle +\frac{256}{35}\langle
  w_2w_6\rangle-\frac{1152}{175}\langle w_8\rangle)
\end{array}$$
modulo terms with univalent vertices.

Now consider the right hand side. Firstly, if we include terms with
univalent vertices we find
\begin{eqnarray*}
\hat{\Omega}(\ell) & = & \hat{1}(\ell)+\frac{1}{48}\widehat{w_2}(\ell) \\
 & = & \ell+\frac{1}{24}\langle w_2\rangle \\
 & = & \ell+\frac{1}{24}\Theta
\end{eqnarray*}
which looks like
$$\twoVchord +\frac{1}{24}\thetachord$$
when written as an element of ${\cal A}(S^1)^{\prime}$. Therefore
$$(\hat{\Omega}(\ell))^{\times k}=(\twoVchord
+\frac{1}{24}\thetachord)^{\times k}$$
and modulo terms with univalent vertices this equals
$$\frac{1}{24^k}(\thetachord)^{\times k}.$$
Written as an element of ${\cal B}^{\prime}$ this is simply
$\frac{1}{24^k}\Theta^k$. Comparing with the left hand side we see
that
$$\begin{array}{lll}
  \Theta   & = & \langle w_2\rangle \\
  \Theta^2 & = & \langle w_2^2\rangle -\frac{4}{5}\langle w_4\rangle \\
  \Theta^3 & = & \langle w_2^3\rangle-\frac{12}{5}\langle
  w_2w_4\rangle +\frac{64}{35}\langle w_6\rangle \\
  \Theta^4 & = & \langle w_2^4\rangle-\frac{24}{5}\langle
  w_2^2w_4\rangle +\frac{48}{25}\langle w_4^2\rangle
  +\frac{256}{35}\langle w_2w_6\rangle-\frac{1152}{175}\langle w_8\rangle
\end{array}$$
etc. Therefore we can conclude:
\begin{prp}
For all $k$, the graph $\Theta^k$ lies in the polywheel subspace. The
formula expressing it as a linear combination of polywheels is given
by taking $48^kk!$ times the $k$th term ${\mathrm Td}^{-1/2}_k$ of
${\mathrm Td}^{-1/2}$, turning the characteristic classes $s_{2m}$
into wheels $w_{2m}$, and then summing over all the ways of joining
the spokes of these wheels.
\end{prp}
To prove Theorem~\ref{thm1} it only remains to take the
Rozansky-Witten invariants corresponding to the above graphs. We get
the following sequence of relations:
$$\begin{array}{lll}
  b_{\Theta}   & = & -s_2 \\
  b_{\Theta^2} & = & s_2^2+\frac{4}{5}s_4 \\
  b_{\Theta^3} & = & -s_2^3-\frac{12}{5}s_2s_4-\frac{64}{35}s_6 \\
  b_{\Theta^4} & = &
  s_2^4+\frac{24}{5}s_2^2s_4+\frac{48}{25}s_4^2+\frac{256}{35}s_2s_6+\frac{1152}{175}s_8
\end{array}$$
etc. By taking the Rozansky-Witten invariant corresponding to the
polywheel $\langle w_{\lambda_1}\cdots w_{\lambda_j}\rangle$ we get
$(-1)^js_{\lambda_1}\cdots s_{\lambda_j}$, so we have effectively
turned the wheels $w_{2m}$ back into characteristic classes. The
overall effect of the sign changes $(-1)^j$ is to change the
generating sequence from ${\mathrm Td}^{-1/2}$ to ${\mathrm
Td}^{1/2}$. Therefore the general case looks like
$$b_{\Theta^k}=48^kk!{\mathrm Td}^{1/2}_k$$
or in other words a generating sequence for
$\frac{1}{48^kk!}b_{\Theta^k}$ is ${\mathrm Td}^{1/2}$, which is
precisely Theorem~\ref{thm1}.

\subsection{The bubbling effect}

\begin{figure}[htpb]
\epsfxsize=120mm
\centerline{\epsfbox{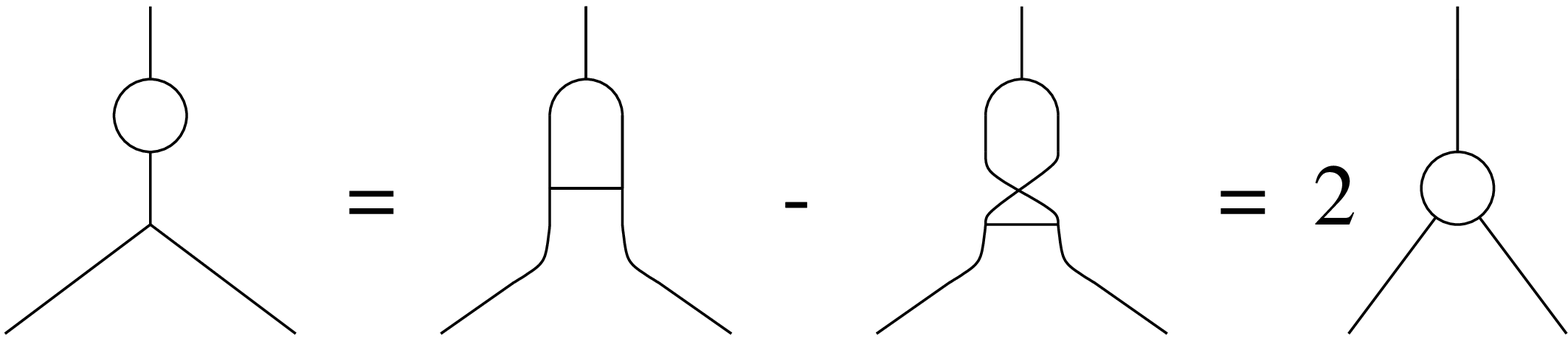}}
\caption{The bubbling effect}
\label{bubble}
\end{figure}

Let $\Gamma$ be a connected trivalent graph with $2k$
vertices. Consider replacing an edge of $\Gamma$ by a two-wheel $w_2$,
which we shall call a {\em bubble\/}, to get a new graph
$\Gamma^{\prime}$ with $2k+2$ vertices. We claim that the graph
homology class of $\Gamma^{\prime}$ is independent of where we add
the bubble. It suffices to show that the bubble commutes past
vertices, and this follows from the IHX relation: a graph with the
bubble near a vertex is equivalent to twice the graph with the bubble
at the vertex (see Figure~\ref{bubble}). Thus we can move the bubble to the
vertex then out the other side, and this proves our claim. We call
this the {\em bubbling effect\/}. Note that if $\Gamma$ has several
connected components then there is some choice as to where we add the
bubble, and different choices may give non-homologous
graphs. Nevertheless, if $\Gamma$ decomposes as
$\gamma_1\cup\cdots\cup\gamma_m$ where each component $\gamma_i$ is
connected with $2k_i$ vertices, then we can define
$$\Gamma^{\prime}=\sum_{i=1}^m\frac{k_i}{k}\gamma_1\cup\cdots\cup\gamma_{i-1}\cup\gamma_i^{\prime}\cup\gamma_{i+1}\cup\cdots\cup\gamma_m.$$

Suppose that $\Gamma$ is in the polywheel subspace, and write
$$\Gamma=\sum a_{(\lambda_1,\ldots,\lambda_j)}\langle
w_{\lambda_1}\cdots w_{\lambda_j}\rangle$$
where the sum is over all even partitions
$(\lambda_1,\ldots,\lambda_j)$ of $2k$. Let us add another two-wheel
$w_2$ to the polywheels. Then either the two spokes of $w_2$ join to
each other and we get a copy of $\Theta$ or they join to the other
wheels and this is equivalent to adding a bubble. There are $k$
choices of where the bubble is added, corresponding to the $k$ pairs
of spokes which are joined, and two ways to connect the bubble at each
location (as we can reverse the ends of $w_2$). If $\Gamma$ decomposes
as $\gamma_1\cup\cdots\cup\gamma_m$ where each component $\gamma_i$ is
connected with $2k_i$ vertices, then $k_i$ of the pairs will be on
the part giving rise to $\gamma_i$. Therefore
\begin{eqnarray*}
\sum a_{(\lambda_1,\ldots,\lambda_j)}\langle w_2w_{\lambda_1}\cdots
w_{\lambda_j}\rangle\hspace*{-20mm} & & \\
 & = & \sum a_{(\lambda_1,\ldots,\lambda_j)}\Theta\langle w_{\lambda_1}\cdots
w_{\lambda_j}\rangle +\sum_{i=1}^m 2k_i\gamma_1\cup\cdots\cup\gamma_i^{\prime}\cup\cdots\cup\gamma_m \\
 & = & \Theta\Gamma +2k\Gamma^{\prime}
\end{eqnarray*}
where we write the disjoint union $\Theta\cup\Gamma$ simply as
$\Theta\Gamma$. We can conclude that $\Theta\Gamma +2k\Gamma^{\prime}$
is also in the polywheel subspace.

For example, we have a generating function ${\mathrm Td}^{-1/2}$
which tells us how to express $\Theta^k$ as a linear combination of
polywheels (for all $k$). Let us write
$$\Theta^k=\sum a^{(k)}_{(\lambda_1,\ldots,\lambda_j)}\langle
w_{\lambda_1}\cdots w_{\lambda_j}\rangle.$$
Then the formula above says that
$$\sum a^{(k)}_{(\lambda_1,\ldots,\lambda_j)}\langle w_2w_{\lambda_1}\cdots
w_{\lambda_j}\rangle =\Theta^{k+1} +2k(\Theta^k)^{\prime}.$$
Now
$$\Theta^{k+1}=\sum a^{(k+1)}_{(\eta_1,\ldots,\eta_j)}\langle
w_{\eta_1}\cdots w_{\eta_j}\rangle$$
summed over all even partitions $(\eta_1,\ldots,\eta_j)$ of $2k+2$,
and $\Theta^{\prime}=\Theta_2$, so that
\begin{eqnarray*}
k(\Theta^k)^{\prime} & = &
\sum_{i=1}^k\Theta\cup\cdots\cup\Theta^{\prime}\cup\cdots\Theta
\\
 & = & k\Theta^{k-1}\Theta_2.
\end{eqnarray*}
Then the above formula tells us that $\Theta^{k-1}\Theta_2$ can also
be expressed as a linear combination of polywheels, and we can
determine this precisely since we know the coefficients $a^{(k)}$ and
$a^{(k+1)}$. For example, if $k=3$ we have
\begin{eqnarray*}
\Theta^2\Theta_2 & = & \frac{1}{6}(\langle w_2^4\rangle-\frac{12}{5}\langle
w_2^2w_4\rangle +\frac{64}{35}\langle w_2w_6\rangle) \\
 & & -\frac{1}{6}(\langle w_2^4\rangle-\frac{24}{5}\langle w_2^2w_4\rangle
+\frac{48}{25}\langle w_4^2\rangle +\frac{256}{35}\langle w_2w_6\rangle-\frac{1152}{175}\langle w_8\rangle) \\
 & = & \frac{2}{5}\langle w_2^2w_4\rangle -\frac{8}{25}\langle
w_4^2\rangle -\frac{32}{35}\langle w_2w_6\rangle
+\frac{192}{175}\langle w_8\rangle.
\end{eqnarray*}

Now we take $\Gamma$ to be $\Theta^{k-2}\Theta_2$, which can be
written
$$\Theta^{k-2}\Theta_2=\sum a^{(k)}_{(\lambda_1,\ldots,\lambda_j)}\langle
w_{\lambda_1}\cdots w_{\lambda_j}\rangle$$
and hence
$$\sum a^{(k)}_{(\lambda_1,\ldots,\lambda_j)}\langle
w_2w_{\lambda_1}\cdots
w_{\lambda_j}\rangle=\Theta^{k-1}\Theta_2+2k(\Theta^{k-2}\Theta_2)^{\prime}.$$
As before, the first term on the right can be written
$$\Theta^{k-1}\Theta_2=\sum a^{(k+1)}_{(\eta_1,\ldots,\eta_j)}\langle
w_{\eta_1}\cdots w_{\eta_j}\rangle$$
and in the second term 
$$k(\Theta^{k-2}\Theta_2)^{\prime}=(k-2)\Theta^{k-3}\Theta_2^2+2\Theta^{k-2}\Theta_3$$
which we can therefore express as a linear combination of
polywheels. In particular $\Theta^{k-2}\Theta_3$ can be expressed as a
linear combination of polywheels and $\Theta^{k-3}\Theta_2^2$.

Next we take $\Gamma$ to be 
\begin{eqnarray*}
(k-1)(\Theta^{k-3}\Theta_2)^{\prime} & = &
(k-3)\Theta^{k-4}\Theta_2^2+2\Theta^{k-3}\Theta_3
\\
 & = & (k-1)\sum a^{(k)}_{(\lambda_1,\ldots,\lambda_j)}\langle
 w_{\lambda_1}\cdots w_{\lambda_j}\rangle
\end{eqnarray*}
and hence
\begin{eqnarray*}
(k-1)\sum a^{(k)}_{(\lambda_1,\ldots,\lambda_j)}\langle
w_2w_{\lambda_1}\cdots w_{\lambda_j}\rangle \hspace*{-50mm} & & \\
 & = & (k-3)\Theta^{k-3}\Theta_2^2+2\Theta^{k-2}\Theta_3+2k((k-3)\Theta^{k-4}\Theta_2^2+2\Theta^{k-3}\Theta_3)^{\prime}
\\
 & = &
 (k-3)\Theta^{k-3}\Theta_2^2+2\Theta^{k-2}\Theta_3+2(k-3)((k-4)\Theta^{k-5}\Theta_2^3+4\Theta^{k-4}\Theta_2\Theta_3) \\
 & & +4((k-2)\Theta^{k-4}\Theta_2\Theta_3+2\Theta^{k-3}\Theta_4) \\
 & = & (k-3)\Theta^{k-3}\Theta_2^2+2\Theta^{k-2}\Theta_3+2(k-3)(k-4)\Theta^{k-5}\Theta_2^3 \\
 & & +4(3k-8)\Theta^{k-4}\Theta_2\Theta_3+8\Theta^{k-3}\Theta_4.
\end{eqnarray*}
It is beginning to become rather difficult to keep track of all the
terms here, but the important point we wish to observe is that we can
express $\Theta^{k-3}\Theta_4$ as a linear combination of polywheels
and graphs made from the disjoint union of copies of the necklace
graphs $\Theta$, $\Theta_2$, and $\Theta_3$. Furthermore, continuing
with these kinds of calculations allows us to prove inductively the
following result. 
\begin{prp}
For $m=1$ to $k$, the graph $\Theta^{k-m}\Theta_m$ with $2k$ vertices
can be expressed as a linear combination of polywheels and graphs made
from the disjoint union of copies of the necklace graphs $\Theta$,
$\Theta_2$, $\ldots$, and $\Theta_{m-1}$.
\label{necklace}
\end{prp}
For $m=1$ this is just the statement that $\Theta^k$ is in the
polywheel subspace. In terms of Rozansky-Witten invariants, this says
that $b_{\Theta^k}$ is a characteristic number. The significance of
this result to the Rozansky-Witten theory for $m>1$ will be discussed
in the next chapter.

%% file: ch5.tex
\section{Calculations}

\subsection{Consequences of our graph homology relations}

In the last two chapters we have derived certain relations in graph
homology. Our primary aim was to express certain graphs or
combinations of graphs as linear combinations of polywheels, because
this implies that the corresponding Rozansky-Witten invariants are
characteristic numbers. For example, we proved in Subsection $4.2$
that a generating sequence for $\frac{1}{48^kk!}b_{\Theta^k}$ is
${\mathrm Td}^{1/2}$. Combining this with our expression for
$b_{\Theta^k}$ from Subsection $2.4$ proves the following theorem.
\begin{thm}
The ${\cal L}^2$-norm of the curvature of an irreducible compact
hyperk{\"a}hler manifold can be expressed in terms of the volume and
characteristic numbers. More specifically, we have
$$\frac{1}{(192\pi^2k)^k}\frac{\|K\|^{2k}}{{\mathrm
vol}(X)^{k-1}}=\int_X{\mathrm Td}^{1/2}_k(X).$$
\end{thm}
When $k=1$ we get
$$\|K\|^2=-4\pi^2 s_2(X)$$
so a K$3$ surface $S$ has $\|K\|^2=192\pi^2$. For $k=2$ we have 
$$\|K\|^4=32\pi^4{\mathrm vol}(X)(s_2^2(X)+\frac{4}{5}s_4(X))$$
from which we conclude that
$$s_2^2(X)+\frac{4}{5}s_4(X)>0.$$
Writing this as
$$\frac{4}{5}(7c_2^2(X)-4c_4(X))$$
and noting that the Todd genus of an irreducible eight-dimensional
hyperk{\"a}hler manifold must equal three, ie.\ 
$$\frac{1}{720}(3c_2^2(X)-c_4(X))=3$$
tells us that the Euler characteristic $c_4(X)$ must be less than
$3024$. Note that a reformulation by Bogolomov of a result of
Verbitsky~\cite{bogomolov96} (see also Beauville~\cite{beauville99})
shows that $c_4(X)\leq 324$, so our bound is quite crude. Actually,
$c_4(S^{[2]})=324$ for the Hilbert scheme of two points on a K$3$
surface $S$, so the upper bound $324$ is sharp. We will discuss these
specific examples of compact hyperk{\"a}hler manifolds in later
subsections.

In general, since $\|K\|$ is positive we get the following bounds on
the characteristic numbers.
\begin{cor}
Suppose $X$ is an irreducible compact hyperk{\"a}hler manifold of
real-dimension $4k$. Then we have
$$\int_X{\mathrm Td}_k^{1/2}(X)>0.$$
\label{bounds}
\end{cor}
Equality would imply that $\|K\|=0$ and hence $X$ is flat, ie.\ a
complex torus $T^k$, but such an $X$ is not irreducible.

We saw in the previous chapter that $\Theta^{k-2}\Theta_2$ also lies
in the polywheel subspace, and therefore the corresponding
Rozansky-Witten invariant $b_{\Theta^{k-2}\Theta_2}$ is also a
characteristic number. Suppose that $X$ is irreducible, and recall
the discussion in Subsection $2.4$. We showed there that the
Rozansky-Witten invariant corresponding to the disjoint union of
several graphs is given by 
$$b_{\gamma\gamma^{\prime}}(X)=(2\pi^2)^{-k}k!\beta_{\gamma}\beta_{\gamma^{\prime}}{\mathrm vol}(X)$$
where $\beta_{\gamma}$ and $\beta_{\gamma^{\prime}}$ are some scalars
depending on $X$ and the graphs $\gamma$ and $\gamma^{\prime}$
respectively. In particular
$$b_{\Theta^k}(X)=(2\pi^2)^{-k}k!\beta_{\Theta}^k{\mathrm vol}(X)$$
and
$$b_{\Theta^{k-2}\Theta_2}(X)=(2\pi^2)^{-k}k!\beta_{\Theta}^{k-2}\beta_{\Theta_2}{\mathrm vol}(X).$$
Now observe that
\begin{eqnarray*}
b_{\Theta^{k-4}\Theta_2^2}(X) & = & (2\pi^2)^{-k}k!\beta_{\Theta}^{k-4}\beta_{\Theta_2}^2{\mathrm vol}(X) \\
 & = & \frac{((2\pi^2)^{-k}k!\beta_{\Theta}^{k-2}\beta_{\Theta_2}{\mathrm
vol}(X))^2}{(2\pi^2)^{-k}k!\beta_{\Theta}^k{\mathrm vol}(X)} \\
 & = & \frac{b_{\Theta^{k-2}\Theta_2}(X)^2}{b_{\Theta^k}(X)}
\end{eqnarray*}
and more generally
$$b_{\Theta^{k-2m}\Theta_2^m}(X)=\frac{b_{\Theta^{k-2}\Theta_2}(X)^m}{b_{\Theta^k}(X)^{m-1}}$$
for $m=1$ to $\lfloor k/2\rfloor$, where $\lfloor k/2\rfloor$ is $k/2$
if $k$ is even, and $(k-1)/2$ if it is odd. Since $b_{\Theta^k}$ and
$b_{\Theta^{k-2}\Theta_2}$ are characteristic numbers, we conclude
that for {\em irreducible\/} hyperk{\"a}hler manifolds
$b_{\Theta^{k-2m}\Theta_2^m}$ can be written as a rational function of
characteristic numbers.

We also saw that $\Theta^{k-3}\Theta_3$ could be expressed as a linear
combination of polywheels and $\Theta^{k-4}\Theta_2^2$, and thus
$b_{\Theta^{k-3}\Theta_3}$ can also be written as a rational function
of characteristic numbers. As before we can write
$$b_{\Theta^{k-5}\Theta_2\Theta_3}(X)=\frac{b_{\Theta^{k-2}\Theta_2}(X)b_{\Theta^{k-3}\Theta_3}(X)}{b_{\Theta^k}(X)}$$
and hence this Rozansky-Witten invariant is also a rational function
of characteristic numbers. In fact, this is true for $b_{\Gamma}$
where $\Gamma$ is any graph made from a disjoint union of copies of
$\Theta$, $\Theta_2$, and $\Theta_3$. Proceeding by induction and
using Proposition~\ref{necklace} proves the following theorem.
\begin{thm}
Let $\Gamma$ be a trivalent graph with $2k$ vertices constructed by
taking a disjoint union of copies of the necklace graphs $\Theta$,
$\Theta_2$, $\ldots,$ and $\Theta_k$. Then for an irreducible
hyperk{\"a}hler manifold $X$, $b_{\Gamma}(X)$ can be expressed as a
rational function of the characteristic numbers of $X$.
\label{rational}
\end{thm}
Note that the only denominators in these expressions are powers of
$b_{\Theta}(X)$, or in other words ${\mathrm Td}^{1/2}(X)$. We have
already noted in Corollary~\ref{bounds} that this is strictly positive
for irreducible $X$, by virtue of the fact that $\|K\|$ is strictly
positive, hence none of these rational functions can be singular.

Recall that for reducible hyperk{\"a}hler manifolds we have a product
formula
$$b_{\Gamma}(X \times Y)=\sum_{\gamma\sqcup{\gamma}^{\prime}=\Gamma}
b_{\gamma}(X) b_{\gamma^{\prime}}(Y).$$ 
More generally, if $X$ decomposes into irreducible factors
$X_1\times\cdots\times X_m$ then 
$$b_{\Gamma}(X_1\times\cdots\times
X_m)=\sum_{\gamma_1\sqcup\cdots\sqcup\gamma_m=\Gamma} 
b_{\gamma_1}(X_1)\cdots b_{\gamma_m}(X_m).$$
If $\Gamma$ is constructed from a disjoint union of necklace graphs
$\Theta$, $\Theta_2$, $\ldots$, and $\Theta_k$, then each of
$\gamma_1,\ldots,\gamma_m$ must be too. Thus $b_{\gamma_i}(X_i)$ can
be determined from the characteristic numbers of $X_i$. So to
determine $b_{\Gamma}(X_1\times\cdots\times X_m)$ we need to know the
characteristic numbers of all the irreducible factors of $X$. These
cannot necessarily be determined from the knowledge of the
characteristic numbers of $X$.

For degree $k\leq 5$ the graphs constructed from the polywheels and
the necklace graphs $\Theta$, $\Theta_2$, $\ldots,$ and $\Theta_5$
span graph homology (see Appendix A). Therefore if we know the
characteristic numbers of $X$ (or the characteristic numbers of the
irreducible factors of $X$, in the case when $X$ is reducible) we can
determine all the Rozansky-Witten invariants of $X$. So for
real-dimension less than or equal to $20$ the Rozansky-Witten
invariants are {\em essentially\/} nothing more than the
characteristic numbers. We say `essentially' because normally we would
only take linear combinations of characteristic numbers, not rational
functions. We will prove in Subsection $5.5$ that some Rozansky-Witten
invariants can certainly not be written as linear combinations of
characteristic numbers.

It is worth considering what happens as $k$ gets larger. We want to
roughly count how many Rozansky-Witten invariants there are and how
many can be expressed as rational functions of characteristic numbers
according to Theorem~\ref{rational}. In particular, we would like to know
how many graphs {\em cannot\/} be expressed as linear combinations of
polywheels and disjoint unions of necklace graphs $\Theta_i$. We will
think of the latter as being part of a `conventional' basis for graph
homology, and the polywheels as being certain linear combinations of
the completed basis. 

Firstly, in degree $k$ our conventional basis contains ${\mathrm
dim}{\cal A}(\emptyset)^k$ graphs, of which ${\mathrm dim}{\cal
A}(\emptyset)^k_{\mathrm conn}$ are connected. The disjoint unions of
necklace graphs $\Theta_i$ account for $p(k)$ of these graphs, where
$p(k)$ is the number of partitions of $k$. This leaves 
$${\mathrm dim}{\cal A}(\emptyset)^k-p(k)$$
graphs. The polywheels span a $p(k)$-dimensional subspace, and so
there are $p(k)$ equations given by writing the polywheels in terms of
our conventional basis. However, we can write $\Theta^k$ and
$\Theta^{k-2}\Theta_2$ in terms of polywheels. We can also write
$\Theta^{k-3}\Theta_3$ in terms of polywheels and
$\Theta^{k-4}\Theta_2^2$. Indeed, for $m=1$ to $k$ we can write 
$\Theta^{k-m}\Theta_m$ in terms of polywheels and disjoint unions of
necklace graphs $\Theta_i$ with $i<m$ (this is precisely
Proposition~\ref{necklace}). Thus there is a subspace of the span of
the $p(k)$ polywheel equations which has dimension at least $k$ and
involves only disjoint unions of necklace graphs $\Theta_i$. So
overall we are left with at most $p(k)-k$ equations involving
polywheels, disjoint unions of necklace graphs $\Theta_i$, and the
remaining 
$${\mathrm dim}{\cal A}(\emptyset)^k-p(k)$$ 
graphs. It follows that if
$$p(k)-k<{\mathrm dim}{\cal A}(\emptyset)^k-p(k)$$
then there are certainly some graphs which cannot be expressed as
linear combinations of polywheels and disjoint unions of necklace
graphs $\Theta_i$. We summarize these values for $k\leq 10$ in the
table below.
$$
\begin{array}{|l|cccccccccc|}  
\hline
k & 1 & 2 & 3 & 4 & 5 & 6 & 7 & 8 & 9 & 10 \\
\hline
{\mathrm dim}{\cal A}(\emptyset)^k_{\mathrm conn} & 1 & 1 & 1 & 2 & 2
& 3 & 4 & 5 & 6 & 8 \\
{\mathrm dim}{\cal A}(\emptyset)^k & 1 & 2 & 3 & 6 & 9 & 16 & 25 & 42
& 65 & 105 \\
p(k) & 1 & 2 & 3 & 5 & 7 & 11 & 15 & 22 & 30 & 42 \\
{\mathrm dim}{\cal A}(\emptyset)^k-p(k) & 0 & 0 & 0 & 1 & 2 & 5 & 10 &
20 & 35 & 63 \\
p(k)-k & 0 & 0 & 0 & 1 & 2 & 5 & 8 & 14 & 21 & 32 \\
\hline
\end{array}
$$
It can be seen that for $k\geq 7$ there are clearly graphs which
cannot be expressed as linear combinations of polywheels and disjoint
unions of necklace graphs $\Theta_i$. We can conclude that these
graphs give rise to Rozansky-Witten invariants which cannot be
determined from characteristic numbers {\em using the methods
described above\/}. There may, however, be other ways to determine
them using simply the knowledge of the characteristic numbers.

\subsection{Deriving Chern numbers from the $\chi_y$-genus}

From the discussion of the previous subsection we know that when $k$
is small the Rozansky-Witten invariants can all be calculated from
knowledge of the Chern numbers . In particular, this works for $k\leq
5$ and possibly for $k=6$. On the other hand, these methods break down
for $k\geq 7$, meaning that it isn't possible to determine {\em all\/}
Rozansky-Witten invariants in this way. 

In the next subsections we will look at some specific examples of compact
hyperk{\"a}hler manifolds, namely the Hilbert schemes of points on a
K$3$ surface and the generalized Kummer varieties. At this stage, it
suffices to say that for these examples the Chern numbers are not well
known, meaning that there are no explicit formulas for calculating
them in all dimensions. We discuss some partial answers to this
question in the next subsections, but first we wish to describe a way
of obtaining some preliminary information about the Chern numbers from
the Hirzebruch $\chi_y$-genus, which is known explicitly (in terms of
generating functions) in all dimensions for the examples of compact
hyperk{\"a}hler manifolds we wish to study. 

Our method is to use the Riemann-Roch formula to express the
$\chi_y$-genus in terms of Chern numbers. In real-dimensions $4$, $8$, and
$12$ we can invert these relations. In dimension $16$ we are left with
an additional unknown variable. Although we will not proceed beyond
dimension $16$, let us note the following: in real-dimension $4k$
there are $p(k)$ Chern numbers. On the other hand, the Hirzebruch
$\chi_y$-genus
$$\chi_y =\chi^0 +\chi^1y+\ldots +\chi^{2k-1}y^{2k-1}+\chi^{2k}y^{2k}$$
contains $2k+1$ terms, though it is symmetric ($\chi^0=\chi^{2k}$,
$\chi^1=\chi^{2k-1}$, $\ldots$). Since we are in the hyperk{\"a}hler
case, it also satisfies
$$\sum_{m=0}^{2k}(-1)^m(6m^2-k(6k+1))\chi^m=0$$
due to a result of Salamon~\cite{salamon96}. So writing the
$\chi_y$-genus in terms of Chern numbers gives us no more than $k$
equations for the Chern numbers. Note that the relations
$\chi^0=\chi^{2k}$, etc.\ and Salamon's result don't give us any
additional relations on the Chern numbers; ie.\ these relations are
tautologous when we write the $\chi^m$ in terms of Chern numbers. So
in general we would expect there to be at least $p(k)-k$ additional
unknown variables when we try to write the Chern numbers in terms of
the $\chi^m$. In other words, for large $k$ the Hirzebruch
$\chi_y$-genus will contain far less information than is required to
determine the Chern numbers.

Let $X$ have Chern roots
$\{\gamma_1,\gamma_2,\ldots,\gamma_{2k-1},\gamma_{2k}\}$. For a
hyperk{\"a}hler manifold they occur in plus/minus pairs, but we shall
explain the theory for a general manifold. The characteristic classes
are given by
$$s_{\lambda}=\gamma_1^{\lambda}+\gamma_2^{\lambda}+\ldots+\gamma_{2k-1}^{\lambda}+\gamma_{2k}^{\lambda}.$$
Let $h^{p,q}$ be the Hodge numbers, defined by ${\mathrm dim}{\mathrm
H}^q(X,\Lambda^pT^*)$. Then the Hirzebruch $\chi_y$-genus is defined 
to be
$$\chi_y(X)=\sum_{p,q=0}^{2k}(-1)^qh^{p,q}y^p$$
so that each individual coefficient is given by the Euler
characteristic of a bundle of forms
\begin{eqnarray*}
\chi^m(X) & = & \chi(\Lambda^mT^*) \\
 & = & \sum_{q=0}^{2k}(-1)^qh^{m,q}. 
\end{eqnarray*}
The Riemann-Roch theorem says
$$\chi^m(X)=\int_X[{\mathrm ch}(\Lambda^mT^*){\mathrm Td}(X)]_{4k}$$
where ${\mathrm ch}(\Lambda^mT^*)$ is the Chern character of
$\Lambda^mT^*$, ${\mathrm Td}(X)$ is the Todd genus of
(the tangent bundle of) $X$, and $[\ldots ]_{4k}$ means that we pick
out the component of degree $4k$. In terms of the Chern roots we can
write these terms as
$${\mathrm
ch}(\Lambda^mT^*)=\sum_{i_1<i_2<\ldots<i_m}e^{-\gamma_{i_1}-\gamma_{i_2}-\ldots
-\gamma_{i_m}}$$
and
$${\mathrm
Td}(X)=\prod_{i=1}^{2k}\frac{\gamma_i}{1-e^{\gamma_i}}.$$
Substituting these into the Riemann-Roch formula allows us to express
each $\chi^m(X)$ as an integral of some product of Chern roots. In
fact the integrand is symmetric in permutations of the roots and hence
can be written in terms of the characteristic classes, and picking the
degree $4k$ component means taking terms $s_{\lambda_1}\cdots
s_{\lambda_j}$ where $(\lambda_1,\ldots,\lambda_j)$ is a partition of
$2k$. Overall we obtain an expression for the $\chi_y$-genus in terms
of the Chern numbers. Note that these calculations are simpler in the
hyperk{\"a}hler case once we pair the Chern roots into plus/minus
pairs. 

These expressions are shown in the Appendix B for $k\leq 4$. In fact a
different basis for the Chern numbers has been used there, namely
$c_{\lambda_1}\cdots c_{\lambda_j}$. It is of course possible to
rewrite these in terms of the `s' basis, but we do this later, after
attempting to invert the relations. We only show $\chi^m$ for $m\leq
k$, as $\chi^{k+1}=\chi^{k-1}$, $\chi^{k+2}=\chi^{k-2}$, etc. Indeed
by Salamon's result $\chi^k$ can also be expressed as $\chi^m$ for
$m<k$. Thus in the subsequent table, where we attempt to invert the
relations, $\chi^k$ does not appear. We observe that for $k=1$, $2$,
and $3$ it is possible to fully invert these relations, whereas for
$k=4$ there is an additional unknown variable $s$. It is these
expressions which we rewrite in the `s' basis. 

In the next two subsections we will look at some particular examples
of compact hyperk{\"a}hler manifolds, for which the $\chi_y$-genus is
known for all $k$. Substituting into the equations derived above gives
us all the Chern numbers in real-dimensions $4$, $8$, and $12$. From
these we can calculate all of the Rozansky-Witten invariants using the
results of the first subsection in this chapter. In dimension $16$
there is still the unknown $s$. However, by directly calculating one
of the Rozansky-Witten invariants we can determine $s$ for these
examples. This will then enable us to calculate all of the Chern
numbers and hence also the remaining Rozansky-Witten invariants in
this dimension.

In particular, in dimension $16$ we can write $b_{\Theta^4}$ in terms
of the Chern numbers. This leads to
\begin{eqnarray*}
b_{\Theta^4} & = & s_2^4+\frac{24}{5}s_2^2s_4+\frac{48}{25}s_4^2+\frac{256}{35}s_2s_6+\frac{1152}{175}s_8 \\
 & = & \frac{1}{175}(-21936s+4396904448\chi^0-7259904\chi^1+2472960\chi^2+278784\chi^3).
\end{eqnarray*}
We will calculate $b_{\Theta^4}$, and hence determine $s$ from the
above formula. This method leads to the complete tables of Chern
numbers and Rozansky-Witten invariants which appear in Appendices D
and E$.1$.

\subsection{The Hilbert scheme of points on a K$3$ surface}

For a long time the only known compact hyperk{\"a}hler manifold was
the K$3$ surface in real-dimension four. Fujiki~\cite{fujiki83} was
the first to discover higher dimensional examples, of dimension eight,
and shortly afterwards Beauville~\cite{beauville83} generalized these
to produce two families of examples in all dimensions. They are
the Hilbert scheme of points on a K$3$ surface and the generalized
Kummer varieties. The latter are also constructed via Hilbert schemes,
which we will discuss shortly. Other examples of compact
hyperk{\"a}hler manifolds were discovered by Mukai~\cite{mukai84}
around the same time by exhibiting a holomorphic symplectic structure on
the moduli space of sheaves on a K$3$ surface or complex torus. These
later proved to be deformations of Beauville's examples. Indeed apart
from these two main families there is just one other known compact
hyperk{\"a}hler manifold, in real-dimension $20$. It was constructed
by O'Grady~\cite{ogrady99} as a desingularization of a certain
singular moduli space of sheaves on a K$3$ surface.

Let $M$ be an algebraic variety, or more generally a scheme (although
the following discussion holds in a more general setting, we will only
be interested in base field $\Bbb C$). The {\em
Hilbert scheme of $k$ points on $M$\/} is the moduli space of
zero-dimensional subschemes of length $k$, and we denote it by
$M^{[k]}$. The generalization to non-algebraic $M$ is known as the
Douady space, which we also denote by $M^{[k]}$. An example of a
length $k$ subscheme is given by a collection of $k$ distinct
unordered points $\{x_1,\ldots,x_k\}$. More generally, we could allow
some of these points to collide, in which case there should be some
additional information at those points. For two points colliding, this
amounts to an element of the projectivization of the tangent space at
the double point, which corresponds to the direction the two points
collided along. In fact this gives a complete description of the
Hilbert scheme of two points as
$$M^{[2]}={\mathrm Blow}_{\Delta}(M\times M)/{{\cal S}_2}$$
where ${\mathrm Blow}_{\Delta}$ denotes the blow-up along the diagonal
and ${\cal S}_2$ the symmetric group on
two elements, acting by interchanging the two factors of $M\times
M$. Unfortunately for $k>2$ it is not possible to describe $M^{[k]}$ in
such a simple way. 

There is a surjective morphism
$$\pi :M^{[k]}\rightarrow {\mathrm Sym}^kM=M^k/{\cal S}_k$$
to the $k$th symmetric product known as the Hilbert-Chow morphism, given by
$$\xi\mapsto\sum_{x\in M}{\mathrm length}(\xi)_x[x].$$
This map is a bijection on the open subset where elements are given
by $k$ distinct points. When the complex-dimension of $M$ is greater
than one, the symmetric product is singular. However, a theorem of
Fogarty~\cite{fogarty68} says that the Hilbert scheme gives a smooth
resolution when $M$ is a complex surface, ie.\ has complex-dimension
two. In particular, $M^{[k]}$ is of complex-dimension $2k$ and
non-singular in this case.

Suppose that the complex surface $M$ has a holomorphic symplectic form
$\omega$. Then on $M^k$ there is a natural holomorphic symplectic form
given by
$$p_1^*\omega+\ldots +p_k^*\omega$$
where $p_i$ is projection onto the $i$th factor. This is clearly
${\cal S}_k$-invariant, and hence we get a holomorphic two-form on
${\mathrm Sym}^kM$. We claim that the pull-back of this two-form to
$M^{[k]}$ is non-degenerate. Consider $M^{[2]}$ for example: the
corresponding holomorphic two-form on
$${\mathrm Blow}_{\Delta}(M\times M)$$
is degenerate on the blown-up diagonal, but a direct calculation shows
that when we take the quotient by ${\cal S}_2$ we get a non-degenerate
two-form on the Hilbert scheme. For general $k$, the same direct
calculation shows that the holomorphic two-form on $M^{[k]}$ is
non-degenerate on the open subset where at most two points
collide. The complement of this open subset has codimension two, hence
by Hartog's theorem the two-form extends to a holomorphic two-form on
all of $M^{[k]}$. Furthermore, if this two-form were to degenerate
anywhere, it would have to degenerate on an entire divisor, ie.\
codimension one submanifold, which is clearly impossible. This is
Beauville's construction of a holomorphic symplectic form on $M^{[k]}$
(see~\cite{beauville83}). More generally, Mukai~\cite{mukai84}
constructed a holomorphic symplectic form on the moduli space of
stable sheaves (of fixed rank and Chern classes) on a K$3$ surface or
abelian surface, ie.\ on compact surfaces with holomorphic symplectic
forms. This is a generalization as the Hilbert scheme $M^{[k]}$ can be
thought of as the moduli space of rank-one torsion-free sheaves on
$M$.

In the compact K{\"a}hler case, the existence of a holomorphic
symplectic form is equivalent to the existence of a hyperk{\"a}hler
metric by Yau's theorem~\cite{yau78} and the results of Bochner and
Yano~\cite{by52}. This means that we should choose $M$ to be a K$3$
surfaces $S$ or an abelian surface $T$, ie.\ a torus. We shall discuss
the latter case in the next subsection; first we wish to look at the
Hilbert scheme of $k$ points on $S$.

To begin with, the Hodge numbers of the Hilbert scheme of points on an
arbitrary smooth projective surface $M$ were calculated by
Cheah~\cite{cheah96}. Writing 
$$h(M^{[k]})=\sum_{p,q=1}^{2k}h^{p,q}(M^{[k]})x^qy^p$$
then in terms of the Hodge numbers $h^{p,q}(M)$ of $M$ we have
$$\sum_{k=0}^{\infty}h(M^{[k]})t^k=\prod_{n=1}^{\infty}\left(\frac{\prod_{p+q\,\mbox{\footnotesize{odd}}}(1+x^{p+n-1}y^{q+n-1}t^n)^{h^{p,q}(M)}}{\prod_{p+q\,\mbox{\footnotesize{even}}}(1+x^{p+n-1}y^{q+n-1}t^n)^{h^{p,q}(M)}}\right).$$
In particular, putting $x=-1$ gives us the Hirzebruch
$\chi_y$-genus. The Hodge diamond of a K$3$ surface $S$ looks like
$$\begin{array}{ccccc}
  &   & 1  &   & \\
  & 0 &    & 0 & \\
1 &   & 20 &   & 1 \\
  & 0 &    & 0 & \\
  &   & 1  &   & 
\end{array}.$$
Using this to calculate the $\chi_y$-genus of $S^{[k]}$ we get the
results in Appendix C for $k\leq 4$. Substituting these values into
our Riemann-Roch formula gives us the results in Appendix D. In
particular, we know all the Chern numbers and hence all the
Rozansky-Witten invariants for $k=1$, $2$, and $3$. Our next aim is to
calculate $b_{\Theta^4}(S^{[4]})$ and hence determine $s$, in order to
complete the tables in the appendices up to $k=4$. In fact we will
derive a formula for $b_{\Theta^k}(S^{[k]})$ for all $k$.

Recall that the Rozansky-Witten invariants are
invariant under rescaling of the holomorphic symplectic form. So let
us normalize $\omega$ on $S$ so that
$$\int_S\omega\bar{\omega}=1.$$
Let us also denote the induced holomorphic symplectic form on
$S^{[k]}$ by $\omega$; if there is likely to be some confusion we
shall specify which space $\omega$ is on with a subscript, eg.\
$\omega_{S^{[k]}}$. It will have the normalization
\begin{eqnarray*}
\int_{S^{[k]}}\omega^k\bar{\omega}^k & = & \int_{{\mathrm
Sym}^kS}(p_1^*\omega+\ldots+p_k^*\omega)^k(p_1^*\bar{\omega}+\ldots+
p_k^*\bar{\omega})^k \\
 & = & \int_{S^k/{\cal S}_k}(k!)^2p_1^*\omega\cdots p_k^*\omega
p_1^*\bar{\omega}\cdots p_k^*\bar{\omega} \\
 & = & \int_{S^k}k!p_1^*(\omega\bar{\omega})\cdots
p_k^*(\omega\bar{\omega}) \\
 & = & k!
\end{eqnarray*}
Recall from Subsection $2.4$ that for irreducible $X$ we can write
$$b_{\Theta^k}(X)=\frac{1}{(8\pi^2)^kk!}\beta_{\Theta}^k\int_X\omega^k\bar{\omega}^k$$ 
where $\beta_{\Theta}$ is a scalar whose value is given by
$$[\Theta(\Phi)]=\beta_{\Theta}[\bar{\omega}]\in{\mathrm
H}^{0,2}_{\bar{\partial}}(X).$$
Actually $\beta_{\Theta}$ is not a canonically defined number, for if
we rescale the holomorphic symplectic form $\omega$ by $\lambda$ then
$\beta_{\Theta}$ rescales by $\lambda^{-2}$. We needn't be
concerned however, as we shall be looking at the Hilbert scheme
$S^{[k]}$ for which we have fixed a normalization of the holomorphic
symplectic form, and we shall assume this throughout. Indeed, this
normalization leads to the formula
$$b_{\Theta^k}(S^{[k]})=\frac{1}{(8\pi^2)^k}\beta_{\Theta}^k$$
with
$$\Theta(\alpha_T)=\beta_{\Theta}[\bar{\omega}]\in{\mathrm
H}^{0,2}_{\bar{\partial}}(S^{[k]}).$$
Note that we have replaced $[\Theta(\Phi)]$ by $\Theta(\alpha_T)$,
where the Atiyah class $\alpha_T$ is represented by $\Phi$ in
Dolbeault cohomology. This is because we wish to use descriptions of
the Atiyah class other than the Dolbeault description. More generally,
if $\Gamma$ decomposes into connected components as
$\gamma_1\cdots\gamma_j$ then
$$b_{\Gamma}(S^{[k]})=\frac{1}{(8\pi^2)^k}\beta_{\gamma_1}\cdots\beta_{\gamma_j}$$
with
$$\gamma_i(\alpha_T)=\beta_{\gamma_i}[\bar{\omega}^{m_i}]\in{\mathrm
H}^{0,2m_i}_{\bar{\partial}}(S^{[k]})$$ 
where $2m_i$ is the number of vertices in $\gamma_i$. In general
$\beta_{\gamma_i}$ will rescale by $\lambda^{-2m_i}$ under a rescaling
of the holomorphic symplectic form, though once again we assume the
normalization of $\omega$ has been fixed, and hence $\beta_{\gamma_i}$
is a scalar which depends only on $S^{[k]}$ (in other words on $k$),
and which can be evaluated according to the formula
\begin{eqnarray*}
k!\beta_{\gamma_i} & = &
\int_{S^{[k]}}\beta_{\gamma_i}[\omega^k][\bar{\omega}^k] \\
 & = &
 \int_{S^{[k]}}\gamma_i(\alpha_T)[\omega^k][\bar{\omega}^{k-m_i}].
\end{eqnarray*}

In fact this approach is not the simplest. Let us drop the subscript
`i' in what follows. Instead of multiplying
$\gamma(\alpha_T)[\omega^m]$ by $[\omega^{k-m}\bar{\omega}^{k-m}]$ and
then integrating over the entire space $S^{[k]}$, it is easier to
integrate the former over a complex submanifold $X$ of
real-dimension $4m$. We can restrict
$$\gamma(\alpha_T)[\omega^m]\in{\mathrm H}^{2m,2m}_{\bar{\partial}}(S^{[k]})$$
to $X$, and there are projections $T^*_{S^{[k]}}\rightarrow T^*_X$ and
$\bar{T}^*_{S^{[k]}}\rightarrow \bar{T}^*_X$, which induce a projection
$${\mathrm
H}^0(X,\Lambda^{2m}T^*_{S^{[k]}}\wedge\Lambda^{2m}\bar{T}^*_{S^{[k]}}|_X)\rightarrow
{\mathrm H}^0(X,\Lambda^{2m}T^*_X\wedge \Lambda^{2m}\bar{T}^*_X).$$
Applying this to $\gamma(\alpha_T)[\omega^m]|_X$ gives us
something we can integrate over $X$. Actually, we shall assume when we
integrate that forms with some normal component to $X$ give zero, so
that the projection is implicitly understood. The integral equals
$$\int_X\gamma(\alpha_T)[\omega^m]|_X=\beta_{\gamma}\int_X\omega^m\bar{\omega}^m|_X$$
from which we can determine $\beta_{\gamma}$.

From a {\v C}ech cohomology point of view, we begin with
$$\gamma(\alpha_T)\in{\mathrm H}^{2m}(S^{[k]},{\cal O}_{S^{[k]}})$$
and restrict this to a cohomology class in ${\mathrm H}^{2m}(X,{\cal
O}_X)$. Note that there is no reason to expect this cohomology group
to be one-dimensional. We also take the section
$$\omega^m\in{\mathrm H}^0(S^{[k]},\Lambda^{2m}T^*_{S^{[k]}})$$
restrict it to $X$, and project to a holomorphic $2m$-form on
$X$. This gives us
$$\omega_X^m\in{\mathrm H}^0(X,\Lambda^{2m}T^*_X)$$
where $\omega_X$ is the restriction and projection of $\omega$ to $X$,
and it is a holomorphic two-form which will most likely be
degenerate. The Serre duality pairing of $\gamma(\alpha_T)$ with
$\omega_X^m$ then gives us a canonical number, which is equal to
$\beta_{\gamma}$ multiplied by
$$\int_X\omega^m\bar{\omega}^m|_X.$$
As before, this integral factor must be included as $\beta_{\gamma}$
on its own is not a canonically defined scalar.

We can also calculate $\beta_{\gamma}$ via the residue
approach. Recall that this involves a meromorphic connection on
$S^{[k]}$ with simple poles along a smooth divisor $D$. The element
$\beta_T$ is the residue of this connection, and we can use it along
with the Atiyah class $\alpha_T$ to construct an element
$$\gamma(\alpha_T,\beta_T)\in{\mathrm H}^{2m-1}(D,L|_D)$$
which we restrict to $D\cap X$. Now we require that $D$ intersects $X$
normally, or in other words $D\cap X$ is of codimension one in $X$,
and hence a divisor in $X$. In particular, this means that if $D$ is
given in local holomorphic coordinates $z_1,\ldots,z_{2k}$ by
$f(z)=0$, then the normal form $df$ in $T^*_{S^{[k]}}$ will project to
a non-zero element in $T^*_X$. Now consider
$$\omega^m\in{\mathrm H}^0(S^{[k]},\Lambda^{2m}T^*_{S^{[k]}})$$
and restrict and project to $X$ as before to get
$$\omega^m_X\in{\mathrm H}^0(X,\Lambda^{2m}T^*_X).$$
Then we can construct
$$\gamma(\alpha_T,\beta_T)\omega^m_X|_{D\cap X}\in{\mathrm
H}^{2m-1}(D\cap X,{\cal K}_X\otimes L|_{D\cap X})$$
and by the adjunction formula the canonical line bundle of $D\cap X$
is 
$${\cal K}_{D\cap X}={\cal K}_X\otimes L|_{D\cap X}.$$
Therefore the above element is really in ${\mathrm H}^{2m-1}(D\cap
X,{\cal K}_{D\cap X})$. In the case that $D\cap X$ is connected this
cohomology group is canonically isomorphic to $\Bbb C$, with the map
to the complex numbers given by contour integration in the case of {\v
C}ech cohomology and integration over all of $D\cap X$ in the case of
Dolbeault cohomology. More generally we need to sum the contributions
coming from the different connected components of $D\cap X$. This
gives us the integral of the element 
$$\gamma(\alpha_T,\beta_T)\omega^m_X|_{D\cap X}$$
which by the residue formula equals $\beta_{\gamma}$, up to the usual
normalization factor 
$$\int_X\omega^m\bar{\omega}^m|_X$$
and the addition $2\pi i$ factor. As before, note that the
Poincar{\'e} residue formula implies that if we write
$$\frac{\omega^m_X}{f}|_{D\cap X}=\frac{df}{f}\wedge\nu_{D\cap X}^{\prime}$$
locally, then the section 
$$\omega^m_X|_{D\cap X}\in{\mathrm H}^0(D\cap X,{\cal K}_X|_{D\cap
X})$$
should be replaced by the section
$$f\nu_{D\cap X}^{\prime}\in{\mathrm H}^0(D\cap X,{\cal K}_{D\cap
X}\otimes L^*|_{D\cap X})$$
in the above integral. 

This approach can be extended to the case when $D$ is a normal
crossing divisor in $S^{[k]}$, although describing the residue
$\beta_T$ in that case is slightly more complicated. On the other
hand, we only need to know $\beta_T|_{D\cap X}$ so if we assume that
$D\cap X$ in $X$ is still a smooth divisor then the calculation is no
different to before. The reason this last approach is so useful is
that it enables us to restrict the calculations to integrals over
submanifolds of small dimension, namely $D\cap X$. We then only need a
description of the Atiyah class in a neighbourhood of these
submanifolds, which potentially allows us to avoid regions in
$S^{[k]}$ where the Atiyah class becomes more difficult to
describe. This reasoning should become clearer in the calculation of
$\beta_{\Theta}$, which we now turn to.

In the case $\gamma=\Theta$ the submanifold $X$ that we should choose
is 
$$\widehat{S}=\{\xi\in
S^{[k]}|\{x_2,x_3,\ldots,x_k\}\subset{\mathrm supp}(\xi)\}$$
where $\{x_2,\ldots,x_k\}$ is an unordered set of $k-1$ distinct fixed
points chosen to lie in generic position in $S$. Let us write
$$\pi(\xi)=[x_1]+[x_2]+[x_3]+\ldots+[x_k]\in{\mathrm Sym}^kS$$
where $\pi$ is the Hilbert-Chow morphism. Here $x_1$ is an arbitrary
point in $S$, which may coincide with one of
$x_2,\ldots,x_k$, and therefore
$$\pi(\widehat{S})\cong S.$$
In $S^{[k]}$ however, a neighbourhood where $x_1$ collides with one of
the fixed points, say $x_2$, looks locally like the product of a
neighbourhood of the blown-up double point $2[x_2]$ in the Hilbert
scheme of two points 
$$S^{[2]}={\mathrm Blow}_{\Delta}(S\times S)/{\cal S}_2$$
and neighbourhoods of the remaining $k-2$ fixed points
$x_3,\ldots,x_k$. The intersection of $\widehat{S}$ with this
neighbourhood is a neighbourhood of the blow-up of $x_2$ in $S$. This
explains our notation, as $\widehat{S}$ is isomorphic to $S$ with the
$k-1$ fixed points blown-up.

From the construction of $\omega$ on $S^{[k]}$ it is clear that
restricting to $\widehat{S}$ gives us the holomorphic symplectic form
on $S$ extended over the exceptional curves $C_2,\ldots,C_k$ of the
blow-ups in some way. Since these curves have codimension one we can
ignore them in calculations of volume, and hence
\begin{eqnarray*}
\int_{\widehat{S}}\omega\bar{\omega}|_{\widehat{S}} & = &
\int_{S\backslash\{x_2,\ldots,x_k\}}\omega_S\bar{\omega}_S \\
 & = & 1.
\end{eqnarray*}
Therefore
$$\beta_{\Theta}=\int_{\widehat{S}}\Theta(\alpha_T)[\omega]|_{\widehat{S}}.$$
We claim that the right hand side is the sum 
$$8\pi^2b_{\Theta}(S)+(k-1)\delta$$
of
$$\int_{S}\Theta(\alpha_{S})[\omega_{S}]=8\pi^2b_{\Theta}(S)$$
where $\alpha_S$ is the Atiyah class of the tangent bundle of $S$, and
$k-1$ additional (identical) contributions $\delta$ to the integral
coming from the neighbourhoods $V_2,\ldots,V_k$ of $C_2,\ldots,C_k$ in
$\widehat{S}$. The point to note is that this expression is linear in
$k$.
\begin{prp}
We know that for some scalar $\beta_{\Theta}$
$$\Theta(\alpha_T)=\beta_{\Theta}[\bar{\omega}]\in{\mathrm
H}^{0,2}_{\bar{\partial}}(S^{[k]}).$$
The dependence of $\beta_{\Theta}$ on $S^{[k]}$ is that it is a linear
expression in $k$.
\end{prp}
\begin{prf}
As with our calculation of $b_{\Theta}(S)$ in Subsection $1.7$ we
prove this result in three different ways, based on the {\v C}ech,
Dolbeault, and residue descriptions of the Atiyah class. First let us
make some general comments. Recall that we had an open cover
$\{U_0,U_1,\ldots,\tilde{U}_{16}\}$ of the Kummer surface $S$. Let
$\Xi$ be the support of $\xi\in S^{[k]}$, which we think of as an
unordered set of $k$ points of $S$ with multiplicities. Define
\begin{eqnarray*}
\Delta & = & \{\xi\in S^{[k]}|\mbox{at least two points of }\Xi\mbox
{ coincide}\} \\
\Delta_2 & = & \{\xi\in S^{[k]}|\mbox{exactly two points of }\Xi\mbox
{ coincide}\}.
\end{eqnarray*}
Then $\Delta$ is a divisor in $S^{[k]}$ which we call the {\em large
diagonal\/}, and $\Delta_2$ is a dense open subset of $\Delta$ which
we call the {\em small diagonal\/}. Let $\bar{\Delta}^{(2)}$ be a
small neighbourhood of the large diagonal in $S^{[2]}$ (which is the
same as the small diagonal in $S^{[2]}$). Define open sets in
$S^{[k]}$ by
\begin{eqnarray*}
U_{\scriptsize{\underbrace{0\cdots 000}_{k}}} & = & \{\xi\in
S^{[k]}|\Xi\cap\Delta=\emptyset,\Xi\subset U_0\} \\
U_{\scriptsize{\underbrace{0\cdots 00}_{k-1}}i} & = & \{\xi\in 
S^{[k]}|\Xi\cap\Delta=\emptyset,\exists x\in\Xi\cap
U_i\mbox{ such that }\Xi\backslash\{x\}\subset U_0\} \\
U_{\scriptsize{\underbrace{0\cdots 00}_{k-1}}\tilde{i}} & = & \{\xi\in
S^{[k]}|\Xi\cap\Delta=\emptyset,\exists x\in\Xi\cap\tilde{U}_i\mbox
{ such that }\Xi\backslash\{x\}\subset U_0\} \\
U_{\scriptsize{\underbrace{0\cdots 0}_{k-2}}\bar{\Delta}} & = &
\{\xi\in S^{[k]}|\exists\xi_2\in\bar{\Delta}^{(2)}\mbox{ such that
}\Xi_2\subset\Xi,\Xi\cap\Delta\subset\Xi_2,\Xi\backslash\Xi_2\subset
U_0\}
\end{eqnarray*}
where $\Xi_2$ obviously means the support of $\xi_2$. In other words,
the first of these is the set where the support of $\xi$ consists of
$k$ distinct points in $U_0$, the second and third are the sets where
the support of $\xi$ consists of $k$ distinct points, $k-1$ of which
are in $U_0$ and the remaining one in $U_i$ or $\tilde{U}_i$
(respectively), and the last is the set where the support of $\xi$
contains two points which either collide or are very close together,
with the remaining $k-2$ being distinct points in $U_0$.

The submanifold $\widehat{S}$ of $S^{[k]}$ is contained in the union
of the above sets. Indeed
\begin{eqnarray*}
U_{0\cdots 0}\cap\widehat{S} & \cong & U_0\backslash\{x_2,\ldots,x_k\}
\\
U_{0\cdots 0i}\cap\widehat{S} & \cong & U_i \\
U_{0\cdots 0\tilde{i}}\cap\widehat{S} & \cong & \tilde{U}_i \\
U_{0\cdots 0\bar{\Delta}}\cap\widehat{S} & \cong & V_2\cup\cdots\cup V_k
\end{eqnarray*}
where the open sets $V_2,\ldots,V_k$ are the small neighbourhoods of
the exceptional curves $C_2,\ldots,C_k$ in $\widehat{S}$, and are
disjoint. Note that by choosing $x_2,\ldots,x_k$ to lie in {\em
generic\/} position in $S$ we mean that they lie in
$$U_0\backslash (U_1\cup\cdots\cup \tilde{U}_{16}).$$
It is clear that the collection
$$\{U_{0\cdots 0},U_{0\cdots 01},\ldots,U_{0\cdots
0\tilde{16}},U_{0\cdots 0\bar{\Delta}}\}$$
can be completed to an open cover of $S^{[k]}$ in such a way that the
additional sets will not intersect with $\widehat{S}$. Since we will
only perform calculations on $\widehat{S}$ (or on divisors in
$\widehat{S}$) it follows that we only need a description of the
Atiyah class over the above sets.

Let us begin with a {\v C}ech cohomology description. Recall that this
is given by taking local holomorphic connections and looking at their
differences on intersections of open sets. Take the product of $k$
copies of the flat connection $\nabla_0$, and quotient by the action
of ${\cal S}_k$. This gives the flat connection $\nabla_{0\cdots 0}$
on $U_{0\cdots 0}$ (as this open set avoids the large
diagonal). Similarly, the flat connections $\nabla_{0\cdots
0i}$ on $U_{0\cdots 0i}$ and $\nabla_{0\cdots 0\tilde{i}}$ on
$U_{0\cdots 0\tilde{i}}$ can be constructed from $k-1$ copies of 
$\nabla_0$ and a copy of the flat connections $\nabla_i$ and
$\tilde{\nabla}_i$ on $U_i$ and $\tilde{U}_i$ respectively. With
$U_{0\cdots 0\bar{\Delta}}$ we encounter a slight problem: there does
not exist a flat (or even holomorphic) connection on
$\bar{\Delta}^{(2)}$ and hence there will not exist such a connection
on $U_{0\cdots 0\bar{\Delta}}$. However, we can cover
$\bar{\Delta}^{(2)}$ with open sets which do admit holomorphic (even
flat) connections. Taking the products of these with $k-2$ copies of
$\nabla_0$, and quotienting by ${\cal S}_k$ as before, will give us
holomorphic (flat) connections on the corresponding open cover of
$U_{0\cdots 0\bar{\Delta}}$. The important feature here is that the
$k-2$ extra factors $\nabla_0$ are all {\em flat\/} connections and
are common to all of the connections over the open sets in the cover
of $U_{0\cdots 0\bar{\Delta}}$.

Next we look at the differences
$$(\alpha_T)_{(0\cdots 0)(0\cdots 0i)}=\nabla_{0\cdots 0}-\nabla_{0\cdots
0i}\in{\mathrm H}^0(U_{0\cdots 0}\cap U_{0\cdots
0i},T^*\otimes{\mathrm End}T)$$
etc. These give a {\v C}ech representative of the Atiyah class
$\alpha_T$ of $S^{[k]}$. Of course we wish to restrict to
$\widehat{S}$ so we need only consider those open sets which intersect
this submanifold. For example, consider the above section over
$U_{0\cdots 0}\cap U_{0\cdots 0i}$, which we can write in local
coordinates on either $U_{0\cdots 0}$ or $U_{0\cdots 0i}$. The
connections $\nabla_{0\cdots 0}$ and $\nabla_{0\cdots 0i}$ both
contain the same $k-1$ flat factors $\nabla_0$ which contribute
nothing to the connection one-forms in either coordinate
patch. Therefore
$$(\alpha_T)_{(0\cdots 0)(0\cdots 0i)}=p_k^*(\nabla_0-\nabla_i)$$
where $p_k$ is the projection to the $k$th factor (this makes sense at
a local level, but not globally as the ${\cal S}_k$ action permutes
the factors). Now
$$\nabla_0-\nabla_i=(\alpha_S)_{0i}\in{\mathrm H}^0(U_0\cap
U_i,T^*_S\otimes{\mathrm End}T_S)$$
gives a {\v C}ech representative of the Atiyah class on $S$. In
particular, restricting the Atiyah class $\alpha_T$ to
$$\widehat{S}\backslash (V_2\cup\cdots\cup V_k)$$
gives us the restriction of the Atiyah class $\alpha_S$ to $S$ with
neighbourhoods of the fixed points $x_2,\ldots,x_k$ removed. As these
points lie in
$$U_0\backslash (U_1\cup\cdots\cup\tilde{U}_{16})$$
their small neighbourhoods will not intersect with any two-fold
intersections of sets on the open cover
$$\{U_0,U_1,\ldots,\tilde{U}_{16}\}$$
of $S$. In other words removing these neighbourhood will not change
the {\v C}ech representative of $\alpha_S$. It follows that
restricting
$$\Theta(\alpha_T)\in{\mathrm H}^2(S^{[k]},{\cal O}_{S^{[k]}})$$
to a cohomology class in ${\mathrm H}^2(\widehat{S},{\cal
O}_{\widehat{S}})$ gives us
$$\Theta(\alpha_S)\in{\mathrm H}^2(S,{\cal O}_S)$$
away from the small neighbourhoods $V_2,\ldots,V_k$ of the exceptional
curves $C_2,\ldots,C_k$. Since the restricted holomorphic
symplectic form $\omega_{\widehat{S}}$ equals $\omega_S$ on
this region, we see that the Serre duality pairing of
$\Theta(\alpha_T)$ with $\omega_{\widehat{S}}$ is equal to the Serre
duality pairing of $\Theta(\alpha_S)$ with $\omega_S$, plus $k-1$
additional terms coming from the neighbourhoods $V_2,\ldots,V_k$. In
other words, $\beta_{\Theta}$ is equal to $8\pi^2b_{\Theta}(S)$ plus
these additional terms.

It is clear by symmetry that these $k-1$ additional terms must be
identical, but we also need them to be independent of $k$. Recall that
we covered $U_{0\cdots 0\bar{\Delta}}$ with open sets such the (flat)
holomorphic connections on these sets all looked locally like a
product of a flat connection on some open set in $\bar{\Delta}^{(2)}$
and $k-2$ copies of the flat connection $\nabla_0$. The flat
connection $\nabla_{0\cdots 0}$ on $U_{0\cdots 0}$ also factorizes
locally into a product of $k$ copies of $\nabla_0$. Therefore the
differences of the connections on the two-fold intersections of these
open sets with $U_{0\cdots 0}$ (or with themselves) will still contain
these $k-2$ flat factors which contribute nothing to the connection
one-forms in either coordinate patch. Note that $U_{0\cdots
0\bar{\Delta}}$ does not intersect $U_{0\cdots 0i}$ or $U_{0\cdots
0\tilde{i}}$ (at least not near $\widehat{S}$). Thus the Atiyah class
on these two-fold intersections of sets does not depend on how many
additional flat factors we add to the local connections. In other
words it does not vary as we change $k$. Now if we restrict to
$\widehat{S}$, we can observe that in the Serre duality pairing of
$\Theta(\alpha_T)$ with $\omega_{\widehat{S}}$ the contributions
which come from the neighbourhoods $V_2,\ldots,V_k$ will be identical
and independent of $k$. If we call these terms $\delta$, it follows
that
$$\beta_{\Theta}=8\pi^2b_{\Theta}(S)+(k-1)\delta$$
is linear in $k$.

We will now repeat the proof using the Dolbeault representative of the
Atiyah class, which is given by the $(1,1)$ part of the curvature of a
smooth global connection of type $(1,0)$ on $S^{[k]}$. Such a
connection may be obtained by patching together local connections of
type $(1,0)$ using a partition of unity. In fact, we can use the local
flat connections 
$$\nabla_{0\cdots 0},\nabla_{0\cdots 01},\ldots,\nabla_{0\cdots
0\tilde{16}}$$
from the {\v C}ech description, plus a connection of type $(1,0)$ on
$U_{0\cdots 0\bar{\Delta}}$. The latter can be constructed by taking
the product of some arbitrary type $(1,0)$ connection on
$\bar{\Delta}^{(2)}$ with $k-2$ copies of the flat connection
$\nabla_0$, and quotienting by the action of ${\cal S}_k$. As before,
we need only consider the collection of open sets
$$\{U_{0\cdots 0},U_{0\cdots 01},\ldots,U_{0\cdots
0\tilde{16}},U_{0\cdots 0\bar{\Delta}}\}$$
which intersect $\widehat{S}$.

Let $\{\psi_0,\psi_1,\ldots,\tilde{\psi}_{16}\}$ be a partition of
unity subordinate to the open cover
$\{U_0,U_1,\ldots,\tilde{U}_{16}\}$ of $S$. We can use this to
construct a partition of unity for the open cover of $S^{[k]}$. For
example, the product of $k$ copies of $\psi_0$, quotiented by the
action of ${\cal S}_k$, will give a partition function $\psi_{0\cdots
0}$ for $U_{0\cdots 0}$. Similarly, on $U_{0\cdots 0i}$ and
$U_{0\cdots 0\tilde{i}}$ we can take the product of $k-1$ copies of
$\psi_0$ with a copy of $\psi_i$ or $\tilde{\psi}_i$, quotient by the
action of ${\cal S}_k$, to get $\psi_{0\cdots 0i}$ and $\psi_{0\cdots
0\tilde{i}}$ respectively. On $U_{0\cdots 0\bar{\Delta}}$ we take the
product of $k-2$ copies of $\psi_0$ with some arbitrary partition
function $\psi_{\bar{\Delta}}$ supported in $\bar{\Delta}^{(2)}$ and
identically one on the diagonal in $S^{[2]}$. This gives
$\psi_{0\cdots 0\bar{\Delta}}$. For these to really be partition
functions on $S^{[k]}$ we may need to modify them slightly in a
neighbourhood of the large diagonal, but nowhere else.

Using this partition of unity we can patch together the local
connections in order to get a global smooth connection of type
$(1,0)$ on $S^{[k]}$. In a neighbourhood of $\widehat{S}$ this will
look like
$$\psi_{0\cdots 0}\nabla_{0\cdots 0}+\psi_{0\cdots 01}\nabla_{0\cdots
01}+\ldots+\psi_{0\cdots 0\tilde{16}}\nabla_{0\cdots
0\tilde{16}}+\psi_{0\cdots 0\bar{\Delta}}\nabla_{0\cdots
0\bar{\Delta}}$$
from which we see that the Atiyah class will look like
$$\alpha_T=A_{0\cdots 0}\bar{\partial}\psi_{0\cdots 0}+A_{0\cdots
01}\bar{\partial}\psi_{0\cdots 01}+\ldots+A_{0\cdots
0\tilde{16}}\bar{\partial}\psi_{0\cdots 0\tilde{16}}+\bar{\partial}(A_{0\cdots
0\bar{\Delta}}\psi_{0\cdots 0\bar{\Delta}})$$
where $A_{0\cdots 0}$ is the connection one-form of the local
connection $\nabla_{0\cdots 0}$, etc. Restrict to $\widehat{S}$ and
ignore the last term which is only supported in the neighbourhoods
$V_2,\ldots,V_k$. Let $\xi$ in $S^{[k]}$ have support
$\{x_1,x_2,\ldots,x_k\}$ and lie on
$$\widehat{S}\backslash (V_2\cup\ldots\cup V_k).$$
In other words $x_1$ lies in $S$ but not in a neighbourhood of any of
the points $x_2,\ldots,x_k$. Then
\begin{eqnarray*}
\psi_{0\cdots 0}(\xi) & = & \psi_0(x_1) \\
\psi_{0\cdots 0i}(\xi) & = & \psi_i(x_1) \\
\psi_{0\cdots 0\tilde{i}}(\xi) & = & \tilde{\psi}_i(x_1)
\end{eqnarray*}
as $\psi_0$ is identically one on the neighbourhoods of the points
$x_2,\ldots,x_k$. Similarly, the connection one forms at $\xi$ are
given by
\begin{eqnarray*}
A_{0\cdots 0}(\xi) & = & p_1^*A_0(x_1) \\
A_{0\cdots 0i}(\xi) & = & p_1^*A_i(x_1) \\
A_{0\cdots 0\tilde{i}}(\xi) & = & p_1^*\tilde{A}_i(x_1)
\end{eqnarray*}
where $p_1$ is the projection to the first factor (which makes sense
locally, though not globally due to the action of the symmetric group
permuting the factors). Therefore
\begin{eqnarray*}
\alpha_T|_{\widehat{S}\backslash (V_2\cup\cdots\cup V_k)} & = &
(p_1^*A_0\bar{\partial}\psi_0+p_1^*A_1\bar{\partial}\psi_1+\ldots+p_1^*A_{\tilde{16}}\bar{\partial}\psi_{\tilde{16}})(x_1)
\\
 & = & p_1^*\alpha_S(x_1)
\end{eqnarray*}
where $\alpha_S$ is the Dolbeault representative of the Atiyah class
on $S$. Since this is supported in a neighbourhood of the sixteen
exceptional curves $D_1,\ldots,D_{16}$ in $S$, removing the
neighbourhoods of the points $x_2,\ldots,x_k$ (where $\alpha_S$
vanishes anyway) will have no effect. Therefore
\begin{eqnarray*}
\int_{\widehat{S}\backslash (V_2\cup\cdots\cup
V_k)}\Theta(\alpha_T)[\omega]|_{\widehat{S}} & = &
\int_S\Theta(\alpha_S)[\omega_S] \\
 & = & 8\pi^2b_{\Theta}(S)
\end{eqnarray*}
as we have already seen that $\omega$ restricted to $\widehat{S}$ is
just $\omega_S$ away from the blow-ups of $x_2,\ldots,x_k$. As before,
this just leaves the terms
$$\int_{V_2}\Theta(\alpha_T)[\omega]|_{\widehat{S}}+\ldots+\int_{V_k}\Theta(\alpha_T)[\omega]|_{\widehat{S}}$$
which we can assume are identical by symmetry. It only remains to
prove that each of these integrals is independent of $k$. As in the
{\v C}ech cohomology case, this follows from the observation that
$\nabla_{0\cdots 0\bar{\Delta}}$ contains $k-2$ copies of the {\em
flat\/} connection $\nabla_0$ as factors, which contribute nothing to
the connection one-form $A_{0\cdots 0\bar{\Delta}}$. Also the
partition function $\psi_{0\cdots 0\bar{\Delta}}$ contains $k-2$
factor equal to $\psi_0$ which are identically one on
$V_2,\ldots,V_k$ (since $\psi_0$ is identically one on the
neighbourhoods of the fixed points $x_2,\ldots,x_k$). In other words,
over $V_2,\ldots,V_k$ the Atiyah class $\alpha_T$ does not vary as we
change $k$, as this merely involves adding `redundant' factors. In
particular
$$\int_{V_j}\Theta(\alpha_T)[\omega]|_{\widehat{S}}=\delta$$
is independent of $k$ (for $j=2,\ldots,k$), and therefore
\begin{eqnarray*}
\beta_{\Theta} & = &
\int_{\widehat{S}}\Theta(\alpha_T)[\omega]|_{\widehat{S}} \\
 & = & \int_{\widehat{S}\backslash (V_2\cup\cdots\cup
V_k)}\Theta(\alpha_T)[\omega]|_{\widehat{S}}+\int_{V_2}\Theta(\alpha_T)[\omega]|_{\widehat{S}}+\ldots+\int_{V_k}\Theta(\alpha_T)[\omega]|_{\widehat{S}}
\\
 & = & 8\pi^2b_{\Theta}(S)+(k-1)\delta
\end{eqnarray*}
is linear in $k$.

Finally, we show how this result can be proved using the residue
approach. Recall that on the Kummer surface $S$ we took the flat
connection $\nabla_0$ on $U_0$ and observed that it could be extended
to a global meromorphic connection with a simple pole along the smooth
divisor
$$D=D_1+\ldots+D_{16}$$
where $D_1,\ldots,D_{16}$ are the sixteen exceptional curves in
$S$. In the case of $S^{[k]}$ we take the flat connection
$\nabla_{0\cdots 0}$ on $U_{0\cdots 0}$. Note that this open set is
dense in $S^{[k]}$ and its complement is the union of the large
diagonal $\Delta$ and the divisors
$$E_i=\{\xi\in S^{[k]}|\exists x\in\Xi\cap D_i\}.$$
In fact, $\nabla_{0\cdots 0}$ can be extended to a global meromorphic
connection on $S^{[k]}$ with a simple pole along the divisor
$$E=E_1+\ldots+E_{16}+\Delta.$$
Although this is a normal crossing divisor, its intersection
$E\cap\widehat{S}$ with $\widehat{S}$ is a smooth divisor in this
submanifold. Indeed
\begin{eqnarray*}
E\cap\widehat{S} & = &
E_1\cap\widehat{S}+\ldots+E_{16}\cap\widehat{S}+\Delta\cap\widehat{S}
\\
 & = & D_1+\ldots+D_{16}+C_2+\ldots+C_k
\end{eqnarray*}
and the terms in the last line are all smooth curves in
$\widehat{S}$. Consider how this connection looks in a neighbourhood
of
$$\xi\in D_i\subset\widehat{S}\subset S^{[k]}.$$
Let the support of $\xi$ be $\Xi=\{x_1,x_2,\ldots,x_k\}$. The points
$x_2,\ldots,x_k$ lie in $U_0$ over which $\nabla_0$ is holomorphic,
whereas $x_1$ lies on $D_i$ along which $\nabla_0$ has a simple
pole. Therefore the residue of $\nabla_{0\cdots 0}$ at $\xi$ is simply
the residue of $\nabla_0$ at $x_1$, as the $k-1$ additional
holomorphic factors do not contribute to the residue. Thus
$$\beta_T|_{D_i}=p_1^*\beta_S|_{D_i}$$
where $p_1$ is the (locally well-defined) projection to the first
factor and $\beta_S$ is the residue of the meromorphic connection
$\nabla_0$ on $S$. Consider also how the connection $\nabla_{0\cdots
0}$ looks in a neighbourhood of
$$\xi\in C_2\subset\widehat{S}\subset S^{[k]}.$$ 
Writing the support of $\xi$ as above, we see that we can still say
that $\nabla_0$ is holomorphic in neighbourhood of the $k-2$ distinct
points $x_3,\ldots,x_k$. Hence these $k-2$ factors do not contribute
to the residue, which therefore does not vary with $k$. Similarly for
neighbourhoods of $\xi$ in $C_3,\ldots,C_k$.

Now the residue calculation of $\beta_{\Theta}$ involves an integral
of
$$\Theta(\alpha_T,\beta_T)\omega_{\widehat{S}}|_{E\cap\widehat{S}}$$
which is either a contour integral if we use {\v C}ech cohomology or
an integral over $E\cap\widehat{S}$ if we use Dolbeault
cohomology. Firstly, from what we know about the forms of the Atiyah
class $\alpha_T$ (in either {\v C}ech or Dolbeault descriptions), the
holomorphic symplectic form $\omega$, and the residue $\beta_T$ when
they are restricted to $\widehat{S}$, we can observe that the
integrand restricted to
$$D=D_1+\ldots+D_{16}$$
is just
$$\Theta(\alpha_S,\beta_S)\omega_S|_D$$
and restricted to
$$\Delta\cap\widehat{S}=C_2+\ldots+C_k$$
gives $k-1$ terms which are independent of $k$, and identical by
symmetry. Therefore, in Dolbeault cohomology for example
\begin{eqnarray*}
\beta_{\Theta} & = & 2\pi
i\int_{E\cap\widehat{S}}\Theta(\alpha_T,\beta_T)\omega_{\widehat{S}}|_{E\cap\widehat{S}}
\\
 & = & 2\pi i\int_D\Theta(\alpha_S,\beta_S)\omega_S|_D+(k-1)\delta \\
 & = & 8\pi^2b_{\Theta}(S)+(k-1)\delta
\end{eqnarray*}
and similarly if we use {\v C}ech cohomology. So again we see that
$\beta_{\Theta}$ is linear in $k$.\end{prf}
Since we know $b_{\Theta}(S)=48$ and $b_{\Theta^2}(S^{[2]})=3600$,
this is enough to determine
$$\beta_{\Theta}=(8\pi^2)12(k+3)$$
and hence
\begin{prp}
For the Hilbert scheme of $k$ points on a K$3$ surface $S$ we have
$$b_{\Theta^k}(S^{[k]})=12^k(k+3)^k$$
and therefore
$$\int_{S^{[k]}}{\mathrm Td}^{1/2}_k(S^{[k]})=\frac{(k+3)^k}{4^kk!}.$$
\end{prp}
We can verify that this gives $b_{\Theta^3}(S^{[3]})=373248$, which is
in agreement with the value we have already obtained from the Chern
numbers. Also, knowing $b_{\Theta^4}(S^{[4]})$ allows us to determine
$s$ and hence complete our tables of Chern numbers and Rozansky-Witten
invariants for $k\leq 4$ (see Appendices D and E$.1$). 

We stated earlier that there are no explicit formulae for the Chern
numbers of these Hilbert schemes in all dimensions. However, there are
some implicit formulae for the Chern numbers of Hilbert schemes of
points on an arbitrary complex surface due to Ellingsrud,
G{\"o}ttsche, and Lehn~\cite{egl99}, and these methods have enabled
G{\"o}ttsche~\cite{gottsche} to calculate the Chern numbers of
$S^{[k]}$ for $k\leq 6$ using the Bott residue formula. It is
reassuring that these values agree with those obtained above for
$k\leq 4$. Furthermore, when $k$ is $5$ and $6$ we can use
G{\"o}ttsche's calculations to compare our values for
$b_{\Theta^5}(S^{[5]})$ and $b_{\Theta^6}(S^{[6]})$ as calculated
above with the value calculated by using the formulae for
$b_{\Theta^5}$ and $b_{\Theta^6}$ in terms of Chern numbers. Once
again we find the values agree.

In order to extend our calculations to higher dimensions we need to
know something about the form of $\beta_{\gamma}$ in the case that
$\gamma$ is a trivalent graph with $2m>2$ vertices, ie.\ not
simply just $\Theta$. In fact there is strong evidence to suggest that
$\beta_{\gamma}$ is linear in $k$ for these cases also. Indeed the
calculations we have already done verify this for $\beta_{\Theta_2}$
with $k\leq 4$. Furthermore, we can calculate a few more values of
Rozansky-Witten invariants by using G{\"o}ttsche's values for the
Chern numbers for $k\leq 6$, and this allows us to verify the
linearity of $\beta_{\gamma}$ in several other cases.

If we attempt to prove the linearity of $\beta_{\gamma}$ using either
the {\v C}ech or Dolbeault descriptions of the Atiyah class things
soon become complicated: the submanifold $X$ in $S^{[k]}$ which we
need to perform calculations on will have real-dimension $4m$ and be
at least as complicated as $S^{[m]}$. In the $\beta_{\Theta}$ case
performing calculations on $\widehat{S}$ was just manageable. For
$m>1$, the approach most likely to succeed would probably involve
meromorphic connections. However, there does not appear to be a way to
reduce the calculation to one over a smooth divisor, so we ultimately
need a good description of the residue of a meromorphic connection
with a simple pole along a normal crossing divisor. Note that in the
case of connections with {\em logarithmic\/} singularities such a
theory is well-developed, and Kapranov has used this approach to
derive quite explicit formulae for the Rozansky-Witten invariants of a
hyperk{\"a}hler manifold which admits such a connection
(see~\cite{kapranov98}). These formulae generalize those of
Ohtsuki~\cite{ohtsuki82} for Chern numbers. Unfortunately, it seems
that such a connection does not exist on the Hilbert scheme of points
on a K$3$ surface (for example, the poles of $\nabla_0$ on $S$ are
certainly not logarithmic).

\subsection{Generalized Kummer varieties}

The Hilbert scheme of $k$ points on a torus $T$ has a
hyperk{\"a}hler metric, but these spaces are not irreducible. There is
a map $p$ to the torus itself given by composing the Hilbert-Chow
morphism with addition (since the torus is an abelian surface, we can
add $k$ points together and the result will lie in $T$). The fibres of
this map are all isomorphic hyperk{\"a}hler manifolds of dimension
$4(k-1)$, and are also irreducible. We denote them by $T^{[[k-1]]}$,
so that $T^{[[k]]}$ would be the $4k$-dimensional space obtained by
taking the fibres in $T^{[k+1]}$. This is the second main family of
examples of compact hyperk{\"a}hler manifolds, and they are known as
generalized Kummer varieties (as $T^{[[1]]}$ is simply the Kummer
model of a K$3$ surface).

The Hodge numbers of these spaces were calculated by G{\"o}ttsche and
Soergel~\cite{gs93}. In particular, the $\chi_y$-genus is given by
$$\chi_y(T^{[[k]]})=(k+1)\sum_{d|(k+1)}d^3(1-y+y^2-\ldots+(-y)^{\frac{k+1}{d}-1})^2(-y)^{k+1-\frac{k+1}{d}}.$$
For $k\leq 4$ these polynomials are written out in Appendix C. Once
again we can substitute these values into the Riemann-Roch formula to
give the results in Appendix D, from which we can determine all the
Rozansky-Witten invariants for $k=1$, $2$, and $3$. The next step is
to compute $b_{\Theta^k}(T^{[[k]]})$, and in particular
$b_{\Theta^4}(T^{[[4]]})$ in order to determine $s$, and hence
complete the tables up to $k=4$. 

To begin with let us normalize the holomorphic symplectic form on $T$
so that 
$$\int_T\omega\bar{\omega}=1.$$
The same calculation as before shows that the induced holomorphic
symplectic form on $T^{[k+1]}$, also denoted $\omega$, will have
normalization
$$\int_{T^{[k+1]}}\omega^{k+1}\bar{\omega}^{k+1}=(k+1)!$$
The generalized Kummer variety $T^{[[k]]}$ sits inside this Hilbert
scheme as a fibre of the fibration over $T$. Denote by $i$ the
inclusion of the fibre; thus we have
$$\begin{array}{ccl}
T^{[[k]]} & \stackrel{i}{\hookrightarrow} & T^{[k+1]} \\
 & & \downarrow p \\
 & & T.
\end{array}$$
We wish to calculate the volume of the generalized Kummer variety, or
rather we wish to know
$$\int_{T^{[[k]]}}\omega^k\bar{\omega}^k.$$
Recall that the holomorphic symplectic form $\omega$ on the Hilbert
scheme $T^{[k+1]}$ is given by
$$\pi^*(p_0^*\omega+\ldots+p_k^*\omega)$$
where $p_j$ is the projection onto the $(j-1)$th factor in $T^{k+1}$,
the sum descends to the symmetric product ${\mathrm Sym}^{k+1}T$, and
then we pull-back to the Hilbert scheme using the Hilbert-Chow
morphism $\pi$. The holomorphic symplectic form on the generalized
Kummer variety $T^{[[k]]}$, also denoted $\omega$, is defined by the
restriction of $\omega_{T^{[k+1]}}$. According to Subsection $2.4$ the
volume form on $T^{[k+1]}$ is given by
\begin{eqnarray*}
V_{T^{[k+1]}} & = &
 \frac{1}{2^{2(k+1)}((k+1)!)^2}\omega^{k+1}\bar{\omega}^{k+1} \\
 & = & \frac{1}{2^{2(k+1)}((k+1)!)^2}\pi^*(p_0^*\omega+\ldots+
p_k^*\omega)^{k+1}\pi^*(p_0^*\bar{\omega}+\ldots+p_k^*\bar{\omega})^{k+1}
\\
 & = & \frac{1}{2^{2(k+1)}}\pi^*(p_0^*\omega\cdots p_k^*\omega
 p_0^*\bar{\omega}\cdots p_k^*\bar{\omega})
\end{eqnarray*}
on $T^{[[k]]}$ by
\begin{eqnarray*}
V_{T^{[[k]]}} & = & \frac{1}{2^{2k}(k!)^2}\omega^k\bar{\omega}^k \\
 & = & \frac{1}{2^{2k}(k!)^2}\pi^*(p_0^*\omega+\ldots+
p_k^*\omega)^k\pi^*(p_0^*\bar{\omega}+\ldots+p_k^*\bar{\omega})^k \\
 & = & \frac{1}{2^{2k}}\pi^*\left(\sum_{j=0}^k
 p_0^*\omega\cdots\widehat{p_j^*\omega}\cdots
 p_k^*\omega\right)\pi^*\left(\sum_{j=0}^k
 p_0^*\bar{\omega}\cdots\widehat{p_j^*\bar{\omega}}\cdots
 p_k^*\bar{\omega}\right)
\end{eqnarray*}
where $\widehat{p_j^*\omega}$ and $\widehat{p_j^*\bar{\omega}}$ mean
that we omit these terms, and on $T$ by
$$V_T=\frac{1}{2^2}\omega\bar{\omega}.$$
We wish to compare the volume form on $T^{[k+1]}$ with the one induced
by the fibration from the volume forms on the fibre $T^{[[k]]}$ and
the base $T$, ie.\ we wish to compare
$$V_{T^{[k+1]}}\qquad\mbox{and}\qquad i_*V_{T^{[[k]]}}\wedge p^*V_T.$$
Consider the map $p_0+\ldots+p_k$ from $T^{k+1}$ to $T$. Clearly it is
symmetric and so descends to a map from ${\mathrm Sym}^{k+1}T$ to $T$.
In fact $p$ is given by composing the Hilbert-Chow morphism $\pi$ with
$p_0+\ldots+p_k$. Now
$$p^*V_T=\frac{1}{2^2}p^*\omega p^*\bar{\omega}$$
and $p^*\omega$ is `almost'
$\pi^*(p_0^*\omega+\ldots+p_k^*\omega)=\omega_{T^{[k+1]}}$ (similarly
for $p^*\bar{\omega}$). Indeed given two vector fields $v$ and $u$ on
$T^{[k+1]}$, we find
\begin{eqnarray*}
p^*\omega(v,u) & = & \omega(p_*v,p_*u) \\
 & = &
 \omega((p_{0*}+\ldots+p_{k*})\pi_*v,(p_{0*}+\ldots+p_{k*})\pi_*u) \\
 & = & \sum_{i,j=0}^k \omega(p_{i*}\pi_*v,p_{j*}\pi_*u) \\
 & = & \sum_{j=0}^k \omega(p_{j*}\pi_*v,p_{j*}\pi_*u)+\sum_{i\neq
 j}\omega(p_{i*}\pi_*v,p_{j*}\pi_*u) \\
 & = & \sum_{j=0}^k \pi^*p_j^*\omega(v,u)+\sum_{i\neq
 j}\omega(p_{i*}\pi_*v,p_{j*}\pi_*u) \\
 & = & \omega_{T^{[k+1]}}(v,u)+\sum_{i\neq
 j}\omega(p_{i*}\pi_*v,p_{j*}\pi_*u) 
\end{eqnarray*}
However, when we calculate $i_*V_{T^{[[k]]}}\wedge p^*V_T$ the second
term above gives zero contribution (this can easily be seen in local
coordinates) and therefore
\begin{eqnarray*}
i_*V_{T^{[[k]]}}\wedge p^*V_T & = &
\frac{1}{2^{2k}(k!)^2}\omega^k\bar{\omega}^k\frac{1}{2^2}\omega\bar{\omega}
\\
 & = & \frac{1}{2^{2(k+1)}(k!)^2}\omega^{k+1}\bar{\omega}^{k+1} \\
 & = & (k+1)^2V_{T^{[k+1]}}.
\end{eqnarray*}
It follows that
\begin{eqnarray*}
\frac{1}{2^2}\frac{1}{2^{2k}(k!)^2}\int_{T^{[[k]]}}\omega^k\bar{\omega}^k
& = & {\mathrm vol}(T){\mathrm vol}(T^{[[k]]}) \\
 & = & (k+1)^2{\mathrm vol}(T^{[k+1]}) \\
 & = &
 (k+1)^2\frac{1}{2^{2(k+1)}((k+1)!)^2}\int_{T^{[k+1]}}\omega^{k+1}\bar{\omega}^{k+1} \\
 & = & \frac{1}{2^2}\frac{1}{2^{2k}(k!)^2}(k+1)!
\end{eqnarray*}
and hence
$$\int_{T^{[[k]]}}\omega^k\bar{\omega}^k=(k+1)!$$

As in the last subsection, if $\Gamma$ decomposes into connected
components as $\gamma_1\cdots\gamma_j$ then
\begin{eqnarray*}
b_{\Gamma}(T^{[[k]]}) & = & \frac{1}{(8\pi^2)^kk!}\beta_{\gamma_1}\cdots\beta_{\gamma_j}\int_{T^{[[k]]}}\omega^k\bar{\omega}^k
\\
 & = & \frac{(k+1)}{(8\pi^2)^k}\beta_{\gamma_1}\cdots\beta_{\gamma_j}
\end{eqnarray*}
with
$$\gamma_i(\alpha_T)=\beta_{\gamma_i}[\bar{\omega}^{m_i}]\in{\mathrm
H}^{0,2m_i}_{\bar{\partial}}(T^{[[k]]})$$
where $\gamma_i$ has $2m_i$ vertices, and the scalars
$\beta_{\gamma_i}$ depend on $T^{[[k]]}$, in other words on $k$. When
$\gamma$ is the graph $\Theta$ we can show that this dependence is
linear.
\begin{prp}
We know that for some scalar $\beta_{\Theta}$
$$\Theta(\alpha_T)=\beta_{\Theta}[\bar{\omega}]\in{\mathrm
H}^{0,2}_{\bar{\partial}}(T^{[[k]]}).$$
The dependence of $\beta_{\Theta}$ on $T^{[[k]]}$ is that it is a
linear expression in $k$.
\end{prp}
\begin{prf}
The proof is similar to the proof of Proposition $18$ and will be
omitted.\end{prf}
From $b_{\Theta}(T^{[[1]]})=48$ and $b_{\Theta^2}(T^{[[2]]})=3888$ we
can determine
$$\beta_{\Theta}=(8\pi^2)12(k+1)$$
and hence
\begin{prp}
For the generalized Kummer variety $T^{[[k]]}$ we have
$$b_{\Theta^k}(T^{[[k]]})=12^k(k+1)^{k+1}$$
and therefore
$$\int_{T^{[[k]]}}{\mathrm Td}^{1/2}_k(T^{[[k]]})=\frac{(k+1)^{k+1}}{4^kk!}.$$
\end{prp}
We can then verify $b_{\Theta^3}(T^{[[3]]})=442368$, and use
$b_{\Theta^4}(T^{[[4]]})$ to determine $s$ and complete our tables of
Chern numbers and Rozansky-Witten invariants for $k\leq 4$ (see
Appendices D and E$.1$).

Unlike for the Hilbert schemes of points on a K$3$ surface, these values
for the Chern numbers of the generalized Kummer varieties appear to be
completely new. In other words, this method of using the
Rozansky-Witten invariants to determine the Chern numbers for $k=4$
appears to give results not previously known.

For trivalent graphs $\gamma$ with $2m>2$ vertices we can probably
also expect $\beta_{\gamma}$ to be linear in $k$ for generalized
Kummer varieties. As with the Hilbert schemes of points on a K$3$
surface there is some evidence to support this, and a proof will most
likely follow from a generalization of the residue approach to
meromorphic connections with simple poles along normal crossing
divisors.

\subsection{Cobordism classes}

Up to now in this chapter we have been dealing only with irreducible
hyperk{\"a}hler manifolds. For reducible manifolds we can use the
product formula. Since we know the values of all the Rozansky-Witten
invariants up to $k=4$ for the Hilbert schemes of points on a K$3$
surface and for the generalized Kummer varieties, we can therefore
calculate the Rozansky-Witten invariants for all hyperk{\"a}hler
manifolds constructed by taking products of these (ie.\ this includes
reducible manifolds of arbitrarily large dimension). Some examples are
given in Appendix E$.2$. It follows that
we can, in principle, calculate the characteristic numbers for these
products as well. Indeed up to real-dimension $20$ we have explicit
formulae for the Chern numbers in terms of Rozansky-Witten invariants
by using the expansions of polywheels as given in Appendix A$.1$.

We know that in dimensions $4$, $8$, and $12$ that the Rozansky-Witten
invariants are all characteristic numbers. However, in
dimension $16$ we only know that $b_{\Theta^4}$ and
$b_{\Theta^2\Theta_2}$ are characteristic numbers. At this stage we
can only write the remaining invariants as rational functions of the
characteristic numbers (for irreducible hyperk{\"a}hler manifolds). In
this subsection we will prove that they are not characteristic
numbers, ie.\ they cannot be written as a {\em linear combinations\/}
of Chern numbers. It suffices to show this for $b_{\Theta^2_2}$, as
the other invariants can be written as linear combinations of
$b_{\Theta^2_2}$ and Chern numbers (see Appendix A$.2$).

All we shall actually do is exhibit an example of two hyperk{\"a}hler
manifolds which have the same Chern numbers but different values of
$b_{\Theta^2_2}$. Recall that in the hyperk{\"a}hler case, the Chern
numbers are equivalent to the Pontryagin numbers. By a theorem of
Thom~\cite{thom54} (see also Hirzebruch's book~\cite{hirzebruch78})
two manifolds have the same Pontryagin numbers if and only if they
represent the same class in the oriented rational cobordism ring
``modulo torsion''. We will ignore subtleties arising from torsion
elements. Two oriented manifolds $X$ and $Y$ of real-dimension $4k$
are cobordant if there exists an oriented manifold $M$ of
real-dimension $4k+1$ with boundary $\partial M=X+(-Y)$, ie.\ the
disjoint sum of $X$ and the manifold $Y$ with its orientation
reversed. The oriented rational cobordism ring is the ring of
cobordism classes over the rational numbers, with multiplication given
by taking products of manifolds and addition given by disjoint sums.
Stated in these terms, we wish to find two cobordant hyperk{\"a}hler
manifolds which are nonetheless distinguished by the Rozansky-Witten
invariants.

The construction is simple enough. We consider all hyperk{\"a}hler
manifolds in dimension $16$ obtained by taking products of Hilbert
schemes of points on a K$3$ surface. They are:
$$S^{[4]}\qquad S\times S^{[3]}\qquad S^{[2]}\times S^{[2]}\qquad
S^2\times S^{[2]}\qquad S^4$$
Now we take a linear combinations of these (ie.\ a disjoint
sum) and choose the rational coefficients in such a way that the
resulting (disconnected, reducible) hyperk{\"a}hler manifold has the
same Chern numbers as the generalized Kummer variety $T^{[[4]]}$. We
find that
$$X=7S^{[4]}-\frac{49}{8}S\times S^{[3]}-3S^{[2]}\times
S^{[2]}+\frac{67}{12}S^2\times S^{[2]}-\frac{21}{16}S^4$$
is the required manifold. However, calculating the value of
$b_{\Theta^2_2}$ we find that
\begin{eqnarray*}
b_{\Theta^2_2}(X) & = & 278784 \\
 & \neq & 288000 \\
 & = & b_{\Theta^2_2}(T^{[[4]]})
\end{eqnarray*}
so this Rozansky-Witten invariant distinguishes the two manifolds and
therefore cannot be a linear combination of Chern numbers. 
\begin{thm}
The Rozansky-Witten invariant $b_{\Theta^2_2}$ is not a characteristic
number on hyperk{\"a}hler manifolds of real-dimension $16$. However,
it can be expressed as a rational function of characteristic numbers
on irreducible hyperk{\"a}hler manifolds.
\end{thm}
Note that the linear combination of manifolds in $X$ involves both
rational and negative coefficients. Instead we can rearrange things so
that we have
$$336S^{[4]}+268S^2\times S^{[2]}\sim
48T^{[[4]]}+294S^{[3]}+144S^{[2]}\times S^{[2]}+63S^4$$
where $\sim$ means the two manifolds are cobordant (in this case, we
have two hyperk{\"a}hler manifolds in the same integral cobordism
class). The value of $b_{\Theta^2_2}$ on each manifold is $19795968$
and $19353600$ respectively, and hence the Rozansky-Witten invariants
distinguish these two manifolds.

%% file: ch6.tex
\section{An invariant of knots}

\subsection{Vector bundles on hyperk{\"a}hler manifolds}

Up to now we have studied exclusively the weight system on graph
homology which gives rise to invariants of three-manifolds. Invariants
of knots and links also arise naturally in the Rozansky-Witten theory.
In this final chapter we wish to construct explicitly a weight system
on chord diagrams from a hyperk{\"a}hler manifold with a
holomorphic vector bundle over it. This weight system is a natural
extension of the Rozansky-Witten weight system on graph homology, and
leads to potentially new finite-type invariants of knots (and
links). We begin in this subsection by reviewing the subject of
holomorphic vector bundles.

Let $E$ be a smooth complex vector bundle of complex-rank $r$ on a
compact hyperk{\"a}hler manifold $X$ of real-dimension $4k$. Once
again, we wish to use the techniques of complex geometry so we will
choose a specific complex structure $I$ from the space of compatible
complex structures on $X$, and regard $X$ as a complex manifold with
respect to this choice. A connection $\nabla$ on $E$ splits into the
sum of two differential operators
$$\partial_A:\Omega^{0,0}(E)\rightarrow\Omega^{1,0}(E)$$
and
$$\bar{\partial}_A:\Omega^{0,0}(E)\rightarrow\Omega^{0,1}(E).$$
Suppose the connection looks like $d+A$ locally, where $A$ is an
${\mathrm End}E$-valued one-form. We can write
$$A=A^{1,0}+A^{0,1}$$
with $A^{p,q}\in\Omega^{p,q}({\mathrm End}E)$, and then
$$\partial_A=\partial +A^{1,0}$$
and
$$\bar{\partial}_A=\bar{\partial}+A^{0,1}$$
locally. The curvature $R=\nabla\circ\nabla$ of the connection is an
${\mathrm End}E$-valued two-form which can also be decomposed into 
$$R=R^{2,0}+R^{1,1}+R^{0,2}$$
with $R^{p,q}\in\Omega^{p,q}({\mathrm End}E)$. Indeed
\begin{eqnarray*}
R^{2,0} & = & \partial_A\circ\partial_A \\
R^{1,1} & = &
\partial_A\circ\bar{\partial}_A+\bar{\partial}_A\circ\partial_A \\
R^{0,2} & = & \bar{\partial}_A\circ\bar{\partial}_A.
\end{eqnarray*}

If $E$ is a holomorphic vector bundle we can choose a connection
$\nabla$ which is compatible with, or preserves, the complex
structure. This means that relative to a local holomorphic
trivialization 
$$\bar{\partial}_A=\bar{\partial}$$
or $A^{0,1}=0$. For such a choice, the $(0,2)$ part of the curvature
vanishes, ie.\
$$R^{0,2}=\bar{\partial}_A\circ\bar{\partial}_A=0.$$
Conversely, if a complex vector bundle admits such a connection then
it is holomorphic. Note that $\bar{\partial}$ gives the holomorphic
structure on $E$ and the connection $\nabla$ is the sum of the
differential operator $\partial_A$ of type $(1,0)$ and this
holomorphic structure.

An {\em Hermitian structure\/} $h$ in a smooth complex vector bundle
$E$ is a smooth field of Hermitian inner products in the fibres of
$E$, and we call $(E,h)$ an {\em Hermitian vector bundle\/}. A
connection $\nabla$ is a $h$-{\em connection\/} if it preserves the
Hermitian structure $h$, or in other words
$$d(h(s,t))=h(\nabla s,t)+h(s,\nabla t)$$
where $s$ and $t$ are (local) sections of $E$. For a holomorphic
vector bundle $E$ with a Hermitian structure $h$ there is a unique
$h$-connection $\nabla$ which is compatible with the complex
structure, and we call this the {\em Hermitian connection\/}. For
example, recall that the Levi-Civita connection is the unique
connection on the (holomorphic) tangent bundle $T$ which preserves the
metric and is torsion-free. In the K{\"a}hler case, the Levi-Civita
connection is precisely the Hermitian connection for the K{\"a}hler
metric. 

For the Hermitian connection, the curvature is purely of type $(1,1)$
$$R\in\Omega^{1,1}({\mathrm End}E)$$
and locally it can be expressed as $\bar{\partial}A^{1,0}$. In local
complex coordinates it has components
$$R^{I}_{\phantom{I}Jk\bar{l}}dz_k\wedge d\bar{z}_l$$
where $I$ and $J$ refer to a local basis of sections of $E$, and
summation over repeated indices is assumed.

Recall from Subsection $1.5$ that the Atiyah class 
$$\alpha_E\in{\mathrm H}^1(X,T^*\otimes{\mathrm End}E)$$
of a holomorphic vector bundle $E$ is the obstruction to the existence
of a global holomorphic connection on $E$. In the Dolbeault model for
cohomology it can be represented by the $(1,1)$ part of the
curvature of a smooth global connection $\nabla$ of type $(1,0)$ on
$E$ (ie. a complex structure preserving connection). Recall that for
such a connection the $(0,2)$ part of the curvature vanishes, and
$$R=R^{2,0}+R^{1,1}.$$
The Bianchi identity tells us that $\nabla R=0$. Now
\begin{eqnarray*}
\nabla R & = & (\partial_A+\bar{\partial})(R^{2,0}+R^{1,1}) \\
         & = & \partial_A R^{2,0}+(\partial_A R^{1,1}+\bar{\partial}R^{2,0})+\bar{\partial}R^{1,1}
\end{eqnarray*}
where the three terms in the final line are of type $(3,0)$, $(2,1)$,
and $(1,2)$ respectively, and hence each must vanish. In particular,
$R^{1,1}$ is $\bar{\partial}$-closed, as we require for it to
represent a Dolbeault cohomology class
$$[R^{1,1}]\in{\mathrm H}_{\bar{\partial}}^{1,1}(X,{\mathrm
End}E)={\mathrm H}_{\bar{\partial}}^{0,1}(X,T^*\otimes{\mathrm
End}E).$$
In the more specialized case that we have an Hermitian structure $h$
on $E$ we know that the curvature $R$ of the (unique) Hermitian
connection is purely of type $(1,1)$, and hence is itself a
representative of the Atiyah class $\alpha_E$.

At the beginning of Chapter $2$ we discussed how to represent (in
Dolbeault cohomology) the Chern character $ch(T)$ of the tangent
bundle to $X$ by taking the trace of powers of the Riemann curvature
tensor. More generally, if $E$ is a smooth complex bundle with a
connection $\nabla$, then its Chern character can be expressed as
\begin{eqnarray*}
ch(E) & = & r+ch_1(E)+ch_2(E)+\ldots +ch_{2k}(E) \\
      & = & \sum_{m=0}^{2k}\frac{(-1)^m}{m!(2\pi i)^m}[{\mathrm
      Tr}(R^m)]
\end{eqnarray*}
where
$$R\in\Omega^2({\mathrm End}E)$$
is the curvature of $\nabla$, its powers are obtained by composing in
${\mathrm End}E$ and taking the wedge product of forms, and the trace
is in ${\mathrm End}E$. Thus the cohomology classes which ${\mathrm
Tr}(R^m)$ represent are independent of the connection chosen. If $E$
is a holomorphic vector bundle with an Hermitian structure $h$, then
$$R\in\Omega^{1,1}({\mathrm End}E)$$
for the Hermitian connection on $E$. Therefore the component $ch_m(E)$
of the Chern character is of pure Hodge type $(m,m)$. Furthermore,
since $R$ represents the Atiyah class, we can replace $[R]$ by
$\alpha_E$ in the above formula. More specifically, ${\mathrm
Tr}(\alpha_E^m)$ is given by composing in ${\mathrm End}E$ and taking
tensor products, and $ch_m(E)$ is given by projecting
$$\frac{(-1)^m}{m!(2\pi i)^m}{\mathrm Tr}(\alpha_E^m)\in{\mathrm
H}^m(X,(T^*)^{\otimes m})$$ 
to
$${\mathrm H}^m(X,\Lambda^m)={\mathrm H}_{\bar{\partial}}^{m,m}(X).$$
Note that because we are in the K{\"a}hler case the topological Chern
character coincides with the one constructed from the Atiyah
class. Indeed we can use any description of the Atiyah class to
construct the Chern character. For example, forget about the Hermitian
structure $h$ and simply choose a connection on $E$ which is
compatible with the complex structure. Then the $(1,1)$ part of the
curvature will represent $\alpha_E$ and hence
$$ch_m(E)=\frac{(-1)^m}{m!(2\pi i)^m}[{\mathrm
Tr}((R^{1,1})^m)]\in{\mathrm H}^{m,m}_{\bar{\partial}}(X).$$
Finally, let us note that if the bundle carries a symplectic structure
then the odd components of the Chern character vanish (as for the
tangent bundle), though this is not always the case.

\subsection{Wilson lines in Chern-Simons theory}

In this subsection we give a brief overview of the process of
obtaining invariants of knots and links in Chern-Simons theory due to
Witten~\cite{witten89}. These comments are intended to serve as
motivation for our construction in the following subsection of a
``Rozansky-Witten'' weight system on chord diagrams, which will be an
analogue of the weight system obtained in perturbative Chern-Simons
theory. 

Let $\cal L$ be a framed oriented link in a three-manifold $M$. Each
component $C_a$ of the link is an embedding of an oriented circle in
$M$, called a {\em Wilson line\/}. Taking the holonomy of the
Chern-Simons connection $A_i$ around $C_a$ gives us an observable of
the theory
$$W_{V_a}(C_a)={\mathrm Tr}_{V_a}{\mathrm Pexp}\int_{C_a}A_idx_i$$
where the Feynman path integral takes values in the gauge group $G$,
we choose a representation $V_a$ of $G$ for each link component, and
we take the trace in the corresponding representation. Including these
in the partition function of the theory gives us an invariant of the
link $\cal L$ in $M$
$$Z(M;{\cal L})=\int D{\cal A}{\mathrm exp}(iL){\prod}_a W_{V_a}(C_a)$$
where $L$ is the Lagrangian, or Chern-Simons action, and the integral
is over the moduli space of gauge equivalence classes of
connections. This partition function of $\cal L$ in $M$ is also known
as the {\em expectation value\/} of $\cal L$, or as the {\em Wilson
correlation function\/}. 

In the case that $M$ is the three-sphere $S^3$, $G={\mathrm SU}(2)$,
and all the $V_a$ are the standard representation on ${\Bbb C}^2$,
then we recover the Jones polynomial of the link $\cal L$ in
$S^3$. Indeed, the motivation behind Witten's work was to find an
intrinsically three-dimensional interpretation of the Jones
polynomial, which would then allow generalizations to links in
arbitrary three-manifolds. More generally, taking $G$ to be ${\mathrm
SU}(N)$ and all the $V_a$ to be the standard representation on ${\Bbb
C}^N$ allows us to obtain the HOMFLY polynomial of a link in $S^3$. On
the other hand, if we choose a representation to be the trivial one,
the corresponding Wilson line vanishes. Choosing all the
representation to be trivial makes the link vanish and we are left
simply with the Chern-Simons three-manifold invariant.

Taking a Feynman diagram expansion of the partition function of $\cal
L$ in $M$ gives us
$$Z(M;{\cal L})=\sum_D c_D({\frak g};V_a)Z_D(M;{\cal L})$$
where the weights $c_D({\frak g};V_a)$ depend only on the Lie algebra
${\frak g}$ of the gauge group and the representations $V_a$, and
$Z_D(M;{\cal L})$ depends on the link ${\cal L}$ embedded in $M$. The
sum is over all Feynman diagrams with external vertices on Wilson
lines. More precisely, $D$ is an oriented unitrivalent graph with
univalent vertices attached to $m$ oriented circles $S^1$ (the
skeleton), where $m$ is the number of components of ${\cal L}$. The
orientation is an equivalence class of cyclic orderings at the
trivalent, or internal, vertices of $D$, with two such being
equivalent if they differ by an even number of changes. We can regard 
the skeleton as part of the graph, so as to obtain a purely trivalent
graph. Then because the skeleton is made up of {\em oriented\/} circles,
there is automatically a cyclic ordering at the external vertices
(those lying on the skeleton) given by taking the order: incoming part
of skeleton, outgoing part of skeleton, remaining edge. Returning to
the unitrivalent graph, each connected component should have at least
one univalent vertex, and hence be attached to the skeleton. We denote
the space of linear combinations of such diagrams modulo the IHX, AS,
and STU relations (as in Chapter $4$) by ${\cal
A}(S^1\sqcup\ldots\sqcup S^1)$, where we take the disjoint union of
$m$ copies of $S^1$. We have already seen ${\cal A}(S^1)$, whose
elements we called chord diagrams. We will continue to use this
terminology when $m>1$. Recall that ${\cal A}(S^1)^{\prime}$ is the
larger space obtained by allowing diagrams with connected components
which are purely trivalent. Let ${\cal A}(S^1\sqcup\ldots\sqcup
S^1)^{\prime}$ be the corresponding space of diagrams with a skeleton
consisting of $m$ circles $S^1$, ie.\ allowing diagrams to have
connected components consisting of purely trivalent graphs.

It is believed that the terms $Z_D(M;{\cal L})$ should be the same as
the coefficients $Z^{\mathrm Kont}_D(M;{\cal L})$ of the framed
Kontsevich integral 
$$Z^{\mathrm Kont}(M;{\cal L})=\sum_D Z^{\mathrm Kont}_D(M;{\cal
L})D\in{\cal A}(S^1\sqcup\ldots\sqcup S^1)$$
of the link $\cal L$ in $M$. In any case, we are more interested in
the weight system $c_D({\frak g};V_a)$.

As with the weight system $c_{\Gamma}({\frak g})$ on trivalent graphs
described in Chapter $3$, this weight system can be constructed for
any Lie algebra ${\frak g}$ with an invariant inner product and
representations $V_a$. Recall that the structure constants of $\frak
g$ (with respect to some basis $\{x_1,\ldots,x_n\}$) give rise to a
skew-symmetric tensor $c_{ijk}$ and the inner product gives us a
symmetric tensor $\sigma^{ij}$. In each representation $V_a$, the
element $x_i$ is mapped to some endomorphism of $V_a$. Thus we get
tensors $(B_a)_{iJ_a}^{\phantom{iJ_a}K_a}$ where $J_a$ and $K_a$ refer
to some basis of $V_a$. 

Suppose we have some chord diagram in ${\cal
A}(S^1\sqcup\ldots\sqcup S^1)$. More generally, these weight systems
are well defined on diagrams $D$ in ${\cal A}(S^1\sqcup\ldots\sqcup
S^1)^{\prime}$. As with $c_{\Gamma}({\frak g})$, we
place a copy of $c_{ijk}$ at each trivalent vertex of $D$ and attach
the indices $ijk$ to the outgoing edges in a way compatible with the
cyclic ordering given by the orientation. If a univalent vertex is
connected to the $a$th circle $S^1$ of the skeleton, then we place a
copy of $(B_a)_{iJ_a}^{\phantom{iJ_a}K_a}$ there. In this way the
representations $V_a$ are attached to the different Wilson lines. If
some circle $S^1$ of the skeleton has no univalent vertices lying on
it, then we simply introduce a factor of ${\mathrm dim}V_b$, where
$V_b$ is the representation associated to that circle. Returning to
the univalent vertex, the index $J_a$ is attached to the incoming part
of the skeleton, the index $K_a$ to the outgoing part, and the index
$i$ to the remaining edge. Next we use $\sigma^{ij}$ to contract along
all edges of $D$. On the skeleton we set adjacent indices to be equal,
ie.\ if a part of the skeleton has been labelled with $J_a$ at one end
and with $K_a$ at the other, with no univalent vertices lying between,
then we set $J_a=K_a$. Equivalently, we contract along this part of
the skeleton with a Kronecker delta $\delta_{J_a}^{K_a}$. The number
that results from this process is the weight $c_D({\frak g};V_a)$.

For example, let $D$ be the chord diagram
$$\twoVchord$$
where the fact that the skeleton has been broken is, of course, only a
notational convenience. Then
$$c_D({\frak
g};R)=\sum_{i,j,K,L}B_{iK}^{\phantom{iK}L}B_{jL}^{\phantom{jL}K}\sigma^{ij}.$$

We chose bases for $\frak g$ and its representations $V_a$ in order to
define these weights but they are in fact independent of these choices.
Furthermore, they satisfy the IHX and AS relations for the same reasons
that the weights $c_{\Gamma}({\frak g})$ do. They also satisfy the STU
relations. This follows from the fact that $V_a$ is a representation
and so
$${\frak g}\rightarrow{\mathrm End}(V_a)$$
is a Lie algebra homomorphism; in particular, the bracket is
preserved. Thus we have a well defined element of the dual space
$({\cal A}(S^1\sqcup\ldots\sqcup S^1)^{\prime})^*$. 

This is the weight system on chord diagrams which arises naturally in
Chern-Simons theory. It has been extensively studied and generalized
by knot theorists, to arbitrary quadratic Lie algebras and Lie
super-algebras with representations. In the next subsection we will
show how an analogous construction can be made with hyperk{\"a}hler
manifolds and vector bundles in Rozansky-Witten theory.

\subsection{A new weight system on chord diagrams}

Let $X$ be a compact hyperk{\"a}hler manifold of real-dimension $4k$
and choose a complex structure $I$ from the family of complex
structures on $X$ compatible with the hyperk{\"a}hler metric. Recall
that the main ingredients in the construction of $b_{\Gamma}(X)$ were
$$\Phi\in\Omega^{0,1}(X,{\mathrm Sym}^3T^*)$$
which was essentially the Riemannian curvature, and the dual of the
holomorphic symplectic form
$$\tilde{\omega}\in{\mathrm H}^0(X,{\Lambda}^2T).$$
In local complex coordinates these have components $\Phi_{ijk\bar{l}}$
and $\omega^{ij}$ respectively.

Let $E_1,\ldots,E_m$ be holomorphic vector bundles over $X$, of ranks
$r_1,\ldots,r_m$ respectively. Suppose we have Hermitian structures
$h_a$ on each of the vector bundles $E_a$, and take the corresponding
(unique) Hermitian connections on the bundles. Let the curvatures of
these connections be
$$R_a\in\Omega^{1,1}({\mathrm End}E_a).$$
In local complex coordinates they have components
$$(R_a)^{I_a}_{\phantom{I_a}J_ak\bar{l}}dz_k\wedge d\bar{z}_l$$
where $I_a$ and $J_a$ refer to a local basis of sections of $E_a$.

Let $D\in{\cal A}(S^1\sqcup\ldots\sqcup S^1)^{\prime}$ be a chord
diagram with skeleton consisting of $m$ copies of $S^1$. We wish to
construct a weight system $b_D(X;E_a)$ which will be the natural analogue of
$c_D({\frak g};V_a)$. The weights $b_{\Gamma}(X)$ were analogous to
$c_{\Gamma}({\frak g})$ in the sense that $\Phi_{ijk\bar{l}}$ replaced
$c_{ijk}$ and $\omega^{ij}$ replaced $\sigma^{ij}$ in the
construction. For $b_D(X;E_a)$ we also wish to replace
$(B_a)_{iJ_a}^{\phantom{iJ_a}K_a}$ in the construction of $c_D({\frak
g};V_a)$ by $(R_a)^{I_a}_{\phantom{I_a}J_ak\bar{l}}$.

The construction should be obvious enough. We take a chord
diagram $D$ whose total number of vertices, trivalent and univalent,
is $2k$. We place a copy of $\Phi$ at each trivalent vertex and attach
the holomorphic indices $i$, $j$, and $k$ to the outgoing edges in the
usual way. We place copies of $R_a$ at the univalent edges, where $a$
is determined by which oriented circle $S^1$ of the skeleton the
univalent edge is attached to. The index $I_a$ is attached to the
incoming part of the skeleton, $J_a$ to the outgoing part, and $k$ to
the remaining outgoing edge. Then we contract along edges of $D$
using the dual of the holomorphic symplectic form $\tilde{\omega}$,
and along segments of the skeleton we contract using a Kronecker delta
$\delta_{I_a}^{J_a}$ as in the Lie algebra case. This gives us a
section of $(\bar{T}^*)^{\otimes 2k}$, which we project to the
exterior product to get
$$D(\Phi;R_a)\in\Omega^{0,2k}(X).$$
Multiplying by the $k$th power of the holomorphic symplectic form
gives us a $4k$-form which we can integrate to get a number. As with
the Lie algebra case, we need some convention to deal with the
situation where some circle $S^1$ of the skeleton has no univalent
vertices lying on it. When this happens, we introduce a factor of
$-{\mathrm rank}E_b$ where $E_b$ is the holomorphic vector bundle
corresponding to that circle. The minus sign is to ensure
compatibility with a formula we shall state in Subsection $6.5$.
\begin{dfn}
We define the weight $b_D(X;E_a)$ corresponding to the chord diagram
$D$ to be
$$\frac{1}{(8\pi^2)^kk!}\int_X D(\Phi;R_a)\omega^k.$$
\end{dfn}

\noindent
In the case that $m=0$, and we have no vector bundles, the chord
diagrams are just trivalent graphs and the weights reduce to the
Rozansky-Witten invariants of the manifold $X$. This explains our
choice of factor in the above definition. We can also say something
about the case of a trivial vector bundle $E_b$, for which a flat
connection can be chosen, and hence the curvature $R_b$ will
vanish. This means that any chord diagram $D$ with a univalent vertex
on the $b$th circle $S^1_b$ of the skeleton will give rise to a
vanishing weight $b_D(X;E_a)$. The remaining chord diagrams $D$ will
give weights which are $-{\mathrm rank}E_b$ times
$$b_{D\backslash S^1_b}(X;E_a\backslash E_b)$$
ie.\ the weight corresponding to $D$ with the $b$th circle $S^1_b$ of
the skeleton removed, where we also remove $E_b$ from the collection
of vector bundles. In other words, except for the $-{\mathrm
rank}E_b$ factor, trivial vector bundles effectively vanish from the
weight system. In some sense this is analogous to the vanishing of
Wilson lines in Chern-Simons theory when the trivial representation
has been associated to them.

So far we have neglected the role of the orientation in the above
construction. Recall that an orientation of the chord diagram $D$ is
an equivalence class of cyclic orderings of the outgoing edges at each
trivalent vertex, with two such being equivalent if the differ at an
even number of vertices. We can regard the skeleton as being part of
the graph, in which case what we have is a purely trivalent graph. In
fact, it is oriented because we already have an equivalence class of
cyclic orderings of the internal vertices and the vertices which lie
on the skeleton (ie.\ the univalent vertices of $D$) have a canonical
cyclic ordering given by: incoming part of skeleton, outgoing part of
skeleton, remaining outgoing edge. We take the corresponding
Rozansky-Witten orientation, ie.\ an equivalence class consisting of
an ordering of the vertices and orientations of the edges. Now the $m$
circles $S^1$ which make up the skeleton already have an orientation,
so the last thing to do is to make sure that the orientations of our
edges agree with this. Since we are dealing with an equivalence class,
we can achieve this by making an even number of changes.

\begin{figure}[htpb]
\epsfxsize=100mm
\centerline{\epsfbox{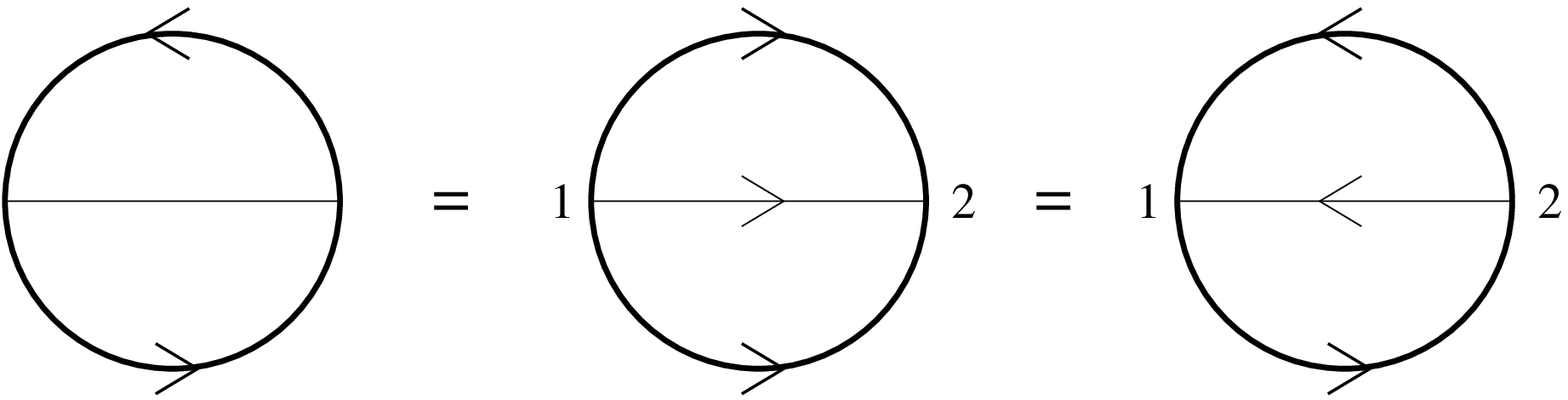}}
\caption{Three equivalent orientations}
\label{equiv_orient}
\end{figure}

For example, if $D\in{\cal A}(S^1)^{\prime}$ is the chord diagram
$$\twoVchord$$
then the orientations of $D$ in Figure~\ref{equiv_orient} are
equivalent. The first of these diagrams has the canonical orientation
given by the anti-clockwise ordering of the outgoing edges at each
vertex. In the second diagram we have converted this to the
corresponding Rozansky-Witten orientation. In the last diagram we have
reversed the orientations of the top and middle edges so that we have
agreement with the original orientation of the skeleton.

Thus an orientation of the chord diagram $D$ is the same as an
equivalence class consisting of an ordering of the vertices (trivalent
and univalent) and orientations of the edges (not including the
skeleton, which has a canonical orientation); if the orderings differ
by a permutation $\pi$ and $n$ edges are oriented in the reverse
manner, then we regard these as equivalent if ${\mathrm
sgn}\pi=(-1)^n$. This is precisely what we need for the weight
$b_D(X;E_a)$ to be well-defined. In particular, the ordering of the
vertices tells us in which order to multiply the copies of $\Phi$ and
$R_a$ associated to the trivalent and univalent vertices
respectively. The orientations of the edges tell us whether to
contract with $\omega^{ij}$ or $\omega^{ji}$ along the edge.

Let us write out $b_D(X;E)$ explicitly when $D$ is
$$\twoVchord$$
and $X$ is a K$3$ surface $S$ with a vector bundle $E$ over it. Using
the third orientation in Figure~\ref{equiv_orient} we see that
$$b{\twoVchord}(S;E)=\frac{1}{8\pi^2}\int_S
R^I_{\phantom{I}Jk_1\bar{l}_1}R^J_{\phantom{J}Ik_2\bar{l}_2}\omega^{k_2k_1}d\bar{l}_1\wedge
d\bar{l}_2\wedge\omega.$$

We can imitate many of the ideas that we used for $b_{\Gamma}(X)$ in
studying the weights $b_D(X;E_a)$. In particular, we will show the
following:
\begin{itemize}
\item the construction of the weights can be formulated in
cohomological terms, as in Kapranov's approach to the Rozansky-Witten
invariants, and this enable us to prove they are independent of the
choice of Hermitian structures $h_a$ on the bundles $E_a$,
\item the cohomological interpretation leads to a proof that the
weights satisfy the STU relations, and hence only depend on the
equivalence class of the chord diagram $D$ in ${\cal
A}(S^1\sqcup\ldots\sqcup S^1)^{\prime}$,
\item the Chern classes of the bundles $E_a$ arise naturally for
certain choices of chord diagrams,
\item the Wheeling Theorem can be used to prove a result which
generalizes Theorem~\ref{thm1} from Chapter $3$.
\end{itemize}

\subsection{The cohomological construction}

Recall that an alternative definition of $b_{\Gamma}(X)$ due to
Kapranov~\cite{kapranov99} replaced differential forms with the
Dolbeault cohomology classes they represent. In particular, $\Phi$
represent the Atiyah class of the tangent bundle
$$\alpha_T=[\Phi]\in{\mathrm H}^{0,1}_{\bar{\partial}}(X,{\mathrm
Sym}^3T^*).$$
We can do the same with $b_D(X;E_a)$, as the curvature $R_a$ of the
Hermitian connection on $E_a$ is $\bar{\partial}$-closed and
represents the Atiyah class of $E_a$
$$\alpha_{E_a}=[R_a]\in{\mathrm
H}^{0,1}_{\bar{\partial}}(X,T^*\otimes{\mathrm End}E_a).$$
We can then construct
$$[D(\Phi;R_a)]\in{\mathrm H}^{0,2k}_{\bar{\partial}}(X)$$
as before, and it only depends on the cohomology classes of $\Phi$ and
$R_a$. Multiplying by $[\omega^k]$ gives an element of ${\mathrm
H}^{2k,2k}_{\bar{\partial}}(X)$ which we integrate to get
$$b_D(X;E_a)=\frac{1}{(8\pi^2)^kk!}\int_X[D(\Phi;R_a)][\omega^k].$$
More generally, we can take arbitrary representatives of the Atiyah
classes $\alpha_T$ and $\alpha_{E_a}$, and construct
$$D(\alpha_T;\alpha_{E_a})\in{\mathrm H}^{2k}(X,{\cal O}_X).$$
This element can be paired with
$$\omega^k\in{\mathrm H}^0(X,{\Lambda}^{2k}T^*)$$
using Serre duality to produce the number $b_D(X;E_a)$, up to the
factor $\frac{1}{(8\pi^2)^kk!}$.

Since the Atiyah class $\alpha_E$ of a holomorphic vector bundle $E$
is well-defined without any reference to a Hermitian structure $h$ on
$E$, the next result follows immediately from the description above.
\begin{prp}
The weights $b_D(X;E_a)$ are independent of the choices of Hermitian
structures $h_a$ on the vector bundles $E_a$.
\end{prp}

It is clear that the weights $b_D(X;E_a)$ satisfy the AS relations,
ie.\ they change sign under a change
of orientation of $D$
$$b_{\bar{D}}(X;E_a)=-b_D(X;E_a).$$
The argument is identical to the one for $b_{\Gamma}(X)$. They also
satisfy the IHX relations: if three chord diagrams $D_I$, $D_H$, and
$D_X$ are identical except inside some small ball (away from the
skeleton) where they look like $I$, $H$, and $X$ respectively, then
$$b_{D_I}(X;E_a)=b_{D_H}(X;E_a)-b_{D_X}(X;E_a).$$
Again this follows from the same argument as for $b_{\Gamma}(X)$,
using the cohomological approach described above. Indeed, the IHX
relations only involve internal vertices, so the fact that the chord
diagram also has a skeleton where we have attached the curvatures of
vector bundles does not interfere with the argument. The STU
relations, on the other hand, do depend on properties of these
curvatures.

Suppose we are given three chord diagrams $D_S$, $D_T$, and $D_X$
which are identical except inside some small ball near the skeleton
where they look like
\begin{center}
\begin{picture}(140,70)(-60,-35)
\put(-110,-20){\vector(1,0){40}} 
\put(-90,10){\line(0,-1){30}} 
\put(-90,10){\line(1,1){20}} \put(-90,10){\line(-1,1){20}} 
%\put(-40,0){$=$}
\put(-30,-20){\vector(1,0){40}}
\put(-18,-20){\line(-1,3){17}} \put(-2,-20){\line(1,3){17}}
%\put(35,0){$-$}
\put(40,0){and}
\put(80,-20){\vector(1,0){40}} 
\put(88,-20){\line(2,3){10}} 
\put(102,1){\line(2,3){20}}
\put(112,-20){\line(-2,3){34}}
\end{picture}
%Figure $1$: The STU relation
\end{center}
respectively. The orientations are induced from the planar embedding
(for internal vertices) and the canonical cyclic orderings at external
vertices induced from the orientation of the skeleton. Converting this
to a Rozansky-Witten orientation, and remembering that we need to
choose an equivalence class for which the orientation on the skeleton
agrees with the given one, we find that we get the following: in $D_S$
let the internal vertex be labelled $t$ and the external vertex $s$,
with the connecting edge oriented from $t$ to $s$. Then in $D_T$ the
left vertex should be labelled $s$ and the right one $t$, while in
$D_U$ the left vertex should be labelled $t$ and the right one
$s$. The orientations of the edges in $D_T$ and $D_U$ are induced from
the orientation of the skeleton. We wish to show that
$$b_{D_S}(X;E_a)=b_{D_T}(X;E_a)-b_{D_U}(X;E_a).$$
Equivalently, we can take $D_T$ with the reverse orientation (ie.\
switch the labels $t$ and $s$ so that they agree with those on $D_U$),
in which case the equivalent relation is
$$b_{D_S}(X;E_a)+b_{\bar{D}_T}(X;E_a)+b_{D_U}(X;E_a)=0.$$

As with the proof of the IHX relation, we will show that
$$D_S(\Phi;R_a)+\bar{D}_T(\Phi;R_a)+D_U(\Phi;R_a)\in\Omega^{0,2k}(X)$$
is $\bar{\partial}$-exact, and thus cohomologous to zero. Since the
three chord diagrams differ only near one of the circles making up the
skeleton, we will temporarily drop the subscript $a$ for the
holomorphic vector bundle $E$ associated to that circle, and for its
curvature $R$. Consider $d^2_AR$, where $d_A$ is the exterior derivative
corresponding to the Hermitian connection on $E$ (corresponding to
some choice of Hermitian structure). By the Bianchi identity $R$ is
$\bar{\partial}$-closed, and hence
$$d^2_AR=\bar{\partial}\partial_AR.$$
On the other hand
$$R\in\Omega^{1,1}(X,{\mathrm End}E)=\Omega^{0,1}(X,T^*\otimes
{\mathrm End}E)$$
so $d^2_A$ gives the curvature of the associated connection on
$T^*\otimes {\mathrm End}E$. This is the curvature of the Levi-Civita
connection on $T^*$ added to the curvature of the connection on
${\mathrm End}E$ induced from the Hermitian connection on $E$. Overall
the curvature acts on $T^*\otimes {\mathrm End}E$ as $K\otimes
1+1\otimes R$, where the second term act on ${\mathrm End}E$ via the
Lie bracket. It is probably easiest to write this out in local
coordinates, in which case we find
$$d^2_AR\in\Omega^{1,2}(X,T^*\otimes {\mathrm
End}E)=\Omega^{0,2}(X,T^*\otimes T^*\otimes {\mathrm End}E)$$
is given by
$$K^i_{\phantom{i}jk\bar{l}}R^A_{\phantom{A}Bi\bar{m}}
+R^A_{\phantom{A}Ck\bar{l}}R^C_{\phantom{C}Bj\bar{m}}
-R^A_{\phantom{A}Cj\bar{m}}R^C_{\phantom{C}Bk\bar{l}}.$$
Just to make things clear, $A$ and $B$ refer to ${\mathrm End}E$, $j$
and $k$ refer to $T^*\otimes T^*$, and $\bar{l}$ and $\bar{m}$ refer
to $\Omega^{0,2}$ (we take a wedge product of forms in order to arrive
in the exterior product). Since we are in this exterior product, we
can interchange the anti-holomorphic indices in the last term,
introducing a sign change
$$K^i_{\phantom{i}jk\bar{l}}R^A_{\phantom{A}Bi\bar{m}}+R^A_{\phantom{A}Ck\bar{l}}R^C_{\phantom{C}Bj\bar{m}}+R^A_{\phantom{A}Cj\bar{l}}R^C_{\phantom{C}Bk\bar{m}}.$$
Note that this is symmetric in $jk$, and hence we have an element
which we write schematically as
\begin{eqnarray}
KR+RR+RR & \in & \Omega^{0,2}(X,{\mathrm Sym}^2T^*\otimes {\mathrm End}E).
\end{eqnarray}
It follows that 
$$\partial_AR\in\Omega^{1,1}(X,T^*\otimes {\mathrm
End}E)=\Omega^{0,1}(X,T^*\otimes T^*\otimes {\mathrm End}E)$$
must have the same symmetry, ie.\ must lie in
$$\Omega^{0,1}(X,{\mathrm Sym}^2T^*\otimes {\mathrm End}E)$$
as
$$\bar{\partial}\partial_AR=KR+RR+RR.$$

Consider the $S$ part of $D_S$, by which we mean the part of $D_S$
which differs from $D_T$ and $D_U$. Let the internal vertex be
labelled with $t$ and the external vertex with $s$, and assume
$t<s$. We have seen that the edge connecting them must be oriented
from $t$ to $s$. When calculating $D_S(\Phi;R_a)$ this part
contributes a copy of $\Phi$ and a copy of $R$ ``joined'' by a copy of
$\tilde{\omega}$, namely the section
$$S(\Phi;R_a)\in C^{\infty}(X,(T^*)^{\otimes 2}\otimes {\mathrm
End}E\otimes (\bar{T}^*)^{\otimes 2})$$
with components
\begin{eqnarray*}
S(\Phi;R_a)^{\phantom{j_tk_t}A_s}_{j_tk_t\phantom{A_s}B_s\bar{l}_t\bar{l}_s}
& = &
\Phi_{i_tj_tk_t\bar{l}_t}\omega^{i_tk_s}R^{A_s}_{\phantom{A_s}B_sk_s\bar{l}_s}
\\
 & = & \omega_{i_tm}K^m_{\phantom{m}j_tk_t\bar{l}_t}\omega^{i_tk_s}R^{A_s}_{\phantom{A_s}B_sk_s\bar{l}_s}
\\
 & = & K^m_{\phantom{m}j_tk_t\bar{l}_t}\delta^{k_s}_mR^{A_s}_{\phantom{A_s}B_sk_s\bar{l}_s}
\\
 & = &
 K^m_{\phantom{m}j_tk_t\bar{l}_t}R^{A_s}_{\phantom{A_s}B_sm\bar{l}_s}.
\end{eqnarray*}
The indices $j_tk_t$ refer to $(T^*)^{\otimes 2}$ and note that
this term is symmetric in $j_tk_t$. Similarly, the $\bar{T}$ part of
$\bar{D}_T$ contributes
$$\bar{T}(\Phi;R_a)\in C^{\infty}(X,(T^*)^{\otimes 2}\otimes {\mathrm
End}E\otimes (\bar{T}^*)^{\otimes 2})$$
with components
$$\bar{T}(\Phi;R_a)^{\phantom{j_tk_t}A_s}_{j_tk_t\phantom{A_s}B_s\bar{l}_t\bar{l}_s}=R^{A_s}_{\phantom{A_s}Cj_t\bar{l}_t}R^{C}_{\phantom{C}B_sk_t\bar{l}_s}$$
to $\bar{D}_T(\Phi;R_a)$, and the $U$ part of $D_U$ contributes
$$U(\Phi;R_a)\in C^{\infty}(X,(T^*)^{\otimes 2}\otimes {\mathrm
End}E\otimes (\bar{T}^*)^{\otimes 2})$$
with components
$$U(\Phi;R_a)^{\phantom{j_tk_t}A_s}_{j_tk_t\phantom{A_s}B_s\bar{l}_t\bar{l}_s}=R^{A_s}_{\phantom{A_s}Ck_t\bar{l}_t}R^{C}_{\phantom{C}B_sj_t\bar{l}_s}$$
to $U(\Phi;R_a)$.
The sum of these three terms will be symmetric in $j_tk_t$, and will
give us an element
$$S(\Phi;R_a)+\bar{T}(\Phi;R_a)+U(\Phi;R_a)\in C^{\infty}(X,{\mathrm
Sym}^2T^*\otimes {\mathrm End}E\otimes (\bar{T}^*)^{\otimes 2}).$$
Projecting to the exterior product $\Omega^{0,2}(X,{\mathrm
Sym}^2T^*\otimes {\mathrm End}E)$ gives us the term (4) which we
described earlier in a schematic way as $KR+RR+RR$ (actually we need
to interchange the second and third terms to get an exact
correspondence). 

It follows that
$$D_S(\Phi;R_a)+\bar{D}_T(\Phi;R_a)+D_U(\Phi;R_a)=D_*(\Phi;R_a,KR+RR+RR)$$
in $\Omega^{0,2k}(X)$ where $D_*$ is the diagram which is identical to
$D_S$ away from the $S$ part (and hence identical to $\bar{D}_T$ and
$D_U$ away from the $\bar{T}$ and $U$ parts respectively), but
contains one bivalent vertex lying on the skeleton instead of the $S$
part. The construct of the right hand side is the same as for a chord
diagram except that at the bivalent vertex we place $KR+RR+RR$. An
orientation of $D_*$ can be defined in such a way as to agree with
those on $D_S$, $\bar{D}_T$, and $D_U$. However, such a choice only
affects the overall sign so is not actually important here (since we
are only interested in showing $\bar{\partial}$-exactness).

Replacing $KR+RR+RR$ by $\bar{\partial}\partial_AR$ and using the fact
that $\Phi$, $R_a$, and $\tilde{\omega}$ are all
$\bar{\partial}$-closed, we get
$$D_*(\Phi;R_a,\bar{\partial}\partial_AR)=\bar{\partial}D_*(\Phi;R_a,\partial_AR)$$
where
$$D_*(\Phi;R_a,\partial_AR)\in\Omega^{0,2k-1}(X).$$
Thus we have shown that in Dolbeault cohomology
$$[D_S(\Phi;R_a)]+[\bar{D}_T(\Phi;R_a)]+[D_U(\Phi;R_a)]=0\in{\mathrm
H}^{0,2k}_{\bar{\partial}}(X),$$
and therefore by the cohomological construction we see that
$$b_{D_S}(X;E_a)+b_{\bar{D}_T}(X;E_a)+b_{D_U}(X;E_a)=0$$
as we wanted.

The fact that the weights $b_D(X;E_a)$ satisfy the STU relations
immediately implies that they also satisfy the IHX and AS relations,
although we have already seen that this is true. Combining these gives
us the following result.
\begin{prp}
The dependence of the weight $b_D(X;E_a)$ on the chord diagram $D$ is
only through its equivalence class in ${\cal A}(S^1\sqcup\ldots\sqcup
S^1)^{\prime}$.
\end{prp}

\subsection{Chern classes}

To obtain Chern numbers from the Rozansky-Witten invariants
$b_{\Gamma}(X)$ we needed to take a graph $\Gamma$ somehow composed of
wheels $w_{\lambda}$. More precisely, we chose $\Gamma$ to be a
polywheel, ie.\ the closure (sum over all the different ways of
joining the spokes) of the disjoint union of a collection of
wheels. Wheels arise naturally in chord diagrams as the skeleton is
made up of a collection of oriented circles $S^1$, each of which can
be regarded as a wheel. Thus we can expect that certain chord diagrams
$D$ will give us weights $b_D(X;E_a)$ which are related to the Chern
classes of the holomorphic vector bundles $E_a$. Indeed this is the
case, and the argument is virtually identical to that given for the
Chern numbers of the tangent bundle in Chapter $2$. Consequently we
will state the results here without proof.

We will denote wheels which consist of a circle $S^1$ which is part of
the skeleton of a chord diagram by a bold ${\bf w}_{\lambda}$. This
looks like 
$$\begin{picture}(120,50)(-60,-20)
\put(-60,-15){\vector(1,0){120}}
\put(-48,-15){\line(0,1){40}}
\put(-36,-15){\line(0,1){40}}
\put(0,0){\ldots}
\put(48,-15){\line(0,1){40}}
\end{picture}$$
where there are $\lambda$ vertical spokes, and we have continued with
our convention of breaking the circle $S^1$ at some point and drawing
it as a directed line.

Suppose that we have a single holomorphic vector bundle $E$ over a
compact hyperk{\"a}hler manifold $X$ of real-dimension $4k$. We know
that the closure $\langle w_{2k}\rangle$ of a $2k$-wheel will give
rise to a Rozansky-Witten invariant which is (up to a factor) the
Chern number given by the top component of the Chern character, ie.\
$$b_{\langle w_{2k}\rangle}(X)=-s_{2k}(X).$$
The same argument shows that
$$b_{\langle{\bf w}_{2k}\rangle}(X;E)=-\int_Xs_{2k}(E)$$
where 
$$s_{2k}(E)=(2k)!ch_{2k}(E)\in\Omega^{2k,2k}(X)$$
is the usual rescaling of a component of the Chern character. More
generally, we can show the following.
\begin{prp}
Suppose we have $m$ holomorphic vector bundles $E_a$ over a compact
hyperk{\"a}hler manifold $X$ of real-dimension $2k$. Let $D$ be the
chord diagram given by the closure of a collection of wheels
$$D=\langle w_{\lambda_1}\cdots w_{\lambda_j}{\bf
w}_{\lambda_{j+1}}\cdots{\bf w}_{\lambda_{j+m}}\rangle\in{\cal A}(S^1\sqcup\ldots\sqcup S^1)^{\prime}$$ 
where $\lambda_1+\ldots+\lambda_{j+m}=2k$ and there are $m$ circles
$S^1$ making up the skeleton. The weight corresponding to this chord
diagram is
$$b_D(X;E_a)=(-1)^{j+m}\int_X s_{\lambda_1}(T)\wedge\cdots\wedge
s_{\lambda_j}(T)\wedge s_{\lambda_{j+1}}(E_1)\wedge\cdots\wedge
s_{\lambda_{j+m}}(E_m)$$
where $T$ is the holomorphic tangent bundle of $X$.
\end{prp}

\noindent
Observe that both sides of this formula are zero unless
$\lambda_1,\ldots,\lambda_j$ are all even. However, when some of 
$\lambda_{j+1},\ldots,\lambda_{j+m}$ are odd we can still get a
non-trivial result. We can also allow some of
$\lambda_{j+1},\ldots,\lambda_{j+m}$ to be zero, in which case ${\bf
w}_0$ will give us a circle $S^1$ of the skeleton which has no
univalent vertices lying on it. Recall that this introduces a factor
of $-{\mathrm rank}E_b$ into $b_D(X;E_a)$, where $E_b$ is the vector
bundle associated to the circle. This agrees with the formula above,
since $s_0(E_b)={\mathrm rank}E_b$ and we have an additional minus
sign in front of the integral coming from this bundle.  

In Chapter $4$ we used the Wheeling Theorem to show that $\Theta^k$ lies
in the subspace generated by polywheels inside graph homology ${\cal
A}(\emptyset)$, thereby proving Theorem~\ref{thm1}. Recall that the
Wheeling Theorem says that the map $\hat{\Omega}$ between ${\cal
B}^{\prime}\cong{\cal A}(S^1)^{\prime}$ and itself is an algebra
isomorphism, where the algebra structure comes from the two different
products $\cup$ and $\times$. Our proof of Theorem~\ref{thm1} used the
consequence
$$\hat{\Omega}(\ell^{\cup k})=(\hat{\Omega}(\ell))^{\times k}$$
though actually we only used the terms on each side of this equation
which have no univalent vertices, ie.\ consist of an element of graph
homology ${\cal A}(\emptyset)$ with a disjoint skeleton. Now we shall
consider the full equation, and the corresponding relations on weights
$b_D(X;E_a)$. In fact, since the skeleton of an element of ${\cal
A}(S^1)^{\prime}$ consists of a single circle $S^1$, we need only one
holomorphic vector bundle $E$ over $X$.

We start with a low degree example, namely $k=2$. Recall that the
first few terms of $\Omega$ are
$$1+\frac{1}{48}w_2+\frac{1}{48^22!}(w_2^2-\frac{4}{5}w_4)+\frac{1}{48^33!}(w_2^3-\frac{12}{5}w_2w_4+\frac{64}{35}w_6)+\ldots$$
and therefore the left hand side looks like
$$\hat{\Omega}(\ell^{\cup 2})=\hat{1}(\ell^{\cup
2})+\frac{1}{48}\widehat{w_2}(\ell^{\cup
2})+\frac{1}{48^22!}(\widehat{w_2^2}-\frac{4}{5}\widehat{w_4})(\ell^{\cup
2}).$$
We wish to rewrite this element of ${\cal B}^{\prime}$ as an element
of ${\cal A}(S^1)^{\prime}$. Recall that the required map is given by
averaging over all the different ways of joining the univalent
vertices to a skeleton $S^1$. Equivalently, on terms with $\lambda$
univalent vertices, we can act with $\widehat{{\bf w}_{\lambda}}$,
where this operator acts in an analogous way to
$\widehat{w_{\lambda}}$. In fact, we also need to divide by $\lambda
!$ in order to get an average and not a sum. Thus as an element of
${\cal A}(S^1)^{\prime}$, we can write $\hat{\Omega}(\ell^{\cup 2})$ as
$$\frac{1}{4!}\widehat{{\bf w}_4}(\ell^{\cup
2})+\frac{1}{2!48}\widehat{w_2{\bf w}_2}(\ell^{\cup
2})+\frac{1}{48^22!}(\widehat{w_2^2{\bf
w}_0}-\frac{4}{5}\widehat{w_4{\bf w}_0})(\ell^{\cup 2}).$$
Now acting on $\ell^{\cup 2}$ is the same as taking the closure of
the above products of wheels, up to a factor of $8=2^22!$ which comes
from the fact that each line $\ell$ can join two given spokes in two
different ways, and the two lines can also be interchanged. Thus we
get
$$8(\frac{1}{4!}\langle{\bf w}_4\rangle+\frac{1}{2!48}\langle w_2{\bf 
w}_2\rangle +\frac{1}{48^22!}(\langle w_2^2{\bf
w}_0\rangle-\frac{4}{5}\langle w_4{\bf w}_0\rangle )).$$
In Chapter $4$ we already saw that the right hand side looks like
\begin{eqnarray*}
(\hat{\Omega}(\ell))^{\times 2} & = &
(\twoVchord+\frac{1}{24}\thetachord)^{\times 2} \\
 & = &
\fourVchord+\frac{2}{24}\thetatwoVchord+\frac{1}{24^2}\thetasqchord
\end{eqnarray*}
as an element of ${\cal A}(S^1)^{\prime}$.

Now suppose that we have a holomorphic vector bundle $E$ over a
compact hyperk{\"a}hler manifold $X$ of real-dimension eight, and
consider the weights $b_D(X;E)$ corresponding to the above
chord diagrams $D\in{\cal A}(S^1)^{\prime}$. Firstly, we have written
the left hand side as a sum of closures of wheels, and we know that
such diagrams give rise to integrals of Chern classes of $E$ and the
tangent bundle $T$. Indeed, the left hand side gives
\begin{eqnarray*}
-8\int_X\frac{s_4(E)}{4!}-\frac{s_2(T)}{48}\wedge\frac{s_2(E)}{2!}+\frac{(s_2^2(T)+\frac{4}{5}s_4(T))}{48^22!}\wedge s_0(E)\hspace*{-90mm} & & \\
 & = & -8\int_X{\mathrm Td}^{1/2}_0(T)\wedge ch_4(E)+{\mathrm
Td}^{1/2}_2(T)\wedge ch_2(E)+{\mathrm Td}^{1/2}_4(T)\wedge ch_0(E) \\
 & = & -8\int_X({\mathrm Td}^{1/2}(T)\wedge ch(E))_4
\end{eqnarray*}
where the subscripts indicate which terms of ${\mathrm Td}^{1/2}(T)$
to take, and in the final line we integrate the fourth term of
${\mathrm Td}^{1/2}(T)\wedge ch(E)$. The right hand side then gives us
a representation of this integral of characteristic numbers in terms
of weights corresponding to simple chord diagrams, namely
$$b\fourVchord(X;E)+\frac{2}{24}b\thetatwoVchord(X;E)+\frac{1}{24^2}b\thetasqchord(X;E).$$

We have just seen the $k=2$ case. The general case of a real-dimension
$4k$ hyperk{\"a}hler manifold $X$ is similar, and produces a formula
expressing
$$-2^kk!\int_X({\mathrm Td}^{1/2}(T)\wedge ch(E))_{2k}$$
as the weight $b_D(X;E)$ corresponding to the chord diagram
$$D=(\twoVchord+\frac{1}{24}\thetachord)^{\times k}.$$
The characteristic class
$$v(E)={\mathrm Td}^{1/2}(T)\wedge ch(E)\in\bigoplus_{m=0}^{2k}{\mathrm
H}^{2m}(X)$$
is known as the Mukai vector of $E$, and warrants a few comments.

The Mukai vector can be defined for any coherent sheaf $\cal E$ on a
smooth manifold. It gives a map from K-theory to the cohomology ring
which satisfies various functorial properties. In the more specific
hyperk{\"a}hler situation, it occurs in connection with moduli spaces
of sheaves on a K$3$ surface $S$ (see the book by Huybrechts and
Lehn~\cite{hl97}). Let ${\cal M}(v)$ and ${\cal M}^s(v)$ be the moduli
spaces of semi-stable and stable sheaves $\cal E$ over $S$,
respectively, with Mukai vectors $v({\cal E})=v$. Provided the sheaves
have rank greater than one, the moduli space ${\cal M}^s(v)$ is
smooth. Its complex-dimension can be defined in terms of the Mukai
vector $v$, as $(v,v)+2$, where the inner product on cohomology is
defined by
$$(v,w)=\int_S -v_0\wedge w_4+v_2\wedge w_2-v_4\wedge w_0.$$
In this context, the Mukai vector $v$ has also been of interest to
physicists (see Dijkgraaf~\cite{dijkgraaf99}), where it occurs as the
{\em charge\/} of a {\em D-brane\/} state, which can be described in
terms of a coherent sheaf $\cal E$ with Mukai vector $v\in{\mathrm
H}^*(S,{\Bbb Z})$. 

Returning to our result above, what we have is a formula for the
integral of the top component of the Mukai vector in terms of weights
of reasonably simple chord diagrams. For example, up to a factor the
last of these chord diagrams is just $\Theta^k$ with a disjoint
skeleton $S^1$.
\begin{rmk}
\begin{enumerate}
\item We know how to interpret the corresponding weight in terms of
classical invariants of $X$ and $E$, namely in terms of the ${\cal
L}^2$-norm of the curvature of the Levi-Civita connection on $X$, the
volume of $X$, and the rank of $E$. Unfortunately, we do not know how
to interpret the weights of the remaining chord diagrams making up $D$
in such a simple way, ie.\ in terms of classical invariants of $X$ and
$E$. 
\item Another goal would be to find a representation of the entire
Mukai vector in terms of weights of simple chord diagrams.
\item More generally, it would be interesting to know how these
weights depend on the holomorphic structures on the vector bundles
$E_a$.
\end{enumerate} 
\end{rmk}

\subsection{Topological quantum field theory}

In this chapter we have described how to construct a weight system
$b_D(X;E_a)$ on chord diagrams $D$ in ${\cal A}(S^1\sqcup\ldots\sqcup
S^1)^{\prime}$ from a compact hyperk{\"a}hler manifold $X$ with a
collection of $m$ holomorphic vector bundles $E_a$ over it, one for
each circle $S^1$ in the skeleton of $D$. Composing this weight system
with the framed Kontsevich integral
$$Z^{\mathrm Kont}(M;{\cal L})=\sum_D Z^{\mathrm Kont}_D(M;{\cal L})D$$
gives us a finite-type invariant 
$$\sum_D b_D(X;E_a)Z^{\mathrm Kont}_D(M;{\cal L})$$
of the $m$-component framed oriented link ${\cal L}$ in the
three-manifold $M$, which is similar in form to the invariant of links
which arises from adding Wilson lines in Chern-Simons theory. Indeed
we can ask: is there some natural extension of the Rozansky-Witten
theory which leads to an invariant of links which, furthermore, gives
rise to the weights $b_D(X;E_a)$ when expanded ``perturbatively''? For
an embedding of a knot in a three-manifold $M$, Rozansky and
Witten~\cite{rw97} constructed an observable for their theory which
depended on a tensor or spinor bundle, but this construction does not
appear to extend to arbitrary holomorphic vector bundles. Let us
consider instead the three-dimensional topological field theory.

In general, such a theory associates to a Riemann surface $\Sigma$ a
Hilbert space ${\cal H}_{\Sigma}$ and to a three-manifold $M$ with
boundary $\partial M=\Sigma$ a vector 
$$v_M\in{\cal H}_{\Sigma}.$$
Suppose we have a knot ${\cal K}$ embedded in a closed three-manifold
$M$. Removing a toroidal neighbourhood $\bar{\cal K}$ of the knot
gives us a three-manifold $M\backslash\bar{\cal K}$ with boundary a
torus $T$, and therefore a vector
$$v_{{\cal K}\subset M}\in{\cal H}_T$$
in the corresponding Hilbert space. If we choose a basis for ${\cal
H}_T$, the components of $v_{{\cal K}\subset M}$ will give us a
scalar-valued invariant of the knot ${\cal K}$ in $M$. Equivalently,
we can choose some vector $w\in{\cal H}_T$ and take the inner product
$$\langle w,v_{{\cal K}\subset M}\rangle$$
with $v_{{\cal K}\subset M}$, and this will also give us a scalar-valued
invariant.

In~\cite{rw97} Rozansky and Witten proposed the following topological
quantum field theory. The Hilbert space associated to a genus $g$
Riemann surface is the direct sum of cohomology groups
$$\bigoplus_{p,q}{\mathrm H}^q(X,(\Lambda^pT^*)^{\otimes g})$$
of the compact hyperk{\"a}hler manifold $X$. Thus for a torus (ie.\
genus $g=1$) we get the cohomology ring 
\begin{eqnarray*}
{\cal H}_T & = & \bigoplus_{p,q}{\mathrm H}^q(X,\Lambda^pT^*) \\
 & = & \bigoplus_{p,q}{\mathrm H}^{p,q}_{\bar{\partial}}(X)
\end{eqnarray*}
of $X$. According to the above argument, a knot ${\cal K}$ embedded in
$M$ will give rise to a cohomology class 
$$v_{{\cal K}\subset M}\in\bigoplus_{p,q}{\mathrm
H}^{p,q}_{\bar{\partial}}(X).$$

Now we must find some way to introduce the holomorphic vector bundle
$E$ into this construction. Of course, $E$ gives rise naturally to 
characteristic classes in ${\cal H}_T$. Indeed we can take the Chern
class, Chern character, or perhaps something like the Mukai vector,
and these will all give us some element $w\in{\cal H}_T$. Pairing this
with $v_{{\cal K}\subset M}$ then gives us a scalar-valued invariant
$$\int_X w\wedge v_{{\cal K}\subset M}$$
of the knot ${\cal K}$ in $M$. Now we may hope that a ``perturbative''
expansion of this knot invariant will give rise to the weights
$b_D(X;E)$ where $D\in{\cal A}(S^1)$. However, the problem here is
that our knot invariant depends on $E$ only through its characteristic
classes, and furthermore only in a linear way. In comparison, we would
expect the weights $b_D(X;E)$ to depend on $E$ in a more subtle
way. Indeed, we saw in Chapter $5$ that some of the Rozansky-Witten
invariants $b_{\Gamma}(X)$ cannot be expressed as linear combinations
of the characteristic numbers of $X$, so presumably a similar
statement should hold for $b_D(X;E)$.

%% file: app1.tex
\section{Graph relations}
\subsection{Expansions of polywheels}

$$
\begin{array}{lll}
{\bf k=1} &   & \\
\langle w_2\rangle        & = & \twoVgraph \\
          &   & \\
{\bf k=2} &   & \\
\langle w_2^2\rangle      & = & \twoVgraph^2+2\fourVgraph \\
\langle w_4\rangle        & = & \frac{5}{2}\fourVgraph \\
          &   & \\
{\bf k=3} &   & \\
\langle w_2^3\rangle      & = & \twoVgraph^3+6\twoVgraph\fourVgraph+8\sixVgraph \\
\langle w_2w_4\rangle     & = & \frac{5}{2}\twoVgraph\fourVgraph+10\sixVgraph \\
\langle w_6\rangle        & = & \frac{35}{4}\sixVgraph \\
          &   & \\
{\bf k=4} &   & \\
\langle w_2^4\rangle      & = &
\twoVgraph^4+12\twoVgraph^2\fourVgraph+32\twoVgraph\sixVgraph+12\fourVgraph^2+48\eightVgraphI
\\
\langle w_2^2w_4\rangle   & = &
\frac{5}{2}\twoVgraph^2\fourVgraph+20\twoVgraph\sixVgraph+5\fourVgraph^2+60\eightVgraphI
\\
\langle w_4^2\rangle      & = &
\frac{25}{4}\fourVgraph^2+48\eightVgraphI+24\eightVgraphII \\
\langle w_2w_6\rangle     & = &
\frac{35}{4}\twoVgraph\sixVgraph+\frac{105}{2}\eightVgraphI \\
\langle w_8\rangle        & = &
\frac{287}{8}\eightVgraphI+7\eightVgraphII \\
          &   & \\
{\bf k=5} &   & \\
\langle w_2^5\rangle      & = &
\twoVgraph^5+20\twoVgraph^3\fourVgraph+60\twoVgraph\fourVgraph^2+80\twoVgraph^2\sixVgraph+160\fourVgraph\sixVgraph
\\
          &   &
\phantom{\twoVgraph}+240\twoVgraph\eightVgraphI+384\tenVgraphI \\
\end{array}
$$

\newpage
$$
\begin{array}{lll}
\langle w_2^3w_4\rangle   & = &
\frac{5}{2}\twoVgraph^3\fourVgraph+15\twoVgraph\fourVgraph^2+30\twoVgraph^2\sixVgraph+80\fourVgraph\sixVgraph
\\
          &   & \phantom{\twoVgraph}+180\twoVgraph\eightVgraphI+480\tenVgraphI \\

\langle w_2w_4^2\rangle   & = &
\frac{25}{4}\twoVgraph\fourVgraph^2+50\fourVgraph\sixVgraph+48\twoVgraph\eightVgraphI+24\twoVgraph\eightVgraphII
\\
          &   & \phantom{\twoVgraph}+384\tenVgraphI+192\tenVgraphII \\
\langle w_2^2w_6\rangle   & = &
\frac{35}{4}\twoVgraph^2\sixVgraph+\frac{35}{2}\fourVgraph\sixVgraph+105\twoVgraph\eightVgraphI+420\tenVgraphI
\\
\langle w_4w_6\rangle     & = &
\frac{175}{8}\fourVgraph\sixVgraph+\frac{483}{2}\tenVgraphI+252\tenVgraphII
\\
\langle w_2w_8\rangle     & = &
\frac{287}{8}\twoVgraph\eightVgraphI+7\twoVgraph\eightVgraphII+287\tenVgraphI+56\tenVgraphII
\\
\langle w_{10}\rangle     & = &
\frac{2541}{16}\tenVgraphI+\frac{231}{2}\tenVgraphII \\ 
\phantom{\bf k=5}         &   & \phantom{\twoVgraph^5+20\twoVgraph^3\fourVgraph+60\twoVgraph\fourVgraph^2+80\twoVgraph^2\sixVgraph+160\fourVgraph\sixVgraph} \\
\end{array}
$$

\subsection{Writing graphs in terms of polywheels}

$$
\begin{array}{lll}
{\bf k=1}      &   & \\
\twoVgraph     & = & \langle w_2\rangle  \\
               &   & \\
{\bf k=2}      &   & \\
\twoVgraph^2   & = & \langle w_2^2\rangle -\frac{4}{5}\langle w_4\rangle  \\
\fourVgraph    & = & \frac{2}{5}\langle w_4\rangle  \\
               &   & \\
{\bf k=3}      &   & \\
\twoVgraph^3   & = & \langle w_2^3\rangle -\frac{12}{5}\langle w_2w_4\rangle +\frac{64}{35}\langle w_6\rangle  \\
\twoVgraph\fourVgraph & = & \frac{2}{5}\langle w_2w_4\rangle -\frac{16}{35}\langle w_6\rangle  \\
\sixVgraph     & = & \frac{4}{35}\langle w_6\rangle  \\
\phantom{\twoVgraph\eightVgraphII} & &
\phantom{\frac{1}{12}\twoVgraph\fourVgraph^2-\fourVgraph\sixVgraph-\frac{1}{75}\langle
w_2w_4^2\rangle +\frac{8}{105}\langle w_4w_6\rangle
+\frac{8}{175}\langle w_2w_8\rangle} \\
\end{array}
$$

\newpage
$$
\begin{array}{lll}
{\bf k=4}      &   & \\
\twoVgraph^4   & = &
\langle w_2^4\rangle -\frac{24}{5}\langle w_2^2w_4\rangle +\frac{48}{25}\langle w_4^2\rangle +\frac{256}{35}\langle w_2w_6\rangle -\frac{1152}{175}\langle w_8\rangle
\\
\twoVgraph^2\fourVgraph & = &
\frac{2}{5}\langle w_2^2w_4\rangle -\frac{8}{25}\langle w_4^2\rangle -\frac{32}{35}\langle w_2w_6\rangle +\frac{192}{175}\langle w_8\rangle 
\\
\twoVgraph\sixVgraph & = &
-\frac{1}{2}\fourVgraph^2+\frac{2}{25}\langle w_4^2\rangle +\frac{4}{35}\langle w_2w_6\rangle -\frac{48}{175}\langle w_8\rangle 
\\
\eightVgraphI  & = &
\frac{1}{12}\fourVgraph^2-\frac{1}{75}\langle w_4^2\rangle +\frac{8}{175}\langle w_8\rangle  \\
\eightVgraphII & = &
-\frac{41}{96}\fourVgraph^2+\frac{41}{600}\langle w_4^2\rangle
 -\frac{16}{175}\langle w_8\rangle \\
               &   & \\
{\bf k=5}      &   & \\
\twoVgraph^5 & = &
\langle w_2^5\rangle -8\langle w_2^3w_4\rangle +\frac{48}{5}\langle
w_2w_4^2\rangle +\frac{128}{7}\langle w_2^2w_6\rangle
-\frac{512}{35}\langle w_4w_6\rangle \\
 & & \phantom{\twoVgraph} -\frac{1152}{35}\langle w_2w_8\rangle +\frac{12288}{385}\langle w_{10}\rangle  \\
\twoVgraph^3\fourVgraph & = &
\frac{2}{5}\langle w_2^3w_4\rangle -\frac{24}{25}\langle w_2w_4^2\rangle -\frac{48}{35}\langle w_2^2w_6\rangle +\frac{64}{35}\langle w_4w_6\rangle +\frac{576}{175}\langle w_2w_8\rangle
\\
 & & \phantom{\twoVgraph} -\frac{1536}{385}\langle w_{10}\rangle  \\
\twoVgraph^2\sixVgraph & = &
-\twoVgraph\fourVgraph^2+\frac{4}{25}\langle w_2w_4^2\rangle
 +\frac{4}{35}\langle w_2^2w_6\rangle -\frac{16}{35}\langle
 w_4w_6\rangle -\frac{96}{175}\langle w_2w_8\rangle\\
 & & \phantom{\twoVgraph} +\frac{384}{385}\langle w_{10}\rangle  \\
\twoVgraph\eightVgraphI & = &
\frac{1}{12}\twoVgraph\fourVgraph^2-\fourVgraph\sixVgraph-\frac{1}{75}\langle
w_2w_4^2\rangle +\frac{8}{105}\langle w_4w_6\rangle +\frac{8}{175}\langle w_2w_8\rangle \\
 & & \phantom{\twoVgraph} -\frac{64}{385}\langle w_{10}\rangle  \\
\twoVgraph\eightVgraphII & = &
-\frac{41}{96}\twoVgraph\fourVgraph^2-\frac{9}{8}\fourVgraph\sixVgraph+\frac{41}{600}\langle w_2w_4^2\rangle -\frac{11}{105}\langle w_4w_6\rangle 
\\
 & & \phantom{\twoVgraph} -\frac{16}{175}\langle w_2w_8\rangle +\frac{184}{1155}\langle w_{10}\rangle  \\
\tenVgraphI & = &
\frac{5}{24}\fourVgraph\sixVgraph-\frac{1}{105}\langle w_4w_6\rangle +\frac{8}{385}\langle w_{10}\rangle 
\\
\tenVgraphII & = & 
-\frac{55}{192}\fourVgraph\sixVgraph+\frac{11}{840}\langle w_4w_6\rangle -\frac{23}{1155}\langle w_{10}\rangle 
\\
\end{array}
$$

%% file: app2.tex
\section{Riemann-Roch formula}
\subsection{$\chi_y$-genus in terms of Chern numbers} 

$$
\begin{array}{|c|c|c|}
\hline
k & \chi^m & \mbox{Chern number expression} \\
\hline\hline
1 & \chi^0 & c_2/12 \\
  & \chi^1 & -10c_2/12 \\
\phantom{k} & & \\
\hline
2 & \chi^0 & (3c_2^2-c_4)/{720} \\
  & \chi^1 & (12c_2^2-124c_4)/{720} \\
  & \chi^2 & (18c_2^2+474c_4)/{720} \\
\phantom{k} & & \\
\hline
3 & \chi^0 & (10c_2^3-9c_2c_4+2c_6)/{60480} \\
  & \chi^1 & (60c_2^3-306c_2c_4-492c_6)/{60480} \\
  & \chi^2 & (150c_2^3-1143c_2c_4+13134c_6)/{60480} \\
  & \chi^3 & (200c_2^3-1692c_2c_4-33224c_6)/{60480} \\
\phantom{k} & & \\
\hline
4 & \chi^0 & (21c_2^4-34c_2^2c_4+5c_4^2+13c_2c_6-3c_8)/{3628800} \\
  & \chi^1 &
(168c_2^4-872c_2^2c_4+1240c_4^2-1816c_2c_6-744c_8)/{3628800} \\
  & \chi^2 &
(588c_2^4-4552c_2^2c_4+7340c_4^2+3964c_2c_6+86316c_8)/{3628800} \\
  & \chi^3 &
(1176c_2^4-10904c_2^2c_4+18280c_4^2+32408c_2c_6-857688c_8)/{3628800} \\
  & \chi^4 &
\hspace*{2mm}(1470c_2^4-14380c_2^2c_4+24350c_4^2+53230c_2c_6+1739310c_8)/{3628800}\hspace*{2mm}
  \\
\phantom{k} & & \\
\hline
\end{array}
$$

\subsection{Chern numbers in terms of $\chi_y$-genus}

$$
\begin{array}{|c|c|c|}
\hline
k & \mbox{Chern no.} & \chi_y\mbox{-genus} \\
\hline\hline
1 & c_2 & 12\chi^0 \\
\phantom{k} & & \\
\hline
2 & c_2^2 & 248\chi^0-2\chi^1 \\
  & c_4 & 24\chi^0-6\chi^1 \\
\phantom{k} & & \\
\hline
3 & c_2^3 & 7272\chi^0-184\chi^1-8\chi^2 \\
  & c_2c_4 & 1368\chi^0-208\chi^1-8\chi^2 \\
  & c_6 & 36\chi^0-16\chi^1+4\chi^2 \\
\phantom{k} & & \\
\hline
4 & c_2^4 & 3s \\
  & c_2^2c_4 & \hspace*{56pt}2s-116032\chi^0-372\chi^1+112\chi^2+12\chi^3\hspace*{56pt} \\
  & c_4^2 & s-74960\chi^0+777\chi^1+332\chi^2+33\chi^3 \\
  & c_2c_6 & 4512\chi^0-1278\chi^1+168\chi^2+18\chi^3 \\
  & c_8 & 48\chi^0-27\chi^1+12\chi^2-3\chi^3 \\
\phantom{k} & & \\
\hline
\end{array}
$$

\newpage
$$
\begin{array}{|c|c|c|}
\hline
k & \mbox{Chern no.} & \chi_y\mbox{-genus} \\
\hline\hline
1 & s_2 & -24\chi^0 \\
\phantom{k} & &
\hspace*{112pt}\phantom{2s-116032\chi^0-372\chi^1+112\chi^2+12\chi^3}
\\
\hline
2 & s_2^2 & 992\chi^0-8\chi^1 \\
  & s_4 & 400\chi^0+20\chi^1 \\
\phantom{k} & & \\
\hline
3 & s_2^3 & -58176\chi^0+1472\chi^1+64\chi^2 \\
  & s_2s_4 & -18144\chi^0-928\chi^1-32\chi^2 \\
  & s_6 & -6552\chi^0-784\chi^1-56\chi^2 \\
\phantom{k} & & \\
\hline
4 & s_2^4 & 48s \\
  & s_2^2s_4 & -8s+1856512\chi^0+5952\chi^1-1792\chi^2-192\chi^3 \\
  & s_4^2 & -4s+657152\chi^0+18384\chi^1+3520\chi^2+336\chi^3 \\
  & s_2s_6 & -12s+1446528\chi^0-10872\chi^1+672\chi^2+72\chi^3 \\
  & s_8 & -6s+664128\chi^0-3924\chi^1+1680\chi^2+204\chi^3 \\
\phantom{k} & & \\
\hline
\end{array}
$$

%% file: app3.tex
\section {Hirzebruch $\chi_y$-genus of $S^{[k]}$ and $T^{[[k]]}$}

$$
\begin{array}{|c|c|c|}
\hline
k & \mbox{Space} & \chi_y\mbox{-genus} \\
\hline\hline
1 & S           & 2-20y+2y^2 \\
2 & S^{[2]}     & 3-42y+234y^2-42y^3+3y^4 \\
3 & S^{[3]}     & 4-64y+508y^2-2048y^3+508y^4-64y^5+4y^6 \\
4 & S^{[4]}     & 5-86y+785y^2-4556y^3+14786y^4 -4556y^5+785y^6-86y^7+5y^8\\
\phantom{k} & & \\
\hline\hline
2 & T^{[[2]]}   & 3-6y+90y^2-6y^3+3y^4 \\
3 & T^{[[3]]}   & 4-8y+44y^2-336y^3+44y^4-8y^5+4y^6 \\
4 & T^{[[4]]}   & 5-10y+15y^2-20y^3+650y^4-20y^5+15y^6-10y^7+5y^8 \\
\phantom{k} & & \\
\hline
\end{array}
$$

%% file: app4.tex
\section {Chern numbers in low degree}

$$
\begin{array}{|c|c|c|c|}
\hline
k & \mbox{Chern number} & S^{[k]} & T^{[[k]]} \\
\hline\hline
1 & s_2         & -48                  & -48 \\
\phantom{k} & & & \\
\hline
2 & s_2^2       & 3312                 & 3024 \\
  & s_4         & 360                  & 1080 \\
\phantom{k} & & & \\
\hline
3 & s_2^3       & -294400              & -241664 \\
  & s_2s_4      & -29440               & -66560 \\
  & s_6         & -4480                & -22400 \\
\phantom{k} & & & \\
\hline 
4 & s_2^4       & 48s                  & 48s \\
  & s_2^2s_4    & -8s+8238720          & -8s+9200000 \\
  & s_4^2       & -4s+2937120          & -4s+3148000 \\
  & s_2s_6      & \hspace*{10mm}-12s+8367120\hspace*{10mm}         & \hspace*{10mm}-12s+7350000\hspace*{10mm} \\
  & s_8         & -6s+4047480          & -6s+3381000 \\
\phantom{k} & & & \\
\hline
\hline
4 & s           & 664080               & 490000 \\
  & s_2^4       & 31875840             & 23520000 \\
  & s_2^2s_4    & 2926080              & 5280000 \\
  & s_4^2       & 280800               & 1188000 \\
  & s_2s_6      & 398160               & 1470000 \\
  & s_8         & 63000                & 441000 \\
\phantom{k} & & & \\
\hline
\end{array}
$$

%% file: app5.tex
\section {Rozansky-Witten invariants in low degree}
\subsection{Irreducible manifolds}

$$
\begin{array}{|c|c|c|c|}
\hline
k & \Gamma     & b_{\Gamma}(S^{[k]}) & b_{\Gamma}(T^{[[k]]}) \\
\hline\hline
1 & \twoVgraph                 & 48           & 48 \\
\phantom{k} & & & \\
\hline
2 & \twoVgraph^2               & 3600         & 3888 \\
  & \fourVgraph                & -144         & -432 \\
\phantom{k} & & & \\
\hline
3 & \twoVgraph^3               & 373248       & 442368 \\
  & \twoVgraph\fourVgraph      & -13824       & -36864 \\
  & \sixVgraph                 & 512          & 2560 \\
\phantom{k} & & & \\
\hline
4 & \twoVgraph^4               & 49787136     & 64800000 \\
  & \twoVgraph^2\fourVgraph    & \hspace*{18mm}-1693440\hspace*{18mm}    & \hspace*{18mm}-4320000\hspace*{18mm} \\
  & \fourVgraph^2              & 57600        & 288000 \\
  & \twoVgraph\sixVgraph       & 56448        & 240000 \\
  & \eightVgraphI              & -1824        & -12000 \\
  & \eightVgraphII             & 348          & -1500 \\
\phantom{k} & & & \\
\hline
\end{array}
$$

\newpage
\subsection{Reducible manifolds}
$$
\begin{array}{|c|c|c|c|c|c|}
\hline
k & \Gamma     & \hspace*{25mm} & \hspace*{25mm} & \hspace*{25mm} &
\hspace*{25mm} \\
\hline\hline
2 &            & b_{\Gamma}(S^2) & & & \\ 
\hline
  & \twoVgraph^2            & 4608     & & & \\
  & \fourVgraph             & 0        & & & \\
\phantom{k} & & & & & \\
\hline
3 &            & b_{\Gamma}(S\times S^{[2]}) & b_{\Gamma}(S\times
  T^{[[2]]}) & b_{\Gamma}(S^3) & \\ 
\hline
  & \twoVgraph^3            & 518400   & 559872   & 663552 & \\
  & \twoVgraph\fourVgraph   & -6912    & -20736   & 0      & \\
  & \sixVgraph              & 0        & 0        & 0      & \\
\phantom{k} & & & & &\\
\hline
4 &            & b_{\Gamma}(S\times S^{[3]}) &
b_{\Gamma}(S^{[2]}\times S^{[2]}) & b_{\Gamma}(S^2\times S^{[2]}) &
b_{\Gamma}(S^4) \\
\hline
  & \twoVgraph^4            & 71663616 & 77760000 & 99532800 & 127401984  \\
  & \twoVgraph^2\fourVgraph & -1327104 & -1036800 & -663552  & 0 \\
  & \fourVgraph^2           & 0        & 41472    & 0        & 0 \\
  & \twoVgraph\sixVgraph    & 24576    & 0        & 0        & 0 \\
  & \eightVgraphI           & 0        & 0        & 0        & 0 \\
  & \eightVgraphII          & 0        & 0        & 0        & 0 \\
\phantom{k} & & & & & \\
\hline
\end{array}
$$

\newpage
\subsection{Virtual manifolds}
$$
\begin{array}{|c|c|c|}
\hline
k & \hspace*{10mm}\Gamma\hspace*{10mm} & \hspace*{47mm}b_{\Gamma}(C_k)\hspace*{47mm} \\
\hline\hline
1 & \twoVgraph                 & -6 \\
\phantom{k} & & \phantom{-7776} \\
\hline
2 & \twoVgraph^2               & 36 \\
  & \fourVgraph                & 12 \\
\phantom{k} & & \\
\hline
3 & \twoVgraph^3               & -216 \\
  & \twoVgraph\fourVgraph      & -72 \\
  & \sixVgraph                 & -24 \\
\phantom{k} & & \\
\hline
4 & \twoVgraph^4               & 1296 \\
  & \twoVgraph^2\fourVgraph    & 432 \\
  & \fourVgraph^2              & 144 \\
  & \twoVgraph\sixVgraph       & 144 \\
  & \eightVgraphI              & 48 \\
  & \eightVgraphII             & 24 \\
\phantom{k} & & \\
\hline
\end{array}
$$

%% file: thesis.bbl
\providecommand{\bysame}{\leavevmode\hbox to3em{\hrulefill}\thinspace}
\providecommand{\MR}{\relax\ifhmode\unskip\space\fi MR }
\providecommand{\MRhref}[2]{%
  \href{http://www.ams.org/mathscinet-getitem?mr=#1}{#2}
}
\providecommand{\href}[2]{#2}
\begin{thebibliography}{10}

\bibitem{atiyah57}
M.~Atiyah, \emph{Complex analytic connections in fibre bundles}, Trans. Am.
  Math. Soc. \textbf{85} (1957), 181--207.

\bibitem{as92}
S.~Axelrod and I.~Singer, \emph{Chern-{S}imons perturbation theory},
  Proceedings of the XXth Conference on Differential Geometric Methods in
  Physics (New York, 1991) (S. Catto and A. Rocha, eds) World Scientific
  (1992), 3--45.

\bibitem{as94}
\bysame, \emph{Chern-{S}imons perturbation theory {II}}, Jour. Diff. Geom.
  \textbf{39} (1994), 173--213.

\bibitem{barnatan95}
D.~Bar-Natan, \emph{On the {V}assiliev knot invariants}, Topology \textbf{34}
  (1995), 423--472.

\bibitem{barnatan95I}
\bysame, \emph{Perturbative {C}hern-{S}imons theory}, Jour. of Knot Theory and
  its Ramifications \textbf{4} (1995), no.~4, 503--548.

\bibitem{bgrt97}
D.~Bar-Natan, S.~Garoufalidis, L.~Rozansky, and D.~Thurston, \emph{The {A}arhus
  integral {I}-{III}}, preprints (1997--1998).

\bibitem{bgrt98}
\bysame, \emph{Wheels, wheeling, and the {K}ontsevich integral of the unknot},
  preprint {\bf q-alg/9703025 v3} (1998).

\bibitem{blt}
D.~Bar-Natan, T.~Le, and D.~Thurston, in preparation.

\bibitem{bpv84}
W.~Barth, C.~Peters, and A.~Van~de Ven, \emph{Compact complex surfaces}, A
  Series of Modern Surveys in Math. {\bf 4}, Springer-Verlag, 1984.

\bibitem{beauville83}
A.~Beauville, \emph{Vari{\'e}t{\'e}s {K}{\"a}hl{\'e}riennes dont la 1{\'e}re
  classe de {C}hern est nulle}, Jour. Diff. Geom. \textbf{18} (1983), 755--782.

\bibitem{beauville99}
\bysame, \emph{Riemannian holonomy and algebraic geometry}, Emmy Noether
  lectures (1999).

\bibitem{besse87}
A.~L. Besse, \emph{Einstein manifolds}, Springer-Verlag, Berlin Heidelberg,
  1987.

\bibitem{by52}
S.~Bochner and K.~Yano, \emph{Tensor-fields in non-symmetric connections}, Ann.
  of Math. \textbf{56} (1952), no.~3, 504--519.

\bibitem{bogomolov96i}
F.~Bogomolov, \emph{On {G}uan's examples of simply connected non-{K}{\"a}hler
  compact complex manifolds}, Amer. Jour. Math. \textbf{118} (1996),
  1037--1046.

\bibitem{bogomolov96}
\bysame, \emph{On the cohomology ring of a simple hyperk{\"a}hler manifold (on
  the results of {V}erbitsky)}, Geom. Funct. Anal. \textbf{6} (1996), 612--618.

\bibitem{bogomolov78}
V.~A. Bogomolov, \emph{Hamiltonian {K}{\"a}hler manifolds}, Dokl. Akad. Nauk
  SSSR (Mat.) \textbf{243} (1978), 1101--1104. English translation: Soviet
  Math. Dokl. {\bf 19} (1978), 1462--1465.

\bibitem{bt94}
R.~Bott and C.~Taubes, \emph{On the self-linking of knots}, Jour. Math. Phys.
  \textbf{35} (1994), no.~10, 5247--5287.

\bibitem{cheah96}
J.~Cheah, \emph{On the cohomology of {H}ilbert schemes of points}, Jour. Alg.
  Geom. \textbf{5} (1996), 479--511.

\bibitem{deligne96}
P.~Deligne, \emph{Letter to {D}. {B}ar-{N}atan {J}an 25th 1996}, available from
  {\bf http://www.ma.huji.ac.il/$\sim$drorbn/Deligne/}.

\bibitem{dijkgraaf99}
R.~Dijkgraaf, \emph{Instanton strings and hyperk{\"a}hler geometry}, Nuclear
  Phys. B \textbf{543} (1999), no.~3, 545--571.

\bibitem{duflo70}
M.~Duflo, \emph{Caract{\`e}res des groupes et des alg{\`e}bres de {L}ie
  r{\'e}solubles}, Ann. scient. {\'E}c. Norm. Sup. \textbf{4(3)} (1970),
  23--74.

\bibitem{egl99}
G.~Ellingsrud, L.~G{\"o}ttsche, and M.~Lehn, \emph{On the cobordism class of
  {H}ilbert schemes of points on surfaces}, preprint {\bf math.AG/9904095}
  (1999).

\bibitem{fogarty68}
J.~Fogarty, \emph{Algebraic families on an algebraic surface}, Amer. Jour.
  Math. \textbf{90} (1968), 511--521.

\bibitem{fujiki83}
A.~Fujiki, \emph{On primitive symplectic compact {K}{\"a}hler {V}-manifolds of
  dimension four}, in {"}Classification of Algebraic and Analytic Manifolds{"},
  K. Ueno (ed.), Progress in Mathematics, Birkh{\"a}user \textbf{39} (1983),
  71--125.

\bibitem{gottsche}
L.~G{\"o}ttsche, private communication (1998).

\bibitem{gs93}
L.~G{\"o}ttsche and W.~Soergel, \emph{Perverse sheaves and the cohomology of
  {H}ilbert schemes of smooth algebraic surfaces}, Math. Ann. \textbf{296}
  (1993), 235--245.

\bibitem{guan95}
D.~Guan, \emph{Examples of compact holomorphic symplectic manifolds which are
  not {K}{\"a}hlerian {II}}, Invent. Math. \textbf{121} (1995), no.~1,
  135--145.

\bibitem{ht99}
N.~Habegger and G.~Thompson, \emph{The universal perturbative quantum
  3-manifold invariant, {R}ozansky-{W}itten invariants, and the generalized
  {C}asson invariant}, preprint (1999).

\bibitem{hirzebruch78}
F.~Hirzebruch, \emph{Topological methods in algebraic geometry, 3rd edn.},
  Springer-Verlag, 1978.

\bibitem{hklr87}
N.~J. Hitchin, A.~Karlhede, U.~Lindstr{\"o}m, and M.~Ro{\v c}ek,
  \emph{Hyperk{\"a}hler metrics and supersymmetry}, Comm. Math. Phys.
  \textbf{108} (1987), 535--589.

\bibitem{hl97}
D.~Huybrechts and M.~Lehn, \emph{The geometry of moduli spaces of sheaves},
  Aspects of Mathematics {\bf E31}, Verlag-Vieweg, 1997.

\bibitem{kapranov98}
M.~Kapranov, \emph{Rozansky-{W}itten invariants via {A}tiyah classes (old
  version)}, privately communicated (1998).

\bibitem{kapranov99}
\bysame, \emph{Rozansky-{W}itten invariants via {A}tiyah classes}, Compositio
  Math. \textbf{115} (1999), 71--113.

\bibitem{kodaira64}
K.~Kodaira, \emph{On the structure of compact complex analytic surfaces {I}},
  Amer. Jour. Math. \textbf{86} (1964), 751--798.

\bibitem{kodaira66}
\bysame, \emph{On the structure of compact complex analytic surfaces {II}},
  Amer. Jour. Math. \textbf{88} (1966), 682--721.

\bibitem{kontsevich93}
M.~Kontsevich, \emph{{V}assiliev's knot invariants}, Advance in Soviet
  Mathematics \textbf{16} (1993), 137--150.

\bibitem{kontsevich94}
\bysame, \emph{Feynman diagrams and low-dimensional topology}, First European
  Congress of Mathematics (Paris, 1992), Vol. II, Progress in Mathematics {\bf
  120}, Birkh{\"a}user (1994), 97--121.

\bibitem{kontsevich97}
\bysame, \emph{Deformation quantization of poisson manifolds}, IHES preprint.
  See also {\bf q-alg/9709040} (1997).

\bibitem{kontsevich99}
\bysame, \emph{Rozansky-{W}itten invariants via formal geometry}, Compositio
  Math. \textbf{115} (1999), 115--127.

\bibitem{lmo98}
T.~T.~Q. Le, J.~Murakami, and T.~Ohtsuki, \emph{On a universal perturbative
  invariant of 3-manifolds}, Topology \textbf{37} (1998), no.~3, 539--574.

\bibitem{mukai84}
S.~Mukai, \emph{Symplectic structure of the moduli space of sheaves on an
  abelian or {K}3 surface}, Invent. Math. \textbf{77} (1984), 101--116.

\bibitem{ogrady99}
K.~O'Grady, \emph{Desingularized moduli spaces of sheaves on a {K}3}, J. Reine
  Angew. Math. \textbf{512} (1999), 49--117.

\bibitem{ohtsuki82}
M.~Ohtsuki, \emph{A residue formula for {C}hern classes associated with
  logarithmic connections}, Tokyo J. Math. \textbf{5} (1982), no.~1, 13--21.

\bibitem{rw97}
L.~Rozansky and E.~Witten, \emph{Hyperk{\"a}hler geometry and invariants of
  three-manifolds}, Selecta Math. \textbf{3} (1997), 401--458.

\bibitem{salamon89}
S.~M. Salamon, \emph{Riemannian geometry and holonomy groups}, Pitman Research
  Notes in Mathematics {\bf 201}, Longman, Harlow, 1989.

\bibitem{salamon96}
\bysame, \emph{On the cohomology of k{\"a}hler and hyperk{\"a}hler manifolds},
  Topology \textbf{35} (1996), no.~1, 137--155.

\bibitem{thom54}
R.~Thom, \emph{Quelques propri{\'e}t{\'e}s globales des vari{\'e}t{\'e}s
  diff{\'e}rentiables}, Comm. Math. Helv. \textbf{28} (1954), 17--86.

\bibitem{tian87}
G.~Tian, \emph{Smoothness of the universal deformation space of compact
  {C}alabi-{Y}au manifolds and its {P}etersson-{W}eil metric}, In: Yau, S. T.
  Mathematical Aspects of String theory. World Scientific (1987), 629--646.

\bibitem{todorov89}
A.~N. Todorov, \emph{The {P}etersson-{W}eil geometry of the moduli space of
  {SU}(n$\geq$3) ({C}alabi-{Y}au) manifolds}, Comm. Math. Phys. \textbf{126}
  (1989), 325--346.

\bibitem{vogel96}
P.~Vogel, \emph{Algebraic structures on modules of diagrams}, Universit{\'e}
  Paris VII preprint, to appear (1996).

\bibitem{witten89}
E.~Witten, \emph{Quantum field theory and the {J}ones polynomial}, Comm. Math.
  Phys. \textbf{121} (1989), 351--399.

\bibitem{yau78}
S.~T. Yau, \emph{On the {R}icci curvature of a compact {K}{\"a}hler manifold
  and the complex {M}onge-{A}mp{\`e}re equation}, Comm. Pure Appl. Math.
  \textbf{31} (1978), 339--411.

\end{thebibliography}
